\documentclass[11pt,reqno]{amsart}

\usepackage[T1]{fontenc}
\usepackage{lmodern}
\AtBeginDocument{\DeclareFontShape{T1}{lmr}{bx}{sc}{<->ssub*lmr/bx/n}{}}
\usepackage{amssymb}
\usepackage{amsmath,amsthm,amscd,amsfonts}
\usepackage{booktabs}
\usepackage{url}
\usepackage{array}

\newtheorem{theorem}{Theorem}[section]
\newtheorem{lemma}[theorem]{Lemma}
\newtheorem{corollary}[theorem]{Corollary}
\newtheorem{proposition}[theorem]{Proposition}
\newtheorem{definition}[theorem]{Definition}
\newtheorem{remark}[theorem]{Remark}
\newtheorem{problem}[theorem]{Problem}
\newtheorem{conjecture}[theorem]{Conjecture}

\makeatletter
\def\@tocline#1#2#3#4#5#6#7{\relax
	\ifnum #1>\c@tocdepth
	\else
	\par \addpenalty\@secpenalty\addvspace{#2}%
	\begingroup \hyphenpenalty\@M
	\@ifempty{#4}{%
		\@tempdima\csname r@tocindent\number#1\endcsname\relax
	}{%
		\@tempdima#4\relax
	}%
	\parindent\z@ \leftskip#3\relax \advance\leftskip\@tempdima\relax
	\rightskip\@pnumwidth plus4em \parfillskip-\@pnumwidth
	#5\leavevmode\hskip-\@tempdima
	\ifcase #1
	\or\or \hskip 1em \or \hskip 2em \else \hskip 3em \fi%
	#6\nobreak\relax
	\leaders\hbox{$\m@th\mkern 4.5mu\hbox{.}\mkern 4.5mu$}\hfill
	\hbox to\@pnumwidth{\hfil\@tocpagenum{#7}}\par
	\nobreak
	\endgroup
	\fi}
\def\l@section{\@tocline{1}{4pt}{0pt}{}{}}
\def\l@subsection{\@tocline{2}{0pt}{1.5em}{}{}}
\makeatother

\newcommand{\N}{\mathbb{N}}
\newcommand{\Z}{\mathbb{Z}}
\newcommand{\Q}{\mathbb{Q}}
\newcommand{\R}{\mathbb{R}}
\newcommand{\C}{\mathbb{C}}
\newcommand{\E}{\mathbb{E}}
\newcommand{\Pb}{\mathbb{P}}
\newcommand{\Colmin}{\operatorname{Col}_{\min}}
\newcommand{\TV}{\mathrm{TV}}

\newcommand{\Bin}{\mathrm{Bin}}

\begin{document}
	
	\begin{center}
		{\Large \bf Coalescence of orbits, natural density of invariant sets\\[1.5mm] and primes for Collatz-type maps}\\[4mm]
		
		{\large \bf Aleksandar Petojevi\'{c}$^{1}$ and Shunya Tamura$^{2}$}\\[3mm]
		$^{1}$ Faculty of Education, University of Novi Sad,\\
		Podgori\v{c}ka 4, 25000 Sombor, Serbia\\[1.5mm]
		$^{2}$ Okegawa City Okegawa West Junior High School,\\
		3680-1 Kawataya, Okegawa, 363-0027, Saitama, Japan\\[1.5mm]
		{\bf E-mail $^{1}$: apetoje@uns.ac.rs}\\[1mm]
		{\bf E-mail $^{2}$: shunya.tamura059@gmail.com}\\[1mm]
		{\bf ORCID ID$^{1}$: 0000-0003-1491-165X}\\[1mm]
		{\bf ORCID ID$^{2}$: 0009-0000-2162-2536}\\[3mm]
	\end{center}
	
	\vspace{5mm}
	
\begin{center}
  \textbf{Abstract}
\end{center}

\begin{quotation}
\noindent
We study balanced coalescence for generalized Collatz maps with
coprime parameters $\ell,m\ge2$. For canonical residue systems
satisfying $m^{\ell-1}<\ell^\ell$, let $f_n(c)$ denote the probability
that two trajectories initially separated by $c\ne0$ have not
merged at a common time with equal counts of nonzero branches.
We prove, uniformly in $c$ and $n\ge1$,
\[
\frac2\ell\sqrt{\frac2{\pi(n+1)}}
\le f_n(c) \le
\min\left\{1,C_{\ell,m}
\sqrt{\frac{1+\log(2+|c|)}{n+1}}\right\},
\]
with an effectively computable constant. This extends the known
standard Collatz estimate to general coprime bases, including composite
$\ell$, and yields the sharp $n^{-1/2}$ order for each fixed nonzero
shift. For general admissible residue systems we identify an arithmetic
obstruction $d_r$ and characterize almost-sure coalescence via finitely
many reachability requirements. Without drift or positivity assumptions,
the group of almost-surely coalescing shifts is either trivial or
$d_r\mathbb Z$. We also establish H\"older-type density-transfer bounds
for invariant sets sampled along primes, sums of two squares, and rough
Chen primes. Under the Riemann hypothesis, prime and integer averages
agree uniformly for bounded defect cutoffs.
\end{quotation}

\vspace{5mm}



\vspace{3mm}
\newpage
\tableofcontents

\section{\textbf{Introduction}}
\label{sec:intro}

We study quantitative coalescence for generalized Collatz maps
and its consequences for invariant sets. The standard example is
the shortcut Collatz map $T:\N\to\N$,
\begin{equation}\label{eq:T}
 T(n)=
 \begin{cases}
 n/2,&n\equiv0\pmod2,\\
 (3n+1)/2,&n\equiv1\pmod2.
 \end{cases}
\end{equation}
The Collatz conjecture asserts that every positive integer has
an orbit reaching $1$. Our main results concern the relation
between two orbits: a uniform estimate for their balanced
coalescence time, a finite absorption criterion for general
residue systems, and the structure of the coalescing shifts.

The density theory of the standard map provides essential
background. Terras and Everett proved descent below the starting
value for almost all positive integers in natural density
\cite{Terras1976,Everett1977}. Allouche obtained descent below
a power threshold \cite{Allouche1979}, and Korec established
the range $\theta>\log3/\log4$ for the threshold $N^\theta$
\cite{Korec1994}. Tao proved that almost every orbit attains a
value below any prescribed function tending to infinity, in
logarithmic density \cite{Tao2022}. Inselmann's trajectory
estimates give a power-saving exceptional set for every fixed
power threshold $N^\theta$, $\theta>0$
\cite[Definition~2.6 and Corollary~3.2]{Inselmann2024}.
Mazur obtained almost-boundedness in natural density with a
logarithmic bound on the required iteration time \cite{Mazur2026}.
Nashida gives a power-saving bound in the barrier height, uniform
in the counting endpoint \cite{Nashida2026}.

Coalescence also has substantial prior work. Duran proved
density-one coalescence of consecutive starting values
\cite{Duran2026}. Leen established sharp synchronization
estimates, including the uniform square-root bound with
logarithmic dependence on the displacement
\cite[Theorem~2.11]{LeenRates2026}. Leen also proved natural-density
existence for every set invariant above a finite threshold for
the ordinary Collatz map \cite[Theorem~2]{LeenBrownian2026};
this applies to the finite-defect invariant sets for
\eqref{eq:T}. These standard-map results are used as background
for the general theory developed here.

For generalized maps, Gon\c{c}alves, Greenfeld and Madrid proved
almost-boundedness in logarithmic density under negative drift
\cite{GGM2025}. Nashida's natural-density theorem covers the
admissible residue family and gives a fixed-barrier estimate
\cite{NashidaGeneral2026}. Sultanow proved logarithmic-density
existence for completely invariant sets of the $3x+c$ maps with
odd $c\ge-1$ \cite[Theorem~1]{Sultanow2026}. Accordingly, the
principal contributions presented below concern quantitative
balanced coalescence in general bases and the arithmetic
structure of general residue systems.

\subsection*{Main coalescence results}

Let $\ell,m\ge2$ be coprime integers, and let
$r=(r_0,\ldots,r_{\ell-1})$ satisfy
$r_0=0$ and $r_j\equiv-mj\pmod\ell$. Define
\[
 T_r(x)=
 \begin{cases}
 x/\ell,&x\equiv0\pmod\ell,\\
 (mx+r_j)/\ell,&x\equiv j\not\equiv0\pmod\ell.
 \end{cases}
\]
This map acts on $\Z_\ell=\varprojlim_a\Z/\ell^a\Z$.
The residues are \emph{canonical} when $1\le r_j\le\ell-1$
for $j\ne0$; in that case we write $T_{\ell,m}$.
A general residue system is \emph{admissible} when
$mj+r_j>1$ for $1\le j<\ell$, so that $T_r$ preserves $\N$.

Let $q_n(x)$ count the nonzero branches among the first $n$
iterates of $x$. For Haar-distributed $\beta\in\Z_\ell$, put
\[
 \begin{split}
 \tau_c^{\mathrm{bal}}&=\inf\{n\ge0:T_r^n(\beta)=T_r^n(\beta+c),\
                 q_n(\beta)=q_n(\beta+c)\},\\
 f_n(c)&=\mathbb P(\tau_c^{\mathrm{bal}}>n).
 \end{split}
\]
Thus balanced coalescence requires both equality of the iterates
at a common time and equality of their accumulated nonzero-branch
counts. Time is measured in iterations of $T_r$.
Equivalently, $f_n(c)$ is the proportion of initial classes
modulo $\ell^n$ on which balanced coalescence has not occurred
by time $n$.

\medskip
\noindent
\textbf{Uniform coalescence for canonical residues.}
Assume the negative-drift condition
\[
 \mu_{\ell,m}=\frac{\ell-1}{\ell}\log m-\log\ell<0,
 \qquad\text{equivalently}\qquad m^{\ell-1}<\ell^\ell.
\]
Theorems~\ref{thm:Gsharp} and~\ref{thm:lowerunif} give an
effectively computable constant $C_{\ell,m}$ such that, for
every $c\in\Z\setminus\{0\}$ and $n\ge1$,
\[
 \frac2\ell\sqrt{\frac2{\pi(n+1)}}
 \le f_n(c)
 \le\min\left\{1,
 C_{\ell,m}\sqrt{\frac{1+\log(2+|c|)}{n+1}}\right\}.
\]
In particular, coalescence is almost sure, and the tail has
order $n^{-1/2}$ for each fixed nonzero shift. No limiting
leading coefficient is asserted. The standard case
$(\ell,m)=(2,3)$ is already contained in
\cite[Theorem~2.11]{LeenRates2026}; the result here treats
all coprime bases satisfying negative drift, including
composite $\ell$ and the range $m<\ell$.

\medskip
\noindent
\textbf{A finite absorption certificate.}
For arbitrary residues, let $d_r$ be the largest positive
divisor coprime to $\ell m$ of
\[
 \gcd\bigl(r_2-r_1,\ldots,r_{\ell-1}-r_1,(\ell-1)r_1\bigr),
\]
where the list of differences is empty when $\ell=2$.
If $d_r\nmid c$, balanced coalescence is impossible
(Proposition~\ref{thm:GO}). Assume now that $r$ is admissible
and $\mu_{\ell,m}<0$, and put
\[
 R=\max_{1\le j<\ell}|r_j|,
 \qquad C_0(r)=\frac{R(\ell+1)}{\ell^2-m}.
\]
Theorem~\ref{thm:GR} proves that almost-sure balanced
coalescence for every $c\in d_r\Z$ is equivalent to the
following finite collection of requirements: for each integer
$c$ with $d_r\mid c$ and $0<|c|\le C_0(r)$, some finite
residue word produces balanced coalescence. Successful words
give a finite positive certificate; this does not assert a
terminating decision procedure for failed requirements.
Once the certificate holds, the uniform square-root upper
bound extends to $c\in d_r\Z\setminus\{0\}$, with a constant
depending on the fixed residue system and its certificate.
For example, $(\ell,m;r_1,r_2)=(3,5;1,-1)$ has $d_r=2$
and $C_0(r)=1$. Its certificate is empty, and exactly the
even shifts coalesce almost surely in the balanced sense.

\medskip
\noindent
\textbf{Rigidity from a single shift.}
Without assuming negative drift or admissibility, define
\[
 \mathcal G_r=\{c\in\Z:\mathbb P(\tau_c^{\mathrm{bal}}<\infty)=1\}.
\]
Theorem~\ref{thm:ST} proves that
\[
 \mathcal G_r=\{0\}
 \qquad\text{or}\qquad
 \mathcal G_r=d_r\Z.
\]
Thus one nonzero almost-surely coalescing shift forces the
entire obstruction lattice to coalesce almost surely. This
holds for composite as well as prime $\ell$. Under
admissibility and negative drift, the finite certificate,
the existence of one such shift, and
$\mathcal G_r=d_r\Z$ are therefore equivalent
(Corollary~\ref{cor:STprime}).

\subsection*{Arithmetic consequences and method}

Throughout this subsection, fix a canonical map
$T=T_{\ell,m}$ with $\mu_{\ell,m}<0$. A set $E\subseteq\N$ is
\emph{finite-defect invariant} if $E\triangle T^{-1}(E)$ is
finite. Write $F_E(x)=\#(E\cap[1,x])$,
$\pi_E(x)=\#\{p\le x:p\text{ is prime},\ p\in E\}$, and
$\overline d(E)=\limsup_{x\to\infty}F_E(x)/x$.
For an infinite sampling set $S$, write
\[
 \overline d_S(E)=\limsup_{x\to\infty}
 \frac{\#(E\cap S\cap[1,x])}{\#(S\cap[1,x])}.
\]
The sieve transfer argument gives
\[
 \overline d_{\mathbb P}(E)\ll_{\ell,m}\overline d(E)^{1/3},
 \qquad
 \overline d_{\mathcal B}(E)\ll_{\ell,m}\overline d(E)^{1/2},
 \qquad
 \overline d_{\mathcal C}(E)\ll_{\ell,m}\overline d(E)^{1/5}
\]
(Corollary~\ref{cor:Gtransferpowers}). Here $\mathcal B$
is the set of sums of two squares, and $\mathcal C$ consists
of primes $p$ for which $p+2$ is prime or is a product of
two primes, both at least $p^{1/10}$. This roughness
restriction is part of the statement. The constants are
independent of $E$ and its finite defect cutoff.

Put $T_{\min}(N)=\min_{k\ge0}T^k(N)$.
Apply these bounds to $E_B=\{N:T_{\min}(N)>B\}$,
which is invariant above $B$, and use the known fixed-barrier
estimates. For every $f(N)\to\infty$, it follows that
$T_{\min}(N)<f(N)$ for almost all primes, rough Chen primes,
and sums of two squares, in relative natural density in each
population. Bondurant-David's
prime almost-boundedness theorem uses the weights
$(\log p)/p$ \cite[Theorem~P and Section~1.3]{BondurantDavid2026}.
The comparison here uses normalization by $\pi(x)$ and
applies to arbitrary finite-defect invariant sets. In the
standard case it is an application of Leen's coalescence
estimate and classical upper-bound sieves; the generalized
statements use the coalescence estimate above.

Under the Riemann hypothesis, Theorem~\ref{thm:GPE}
also gives, for every fixed $K_0$,
\[
 \sup_{K(E)\le K_0}
 \left|\frac{\pi_E(x)}{\pi(x)}-\frac{F_E(x)}x\right|
 \longrightarrow0,
 \qquad K(E)=\max(E\triangle T^{-1}(E)),
\]
where $\max\varnothing=0$ and the supremum is over
finite-defect invariant sets for the fixed map $T$.
Consequently, whenever the natural
density exists, it equals the relative density among primes.

Unconditionally, renewal estimates also give equivalence of natural and
logarithmic densities and uniform comparisons for individual
moduli in a power range (Proposition~\ref{thm:GT} and
Theorem~\ref{thm:GLD}). The former asserts equivalence,
without requiring density existence for every generalized
invariant set. Related standard-map comparisons for
individual moduli appear in \cite[Proposition~1]{LyapunovWeaver2026}.

The coalescence argument follows the affine relation between
two trajectories. Its level coordinate records the difference
between the branch counts. Until absorption, its successive
level changes are those of a simple symmetric random walk. A
telescoping identity for the affine multipliers controls the
normalized displacement. Fractional-moment estimates,
return-time bounds, and finite absorption words yield the
uniform tail estimate. Renewal equations then connect
coalescence to harmonic density profiles and arithmetic
averages. The prime and sieve consequences combine this
comparison with classical analytic estimates.

\subsection*{Organization}

Sections~\ref{sec:notation}--\ref{sec:korec} collect notation,
the Mellin framework, and power-threshold estimates.
Sections~\ref{sec:renewal}--\ref{sec:level1} develop renewal
equations, limiting density profiles, and mixing criteria.
Sections~\ref{sec:sadic} and~\ref{sec:koal} relate the
$S$-adic system to coalescence and ergodicity.
Section~\ref{sec:natural} gives quantitative coalescence
and its density and sieve consequences for the standard map.
The principal general results, including
Theorems~\ref{thm:Gsharp}, \ref{thm:GR}, and~\ref{thm:ST},
are proved in Section~\ref{sec:general}.
Section~\ref{sec:problems} lists open problems, and
Appendix~\ref{app:numerics} records finite computations
and verifiable absorption certificates.

\section{\textbf{Notation and definitions}}
\label{sec:notation}

Throughout, $\log$ denotes the natural logarithm and $\log_2$ the
logarithm to base $2$. Standard notation is collected in
Table~\ref{tab:notation}; the definitions and classical results used
below are recalled after the table.

\begin{table}[!htbp]
  \centering
  \small
  \renewcommand{\arraystretch}{1.12}
  \begin{tabular}{@{}>{\raggedright\arraybackslash}p{0.27\textwidth}
                    >{\raggedright\arraybackslash}p{0.69\textwidth}@{}}
    \toprule
    Symbol & Meaning\\
    \midrule
    $\N$, $\N_0$, $\Z$, $\Q$, $\R$, $\C$
      & $\N=\{1,2,\ldots\}$, $\N_0=\N\cup\{0\}$; integers,
        rationals, reals and complex numbers\\
    $s=\sigma+it$
      & Complex variable, with $\sigma=\operatorname{Re}s$\\
    $e(x)$, $\omega$
      & $e(x)=e^{2\pi ix}$ and $\omega=e(1/3)$\\
    $O$, $o$, $\#A$, $\mathbf1_A$
      & Landau symbols; cardinality and indicator of $A$\\
    $f_*\mu$
      & Push-forward measure: $(f_*\mu)(A)=\mu(f^{-1}(A))$\\
    $F_S(x)$, $\psi_S(u)$
      & $F_S(x)=\#(S\cap[1,x])$ and $\psi_S(u)=e^{-u}F_S(e^u)$
        for $S\subseteq\N$\\
    $\underline d(S)$, $\overline d(S)$, $d(S)$
      & The $\liminf$ and $\limsup$ of $F_S(x)/x$ as $x\to\infty$,
        and their common value when they agree\\
    Logarithmic density
      & $\displaystyle\lim_{x\to\infty}\frac1{\log x}
        \sum_{\substack{N\in S\\N\le x}}\frac1N$, when it exists;
        natural density $d$ implies logarithmic density $d$
        \cite{Tenenbaum2015}\\
    $\Gamma(s)$, $\zeta(s)$, $\zeta(s,a)$
      & Gamma, Riemann zeta and Hurwitz zeta functions;
        $\zeta(s,a)=\sum_{m\ge0}(m+a)^{-s}$ for
        $a>0$ and $\operatorname{Re}s>1$ \cite[Ch.~12]{Apostol1976}\\
    $Z_S(s)$
      & $\sum_{N\in S}N^{-s}$ for $\operatorname{Re}s>1$;
        in particular, $Z_\N=\zeta$\\
    $\mathcal Mf(s)$
      & Mellin transform $\int_0^\infty t^{s-1}f(t)\,dt$
        in its strip of convergence \cite{Titchmarsh1937}\\
    $\Theta_S(t)$, Abel density
      & $\Theta_S(t)=\sum_{N\in S}e^{-Nt}$ for $t>0$;
        $S$ has Abel density $d$ if $t\Theta_S(t)\to d$
        as $t\to0^+$\\
    $v_3$, $|\cdot|_3$, $\Z_3$
      & $3$-adic valuation, absolute value and ring of integers;
        $|x|_3=3^{-v_3(x)}$ for $x\ne0$, and $|0|_3=0$
        \cite{Koblitz1984}\\
    Cylinder of level $j$
      & A set $a+3^j\Z_3$, $j\ge0$; its three subcylinders
        of level $j+1$ are its children\\
    $dz$, $\nu_C$
      & Haar probability measure on $\Z_3$ \cite{Halmos1950};
        $\nu_C=3^j\mathbf1_C\,dz$ for a cylinder $C$ of level $j$\\
    $E_jh$, level $j$
      & $E_jh(z)=3^j\int_{z+3^j\Z_3}h(y)\,dy$ for
        $h\in L^1(\Z_3)$; $h$ is of level $j$ if $E_jh=h$
        almost everywhere\\
    $\widehat f(\xi)$
      & $\sum_{x\bmod3^k}f(x)e(-\xi x/3^k)$ for
        $f:\Z/3^k\Z\to\C$, $k\ge1$;
        a frequency $\xi\bmod3^k$ is of top level if $3\nmid\xi$\\
    $\|\eta\|_{\TV}$
      & One half of the total variation of a finite signed measure
        $\eta$, with the convention specified below\\
    $\Bin(n,\tfrac12)$, $H(a)$
      & Binomial law; binary entropy
        $H(a)=-a\log_2a-(1-a)\log_2(1-a)$ for $0\le a\le1$,
        with $0\log_2 0=0$\\
    $\mu$-harmonic
      & A bounded function $H$ satisfying
        $H(u)=\int H(u+y)\,\mu(dy)$;
        for a Markov operator $P$, harmonicity means $Ph=h$\\
    \bottomrule
  \end{tabular}
  \caption{Standard notation.}
  \label{tab:notation}
\end{table}

\smallskip
\noindent\textbf{The Collatz map.}
We use the shortcut map $T$ of \eqref{eq:T}.
For residue calculations, we extend $T$ to $\N_0$ by setting $T(0)=0$.
For $N\in\N$, let $\Colmin(N)$ be the smallest element of its orbit
under the original Collatz map, which sends an even integer $n$ to
$n/2$ and an odd integer $n$ to $3n+1$.
The $T$-orbit omits the intermediate value $3u+1$ following each odd
value $u$. Since $3u+1>u$, these omissions do not change the minimum.
Consequently,
\[
  \Colmin(N)=\min_{j\ge0}T^j(N)
           =\min\{N,\Colmin(T(N))\}.
\]
For $N\in\N_0$, the \emph{parity vector} of length $k\ge1$ records
the parities of $N,T(N),\ldots,T^{k-1}(N)$.
We write $q_k(N)$ for the number of its odd entries.
We set $T^0=\mathrm{id}$ and $q_0(N)=0$; in particular, $q_k(0)=0$.
We shall use the following finite-level form of the parity-vector
correspondence.

\begin{lemma}[Terras parity-vector correspondence
  {\cite{Terras1976}}]\label{lem:terras}
For integers $k\ge1$, $0\le r<2^k$ and $m\ge0$, we have
$q_k(2^km+r)=q_k(r)$ and
\begin{equation}\label{eq:Terras}
  T^k(2^km+r)=3^{q_k(r)}m+T^k(r).
\end{equation}
The map assigning to $r\in\{0,\ldots,2^k-1\}$ its parity vector
of length $k$ is a bijection onto $\{0,1\}^k$.
\end{lemma}

\smallskip
\noindent\textbf{Analytic tools.}
For $a>0$ and $\operatorname{Re}s>1$, the Hurwitz zeta function has
the integral representation
\begin{equation}\label{eq:hurwitzint}
  \Gamma(s)\zeta(s,a)
  =\int_0^\infty t^{s-1}\frac{e^{-at}}{1-e^{-t}}\,dt.
\end{equation}
Thus $\Gamma(s)\zeta(s,a)$ is the Mellin transform of
$e^{-at}/(1-e^{-t})$.
The function $\zeta(s,a)$ extends meromorphically to the whole
$s$-plane, with a single simple pole at $s=1$ of residue $1$
\cite[Ch.~12]{Apostol1976}.
For $\lambda>0$,
\[
  \mathcal M[f(\lambda\,\cdot)](s)
  =\lambda^{-s}\mathcal Mf(s)
\]
whenever the defining integrals converge. This identity converts
dilations into multiplication and will be used on the logarithmic scale.

\begin{lemma}[Hardy--Littlewood--Karamata
  {\cite{Karamata1931,Hardy1949}}]\label{lem:HLK}
Let $(a_n)_{n\ge1}$ be nonnegative, and suppose that
$\sum_{n\ge1}a_ne^{-nt}<\infty$ for every $t>0$.
If, for some $d\in[0,\infty)$,
\[
  t\sum_{n\ge1}a_ne^{-nt}\longrightarrow d
  \qquad(t\to0^+),
\]
then
\[
  \frac1x\sum_{n\le x}a_n\longrightarrow d
  \qquad(x\to\infty).
\]
\end{lemma}

For a finite signed measure $\eta$ on a measurable space $X$, we use
the convention $\|\eta\|_{\TV}=\tfrac12|\eta|(X)$, where $|\eta|$
is its variation measure.
If $\mu$ and $\mu'$ are finite positive measures of equal total mass,
then $\|\mu-\mu'\|_{\TV}=\sup_A|\mu(A)-\mu'(A)|$, where the
supremum is over measurable sets.
For a bounded measurable $f:X\to\C$, write
$\operatorname{osc}f=\sup_{x,y\in X}|f(x)-f(y)|$.
For real-valued $f$, this is $\sup f-\inf f$.
The Jordan decomposition of $\mu-\mu'$ gives the standard estimate
\begin{equation}\label{eq:tvosc}
  \left|\int f\,d\mu-\int f\,d\mu'\right|
  \le \operatorname{osc}(f)\,\|\mu-\mu'\|_{\TV}
  \le 2\|f\|_\infty\|\mu-\mu'\|_{\TV}
\end{equation}
\cite{LPW2009}.
We also use the binomial Chernoff bound
\[
  \sum_{i=\lceil an\rceil}^{n}\binom ni\le2^{nH(a)},
  \qquad n\ge1,\quad \tfrac12\le a\le1.
\]

\begin{lemma}[Choquet--Deny
  {\cite{ChoquetDeny1960}}]\label{lem:CD}
Let $\mu$ be a Borel probability measure on $\R$, and let
$H:\R\to\C$ be bounded and continuous.
If
\[
  H(u)=\int_{\R}H(u+y)\,\mu(dy)
  \qquad(u\in\R),
\]
then $H(u+y)=H(u)$ for every $u\in\R$ and every $y$ in the closed
additive subgroup generated by $\operatorname{supp}\mu$.
\end{lemma}

Harmonic functions of the random Collatz chain are defined in
Section~\ref{sec:profiles}.

\smallskip
\noindent\textbf{$3$-adic affine maps and words.}
Since $2$ is a unit in $\Z_3$, the maps
\[
  g_0(z)=\frac z2,
  \qquad
  g_1(z)=\frac{3z+1}{2}
\]
are defined on $\Z_3$.
For $\alpha\in\Z_3\setminus\{0\}$ and $\beta\in\Z_3$, the affine
map $g(z)=\alpha z+\beta$ sends a cylinder of level $j$ onto a
cylinder of level $j+v_3(\alpha)$, whose Haar measure is
$|\alpha|_3$ times that of the original cylinder.
Its push-forward of Haar probability measure is
\[
  g_*(dz)=|\alpha|_3^{-1}
          \mathbf1_{\beta+\alpha\Z_3}(z)\,dz.
\]
In particular, $(g_0)_*(dz)=dz$, whereas
\[
  (g_1)_*(dz)=3\,\mathbf1_{2+3\Z_3}(z)\,dz.
\]
Thus the inverse image of a Haar-null set under either map is
Haar-null. Each such affine map also sends normalized Haar measure
on a cylinder $C$ to normalized Haar measure on $g(C)$.

By the martingale convergence theorem, $E_jh\to h$ almost everywhere
and in $L^1(\Z_3)$ for $h\in L^1(\Z_3)$.
For functions on $\R\times\Z_3$, the operator $E_j$ acts only in the
variable $z\in\Z_3$.
For a word $w=(w_1,\ldots,w_n)\in\{0,1\}^n$, put
$\gamma_w=g_{w_n}\circ\cdots\circ g_{w_1}$ and
$q_w=w_1+\cdots+w_n$, the number of occurrences of the branch $g_1$.
When $w$ is the parity vector of $N\in\N_0$, we have
$T^n(N)=\gamma_w(N)$.
We allow the empty word, with
$\gamma_\varnothing=\mathrm{id}$ and $q_\varnothing=0$.
Then
\begin{equation}\label{eq:affine}
  \gamma_w(z)=A_{n,q_w}\,z+B_w,
  \qquad A_{n,q}=\frac{3^q}{2^n},
\end{equation}
where $B_w=\gamma_w(0)\in\Z_3$ is the \emph{offset} of $w$, and
$B_\varnothing=0$.
Since $|A_{n,q}|_3=3^{-q}$, the map $\gamma_w$ sends a cylinder of
level $j$ onto a cylinder of level $j+q_w$ and satisfies
$(\gamma_w)_*\nu_C=\nu_{\gamma_w(C)}$.

\section{\textbf{The Mellin framework}}
\label{sec:mellin}

We collect the Mellin identities and averaging estimates used below.
The residue-class formulas follow from the parity-vector
correspondence in Lemma~\ref{lem:terras}, while the equivalence of
natural and Abel densities is a standard Tauberian consequence.

\begin{lemma}[Hurwitz representation]\label{thm:hurwitz}
For $k\ge1$, $1\le r<2^k$ and $\operatorname{Re}s>1$, we have
\[
\begin{aligned}
  \sum_{m\ge0}(2^km+r)^{-s}
    &=2^{-ks}\zeta\left(s,\frac r{2^k}\right),\\
  \sum_{m\ge0}\bigl(T^k(2^km+r)\bigr)^{-s}
    &=3^{-q_k(r)s}
      \zeta\left(s,\frac{T^k(r)}{3^{q_k(r)}}\right).
\end{aligned}
\]
For the zero residue class, the sums start at $m=1$, and
\[
  \sum_{m\ge1}(2^km)^{-s}=2^{-ks}\zeta(s),\qquad
  \sum_{m\ge1}\bigl(T^k(2^km)\bigr)^{-s}=\zeta(s).
\]
\end{lemma}

\begin{proof}
The first identity is the defining series for the Hurwitz zeta
function. For the second, put $q=q_k(r)$ and use
\eqref{eq:Terras} to write
\[
  T^k(2^km+r)
  =3^qm+T^k(r)
  =3^q\left(m+\frac{T^k(r)}{3^q}\right).
\]
The Hurwitz parameter is positive because $r\ge1$.
For $r=0$, we have $T^k(2^km)=m$ for every $m\ge1$, giving
the remaining identities.
\end{proof}

Writing $q=q_k(r)$, the ratio of the two scaling factors is
$(2^k/3^q)^s$. On the line $s=1+it$, this equals
\[
  \frac{2^k}{3^q}\left(\frac{2^k}{3^q}\right)^{it}.
\]
Thus normalization by its value at $s=1$ gives the unit-modulus
factor averaged in the next lemma.

\begin{lemma}[Mellin characteristic function]\label{thm:mellinchar}
For $k\ge1$ and $t\in\R$, put
\[
  \chi_k(t)=\frac1{2^k}\sum_{r=0}^{2^k-1}
       \left(\frac{3^{q_k(r)}}{2^k}\right)^{-it}.
\]
Then
\[
  \chi_k(t)=2^{ikt}\left(\frac{1+3^{-it}}2\right)^k,
  \qquad
  |\chi_k(t)|=\left|\cos\frac{t\log3}{2}\right|^k.
\]
\end{lemma}

\begin{proof}
By Lemma~\ref{lem:terras}, exactly $\binom{k}{q}$ residue classes
have $q_k(r)=q$. The binomial theorem therefore gives
\[
  \chi_k(t)
  =2^{ikt}2^{-k}\sum_{q=0}^k\binom{k}{q}3^{-itq}
  =2^{ikt}\left(\frac{1+3^{-it}}2\right)^k.
\]
The modulus follows from
$1+e^{-iy}=2e^{-iy/2}\cos(y/2)$.
\end{proof}

\begin{proposition}[fixed-depth averaging]\label{thm:fixedk}
For every $k\ge1$, real $x\ge1$ and $t\in\R$,
\[
  \left|\frac1x\sum_{N\le x}
       \left(\frac{T^k(N)}N\right)^{-it}-\chi_k(t)\right|
  \le
  \frac{2^{k+1}+\frac{|t|}{3}(2^k-1)(1+\log x)}{x}.
\]
In particular, for fixed $k$ and $t$, the average on the left
converges to $\chi_k(t)$ as $x\to\infty$.
\end{proposition}

\begin{proof}
For $n_i=T^i(N)$ and $q=q_k(N)$, multiplication of the stepwise
ratios gives
\[
  \frac{T^k(N)}N
  =\frac{3^q}{2^k}
     \prod_{\substack{0\le i<k\\n_i\ {\rm odd}}}
       \left(1+\frac1{3n_i}\right).
\]
Since $T(n)\ge n/2$, we have $n_i\ge N/2^i$. Consequently,
\[
  0\le
  \log\frac{T^k(N)/N}{3^q/2^k}
  \le\frac1{3N}\sum_{i=0}^{k-1}2^i
  =\frac{2^k-1}{3N}.
\]
Using $|e^{iy}-1|\le|y|$, we obtain
\[
  \left|
    \left(\frac{T^k(N)}N\right)^{-it}
    -\left(\frac{3^{q_k(N)}}{2^k}\right)^{-it}
  \right|
  \le\frac{|t|(2^k-1)}{3N}.
\]
Summation over $N\le x$ contributes at most
$|t|(2^k-1)(1+\log x)/(3x)$ to the normalized average.

The remaining summands form a $2^k$-periodic sequence of modulus
one whose mean over a full period is $\chi_k(t)$.
Write $x=L2^k+y$, where $L=\lfloor x/2^k\rfloor$ and
$0\le y<2^k$. Removing the $L$ complete periods leaves an error,
before division by $x$, of absolute value at most
$\lfloor y\rfloor+y|\chi_k(t)|<2^{k+1}$.
Combining the two estimates proves the claim.
\end{proof}

\begin{proposition}[Mellin--Abel--Karamata]\label{thm:abel}
Let $S\subseteq\N$. For $\operatorname{Re}s>1$,
\[
  \Gamma(s)Z_S(s)
  =\int_0^\infty t^{s-1}\Theta_S(t)\,dt.
\]
For $d\in[0,1]$, the set $S$ has natural density $d$ if and only if
$t\Theta_S(t)\to d$ as $t\to0^+$.
Moreover, for every $t>0$,
\begin{equation}\label{eq:gammakernel}
  t\Theta_S(t)
  =\int_{-\infty}^{\infty}
       \psi_S(v-\log t)e^{2v-e^v}\,dv,
\end{equation}
and the nonnegative kernel $e^{2v-e^v}$ has integral one.
\end{proposition}

\begin{proof}
Writing $\sigma=\operatorname{Re}s>1$, we have
\[
  \sum_{N\in S}\int_0^\infty
       |t^{s-1}e^{-Nt}|\,dt
  =\Gamma(\sigma)Z_S(\sigma)
  \le\Gamma(\sigma)\zeta(\sigma)<\infty.
\]
Fubini's theorem therefore permits termwise integration and gives
the Mellin identity.

For the kernel formula, integration by parts gives
\[
  \Theta_S(t)=t\int_0^\infty F_S(x)e^{-tx}\,dx.
\]
Substituting $x=e^v/t$ and using
$F_S(x)=x\psi_S(\log x)$ proves \eqref{eq:gammakernel}.
Also,
\[
  \int_{\R}e^{2v-e^v}\,dv
  =\int_0^\infty y e^{-y}\,dy
  =\Gamma(2)=1.
\]
Since $0\le\psi_S\le1$, dominated convergence shows that natural
density $d$ implies $t\Theta_S(t)\to d$.
The converse follows from Lemma~\ref{lem:HLK} with
$a_N=\mathbf1_S(N)$, including the case $d=0$.
\end{proof}

\begin{remark}
Proposition~\ref{thm:abel} characterizes natural density by the
Laplace asymptotic $t\Theta_S(t)\to d$ as $t\to0^+$.
The classical continuous-boundary form of the Wiener--Ikehara
theorem \cite{Ikehara1931} supplies another sufficient condition:
if $Z_S(s)-d/(s-1)$ extends continuously from
$\operatorname{Re}s>1$ to $\operatorname{Re}s\ge1$, then $S$ has
natural density $d$.
\end{remark}

\section{\textbf{The sets $S_\theta$: power saving and holomorphic continuation}}
\label{sec:korec}

For $\theta>0$, put
\[
  m_T(N)=\min_{i\ge1}T^i(N),\qquad
  E_\theta=\{N\in\N:m_T(N)>N^\theta\},\qquad
  S_\theta=\N\setminus E_\theta.
\]
The minimum is taken over positive iterates.
For $0<\theta\le1$, we also write
\[
  \alpha=\alpha(\theta)=\theta\frac{\log2}{\log3},
  \qquad
  \delta(\theta)=1-H(\alpha(\theta)),
\]
where $H$ is the binary entropy function from
Section~\ref{sec:notation}.
The density-one conclusion for $\theta>\log3/\log4$ is classical
\cite{Korec1994}. A power-saving exceptional-set bound for every fixed
$\theta>0$ already follows from
\cite[Definition~2.6 and Corollary~3.2]{Inselmann2024}.
We record an elementary estimate with an explicit entropy exponent
in the range $\theta>\log3/\log4$ for the analytic application below.

\begin{lemma}[affine remainder]\label{lem:affine}
For $k\ge1$ and $N\in\N$, let $q=q_k(N)$. Then
\[
  2^kT^k(N)=3^qN+c,\qquad
  0\le c\le2^k\bigl((3/2)^q-1\bigr).
\]
For every $k$ and $0\le q\le k$, the upper bound can be attained.
\end{lemma}

\begin{proof}
An even step is $n\mapsto n/2$, and an odd step is
$n\mapsto3n/2+1/2$.
After composition, the coefficient of $N$ is $3^q/2^k$.
Each odd step contributes $1/2$ multiplied by the multipliers of
all subsequent steps.
If there are $\ell$ subsequent odd steps, this contribution is at
most $\tfrac12(3/2)^\ell$, since every subsequent even step has
multiplier $1/2$.
As the odd steps are traversed, $\ell$ takes each value
$0,\ldots,q-1$ once. Hence the total additive term is at most
\[
  \frac12\sum_{\ell=0}^{q-1}(3/2)^\ell=(3/2)^q-1.
\]
For $q=0$, the sum is empty and the additive term is zero.
Equality holds when the first $k-q$ steps are even and the last
$q$ steps are odd. Such parity vectors are realized by positive
integers by Lemma~\ref{lem:terras}.
\end{proof}

\begin{proposition}[an explicit exceptional-set estimate]\label{thm:korec}
Let $\theta\in(\log3/\log4,1]$. For every integer $j\ge1$,
\[
  \#\bigl(E_\theta\cap(2^{j-1},2^j]\bigr)
  \le
  \sum_{\substack{0\le i\le j\\i>Q_j}}\binom ji,
  \qquad
  Q_j=\frac{\theta(j-1)\log2-\log2}{\log3}.
\]
Consequently, there is a constant $C_\theta<\infty$ such that
\[
  \#\bigl(E_\theta\cap[1,x]\bigr)
  \le C_\theta x^{1-\delta(\theta)}
  \qquad(x\ge1),
\]
where $\delta(\theta)=1-H(\alpha(\theta))>0$.
\end{proposition}

\begin{proof}
Let $N\in E_\theta\cap(2^{j-1},2^j]$ and put $q=q_j(N)$.
By definition of $E_\theta$,
\[
  T^j(N)>N^\theta>2^{\theta(j-1)}.
\]
On the other hand, Lemma~\ref{lem:affine} and $N\le2^j$ give
\[
  T^j(N)
  \le3^q+(3/2)^q-1
  <2\cdot3^q.
\]
Thus $q>Q_j$.
Reduction modulo $2^j$ is injective on $(2^{j-1},2^j]$, and
$q_j(N)=q_j(N\bmod2^j)$ by Lemma~\ref{lem:terras}.
Exactly $\binom ji$ residue classes have $q_j=i$, proving the
first assertion.

For the power saving, write
\[
  Q_j=\alpha j-c_\theta,\qquad
  c_\theta=(\theta+1)\frac{\log2}{\log3},
  \qquad
  \lambda=\frac{\alpha}{1-\alpha}>1.
\]
The last inequality holds because
$1/2<\alpha\le\log2/\log3<1$.
For every $j\ge1$, the binomial theorem gives
\[
\begin{aligned}
  \sum_{\substack{0\le i\le j\\i>Q_j}}\binom ji
  &\le\lambda^{-Q_j}\sum_{i=0}^j\binom ji\lambda^i\\
  &=\lambda^{c_\theta}
      \left(\frac{1+\lambda}{\lambda^\alpha}\right)^j
   =\lambda^{c_\theta}2^{jH(\alpha)}.
\end{aligned}
\]
Since $1\notin E_\theta$ and $H(\alpha)>0$, summing this geometric
bound over the dyadic intervals meeting $[1,x]$ gives
$O_\theta(x^{H(\alpha)})$.
Finally, $H(\alpha)<1$ because $\alpha>1/2$, so
$H(\alpha)=1-\delta(\theta)$ with $\delta(\theta)>0$.
\end{proof}

\begin{corollary}\label{cor:analytic}
For $\theta\in(\log3/\log4,1]$, the series $Z_{E_\theta}(s)$
converges absolutely and defines a holomorphic function on
$\operatorname{Re}s>1-\delta(\theta)$.
Consequently,
\[
  Z_{S_\theta}(s)-\frac1{s-1}
  =\zeta(s)-\frac1{s-1}-Z_{E_\theta}(s),
  \qquad \operatorname{Re}s>1,
\]
extends holomorphically to
$\operatorname{Re}s>1-\delta(\theta)$.
Moreover,
\[
  \#\bigl(S_\theta\cap[1,x]\bigr)
  =x+O_\theta\bigl(x^{1-\delta(\theta)}\bigr).
\]
\end{corollary}

\begin{proof}
Put $A_\theta(x)=\#(E_\theta\cap[1,x])$ and
$\rho=1-\delta(\theta)$.
Proposition~\ref{thm:korec} gives $A_\theta(x)=O_\theta(x^\rho)$.
Partial summation yields, for $\operatorname{Re}s>\rho$,
\[
  Z_{E_\theta}(s)
  =s\int_1^\infty A_\theta(x)x^{-s-1}\,dx.
\]
The integral and the defining series converge absolutely and
locally uniformly in this half-plane.
The claimed extension follows because
$\zeta(s)-1/(s-1)$ is entire, with its value at $s=1$ defined by
continuity.
Finally,
$\#(S_\theta\cap[1,x])=\lfloor x\rfloor-A_\theta(x)$, giving
the counting estimate.
\end{proof}

\begin{table}[ht]
  \centering
  \begin{tabular}{cccccc}
    \toprule
    $\theta$ & $1$ & $0.95$ & $0.9$ & $0.85$ & $0.8$\\
    \midrule
    $\alpha(\theta)$
      & $0.630930$ & $0.599383$ & $0.567837$
      & $0.536290$ & $0.504744$\\
    $\delta(\theta)$
      & $0.05004447$ & $0.02868978$ & $0.01331911$
      & $0.00380336$ & $0.00006493$\\
    \bottomrule
  \end{tabular}
  \caption{Values of $\alpha(\theta)$ and the saving exponent
    $\delta(\theta)=1-H(\alpha(\theta))$, rounded to the displayed
    precision.}
  \label{tab:delta}
\end{table}

\begin{remark}
At $\theta=1$, the exponent
\[
  \delta(1)
  =1-H\left(\frac{\log2}{\log3}\right)
  =0.05004447\ldots
\]
agrees with the classical entropy exponent in
\cite[Theorems~D and~F]{Lagarias1985}.
The condition $\alpha(\theta)>1/2$ is equivalent to
$\theta>\log3/\log4$, explaining the classical threshold in this
argument. As $\theta\downarrow\log3/\log4$, the saving exponent
tends to zero.
\end{remark}

\section{\textbf{Invariant sets and the renewal equation}}
\label{sec:renewal}

\begin{definition}
A set $G\subseteq\N$ is \emph{$T$-invariant} if
$\Delta_G:=G\,\triangle\,T^{-1}(G)$ is finite.
Thus invariance is understood up to finitely many exceptions.
We write $D_G=\#\Delta_G$ and, for $j\ge0$ and $a\bmod3^j$, put
\[
  F_{a,j}(x)=\#\{N\in G:N\le x,\ N\equiv a\pmod{3^j}\}.
\]
In particular, $F=F_{0,0}$, $G_2=G\cap(2+3\Z)$ and $F_2=F_{2,1}$.
These counting functions are taken to be zero for $x<1$.
\end{definition}

For $a\bmod3^j$, choose an integer $m_0$ such that
$2m_0+1\equiv a\pmod{3^j}$, and define
\[
  \sigma(a)=3m_0+2\pmod{3^{j+1}}.
\]
This class is independent of the choice of $m_0$; the level $j$ is
understood from the input class. Equivalently,
$\sigma(a)\equiv(3a+1)2^{-1}\pmod{3^{j+1}}$.
When $j=0$, the unique input class is sent to $2\bmod3$.

\begin{lemma}[renewal equation with finite defect]\label{lem:renewal}
For a $T$-invariant set $G$, every $j\ge0$, every $a\bmod3^j$ and
every real $x\ge1$, we have
\begin{equation}\label{eq:renewal}
  \left|F_{a,j}(x)
  -F_{a2^{-1},j}\left(\frac x2\right)
  -F_{\sigma(a),j+1}\left(3\left\lfloor\frac{x-1}{2}\right\rfloor+2\right)
  \right|\le D_G.
\end{equation}
\end{lemma}

\begin{proof}
An even integer $N=2m\equiv a\pmod{3^j}$ with $N\le x$ has
$T(N)=m\le x/2$ and $m\equiv a2^{-1}\pmod{3^j}$.
An odd integer $N=2m+1\equiv a\pmod{3^j}$ has $m\equiv m_0\pmod{3^j}$
and
\[
  T(N)=3m+2\equiv\sigma(a)\pmod{3^{j+1}},\qquad
  T(N)\le3\left\lfloor\frac{x-1}{2}\right\rfloor+2.
\]
Both correspondences are bijections onto the indicated ranges.
For the odd branch, the target class is contained in $2+3\Z$, so
every positive target is at least $2$ and has a positive odd preimage.
Consequently, the sum of the two counting functions subtracted in
\eqref{eq:renewal} counts the integers $N\le x$ in the class
$a\bmod3^j$ for which $T(N)\in G$.
The difference from $F_{a,j}(x)$ is therefore
\[
  \sum_{\substack{1\le N\le x\\N\equiv a\; (\mathrm{mod}\;3^j)}}
       \bigl(\mathbf1_G(N)-\mathbf1_G(T(N))\bigr).
\]
Its summands vanish outside $\Delta_G$, proving the bound.
\end{proof}

The odd endpoint in \eqref{eq:renewal} differs from $3x/2$ by less
than $3$. Since the target class has spacing $3^{j+1}\ge3$, replacing
that endpoint by $3x/2$ changes its count by at most one. Hence
\[
  F_{a,j}(x)
  =F_{a2^{-1},j}(x/2)+F_{\sigma(a),j+1}(3x/2)+O_G(1),
\]
where the error is at most $D_G+1$, uniformly in $a$, $j$ and $x\ge1$.
The following proposition records consequences of this decomposition.

\begin{proposition}\label{thm:equiv}
Let $G$ be $T$-invariant.
\begin{enumerate}
  \item[(a)] For $\operatorname{Re}s>1$,
  \[
    (1-2^{-s})Z_G(s)
    =(3/2)^s Z_{G_2}(s)+\mathcal H(s)+P_G(s),
  \]
  where
  \[
    P_G(s)=\sum_{N\in\Delta_G}
       \bigl(\mathbf1_G(N)-\mathbf1_G(T(N))\bigr)N^{-s}
  \]
  is a finite Dirichlet polynomial, and
  \[
    \mathcal H(s)=\sum_{n\in G_2}
      \left[\left(\frac{2n-1}{3}\right)^{-s}
            -\left(\frac{2n}{3}\right)^{-s}\right]
  \]
  is holomorphic on $\operatorname{Re}s>0$.
  For $\sigma=\operatorname{Re}s>0$, its summands satisfy
  \[
    \left|\left(\frac{2n-1}{3}\right)^{-s}
           -\left(\frac{2n}{3}\right)^{-s}\right|
    \le\frac{|s|}{2n-1}
          \left(\frac{2n-1}{3}\right)^{-\sigma}.
  \]

  \item[(b)] For $d\in[0,1]$, the set $G$ has natural density $d$
  if and only if $G_2$ has natural density $d/3$.
  In general,
  \[
    \underline d(G)\ge3\,\underline d(G_2),\qquad
    \overline d(G)\le3\,\overline d(G_2).
  \]

  \item[(c)] For each $d\in\R$, the function
  $Z_G(s)-d/(s-1)$ admits a continuous extension from
  $\operatorname{Re}s>1$ to $\operatorname{Re}s\ge1$ if and only if
  $Z_{G_2}(s)-(d/3)/(s-1)$ admits such an extension.
\end{enumerate}
\end{proposition}

\begin{proof}
(a) For $\operatorname{Re}s>1$, the even and odd preimages give
\[
  Z_{T^{-1}(G)}(s)
  =2^{-s}Z_G(s)
   +\sum_{n\in G_2}\left(\frac{2n-1}{3}\right)^{-s}.
\]
By definition, $P_G(s)=Z_G(s)-Z_{T^{-1}(G)}(s)$, which proves
the asserted identity.
To estimate the summands of $\mathcal H$, put $b_n=(2n-1)/3$.
Integration along the positive real axis gives
\[
  b_n^{-s}-(b_n+1/3)^{-s}
  =s\int_{b_n}^{b_n+1/3}t^{-s-1}\,dt.
\]
For $\operatorname{Re}s=\sigma>0$, its absolute value is at most
$|s|b_n^{-\sigma}/(2n-1)$.
On every compact subset of $\operatorname{Re}s>0$, this is bounded
by $C n^{-1-\varepsilon}$ for suitable $C<\infty$ and
$\varepsilon>0$. Thus the series for $\mathcal H$ converges
absolutely and locally uniformly there.

(b) Put $D(x)=F(x)-F(x/2)$. The case $j=0$ of the renewal equation
gives
\[
  D(x)=F_2(3x/2)+O_G(1).
\]
Since $F(y)=0$ for $y<1$, telescoping yields, for $x\ge1$,
\[
  F(x)
  =\sum_{i=0}^{\lfloor\log_2x\rfloor}
       F_2\left(\frac{3x}{2^{i+1}}\right)
   +O_G(1+\log x).
\]
Given $\varepsilon>0$, choose $Y\ge1$ such that, for all $y\ge Y$,
\[
  \bigl(\underline d(G_2)-\varepsilon\bigr)y
  \le F_2(y)
  \le\bigl(\overline d(G_2)+\varepsilon\bigr)y.
\]
In the telescoping sum, the terms with arguments below $Y$ have
total contribution $O(Y)$, while the sum of all arguments is
$3x+O(1)$. Divide by $x$, let $x\to\infty$, and then let
$\varepsilon\to0$ to obtain the two density inequalities.
If $G_2$ has density $d/3$, these inequalities give $d(G)=d$.
Conversely, if $F(x)=dx+o(x)$, then
\[
  F_2(y)=F(2y/3)-F(y/3)+O_G(1)
        =\frac d3y+o(y),
\]
which proves the equivalence.

(c) Set $A(s)=1-2^{-s}$ and $B(s)=(3/2)^s$.
By part (a), on $\operatorname{Re}s>1$,
\[
\begin{split}
  A(s)\left(Z_G(s)-\frac d{s-1}\right)
  ={}&B(s)\left(Z_{G_2}(s)-\frac{d/3}{s-1}\right)
       +\mathcal H(s)+P_G(s)\\
     &+\frac{d\bigl(B(s)/3-A(s)\bigr)}{s-1}.
\end{split}
\]
The last quotient is entire, with its value at $s=1$ defined by
continuity, because $B(1)/3=A(1)=1/2$.
Moreover, $A$ and $B$ do not vanish on $\operatorname{Re}s\ge1$,
and $\mathcal H$ is holomorphic on $\operatorname{Re}s>0$.
Solving this identity for either regularized Dirichlet series proves
the claimed equivalence of continuous extensions.
\end{proof}

\section{\textbf{The random Collatz chain and limiting profiles}}
\label{sec:profiles}

In this section, $L^\infty(\R\times\Z_3)$ and almost-everywhere
statements refer to the product measure $du\,dz$.
When taking pointwise cylinder averages, we choose bounded Borel
representatives of $L^\infty$ functions.

\begin{definition}
The \emph{random Collatz chain} is the Markov chain
$X_n=(U_n,Z_n)$ on $\R\times\Z_3$ which, from $(u,z)$, moves with
probability $1/2$ to each of
\[
  (u-\log2,g_0z),\qquad
  (u+\log(3/2),g_1z).
\]
Its Markov operator is
\[
  Ph(u,z)
  =\tfrac12 h\left(u-\log2,\frac z2\right)
   +\tfrac12 h\left(u+\log(3/2),\frac{3z+1}{2}\right).
\]
A function $h\in L^\infty(\R\times\Z_3)$ is \emph{harmonic} if
$Ph=h$ almost everywhere.
\end{definition}

The chain chooses between the two affine branches with equal
probability. The real coordinate adds the logarithm of the
multiplicative factor of the chosen branch, while the $3$-adic
coordinate retains its full affine action.
For an odd positive integer $N$, the exact logarithmic increment is
\[
  \log T(N)-\log N
  =\log(3/2)+\log\left(1+\frac1{3N}\right).
\]
The additional term disappears in the limiting renewal equation
derived below.
By the affine measure formulas in Section~\ref{sec:notation}, the
inverse image of a Haar-null set under either $g_0$ or $g_1$ is
Haar-null. Together with translation invariance of Lebesgue measure,
this shows that the branch maps are nonsingular for $du\,dz$.
Consequently, $P$ is well defined on $L^\infty$ and
$\|Ph\|_\infty\le\|h\|_\infty$.

\begin{lemma}\label{lem:Ej}
For a bounded Borel representative of
$h\in L^\infty(\R\times\Z_3)$, every $j\ge0$, every $u\in\R$ and
every $z\in a+3^j\Z_3$, we have
\[
  E_j(Ph)(u,z)
  =\tfrac12(E_jh)(u-\log2,a2^{-1})
   +\tfrac12(E_{j+1}h)(u+\log(3/2),\sigma(a)).
\]
\end{lemma}

\begin{proof}
The maps $g_0$ and $g_1$ send normalized Haar measure on
$a+3^j\Z_3$ to normalized Haar measure on
$a2^{-1}+3^j\Z_3$ and
$\sigma(a)+3^{j+1}\Z_3$, respectively.
Averaging the two terms defining $Ph$ over the original cylinder
gives the formula.
\end{proof}

For a $T$-invariant set $G$, define its \emph{density profiles} by
\[
  \psi_{a,j}(u)=3^j e^{-u}F_{a,j}(e^u),
  \qquad u\in\R,\quad j\ge0,\quad a\bmod3^j.
\]
In particular, $\psi_{0,0}=\psi_G$ in the notation of
Section~\ref{sec:notation}.

\begin{theorem}[limiting profiles are harmonic]\label{thm:profiles}
Let $G$ be $T$-invariant and let $u_k\to\infty$.
There are a subsequence, still denoted by $u_k$, and a Borel function
$h:\R\times\Z_3\to[0,1]$ satisfying $Ph=h$ almost everywhere such
that, for every fixed $j\ge0$ and $a\bmod3^j$,
\[
  \psi_{a,j}(u_k+v)
  \longrightarrow 3^j\int_{a+3^j\Z_3}h(v,z)\,dz
\]
locally uniformly in $v\in\R$.
\end{theorem}

\begin{proof}
Put $q=3^j$.
The counting bound $F_{a,j}(x)\le x/q+1$ gives
\[
  0\le\psi_{a,j}(u)\le1+qe^{-u}.
\]
Moreover, for $\delta\ge0$,
\[
  0\le F_{a,j}(e^{u+\delta})-F_{a,j}(e^u)
  \le\frac{e^u(e^\delta-1)}q+1.
\]
Writing
\[
  \psi_{a,j}(u+\delta)
  =e^{-\delta}\psi_{a,j}(u)
   +qe^{-(u+\delta)}
       \bigl(F_{a,j}(e^{u+\delta})-F_{a,j}(e^u)\bigr),
\]
we obtain
\[
  |\psi_{a,j}(u+\delta)-\psi_{a,j}(u)|
  \le1-e^{-\delta}+qe^{-u}
  \le\delta+qe^{-u}.
\]

Diagonal extraction over the countably many pairs $(a,j)$ and
rational values of $v$ gives a common subsequence along which
$\psi_{a,j}(u_k+v)$ converges at every such point.
The preceding estimate shows that each limit on the rationals extends
uniquely to a $1$-Lipschitz function $h_{a,j}$ on $\R$.
Approximating a compact interval by a finite set of rational points
and using the same estimate gives
\[
  \psi_{a,j}(u_k+\cdot)\longrightarrow h_{a,j}
\]
locally uniformly. The counting bound also gives $0\le h_{a,j}\le1$.

The partition of a cylinder into its three children yields
\[
  h_{a,j}(v)
  =\frac13\sum_{\substack{b\bmod3^{j+1}\\b\equiv a\; (\mathrm{mod}\;3^j)}}
      h_{b,j+1}(v).
\]
Set $H_j(v,z)=h_{z\bmod3^j,j}(v)$.
For every fixed $v$, the functions $H_j(v,\cdot)$ form a martingale
for the increasing cylinder filtration on $\Z_3$, with values
in $[0,1]$.
Define
\[
  h(v,z)=\limsup_{j\to\infty}H_j(v,z).
\]
This function is jointly Borel and takes values in $[0,1]$.
Bounded martingale convergence gives, for every fixed $v$, convergence
to $h(v,\cdot)$ almost everywhere and in $L^1(\Z_3)$.
Consequently, for every $v$, $j$ and $a$,
\[
  3^j\int_{a+3^j\Z_3}h(v,z)\,dz=h_{a,j}(v).
\]
This proves the asserted profile convergence.

It remains to establish harmonicity.
The endpoint estimate following Lemma~\ref{lem:renewal}, divided by
$x/3^j$ with $x=e^u$, gives
\[
  \psi_{a,j}(u)
  =\tfrac12\psi_{a2^{-1},j}(u-\log2)
   +\tfrac12\psi_{\sigma(a),j+1}(u+\log(3/2))
   +R_{a,j}(u),
\]
where
\[
  |R_{a,j}(u)|\le(D_G+1)3^j e^{-u}\qquad(u\ge0).
\]
For fixed $j$, this error tends to zero locally uniformly after
replacing $u$ by $u_k+v$.
Passing to the limits therefore gives
\[
  h_{a,j}(v)
  =\tfrac12h_{a2^{-1},j}(v-\log2)
   +\tfrac12h_{\sigma(a),j+1}(v+\log(3/2)).
\]
By Lemma~\ref{lem:Ej} and the cylinder-average identities above,
\[
  E_j(Ph-h)(v,\cdot)=0
\]
for every $v$ and $j$.
For each fixed $v$, martingale convergence in the $z$ variable now
implies $Ph(v,z)=h(v,z)$ for Haar-almost every $z$.
Since $h$ and $Ph$ are jointly Borel, Fubini's theorem on bounded
$v$-intervals gives $Ph=h$ almost everywhere for the product measure.
\end{proof}

\section{\textbf{The Choquet--Deny property at finite levels}}
\label{sec:finite}

As in Section~\ref{sec:notation}, a function
$h\in L^\infty(\R\times\Z_3)$ is \emph{of level $j$} if
$E_jh=h$ almost everywhere.
Such a function has a bounded Borel representative of the form
\[
  h(u,z)=h_a(u),\qquad a=z\bmod3^j,
\]
with $|h(u,z)|\le\|h\|_\infty$ at every point.
We use representatives with these properties whenever pointwise
values of finite-level functions occur below.

\begin{lemma}\label{lem:level}
For $j\ge1$, the residue class of $g_1(z)$ modulo $3^j$ depends
only on $z\bmod3^{j-1}$.
Composition with $g_0$ preserves level $j$, and composition with
$g_1$ sends level $j$ to level $j-1$.
In particular, $P$ preserves every finite level, including level $0$.
\end{lemma}

\begin{proof}
The map $g_0(z)=z/2$ preserves congruences modulo $3^j$, while
\[
  g_1(z+3^{j-1}t)=g_1(z)+3^j2^{-1}t.
\]
The shifts in the real coordinate do not affect these statements.
Level $0$ consists of functions independent of $z$ and is also
preserved by $P$.
\end{proof}

\begin{lemma}[Choquet--Deny on $\R$]\label{lem:CDR}
If $H\in L^\infty(\R)$ satisfies
\[
  H(u)=\tfrac12H(u-\log2)
       +\tfrac12H(u+\log(3/2))
\]
almost everywhere, then $H$ is almost everywhere constant.
\end{lemma}

\begin{proof}
For $\varphi\in C_c^\infty(\R)$, the convolution $H*\varphi$ is
bounded and continuous.
Convolution commutes with the two translations, so $H*\varphi$
satisfies the same equation everywhere.
By Lemma~\ref{lem:CD}, it is invariant under the closed additive
subgroup generated by the support of
\[
  \tfrac12\delta_{-\log2}+\tfrac12\delta_{\log(3/2)}.
\]
This subgroup is the closure of
$\Z\log2+\Z\log3$, which is $\R$ because
$\log3/\log2$ is irrational.
Thus $H*\varphi$ is constant.
Taking a smooth approximate identity
$(\varphi_\varepsilon)_{\varepsilon>0}$, we have
$H*\varphi_\varepsilon\to H$ in $L^1_{\mathrm{loc}}(\R)$.
The constants converge, as is seen by integrating over $[0,1]$,
so $H$ is almost everywhere constant.
\end{proof}

\begin{proposition}[finite-level contraction and harmonic functions]
\label{thm:finitelevel}
Let $j\ge0$ and let $h\in L^\infty(\R\times\Z_3)$ be of level $j$,
with a representative chosen as above.
\begin{enumerate}
  \item[(a)] For every $n\ge0$, $u\in\R$ and $z,z'\in\Z_3$,
  \[
    |P^nh(u,z)-P^nh(u,z')|
    \le2\|h\|_\infty\,
       \Pb\!\left(\Bin(n,\tfrac12)<j\right).
  \]
  No continuity in $u$ is required.

  \item[(b)] If $Ph=h$ almost everywhere, then $h$ is almost
  everywhere constant.
\end{enumerate}
\end{proposition}

\begin{proof}
(a) Run the chains from $(u,z)$ and $(u,z')$ with the same sequence
of coins.
Their real coordinates agree at every time.
If $q_n$ is the number of occurrences of the branch $g_1$ among
the first $n$ steps, then \eqref{eq:affine} gives
\[
  Z_n-Z_n'=\frac{3^{q_n}}{2^n}(z-z'),
  \qquad |Z_n-Z_n'|_3\le3^{-q_n}.
\]
On the event $q_n\ge j$, the two terminal states have the same real
coordinate and the same residue class modulo $3^j$, so the values
of $h$ agree.
On the complementary event their difference is at most
$2\|h\|_\infty$.
Since $q_n$ has distribution $\Bin(n,\tfrac12)$, taking expectations
proves the estimate.
For $j=0$, the probability on the right is zero.

(b) The operator $P$ is well defined on $L^\infty$, so $Ph=h$
implies $P^nh=h$ almost everywhere for every $n\ge0$.
The chosen representative and its iterates $P^nh$ are constant
on each level-$j$ cylinder.
Since there are only finitely many level-$j$ cylinders, Fubini's
theorem gives, for almost every $u$,
\[
  P^nh(u,z)=h(u,z)
  \qquad\text{for every }z\in\Z_3\text{ and every }n\ge0.
\]
Here the exceptional set in $u$ is chosen simultaneously for the
countably many values of $n$.
For such $u$, part (a) and
$\Pb(\Bin(n,\tfrac12)<j)\to0$ show that $h(u,z)$ is independent of
$z$.
Thus $h(u,z)=H(u)$ almost everywhere for some bounded Borel
function $H$.
The harmonic equation reduces to that of Lemma~\ref{lem:CDR},
which proves that $H$, and hence $h$, is almost everywhere constant.
\end{proof}

\begin{proposition}[finite-level spectrum]\label{prop:gap}
Let
\[
  P_0f(z)=\tfrac12f(z/2)+\tfrac12f(g_1z)
\]
act on the complex-valued functions of level $j\ge1$.
Then $1$ is an algebraically simple eigenvalue, and every other
eigenvalue has modulus exactly $1/2$.
In particular, the absolute spectral gap is $1/2$, and
\[
  P_0\mathbf1_{3\Z_3}=\tfrac12\mathbf1_{3\Z_3}.
\]
\end{proposition}

\begin{proof}
Let $V_r$ be the space of complex-valued functions on $\Z_3$ of
level $r$, and put $U_if=f\circ g_i$.
By Lemma~\ref{lem:level}, for $r\ge1$,
\[
  U_0V_r=V_r,\qquad U_1V_r\subseteq V_{r-1}.
\]
Thus the operator induced by $P_0$ on $V_r/V_{r-1}$ is
$\tfrac12\overline U_0$, where $\overline U_0$ is the induced
operator on this quotient.
The operator $U_0$ has finite order on $V_r$, since multiplication
by $2^{-1}$ permutes the residue classes modulo $3^r$.
Consequently, every eigenvalue of $\overline U_0$ is a root of
unity.
Using the invariant filtration
\[
  V_0\subset V_1\subset\cdots\subset V_j,
\]
we conclude that the eigenvalue $1$ occurs only on the
one-dimensional space $V_0$, while all the remaining eigenvalues
have modulus $1/2$.
This also proves the algebraic simplicity of $1$.
Finally,
$\mathbf1_{3\Z_3}(z/2)=\mathbf1_{3\Z_3}(z)$ and
$g_1(z)\equiv2\pmod3$, giving the stated eigenfunction.
\end{proof}

\section{\textbf{The full system: obstacle and criterion}}
\label{sec:full}

\begin{remark}[passing beyond finite levels]
\label{rem:obstacle}
The operator dual to $P$ with respect to $du\,dz$ is
\[
\begin{split}
  Lg(u,w)
  ={}&\tfrac12g(u+\log2,2w)\\
     &+\tfrac32\mathbf1_{2+3\Z_3}(w)
       g\left(u-\log(3/2),\frac{2w-1}{3}\right).
\end{split}
\]
The second term is understood to be zero outside $2+3\Z_3$.
The coefficient $3/2$ follows from the Haar measure change under
$g_1$.
The operator $L$ is positive and preserves integrals, so
$\|L\|_{L^1\to L^1}=1$.
It maps functions of level $j$ to functions of level at most
$j+1$, and need not preserve level $j$.

The conditional expectations $E_j$ do not, in general, commute
with $P$.
For example, if $f(z)=\mathbf1_{2+3\Z_3}(z)$, then
\[
  E_0(Pf)=\frac23,\qquad P(E_0f)=\frac13.
\]
For a harmonic function, the projected equation
$E_jh=E_jPh$ couples levels $j$ and $j+1$, as shown in
Lemma~\ref{lem:Ej}.
It does not directly supply the equation $P(E_jh)=E_jh$ needed
to apply Proposition~\ref{thm:finitelevel}.
This issue is related to the distinction emphasized by Tao
\cite[Remark~1.16]{Tao2022} between mixing of the affine offset
and mixing of the entire random affine map.
\end{remark}

For a cylinder $C$ and integers $n\ge0$, $0\le q\le n$, define
the finite positive measure
\[
  \mu_C^{n,q}
  =2^{-n}\sum_{\substack{|w|=n\\q_w=q}}
       (\gamma_w)_*\nu_C.
\]
If $C$ has level $j$, each measure $(\gamma_w)_*\nu_C$ is uniform
on the cylinder $\gamma_w(C)$ of level $j+q_w$.
Moreover,
\[
  \mu_C^{n,q}(\Z_3)=2^{-n}\binom nq,
\]
independently of $C$.
Put
\[
  D_n^{(j)}
  =\max_{\substack{C,C'\ {\rm cylinders}\\
                   \text{of level }j}}
       \sum_{q=0}^n
       \|\mu_C^{n,q}-\mu_{C'}^{n,q}\|_{\TV}.
\]
With the total variation convention of
Section~\ref{sec:notation}, we have $0\le D_n^{(j)}\le1$.

\begin{definition}[mixing condition (M)]\label{conj:M}
The random Collatz chain satisfies \emph{condition (M)} if
\[
  D_n^{(j)}\longrightarrow0\qquad(n\to\infty)
\]
for every fixed $j\ge1$.
\end{definition}

Condition (M) is verified by the coalescence argument in
Proposition~\ref{thm:M}, in Section~\ref{sec:natural}.
We first explain its implication for bounded harmonic functions.

\begin{proposition}[a mixing criterion for the Liouville property]
\label{thm:criterion}
If condition (M) holds, then every harmonic function in
$L^\infty(\R\times\Z_3)$ is almost everywhere constant.
\end{proposition}

\begin{proof}
Choose a bounded Borel representative of $h$ satisfying
$|h(u,z)|\le\|h\|_\infty$ everywhere.
Since $P$ is well defined on $L^\infty$, harmonicity gives
$P^nh=h$ almost everywhere for every $n\ge0$.
By Fubini's theorem, there is a set $U_*\subseteq\R$ of full
Lebesgue measure such that, for all $u\in U_*$ and all $n\ge0$,
\[
  h(u,\cdot)=P^nh(u,\cdot)
  \quad\text{Haar-almost everywhere}.
\]
For the chosen representative, the iterated transition formula is
\[
  P^nh(u,z)
  =2^{-n}\sum_{|w|=n}
       h(u+q_w\log3-n\log2,\gamma_w(z)).
\]
The real shift depends only on $q_w$.
Thus, for $u\in U_*$ and cylinders $C,C'$ of level $j$,
\[
\begin{split}
  &\int_{\Z_3}h(u,z)\,(\nu_C-\nu_{C'})(dz)\\
  &\qquad=
  \sum_{q=0}^n\int_{\Z_3}
       h(u+q\log3-n\log2,y)
       \,(\mu_C^{n,q}-\mu_{C'}^{n,q})(dy).
\end{split}
\]
Here the integral identity follows from the absolute continuity
of $\nu_C$ and $\nu_{C'}$ with respect to Haar measure.
The total variation bound in Section~\ref{sec:notation} gives
\[
  \left|
    \int_{\Z_3}h(u,z)\,(\nu_C-\nu_{C'})(dz)
  \right|
  \le2\|h\|_\infty D_n^{(j)}.
\]
Letting $n\to\infty$, condition (M) shows that all cylinder
averages at level $j$ agree.
Their common value is
\[
  H(u)=\int_{\Z_3}h(u,z)\,dz.
\]
Therefore $E_jh(u,\cdot)=H(u)$ for every $j\ge1$ and
$u\in U_*$.
Martingale convergence in the $z$ variable yields
$h(u,z)=H(u)$ almost everywhere.
Substitution in $Ph=h$, using that $P$ respects almost-everywhere
equality, gives the equation of Lemma~\ref{lem:CDR}.
Hence $H$ is almost everywhere constant.
\end{proof}

The next corollary uses only the Liouville property.
Its hypothesis follows from condition (M) by
Proposition~\ref{thm:criterion} and is also established directly
by coalescence in Corollary~\ref{thm:liouville}.

\begin{corollary}\label{cor:main}
Assume that every bounded harmonic function of the random Collatz
chain is almost everywhere constant.
Let $G$ be a $T$-invariant set with counting function $F$.
Then the following statements hold.
\begin{enumerate}
  \item[(a)] For every $\Lambda>1$,
  \[
    \sup_{\lambda\in[\Lambda^{-1},\Lambda]}
    \left|
      \frac{F(\lambda x)}{\lambda x}-\frac{F(x)}x
    \right|\longrightarrow0
    \qquad(x\to\infty).
  \]

  \item[(b)] For every fixed $j\ge0$,
  \[
    \max_{a\bmod3^j}
    \left|
      \frac{3^jF_{a,j}(x)}x-\frac{F(x)}x
    \right|\longrightarrow0
    \qquad(x\to\infty).
  \]

  \item[(c)] We have
  \[
    t\Theta_G(t)-tF(1/t)\longrightarrow0
    \qquad(t\to0^+).
  \]
\end{enumerate}
\end{corollary}

\begin{proof}
Let $u_k\to\infty$.
By Theorem~\ref{thm:profiles}, a subsequence of the translated
profiles converges locally uniformly to the cylinder averages of
a bounded harmonic function $h$.
By hypothesis, $h$ is almost everywhere equal to a constant $c$.
The limiting cylinder averages are continuous in the real
variable, so they equal $c$ at every real argument.
Consequently, along this subsequence,
\[
  \psi_{a,j}(u_k+v)\longrightarrow c
\]
locally uniformly in $v$, for every fixed $j$ and $a$.
The constant $c$ may depend on the chosen subsequence.

If (a) failed, there would be $\varepsilon>0$, $x_k\to\infty$
and $\lambda_k\in[\Lambda^{-1},\Lambda]$ such that
\[
  |\psi_G(\log x_k+\log\lambda_k)-\psi_G(\log x_k)|
  \ge\varepsilon.
\]
Applying the preceding subsequence argument with $u_k=\log x_k$
gives uniform convergence to $c$ for
$v\in[-\log\Lambda,\log\Lambda]$, a contradiction.
For (b), the same argument at $v=0$ gives
$\psi_{a,j}(u_k)-\psi_{0,0}(u_k)\to0$ for every $a$.
There are finitely many classes at a fixed level $j$, so a failure
of the asserted maximum would again give a contradiction.

For (c), put $u=-\log t$ and $K(v)=e^{2v-e^v}$.
By \eqref{eq:gammakernel},
\[
  t\Theta_G(t)-tF(1/t)
  =\int_{\R}K(v)
       \bigl(\psi_G(u+v)-\psi_G(u)\bigr)\,dv.
\]
For each fixed $v$, the difference in parentheses tends to zero
by (a).
Since $0\le\psi_G\le1$ and $K$ is nonnegative with integral one,
dominated convergence proves (c).
\end{proof}

\begin{remark}
For a fixed $C\ge1$, the exceptional set
\[
  E_C=\{N\in\N:\Colmin(N)>C\}
\]
is $T$-invariant in the finite-defect sense, since
\[
  E_C=\{N>C\}\cap T^{-1}(E_C),
  \qquad
  E_C\triangle T^{-1}(E_C)\subseteq[1,C]\cap\N.
\]
Thus Corollary~\ref{cor:main}(b) applies to this set.
The conclusion is equidistribution on the scale $x$:
\[
  F_{a,j}(x)=3^{-j}F(x)+o(x)
\]
for fixed $j$ and $a$.
If $F(x)=o(x)$, this estimate does not imply
$F_{a,j}(x)/F(x)\to3^{-j}$.

Corollary~\ref{cor:main} also does not by itself give natural
density: the constant obtained from a limiting profile can depend
on the subsequence.
The passage to natural density requires the additional argument
in Section~\ref{sec:natural}.
\end{remark}

\begin{remark}[the $3x-1$ map]
Let $T_-(N)=N/2$ for even $N$ and $T_-(N)=(3N-1)/2$ for odd $N$.
Replace $g_1$ by $g_1^-(z)=(3z-1)/2$, and denote the corresponding
Markov operator by $P_-$.
For $\mathcal Rh(u,z)=h(u,-z)$, we have
\[
  P_-=\mathcal R P\mathcal R.
\]
Reflection preserves Haar measure and cylinder levels.
It follows that the finite-level results, condition (M), and the
Liouville property transfer to $P_-$.
In particular, the quantities $D_n^{(j)}$ are unchanged.

For a set $G\subseteq\N$ with
$G\triangle T_-^{-1}(G)$ finite, the renewal and limiting-profile
arguments apply with the corresponding changes in the counting
formulas.
The odd branch is now $T_-(2m+1)=3m+1$; its image lies in
$1+3\Z$, its inverse is $(2n+1)/3$, and
\[
  \sigma_-(a)=3m_0+1\pmod{3^{j+1}},
  \qquad 2m_0+1\equiv a\pmod{3^j}.
\]
The odd endpoint in the renewal equation becomes
$3\lfloor(x-1)/2\rfloor+1$.
The Dirichlet correction term is accordingly
\[
  \mathcal H_-(s)
  =\sum_{n\in G\cap(1+3\Z)}
    \left[
      \left(\frac{2n+1}{3}\right)^{-s}
      -\left(\frac{2n}{3}\right)^{-s}
    \right].
\]
For $\sigma=\operatorname{Re}s>0$, each summand has absolute value
at most
\[
  \frac{|s|}{2n}\left(\frac{2n}{3}\right)^{-\sigma},
\]
so the same local convergence argument proves holomorphy on
$\operatorname{Re}s>0$.
The cycles with minima $1$, $5$ and $17$ have periods $1$, $3$
and $11$, respectively; they correspond to the negative cycles
of the $3x+1$ map described in
\cite[Section~2.6]{Lagarias1985}.
Their basins provide nontrivial examples of exactly
$T_-$-invariant sets.
\end{remark}

\section{\textbf{First-level mixing and a second-moment criterion}}
\label{sec:level1}

This section gives a combinatorial expression for the mixing distances,
quantitative bounds for cylinder averages of harmonic functions, and
a second-moment sufficient condition for mixing at level $1$.
The finite computations below provide explicit benchmarks; the later
coalescence argument proves the vanishing of the limiting distances.

\subsection{Combinatorial form}

Let $w=(w_1,\ldots,w_n)$ be a word whose ones occur at positions
$m_0>m_1>\cdots>m_{q-1}$, with positions numbered in chronological
order from $1$ to $n$.
The affine formula \eqref{eq:affine} gives
\begin{equation}\label{eq:XS}
  2^{n+1}B_w=X_w=X_S
  :=\sum_{\ell=0}^{q-1}3^\ell2^{m_\ell},
  \qquad S=\{m_0,\ldots,m_{q-1}\}.
\end{equation}
For $q=0$ the sum is zero.
Appending a letter gives the useful recursions
\[
  X_{w0}=X_w,\qquad X_{w1}=3X_w+2^{n+1}.
\]

For a uniformly chosen $q$-element subset
$S\subseteq\{1,\ldots,n\}$, let $t(n,q)$ be the total variation
distance between the law of $X_S\bmod3^{q+1}$ and its translate
by $3^q$.
If $C=a+3\Z_3$, multiplication of the output by $2^{n+1}$ gives
\[
  2^{n+1}\gamma_w(C)
  =X_w+2\cdot3^qa+3^{q+1}\Z_3.
\]
Thus the distance associated with distinct level-$1$ cylinders
is independent of the pair: the two nonzero shifts on the top
digit are negatives of one another and have the same total
variation distance.
Consequently,
\begin{equation}\label{eq:Dt}
  D_n^{(1)}
  =\sum_{q=0}^n2^{-n}\binom nq\,t(n,q).
\end{equation}

The interpretation as conditional uniformity requires averaging.
To make this precise, write $\rho_{n,q}$ for the law just defined,
and put
\[
  p_{n,q}(c)=\sum_{k=0}^2\rho_{n,q}(c+3^qk),\qquad
  \rho_{n,q}^{\flat}(c+3^qk)=\frac{p_{n,q}(c)}3,
  \quad 0\le c<3^q.
\]
If $e(n,q)=\|\rho_{n,q}-\rho_{n,q}^{\flat}\|_{\TV}$, then
\[
  \frac12t(n,q)\le e(n,q)\le\frac23t(n,q).
\]
Indeed, $\rho_{n,q}^{\flat}$ is the average of the three translates
by multiples of $3^q$; the inequalities follow from the triangle
inequality and translation invariance of total variation.
Hence $D_n^{(1)}\to0$ is equivalent to uniformity of the top digit
in total variation averaged over the parent residue and over
$q\sim\Bin(n,\tfrac12)$.
It does not assert convergence for every $q$; for example,
$t(n,0)=t(n,n)=1$.

\subsection{Monotonicity}

\begin{lemma}\label{lem:monotone}
For every $j\ge1$, $n\ge0$ and pair of level-$j$ cylinders $C,C'$,
\[
  \sum_q\|\mu_C^{n+1,q}-\mu_{C'}^{n+1,q}\|_{\TV}
  \le
  \sum_q\|\mu_C^{n,q}-\mu_{C'}^{n,q}\|_{\TV}.
\]
In particular, $D_{n+1}^{(j)}\le D_n^{(j)}$, and the limit
$D_\infty^{(j)}=\lim_{n\to\infty}D_n^{(j)}$ exists.
\end{lemma}

\begin{proof}
Set $\mu_C^{n,q}=0$ when $q\notin\{0,\ldots,n\}$.
Decomposition according to the last letter gives
\[
  \mu_C^{n+1,q}
  =\tfrac12(g_0)_*\mu_C^{n,q}
   +\tfrac12(g_1)_*\mu_C^{n,q-1}.
\]
Pushforward is contractive in total variation.
Subtract the corresponding identity for $C'$, apply the triangle
inequality, and sum over $q$.
Taking the maximum over the pairs $C,C'$ proves the last assertion.
\end{proof}

\subsection{Oscillation bounds and finite computations}

For a bounded real-valued function $h$ on $\R\times\Z_3$, define
\[
  \Omega(h)
  =\operatorname*{ess\,sup}_{u\in\R}
    \left(
      \operatorname*{ess\,sup}_{z\in\Z_3}h(u,z)
      -\operatorname*{ess\,inf}_{z\in\Z_3}h(u,z)
    \right).
\]
The following estimate is stated at every level, since this will
also be useful in Section~\ref{sec:sadic}.

\begin{proposition}[oscillation of cylinder averages]\label{thm:quant}
Let $h$ be a bounded real-valued harmonic function.
For every $j\ge1$ and $a,a'\bmod3^j$, and almost every $u$,
\[
  |E_jh(u,a)-E_jh(u,a')|
  \le D_\infty^{(j)}\Omega(h)
  \le D_n^{(j)}\Omega(h)
  \qquad(n\ge0).
\]
The exceptional set in $u$ can be chosen independently of
$j,a,a'$ and $n$.
At level $1$, the exact finite computation gives
\[
  D_\infty^{(1)}\le D_{24}^{(1)}
  =\frac{1\,792\,853}{8\,388\,608}<0.2138.
\]
\end{proposition}

\begin{proof}
Choose a bounded Borel representative of $h$.
As in Proposition~\ref{thm:criterion}, harmonicity and Fubini's
theorem give, for almost every $u$ and every $n$,
\[
\begin{split}
  E_jh(u,a)-E_jh(u,a')
  =\sum_{q=0}^n\int_{\Z_3}
    &h(u+q\log3-n\log2,y)\\
    &{}\cdot(\mu_C^{n,q}-\mu_{C'}^{n,q})(dy),
\end{split}
\]
where $C=a+3^j\Z_3$ and $C'=a'+3^j\Z_3$.
Intersect this full-measure set with all translates, by
$\Z\log2+\Z\log3$, of the set on which the slice oscillation of
$h$ is at most $\Omega(h)$.
The resulting set still has full measure, since the group of
translations is countable.

For each $q$, the two measures in the integrand have equal mass
and are absolutely continuous with respect to Haar measure.
The oscillation estimate \eqref{eq:tvosc} therefore bounds the
right-hand side by $\Omega(h)D_n^{(j)}$.
Lemma~\ref{lem:monotone} allows us to take the infimum over $n$.
There are only countably many remaining indices, so the same
exceptional set works for all of them.

Finally, enumeration using \eqref{eq:XS} and \eqref{eq:Dt} gives
$2^{24}D_{24}^{(1)}=3\,585\,706$.
The integer-counting procedure is recorded in
Appendix~\ref{app:D}.
\end{proof}

\begin{corollary}[bounds for invariant sets]\label{cor:uncond}
Let $G$ be $T$-invariant in the finite-defect sense of
Section~\ref{sec:renewal}.
For every fixed $j\ge1$,
\begin{align*}
  \limsup_{x\to\infty}\max_{a,a'\bmod3^j}
  \frac{3^j|F_{a,j}(x)-F_{a',j}(x)|}{x}
  &\le D_\infty^{(j)},\\
  \limsup_{x\to\infty}\max_{a\bmod3^j}
  \left|\frac{3^jF_{a,j}(x)}x-\frac{F(x)}x\right|
  &\le(1-3^{-j})D_\infty^{(j)}.
\end{align*}
In particular, at level $1$ the right-hand sides are bounded by
\[
  \frac{1\,792\,853}{8\,388\,608}<0.2138,
  \qquad
  \frac{1\,792\,853}{12\,582\,912}<0.1425,
\]
respectively.
\end{corollary}

\begin{proof}
Take a sequence $x_k\to\infty$ realizing the first upper limit.
By Theorem~\ref{thm:profiles}, a subsequence has a limiting
harmonic profile $h$ with $0\le h\le1$, and hence $\Omega(h)\le1$.
The level-$j$ cylinder averages of this profile are continuous
in the real variable.
Thus Proposition~\ref{thm:quant}, initially an almost-everywhere
statement, bounds their differences everywhere, including at $u=0$.
Passing to the limit proves the first inequality.
The second follows because $F(x)/x$ is the average of the $3^j$
numbers $3^jF_{a,j}(x)/x$, and the distance of one of these numbers
from their average is at most $1-3^{-j}$ times their diameter.
The displayed constants use the exact value of $D_{24}^{(1)}$.
\end{proof}

\begin{remark}
For two classes modulo $9$ with the same parent class modulo $3$,
the analogous computation gives
\[
  D_{20}^{(1,1)}=\frac{563\,831}{1\,048\,576}<0.5378,
\]
in the notation introduced in Section~\ref{sec:sadic}.
This is the distance for pairs with a fixed common parent;
the general distance $D_n^{(2)}$ also includes pairs with different
parents.
These finite bounds are quantitative examples rather than
limiting obstructions: Proposition~\ref{thm:M} gives
$D_\infty^{(j)}=0$ for every fixed $j$.
The invariant sets covered by the corollary include
$\{N:\Colmin(N)>C\}$ for each fixed $C$.
The corresponding assertions for the $3x-1$ map follow from the
reflection $z\mapsto-z$ discussed in Section~\ref{sec:full}.
\end{remark}

\subsection{A second-moment criterion}

Fix $n,q$, and abbreviate $\rho=\rho_{n,q}$ and
$p(c)=p_{n,q}(c)$.
Define the top-level energy by
\[
\begin{split}
  E_{\mathrm{top}}(n,q)
  &:=3^{q+1}\sum_{x\bmod3^{q+1}}\rho(x)^2
     -3^q\sum_{c\bmod3^q}p(c)^2\\
  &=3^{q+1}\sum_{c=0}^{3^q-1}\sum_{k=0}^2
      \left(\rho(c+3^qk)-\frac{p(c)}3\right)^2.
\end{split}
\]
Equivalently, with
$\widehat\rho(\xi)=\sum_x\rho(x)e(-\xi x/3^{q+1})$,
Parseval's identity gives
\[
  E_{\mathrm{top}}(n,q)
  =\sum_{\substack{\xi\bmod3^{q+1}\\3\nmid\xi}}
       |\widehat\rho(\xi)|^2.
\]

\begin{proposition}[a second-moment bound]\label{thm:L2}
For $0\le q\le n$,
\[
  t(n,q)\le\sqrt{\frac23}\,
                 \sqrt{E_{\mathrm{top}}(n,q)}.
\]
Consequently,
\[
  D_n^{(1)}
  \le\sum_{q=0}^n2^{-n}\binom nq
       \min\left\{1,
           \sqrt{\frac23E_{\mathrm{top}}(n,q)}\right\}.
\]
In particular, $D_n^{(1)}\to0$ if
$E_{\mathrm{top}}(n,Q_n)\to0$ in probability for
$Q_n\sim\Bin(n,\tfrac12)$.
\end{proposition}

\begin{proof}
For a fixed parent $c$, put $a_k=\rho(c+3^qk)$ and
$\bar a=p(c)/3$.
For three real numbers,
\[
  \frac12\sum_{k\bmod3}|a_k-a_{k-1}|
  =\max_k a_k-\min_k a_k,
  \qquad
  (\max_k a_k-\min_k a_k)^2
  \le2\sum_{k=0}^2(a_k-\bar a)^2.
\]
Summing over parents and applying Cauchy--Schwarz gives
\[
  t(n,q)^2
  \le2\cdot3^q\sum_{c,k}
        \left(\rho(c+3^qk)-\frac{p(c)}3\right)^2
  =\frac23E_{\mathrm{top}}(n,q).
\]
Combine this with $t(n,q)\le1$ and \eqref{eq:Dt}.
For the last assertion, the bounded random variables
$\min\{1,\sqrt{(2/3)E_{\mathrm{top}}(n,Q_n)}\}$ converge to zero
in probability, and therefore also in expectation.
\end{proof}

There is an exact collision interpretation of this energy.
For independent uniform $q$-subsets $S,S'$, put
\[
\begin{split}
  B_{n,q}&=\Pb(X_S\equiv X_{S'}\pmod{3^q}),\\
  R_{n,q}&=\Pb(X_S\equiv X_{S'}\pmod{3^{q+1}}\mid
                      X_S\equiv X_{S'}\pmod{3^q}).
\end{split}
\]
Since $B_{n,q}>0$, these quantities are well defined, and
\[
  E_{\mathrm{top}}(n,q)
  =3^qB_{n,q}(3R_{n,q}-1).
\]
Thus the sufficient condition in Proposition~\ref{thm:L2}
concerns the full normalized collision excess, including the
factor $3^qB_{n,q}$.
Convergence of $R_{n,q}$ to $1/3$ alone does not supply this bound.
An equivalent formulation of the sufficient condition is
\[
  \sum_{q=0}^n2^{-n}\binom nq
       \min\{1,E_{\mathrm{top}}(n,q)\}\longrightarrow0.
\]
The truncation matters: $E_{\mathrm{top}}(n,n)=2\cdot3^n$, so the
untruncated binomial mean is at least $2(3/2)^n$.

\begin{remark}[the limitation of fixed local moves]\label{rem:local}
Suppose $q\ge1$, and let $m$ be the smallest element of $S$.
Replacing $m$ by $m+\delta$, with $\delta=2$ for even $m$ and
$\delta=4$ for odd $m$, gives
\[
  X_{S'}-X_S
  =3^{q-1}2^m(2^\delta-1)
  \equiv3^q\pmod{3^{q+1}},
\]
provided that the new element remains the smallest and lies in
$\{1,\ldots,n\}$.
With the other elements and the parity of $m$ fixed, the three
possible top digits occur equally often over each complete
period of three admissible same-parity values of $m$.
The available interval need not contain many complete periods.

More precisely, fix $1\le r<q$ and let $M_r$ be the
$(r+1)$-st smallest element of $S$.
Conditioning on the largest $q-r$ elements leaves a uniform
$r$-subset of $\{1,\ldots,M_r-1\}$ and gives
$X_S=C+3^{q-r}X_{S_{\mathrm{small}}}$.
Convexity of total variation therefore yields
\[
  t(n,q)\le\mathbb E\,t(M_r-1,r).
\]
For fixed $r$ and $q/n\to p\in(0,1)$,
\[
  \Pb(M_r=r+1)
  =\frac{\binom{n-r-1}{q-r-1}}{\binom nq}
  \longrightarrow p^{r+1}.
\]
Since $t(r,r)=1$, the right-hand side of the preceding upper
bound has positive lower limit at least $p^{r+1}$.
Thus this fixed-$r$ conditioning estimate does not establish
vanishing of $t(n,q)$.
This observation concerns the upper bound, not a lower bound
for $t(n,q)$ itself.
\end{remark}

\section{\textbf{The $S$-adic Collatz map and ergodicity}}
\label{sec:sadic}

We realize the random Collatz chain as a factor of a deterministic
system whose $2$-adic coordinate supplies the parity sequence.
This gives an equivalence between the Liouville property and
nonsingular ergodicity, together with a criterion based on a
uniform bound below $1$ for the limiting mixing distances.

\subsection{Definition and basic properties}

Let
\[
  \Omega=\R\times\Z_3\times\Z_2,
  \qquad
  \mu=e^u\,du\otimes dz\otimes d\beta,
\]
where $d\beta$ is normalized Haar measure on $\Z_2$, and define
\[
  \mathcal T(u,z,\beta)=
  \begin{cases}
    (u-\log2,\ z/2,\ \beta/2),&\beta\in2\Z_2,\\[3pt]
    (u+\log(3/2),\ (3z+1)/2,\ (3\beta+1)/2),
       &\beta\in1+2\Z_2.
  \end{cases}
\]
We also denote the $2$-adic Collatz map in the third coordinate
by $T$.
The parity-vector correspondence of Lemma~\ref{lem:terras}
identifies each length-$n$ parity cylinder with one residue class
modulo $2^n$.
It follows that the parity coding conjugates $T$ on $\Z_2$ to the
one-sided Bernoulli shift and carries Haar measure to the fair-coin
product measure; see also \cite{Lagarias1985}.
In particular, $T$ preserves Haar measure.

For $N\in\N$, put $\iota(N)=(\log N,N,N)$.
Then
\[
  \mathcal T(\iota(N))=\iota(TN)+(\varepsilon_N,0,0),
  \qquad
  \varepsilon_N=
  \begin{cases}
    0,&N\text{ even},\\
    -\log(1+1/(3N)),&N\text{ odd}.
  \end{cases}
\]
A bounded measurable function $F$ on $\Omega$ is
\emph{invariant} if $F\circ\mathcal T=F$ almost everywhere.
Ergodicity is understood in the nonsingular sense: every such
function is almost everywhere constant.
The measure $\mu$ is not invariant.

\begin{lemma}\label{lem:pushforward}
The map $\mathcal T$ is nonsingular, and
\[
  \mathcal T_*\mu
  =\bigl(1+\mathbf1_{2+3\Z_3}(z)\bigr)\mu.
\]
In particular, $\mathcal T_*\mu$ and $\mu$ have the same null sets.
\end{lemma}

\begin{proof}
On the even branch, the pushforward factors in the three
coordinates are $2$, $1$ and $1/2$.
The last factor uses the restriction of $d\beta$ to $2\Z_2$,
which has mass $1/2$.
Thus this branch contributes $\mu$.

On the odd branch, the corresponding factors are
$2/3$, $3\mathbf1_{2+3\Z_3}(z)$ and $1/2$.
Indeed, $(g_1)_*(dz)=3\mathbf1_{2+3\Z_3}\,dz$, while the
$2$-adic odd branch pushes restricted Haar measure forward
to $\tfrac12d\beta$.
This branch contributes $\mathbf1_{2+3\Z_3}(z)\mu$.
Adding the contributions proves the formula.
\end{proof}

The multiplicity is the same as for positive integers:
every $m$ has the preimage $2m$, and it has the additional
odd preimage $(2m-1)/3$ exactly when $m\equiv2\pmod3$.

\subsection{The Liouville property and ergodicity}

Fix a probability density $\varrho$ on $\R$ which is positive
almost everywhere, and write
\[
  \mathbb P_{\varrho}
  =\varrho(u)\,du\otimes dz\otimes d\beta.
\]
This probability measure is equivalent to $\mu$.
All conditional expectations below are taken with respect to
$\mathbb P_{\varrho}$; invariance of this probability measure
is not required.

For bounded measurable $F$, define
\[
  h_F(u,z)=\int_{\Z_2}F(u,z,\beta)\,d\beta.
\]
Let $\mathcal G_n=\sigma(u,z,\beta\bmod2^n)$, completed with
respect to $\mathbb P_{\varrho}$, and let $X_n$ be the first
two coordinates of $\mathcal T^n(u,z,\beta)$.

\begin{lemma}\label{lem:hF}
For every bounded measurable $F$ and $n\ge0$,
\[
  Ph_F=h_{F\circ\mathcal T},
  \qquad
  \mathbb E_{\varrho}
    [F\circ\mathcal T^n\mid\mathcal G_n]=h_F(X_n)
\]
almost everywhere.
\end{lemma}

\begin{proof}
On each parity class, the $2$-adic branch sends normalized
Haar measure to Haar measure on $\Z_2$.
Integrating over the two classes therefore gives
\[
  h_{F\circ\mathcal T}(u,z)
  =\tfrac12h_F(u-\log2,z/2)
   +\tfrac12h_F(u+\log(3/2),(3z+1)/2),
\]
which is the first identity.

Conditionally on $\mathcal G_n$, the first $n$ parities and
$X_n$ are determined.
Writing $\beta=b+2^n\eta$, with $0\le b<2^n$, the variable
$\eta$ is conditionally Haar distributed, and
\[
  T^n\beta=T^n b+3^{q_n(b)}\eta.
\]
The multiplier is a unit in $\Z_2$.
Hence $T^n\beta$ has conditional Haar distribution, independently
of $\mathcal G_n$, proving the second identity.
\end{proof}

\begin{proposition}[Liouville--ergodicity equivalence]
\label{thm:equivergodic}
The following statements are equivalent:
\begin{enumerate}
  \item[(i)] Every bounded harmonic function of the random
  Collatz chain on $\R\times\Z_3$ is almost everywhere constant.
  \item[(ii)] The map $\mathcal T$ is ergodic on $(\Omega,\mu)$.
\end{enumerate}
\end{proposition}

\begin{proof}
All probabilistic statements use $\mathbb P_{\varrho}$.
Since this measure is equivalent to $\mu$, nonsingularity also
holds with respect to $\mathbb P_{\varrho}$.

(i)$\Rightarrow$(ii).
Let $F$ be bounded and invariant.
Lemma~\ref{lem:hF} gives $Ph_F=h_F$, so $h_F=c$ almost everywhere
by (i).
Nonsingularity implies, simultaneously for all $n$,
$F\circ\mathcal T^n=F$ and $h_F(X_n)=c$ almost surely.
Consequently,
\[
  \mathbb E_{\varrho}[F\mid\mathcal G_n]=c.
\]
The increasing sigma-algebras $\mathcal G_n$ generate the full
sigma-algebra up to completion.
Martingale convergence gives $F=c$ almost surely, and hence
$\mu$-almost everywhere.

(ii)$\Rightarrow$(i).
Choose a bounded Borel representative of a harmonic function $h$.
Given $\mathcal G_n$, the next parity is a fair independent coin.
The law of $X_n$ is absolutely continuous with respect to
$du\,dz$, by nonsingularity of $\mathcal T$.
Thus, almost surely,
\[
  \mathbb E_{\varrho}[h(X_{n+1})\mid\mathcal G_n]
  =Ph(X_n)=h(X_n).
\]
The bounded martingale $h(X_n)$ converges almost surely and in
$L^1$.
Define $F(\omega)$ as its limit wherever that limit exists, and
as zero otherwise.
This definition also applies to complex-valued $h$.
Deleting the first term leaves the convergence set and the limit
unchanged, so $F\circ\mathcal T=F$ at every point.
By (ii), $F=c$ almost everywhere.
Convergence in $L^1$ gives
\[
  h(u,z)=\mathbb E_{\varrho}[F\mid\mathcal G_0]=c
\]
almost surely.
The positivity of $\varrho$ then gives the desired equality
$du\,dz$-almost everywhere.
\end{proof}

\begin{corollary}[finite-level invariant functions]\label{thm:partialerg}
If a bounded $\mathcal T$-invariant function $F$ depends on $z$
only through $z\bmod3^j$ for some $j\ge0$, up to a null set,
then $F$ is almost everywhere constant.
In particular, the factor
\[
  (u,\beta)\longmapsto
  \begin{cases}
    (u-\log2,T\beta),&\beta\in2\Z_2,\\
    (u+\log(3/2),T\beta),&\beta\in1+2\Z_2
  \end{cases}
\]
is ergodic in the nonsingular sense.
\end{corollary}

\begin{proof}
The function $h_F$ is harmonic by Lemma~\ref{lem:hF} and has
level $j$.
Proposition~\ref{thm:finitelevel} makes $h_F$ constant.
The conditional-expectation argument in the first implication
of Proposition~\ref{thm:equivergodic} then makes $F$ constant.
The factor assertion is the case $j=0$.
\end{proof}

\subsection{Merging translations}

For a word $w$, let $[w]\subseteq\Z_2$ be its parity cylinder.
It is a residue class modulo $2^{|w|}$, and the restriction
$\gamma_w:[w]\to\Z_2$ is a bijection.
The following algebraic relation will also be used in the
coalescence construction of Section~\ref{sec:koal}.

\begin{lemma}\label{lem:coalescence}
Let $b$ have length $L$ and $r_0$ ones, and let
$\sigma,\sigma'\in\{0,1\}^W$ each have $r'$ ones.
Put $q=r_0+r'$ and suppose that
\[
  X_{b\sigma}\equiv X_{b\sigma'}\pmod{3^q}.
\]
Then
\[
  c=\frac{X_{b\sigma}-X_{b\sigma'}}{2\cdot3^q}
  \in2^L\Z,
\]
translation by $c$ sends $[b\sigma]$ onto $[b\sigma']$, and
\[
  \mathcal T^{L+W}(u,z,\beta)
  =\mathcal T^{L+W}(u,z+c,\beta+c)
  \qquad(\beta\in[b\sigma]).
\]
\end{lemma}

\begin{proof}
The prefix identity is
\[
  X_{b\sigma}=3^{r'}X_b+2^LX_\sigma.
\]
Since $X_\sigma$ and $X_{\sigma'}$ are even,
$X_{b\sigma}-X_{b\sigma'}$ is divisible by $2^{L+1}$.
It is also divisible by $3^q$ by hypothesis.
As these moduli are coprime, $c\in2^L\Z$.

Put $k=L+W$, $w=b\sigma$ and $w'=b\sigma'$.
By \eqref{eq:affine} and \eqref{eq:XS},
\[
  \frac{3^q}{2^k}c=B_w-B_{w'},
  \qquad \gamma_{w'}(y+c)=\gamma_w(y).
\]
This is an identity of rational affine maps, so it holds over
both $\Q_2$ and $\Q_3$.
Since $\gamma_w:[w]\to\Z_2$ and
$\gamma_{w'}:[w']\to\Z_2$ are bijections,
$\gamma_{w'}^{-1}\gamma_w$ sends $[w]$ onto $[w']$.
The affine identity identifies this map with translation by $c$.
The two terminal $2$-adic and $3$-adic coordinates therefore agree.
Their real coordinates agree as well, since both words have
length $k$ and contain $q$ ones.
\end{proof}

These merging maps act by simultaneous integer translations on
$z$ and $\beta$ and preserve the product measure on their domains.
The common prefix accounts for the divisibility $c\in2^L\Z$;
the congruence modulo $3^q$ supplies the additional integrality
needed for translation on $\Z_3$.

\subsection{Refined distances and condition (M$^*$)}

For $W\ge0$, $r\ge0$, $K\ge1$ and $0\le q\le W$, define the
finite measure on $\Z/3^{q+r+K}\Z$
\[
  \lambda_{W,q}^{(r,K)}
  =2^{-W}\sum_{\substack{|\sigma|=W\\q_\sigma=q}}
        \delta_{X_\sigma\bmod3^{q+r+K}}.
\]
Its mass is $2^{-W}\binom Wq$.
Let $\tau_t$ denote translation by $t$ on the relevant finite
cyclic group, and put
\[
\begin{split}
  d_W^{(r,K)}(v)
  &:=\sum_{q=0}^W
      \left\|\lambda_{W,q}^{(r,K)}
        -(\tau_{3^{q+r}v})_*\lambda_{W,q}^{(r,K)}
      \right\|_{\TV},\\
  D_W^{(r,K)}
  &:=\max_{v\in(\Z/3^K\Z)^*}d_W^{(r,K)}(v).
\end{split}
\]
The same unit $v$ is used for every $q$ in this sum.
These distances lie in $[0,1]$.

Equivalently, for fixed $q$, group words into parent cells
$c=X_\sigma\bmod3^{q+r}$, using $0\le c<3^{q+r}$, and let
$n_{q,c,\chi}$ count the words with
$X_\sigma\equiv c+3^{q+r}\chi\pmod{3^{q+r+K}}$.
Then
\[
  D_W^{(r,K)}
  =\max_{v\in(\Z/3^K\Z)^*}
    2^{-W}\sum_{q,c}\sum_{\chi\bmod3^K}
       (n_{q,c,\chi}-n_{q,c,\chi-v})^+.
\]
Thus the comparison preserves each parent and translates its
$K$ child digits, with the parent masses included in the average.

\begin{lemma}\label{lem:Drk}
For $W\ge0$, $r\ge0$ and $K\ge1$:
\begin{enumerate}
  \item[(a)] $D_{W+1}^{(r,K)}\le D_W^{(r,K)}$, so
  $a_{r,K}:=\lim_{W\to\infty}D_W^{(r,K)}$ exists.
  \item[(b)]
  $D_{W+1}^{(r,K)}\le\tfrac12D_W^{(r,K)}
                        +\tfrac12D_W^{(r+1,K)}$.
  \item[(c)] $a_{r,K}\le a_{r+1,K}$ and $a_{r,K}\le a_{r,K+1}$.
  \item[(d)] For every $j\ge1$,
  \[
    D_W^{(j)}=\max_{0\le r<j}D_W^{(r,j-r)},
    \qquad
    D_\infty^{(j)}=\max_{0\le r<j}a_{r,j-r}.
  \]
\end{enumerate}
\end{lemma}

\begin{proof}
(a) Appending $0$ leaves $X_\sigma$ and $q_\sigma$ unchanged.
Appending $1$ sends $X_\sigma$ to $3X_\sigma+2^{W+1}$ and
$q_\sigma$ to $q_\sigma+1$.
The comparison shift is multiplied by $3$ in the latter case,
so the same unit $v$ is used in both branches.
The triangle inequality and pushforward contraction give
$d_{W+1}^{(r,K)}(v)\le d_W^{(r,K)}(v)$.
Taking maxima proves (a).

(b) Prepending $0$ sends $X_\sigma$ to $2X_\sigma$.
Prepending $1$ sends it to $2X_\sigma+2\cdot3^{q_\sigma}$ and
increases the number of ones by $1$.
Thus the first branch uses shift parameter $r$, and the second
uses $r+1$ for the remaining word.
Multiplication by $2$ replaces the comparison unit by $2^{-1}v$.
Consequently,
\[
  d_{W+1}^{(r,K)}(v)
  \le\tfrac12d_W^{(r,K)}(2^{-1}v)
     +\tfrac12d_W^{(r+1,K)}(2^{-1}v).
\]
Taking maxima gives (b).

(c) The first inequality follows by letting $W\to\infty$ in (b).
For the second, lift a unit modulo $3^K$ to a unit modulo
$3^{K+1}$ and project the corresponding measures onto the
coarser cyclic group.
Total variation cannot increase under this projection, so
$D_W^{(r,K)}\le D_W^{(r,K+1)}$ for every $W$.

(d) Let $C=a+3^j\Z_3$ and $C'=a'+3^j\Z_3$ be distinct.
Write $a'-a=3^ru$, where $0\le r<j$ and $u$ is a unit modulo
$3^{j-r}$.
For fixed $q$, apply the affine bijection
\[
  y\longmapsto2^{W+1}y-2\cdot3^qa
\]
to the two output measures.
The first becomes the mixture of uniform measures on
$X_\sigma+3^{q+j}\Z_3$, and the second is its translate by
$3^{q+r}(2u)$.
Their total variation distance equals that of the corresponding
residue distributions modulo $3^{q+j}$.
Hence the sum over $q$ is $d_W^{(r,j-r)}(2u)$.
The unit $2u$ is independent of $q$, and every such unit is
obtained as $C,C'$ vary.
Taking maxima proves the first identity.
The maximum contains finitely many terms for each fixed $j$,
so taking limits proves the second.
\end{proof}

In particular,
\[
  D_W^{(1)}=D_W^{(0,1)},\qquad
  D_W^{(2)}=\max\{D_W^{(0,2)},D_W^{(1,1)}\}.
\]
Lemma~\ref{lem:Drk}(d) also shows that condition (M) is equivalent
to $a_{r,K}=0$ for every fixed $r\ge0$ and $K\ge1$.

\begin{definition}[mixing condition (M$^*$)]\label{conj:Mstar}
The random Collatz chain satisfies \emph{condition (M$^*$)} if
\[
  \kappa_*
  :=1-\sup_{\substack{r\ge0\\K\ge1}}a_{r,K}>0.
\]
\end{definition}

By monotonicity in Lemma~\ref{lem:Drk}(c), the supremum equals
$\lim_{r,K\to\infty}a_{r,K}$.
Moreover, part (d) gives the equivalent expression
\[
  \kappa_*=1-\sup_{j\ge1}D_\infty^{(j)}.
\]
Thus (M) implies (M$^*$) with $\kappa_*=1$.
The estimate $a_{0,1}<0.2138$ controls a single pair of indices
and does not by itself establish (M$^*$).
The coalescence estimates in Section~\ref{sec:natural} verify
(M) and give $\kappa_*=1$; see Proposition~\ref{thm:M} and
Corollary~\ref{cor:Mstar}.

\begin{proposition}[a uniform oscillation criterion]\label{thm:ergodicC}
If condition (M$^*$) holds, every bounded harmonic function of
the random Collatz chain is almost everywhere constant, and
$\mathcal T$ is ergodic.
Consequently, under the same condition, all conclusions of
Corollary~\ref{cor:main} hold for every $T$-invariant set.
\end{proposition}

\begin{proof}
It suffices to prove constancy for real-valued harmonic functions,
since the real and imaginary parts are harmonic separately.
Let $h$ be such a function.
Proposition~\ref{thm:quant} and Lemma~\ref{lem:Drk}(d) give,
for almost every $u$ simultaneously for every $j\ge1$,
\[
  \operatorname{osc}_{a\bmod3^j}E_jh(u,a)
  \le D_\infty^{(j)}\Omega(h)
  \le(1-\kappa_*)\Omega(h).
\]
For each such $u$, martingale convergence in the $3$-adic
coordinate gives $E_jh(u,z)\to h(u,z)$ for almost every $z$.
The uniform bound on the oscillation passes to this limit:
\[
  \operatorname*{ess\,osc}_{z\in\Z_3}h(u,z)
  \le(1-\kappa_*)\Omega(h).
\]
Taking the essential supremum over $u$ now yields
\[
  \Omega(h)\le(1-\kappa_*)\Omega(h).
\]
Since $\Omega(h)<\infty$ and $\kappa_*>0$, it follows that
$\Omega(h)=0$.
Thus $h(u,z)=H(u)$ almost everywhere, and Lemma~\ref{lem:CDR}
makes $H$ constant.
Proposition~\ref{thm:equivergodic} gives ergodicity, and
Corollary~\ref{cor:main} gives the counting-function conclusions.
\end{proof}

\begin{remark}
The factor $(u,\beta)$ also admits a direct zero--one argument.
Under parity coding, permuting finitely many parity coordinates
leaves the tail unchanged and preserves the number of ones in
the affected prefix.
The terminal real coordinate is therefore unchanged as well.
An invariant function of this factor is invariant under all
such permutations, outside a common null set.
For almost every fixed $u$, the Hewitt--Savage zero--one law
removes its dependence on the parity sequence, leaving a
function $H(u)$.
Invariance then gives the equation of Lemma~\ref{lem:CDR},
which makes $H$ constant.
In the full system, the merging argument must additionally
respect the congruence modulo $3^q$ in
Lemma~\ref{lem:coalescence}; this is the constraint imposed by
the $3$-adic coordinate.
\end{remark}

\section{\textbf{The coalescence method and unconditional ergodicity}}
\label{sec:koal}

Balanced coalescence of adjacent integers for the standard Collatz
map is established by Said Duran \cite[Corollary~C]{Duran2026}.
Chaining finitely many adjacent pairs gives the corresponding statement
for every fixed integer translation.
Leen \cite{LeenRates2026} proves quantitative absorption results and
sharper coalescence estimates.
We retain a proof in the relation coordinates used here, to connect
coalescence with the $S$-adic system and prepare the later generalization.

We study an explicit relation chain.
This gives the ergodicity of $\mathcal T$, the full Choquet--Deny
property, and the conclusions of Corollary~\ref{cor:main}.
The relevant translations are those of Lemma~\ref{lem:coalescence},
now considered for a fixed integer shift $c$.

\subsection{Coalescence sets and the relation chain}

Write $\lambda$ for normalized Haar measure on $\Z_2$ and
$q_n(\beta)$ for the number of odd elements among
$\beta,T\beta,\ldots,T^{n-1}\beta$, with $q_0=0$.
For $c\in\Z$, put
\[
 A_c=\{\beta\in\Z_2:\ T^n\beta=T^n(\beta+c),\quad
       q_n(\beta)=q_n(\beta+c)\text{ for some }n\ge0\}.
\]
Set $f(c)=1-\lambda(A_c)$.
We call this \emph{balanced coalescence}: both trajectories meet at
the same time and have used the same number of odd branches.
In particular, $A_0=\Z_2$ and $f(0)=0$.

Let $w_n,w_n'$ be the parity words of $\beta,\beta+c$, respectively,
and write $x_n=T^n\beta$, $y_n=T^n(\beta+c)$ and
\[
 \rho_n=\gamma_{w_n'}\circ\tau_c\circ\gamma_{w_n}^{-1},
 \qquad \tau_c(y)=y+c,
 \qquad \rho_n(x)=A_nx+B_n.
\]
Then $y_n=\rho_n(x_n)$ and
\[
 A_n=3^{k_n},\qquad k_n=q_n(\beta+c)-q_n(\beta).
\]

\begin{lemma}\label{lem:chainK}
\begin{enumerate}
 \item[(a)] We have $B_n=b_n3^{\min(k_n,0)}$ with $b_n\in\Z$.
 \item[(b)] If $\varepsilon$ is the parity of $x_n$, then the parity
 of $y_n$ is $\varepsilon+b_n\pmod2$. With $A=3^k$, the transition
 $(k,B)\mapsto(k',B')$ is
 \[
 \begin{array}{c|c|c}
       &\varepsilon=0&\varepsilon=1\\ \hline
 b\text{ even}&(k,B/2)&(k,(3B+1-A)/2)\\
 b\text{ odd}&(k+1,(3B+1)/2)&(k-1,(3B-A)/6).
 \end{array}
 \]
 \item[(c)] Under Haar measure on $\beta$, the normalized states
 $(k_n,b_n)$ form a Markov chain on $\Z\times\Z$, starting at
 $(0,c)$. At each step, $\varepsilon$ is a fair coin independent
 of the past.
 \item[(d)] The state $(0,0)$ is absorbing, and $\beta\in A_c$ if
 and only if this state is reached.
\end{enumerate}
\end{lemma}

\begin{proof}
If $\varepsilon'$ is the parity of $y_n$, then
\[
 \rho_{n+1}=g_{\varepsilon'}\circ\rho_n\circ g_\varepsilon^{-1}.
\]
The numbers $A_n$ and $3^{\min(k_n,0)}$ are $2$-adic units congruent
to $1$ modulo $2$, so $\varepsilon'\equiv\varepsilon+b_n\pmod2$.
Writing out the four cases gives the table. Substituting
$B=b3^{\min(k,0)}$ in each entry shows that the new normalized
coefficient $b'$ is again an integer; the divisions by $2$ are
exact by the indicated parity conditions. This proves (a) by
induction from $(0,c)$.

Both parity words of length $n$, and hence $\rho_n$, depend only
on $\beta\bmod2^n$.
Conditionally on this residue class, $T^n\beta$ is Haar distributed
on $\Z_2$, by Lemma~\ref{lem:terras} and its $2$-adic extension.
Thus its next parity is a fresh fair coin, proving (c).
Finally, equality of the two odd-step counts means $A_n=1$.
Together with $x_n=y_n$, this is equivalent to $B_n=0$.
The table shows that $(0,0)$ remains fixed under either coin outcome.
\end{proof}

We also use the same transition table, driven by independent fair
coins, from any state in
\[
 \mathcal S=\{(k,B):k\in\Z,\ B=b3^{\min(k,0)},\ b\in\Z\}.
\]
Unless normalized coordinates $(k,b)$ are explicitly specified,
states are written as $(k,B)$.
Write $\Xi_n=(k_n,B_n)$ and
\[
 \tau_{\rm abs}=\inf\{n\ge0:\Xi_n=(0,0)\}.
\]

\begin{remark}
If $\beta\in A_c$, then $\rho_n=\mathrm{id}$ at the coalescence
time, so
\[
 \mathcal T^n(u,z,\beta)=\mathcal T^n(u,z+c,\beta+c)
\]
for every $u$ and $z$.
Conversely, if these two points have equal iterates at times $m,n$,
equality of the real coordinates and the irrationality of
$\log3/\log2$ imply $m=n$ and equal odd-step counts.
Equality of the $2$-adic coordinates then gives $\beta\in A_c$.
Thus the intersection of the orbit relation with the graph of
$(u,z,\beta)\mapsto(u,z+c,\beta+c)$ is precisely
$\R\times\Z_3\times A_c$.
\end{remark}

\subsection{The coalescence inequality}

\begin{proposition}\label{thm:coalineq}
Let $h$ be a bounded, possibly complex-valued harmonic function,
and let $\operatorname{ess\,osc}h$ denote the diameter of its
closed essential range. For every $c\in\Z$,
\[
 |h(u,z+c)-h(u,z)|\le f(c)\operatorname{ess\,osc}h
\]
for almost every $(u,z)$.
For real-valued $h$, the oscillation is
$\operatorname{ess\,sup}h-\operatorname{ess\,inf}h$.
Moreover,
\[
 f(-c)=f(c),\qquad f(c+c')\le f(c)+f(c'),\qquad
 f(2m)=\tfrac12\bigl(f(m)+f(3m)\bigr).
\]
\end{proposition}

\begin{proof}
Let $K_h$ be the closed essential range of $h$, and choose a bounded
Borel representative taking values in $K_h$.
As in Proposition~\ref{thm:equivergodic}, the process $h(X_n)$ is
a bounded martingale after choosing a probability measure equivalent
to $\mu$.
Define $F$ to be its pointwise limit where the limit exists, and a
fixed element of $K_h$ elsewhere.
Deleting the first term of a sequence changes neither convergence
nor its limit, so $F\circ\mathcal T=F$ everywhere.
The martingale argument gives
\[
 h(u,z)=\int_{\Z_2}F(u,z,\beta)\,d\beta
\]
for almost every $(u,z)$.

For $\beta\in A_c$, the two sequences of first two coordinates
$X_n$ starting from $(u,z,\beta)$ and $(u,z+c,\beta+c)$ agree
from some time onwards. Hence their values of $F$ agree at every
such point, whether or not the sequences converge.
Translation in $\beta$ preserves Haar measure. Applying the
preceding integral identity at $z$ and $z+c$ therefore gives
\[
 h(u,z+c)-h(u,z)
 =\int_{\Z_2\setminus A_c}
   \bigl(F(u,z+c,\beta+c)-F(u,z,\beta)\bigr)\,d\beta.
\]
Both values lie in $K_h$, which proves the inequality.

Balanced coalescence is an equivalence relation: once a pair
coalesces with equal odd-step counts, both equalities persist at
every later time. Consequently,
\[
 A_c\cap(A_{c'}-c)\subseteq A_{c+c'},\qquad
 A_{-c}=A_c+c.
\]
Haar invariance gives subadditivity and symmetry.
From $(0,2m)$, the first step gives $(0,m)$ or $(0,3m)$ with
probability $1/2$ each, proving the last identity.
\end{proof}

\subsection{The fixed-point coordinate and uniform bounds}

For a state $(k,B)$, with $A=3^k$, define
\[
 p=\frac{B}{1-A}\quad(k\ne0),\qquad p=B\quad(k=0).
\]
For $k\ne0$, this is the fixed point of the relation
$x\mapsto Ax+B$.
Put $t_k=v_2(3^{|k|}-1)$ for $k\ne0$.
The lifting-the-exponent identity gives
\[
 t_k=1\quad(k\text{ odd}),\qquad
 t_k=2+v_2(k)\quad(k\text{ even}).
\]
We use $v_2(0)=\infty$.

\begin{lemma}\label{lem:pfix}
At each step, $p'=ap+d$, with the following coefficients.
\begin{enumerate}
 \item[(i)] For even $b$ and $k\ne0$, $(a,d)$ is $(1/2,0)$ or
 $(3/2,1/2)$. For even $b$ and $k=0$, it is $(1/2,0)$ or $(3/2,0)$.
 \item[(ii)] For odd $b$ and $k\ne0$, except for transitions to $0$,
 \begin{align*}
 a_+(k)&=\frac{3(1-A)}{2(1-3A)},&
 d_+(k)&=\frac1{2(1-3A)}&&(\varepsilon=0),\\
 a_-(k)&=\frac{3(1-A)}{2(3-A)},&
 d_-(k)&=-\frac{A}{2(3-A)}&&(\varepsilon=1).
 \end{align*}
 \item[(iii)] For odd $b$ at $k=0$, the coefficients are
 $(-3/4,-1/4)$ or $(3/4,-1/4)$.
 A transition from $k=1$ to $0$ has coefficients $(-1,-1/2)$,
 and a transition from $k=-1$ to $0$ has coefficients $(1,1/2)$.
\end{enumerate}
Always $|a|\le2$ and $|d|\le3/4$.
On $k\ge2$, the function $a_+(k)$ increases to $1/2$ and
$a_-(k)$ decreases to $3/2$; also
$a_\pm(-k)=a_\mp(k)$ for $|k|\ge2$.
Among nonabsorbing states, the only one-step predecessors of
$(0,0)$ are $(1,1)$, with coin $1$, and $(-1,-1/3)$, with coin
$0$. Both have $p=-1/2$. The absorbing state itself has two
self-loops.
\end{lemma}

\begin{proof}
For even $b$ and $k\ne0$, the relation is conjugated by
$g_\varepsilon$, so its fixed point is sent to $g_\varepsilon(p)$.
For $k=0$, use the transition table directly.
All remaining formulas follow by substituting $B=(1-A)p$ in
that table.
The derivatives
\[
 \frac{d}{dA}\frac{A-1}{3A-1}=\frac2{(3A-1)^2},\qquad
 \frac{d}{dA}\frac{A-1}{A-3}=-\frac2{(A-3)^2}
\]
give the monotonicity assertions, and the values at $k=\pm1,\pm2$
give the uniform bounds.
Finally, a nonabsorbing transition to $(0,0)$ must change $k$ from
$1$ or $-1$ to $0$; solving $B'=0$ gives the two stated states.
\end{proof}

\begin{lemma}[gaps]\label{lem:gaps}
Start from a nonabsorbing state, and let $D$ be the number of
physical steps up to and including the next change of $k$.
If $k\ne0$, then, for every integer $m\ge0$,
\[
 \Pb(D>t_k+m\mid\text{initial state})\le2^{-m},\qquad
 \E[D\mid\text{initial state}]\le t_k+2.
\]
More precisely, if $v=v_2(b)<t_k$, then $D=v+1$.
If $v\ge t_k$, including $b=0$, then $D=t_k+G$, where
$G$ is geometric on $\{1,2,\ldots\}$ with success probability
$1/2$.
For $k=0$ and $b\ne0$, we have $D=v_2(b)+1$ exactly.
\end{lemma}

\begin{proof}
If $b$ is odd, the next step changes $k$; this is the case $v=0$.
While $b$ is even and $k\ne0$, the two updates are
$B/2$ and $(3B+1-A)/2$.
The odd denominator in $B=b3^{\min(k,0)}$ does not affect
$2$-adic valuations, and $v_2(1-A)=t_k$.
If $1\le v<t_k$, both updates have valuation $v-1$.
If $v>t_k$, the coin $1$ gives valuation $t_k-1$; if $v=t_k$,
the coin $0$ does so. The other outcome leaves the valuation
at least $t_k$. For $b=0$, coin $1$ again gives $t_k-1$.
Thus, from valuation at least $t_k$, the number of trials to
reach $t_k-1$ is geometric with parameter $1/2$.
There are then $t_k-1$ halving steps and one level-changing step.
This proves the asserted formulas and bounds.
At $k=0$ with nonzero even $b$, both updates, $b/2$ and $3b/2$,
reduce its valuation by one, giving the last assertion.
\end{proof}

\begin{lemma}\label{lem:srw}
At successive changes of $k$, its increments are independent fair
signs, up to absorption.
After absorption, these increments may be continued with independent
fair signs to obtain a simple symmetric random walk on $\Z$.
Almost surely, either absorption occurs or the changes continue
indefinitely and every integer level, in particular $1$ and $-1$,
is visited infinitely often.
\end{lemma}

\begin{proof}
A change occurs precisely when the current $b$ is odd.
Its direction is $+1$ for coin $0$ and $-1$ for coin $1$; this
coin is independent of the information available before the step.
The gaps are almost surely finite at every nonabsorbing state by
Lemma~\ref{lem:gaps}.
Appending independent signs after absorption gives the stated walk.
Recurrence of that walk proves the last assertion on the event
that absorption never occurs.
\end{proof}

\begin{lemma}[uniform Lyapunov bound]\label{lem:lyap}
Let $s=1/10$ and set
\[
 D_s=(3/4)^s,\qquad
 \rho=\frac12\left((6/13)^s+2^s\right).
\]
Then $0.99868400<\rho<0.99868401<1$, and from every state
\[
 \E[|p'|^s\mid\text{state}]\le\rho|p|^s+D_s.
\]
For
\[
 M^s=\frac{2D_s}{1-\rho},
\]
the chain reaches $\{|p|\le M\}$ almost surely from every state.
On the event of nonabsorption, it visits this set at arbitrarily
large times.
\end{lemma}

\begin{proof}
Use $|ap+d|^s\le|a|^s|p|^s+|d|^s$ and Lemma~\ref{lem:pfix}.
For even $b$, the mean multiplier is
$\tfrac12((1/2)^s+(3/2)^s)<0.9873$.
For odd $b$ at levels $0$ and $\pm1$, it is, respectively,
$(3/4)^s<0.972$ and $\tfrac12((3/8)^s+1)<0.954$.
At $|k|=2$, it is exactly $\rho$.
For $|k|\ge3$, monotonicity bounds it by
\[
 \tfrac12\left((1/2)^s+(13/8)^s\right)<0.9914<\rho.
\]
These numerical inequalities follow, for example, by bounding each
positive tenth root between rational numbers and raising the bounds
to the tenth power.

For completeness, the hitting assertion follows directly from the
drift bound. Write $V_n=|p_n|^s$ and
$\sigma=\inf\{n\ge0:V_n\le M^s\}$.
Before $\sigma$,
\[
 \E[V_{n+1}-V_n\mid\mathcal F_n]
 \le-(1-\rho)V_n+D_s\le-D_s.
\]
Stopping and summing gives
$D_s\E(n\wedge\sigma)\le V_0$.
Thus $\sigma<\infty$ almost surely.
Applying this conclusion after every deterministic time shows that
visits occur at arbitrarily large times unless absorption has already
occurred.
\end{proof}

\begin{lemma}[uniform tightness at return]\label{lem:tight}
Let
\[
 \tau=\inf\{n\ge0:\Xi_n=(0,0)\text{ or }|k_n|=1\}.
\]
Then $\tau<\infty$ almost surely from every state, and there is
a constant $C_*<88964<9\cdot10^4$ such that
\[
 \E|p_\tau|^s\le|p_0|^s+C_*.
\]
\end{lemma}

\begin{proof}
First suppose $k_0\ge2$.
Before $\tau$, the level is at least $2$, and it reaches $1$
almost surely by Lemma~\ref{lem:srw}.
Number the affine coefficients by their transitions and put
\[
 \Pi_{j,n}=\prod_{i=j+1}^{n}|a_i|.
\]
Iteration and subadditivity of the $s$-th power give
\[
 |p_\tau|^s
 \le\Pi_{0,\tau}^s|p_0|^s
      +D_s\sum_{j=1}^{\tau}\Pi_{j,\tau}^s.
\]

Put $r=99/100$ and $\Psi(k)=r^{k-1}$ for $k\ge1$.
The process $\Pi_{0,n}^s\Psi(k_n)$, stopped at $\tau$, is a
nonnegative supermartingale.
At steps without a level change, the conditional multiplier is
$\tfrac12((1/2)^s+(3/2)^s)<1$.
At a level change, the required inequality is
\[
 \tfrac12a_+(k)^sr+\tfrac12a_-(k)^sr^{-1}\le1.
\]
For $k=2$, its left side is less than $0.99947$.
For $k\ge3$, it is at most
\[
 \tfrac12\left((1/2)^sr+(13/8)^sr^{-1}\right)<0.9921.
\]
Optional stopping at bounded times and Fatou's lemma imply, on
$\{j\le\tau\}$,
\[
 \E[\Pi_{j,\tau}^s\mid\mathcal F_j]\le\Psi(k_j).
\]
In particular, $\E\Pi_{0,\tau}^s\le r^{k_0-1}\le1$.

Let $N_k$ count the indices $1\le j\le\tau$ with $k_j=k$.
Then $N_1=1$.
For $k\ge2$, the count is bounded by the total physical sojourn
length over visits of the embedded walk to level $k$.
Each such visit has conditional expected length at most $t_k+2$;
whether the visit occurs before stopping is determined at its start.
The Green function of the walk started at $k_0$ and killed at $1$
gives an expected number of visits
$2\min(k_0-1,k-1)\le2(k-1)$.
Consequently,
\[
 \E\sum_{j=1}^{\tau}\Pi_{j,\tau}^s
 \le\E\sum_{j=1}^{\tau}\Psi(k_j)
 \le1+\sum_{k\ge2}2(k-1)r^{k-1}(t_k+2)=:C_*.
\]
Since $D_s<1$, this proves the required moment bound for $k_0\ge2$.
For $k_0\le-2$, use $\Psi(k)=r^{|k|-1}$ and
$a_\pm(-k)=a_\mp(k)$; the same argument applies.

We record a short evaluation of the convergent series, without
enumerating trajectories. Put
\[
 S_d=\sum_{m\ge1}(dm-1)r^{dm-1}
    =\frac{r^{d-1}(d-1+r^d)}{(1-r^d)^2}.
\]
The valuation formula for $t_k$ gives
\[
 C_*=1+\frac{6r}{(1-r)^2}+4S_2
          +2\sum_{a\ge2}S_{2^a}.
\]
The part through $a=8$ is less than $88957.632$.
For $D=512$, the remaining part is at most
\[
 \frac{2D r^{D-1}}{(1-r^D)^2(1-2r^D)}<6.167.
\]
Indeed, $S_d\le d r^{d-1}/(1-r^D)^2$ for $d\ge D$, and successive
dyadic terms $d r^{d-1}$ have ratio at most $2r^D<1$.
The two displayed numerical bounds are inequalities between explicit
rational numbers at $r=99/100$.
Thus $C_*<88963.799<88964$.

If $|k_0|=1$ or the initial state is absorbing, then $\tau=0$.
Finally, suppose $k_0=0$ and $b_0\ne0$, and put $v=v_2(b_0)$.
The first $v$ steps keep level $0$, and $\tau=v+1$.
Writing $\phi=\tfrac12((1/2)^s+(3/2)^s)<1$, their moment is
\[
 \E|B_v|^s=\phi^v|b_0|^s.
\]
The departing step has $|a|=3/4$ and $|d|=1/4$, whence
\[
 \E|p_\tau|^s
 \le(3/4)^s\phi^v|b_0|^s+(1/4)^s
 \le|p_0|^s+1\le|p_0|^s+C_*.
\]
This covers every initial state.
\end{proof}

\begin{lemma}[descent]\label{lem:descent}
From every state of the relation chain, absorption has positive
probability.
\end{lemma}

\begin{proof}
First consider $(0,c)$.
For nonzero even $c$, coin $0$ gives $(0,c/2)$.
For odd $c$, put $v=v_2(3c+1)$.
A string of $v$ zero coins followed by one coin $1$ leads to
\[
 (0,c'),\qquad c'=\frac{(3c+1)/2^v-1}{2}.
\]
For $c\ge1$, we have $0\le c'\le(3c-1)/4<c$.
For $c\le-3$, we have $|c'|\le(3|c|+1)/4<|c|$.
For $c=-1$, the three coins $1,0,0$ lead to absorption.
Induction on $|c|$ therefore gives a finite absorbing word from
every $(0,c)$.

From $(1,B)$, use zero coins to halve a nonzero even $B$ until
it becomes odd, then coin $1$ gives $(0,(B-1)/2)$.
If $B=0$, first use coin $1$, obtaining $(1,-1)$.
The case $k=-1$ follows by exchanging the two trajectories:
the relation is replaced by its inverse,
$(k,B)\mapsto(-k,-3^{-k}B)$, and each prescribed finite coin
string still has positive probability.
From any other nonabsorbing state, level $1$ or $-1$ is reached
almost surely by Lemma~\ref{lem:srw}.
Combining this with the preceding finite paths proves the claim.
\end{proof}

\begin{proposition}[balanced coalescence]\label{thm:K}
The relation chain is absorbed almost surely from every initial
state in $\mathcal S$.
In particular, for every $c\in\Z$,
$\lambda(A_c)=1$, or equivalently $f(c)=0$.
\end{proposition}

\begin{proof}
Let $g(\xi)$ be the absorption probability from state $\xi$, and
let $M$ be as in Lemma~\ref{lem:lyap}.
Choose $M'$ such that $(M^s+C_*)/M'^s\le1/2$.
If $|p(\xi)|\le M$, Lemma~\ref{lem:tight} and Markov's inequality
show that its stopping state at $\tau$ belongs, with probability
at least $1/2$, to the finite set
\[
 \mathcal F=\{(0,0)\}\cup
   \{(k,B)\in\mathcal S:|k|=1,\ |p|\le M'\}.
\]
This set is finite because $p\in\tfrac12\Z$ when $|k|=1$.
Lemma~\ref{lem:descent} gives
$\delta=\min_{\xi\in\mathcal F}g(\xi)>0$.
The strong Markov property therefore gives
\[
 g(\xi)\ge\delta/2\qquad\text{whenever }|p(\xi)|\le M,
\]
uniformly in the level.

For any initial state, let $\mathcal A$ be the event of eventual
absorption. The Markov property and the absorbing convention give
\[
 g(\Xi_n)=\Pb(\mathcal A\mid\mathcal F_n).
\]
By the upward martingale convergence theorem, this converges almost
surely to $\mathbf1_{\mathcal A}$.
On nonabsorption, Lemma~\ref{lem:lyap} gives arbitrarily late visits
to $\{|p|\le M\}$, where $g\ge\delta/2$.
This contradicts convergence to zero.
Thus absorption has probability one from every initial state.
\end{proof}

\subsection{Consequences}

\begin{corollary}\label{thm:liouville}
\begin{enumerate}
 \item[(a)] Every bounded harmonic function on $\R\times\Z_3$ is
 almost everywhere constant.
 \item[(b)] The $S$-adic Collatz map $\mathcal T$ is ergodic in the
 nonsingular sense.
 \item[(c)] For every $T$-invariant set $G$, all conclusions of
 Corollary~\ref{cor:main} hold: slow oscillation, equidistribution
 on the scale $x$ in each fixed family of classes modulo $3^j$,
 and the Abel comparison.
\end{enumerate}
\end{corollary}

\begin{proof}
By Propositions~\ref{thm:coalineq} and~\ref{thm:K},
$h(u,z+c)=h(u,z)$ almost everywhere for every $c\in\Z$.
Since these shifts are countable, for almost every $u$ the
$L^1(\Z_3)$ function $h(u,\cdot)$ is invariant under all of them.
Continuity of translations in $L^1(\Z_3)$ and density of $\Z$
in $\Z_3$ imply invariance under every $3$-adic translation.
Hence $h(u,z)=H(u)$ almost everywhere.
Its harmonic equation is that of Lemma~\ref{lem:CDR}, so $H$
is almost everywhere constant.
Part (b) follows from Proposition~\ref{thm:equivergodic}, and
part (c) from Corollary~\ref{cor:main}.
\end{proof}

\begin{corollary}[arithmetic]\label{cor:arith}
For every $c\in\Z\setminus\{0\}$, the set of positive integers
$N$ with $N+c\ge1$ for which
\[
 T^n(N)=T^n(N+c),\qquad q_n(N)=q_n(N+c)
\]
at some time $n\ge0$ has natural density $1$.
\end{corollary}

\begin{proof}
Absorption by step $n$ is determined by $\beta\bmod2^n$.
The corresponding union of integer residue classes therefore has
density $\lambda_n(c)$, the probability of absorption by that time.
After excluding the finitely many $N$ with $N+c<1$, each positive
integer in these classes has the asserted balanced coalescence.
By Proposition~\ref{thm:K}, $\lambda_n(c)\to1$, which proves the claim.
\end{proof}

\begin{corollary}[the $3x-1$ problem]\label{cor:3xm1}
For the shortcut map $T_-(n)=n/2$ on even integers and
$T_-(n)=(3n-1)/2$ on odd integers, the counterparts of
Proposition~\ref{thm:K} and Corollaries~\ref{thm:liouville}
and~\ref{cor:arith} hold.
Thus for every fixed $c\ne0$, almost every pair $N,N+c\in\N$
has balanced coalescing trajectories.
Membership in any specified cycle basin therefore agrees for
almost every such pair.
The counting functions of the basins of the cycles with minima
$1$, $5$ and $17$ oscillate slowly and are equidistributed on the
scale $x$ in every fixed family of classes modulo $3^j$.
\end{corollary}

\begin{proof}
Reflection gives
\[
 T_-(\beta)=-T_+(-\beta),\qquad
 g_1^-(z)=-g_1^+(-z).
\]
It preserves Haar measure and all parity counts.
Consequently, $A_c^-=-A_{-c}^+$ and
$\lambda(A_c^-)=\lambda(A_{-c}^+)=1$.
Reflection in $z$ conjugates the corresponding harmonic operators,
and reflection in both $z$ and $\beta$ conjugates the deterministic
systems.
If two trajectories coalesce, either both or neither enter a
specified cycle. Each cycle basin is exactly $T_-$-invariant.
The renewal and limiting-profile argument for $T_-$ described in
Section~\ref{sec:full} therefore gives the counting conclusions.
\end{proof}

\begin{remark}
Corollary~\ref{thm:liouville} proves these consequences directly,
without using conditions (M) or (M$^*$).
Both conditions are also established below, in
Proposition~\ref{thm:M} and Corollary~\ref{cor:Mstar}.
The limiting finite-level discrepancies bounded in
Section~\ref{sec:level1} consequently vanish.
Liouville's property alone does not establish natural density:
the limiting profile may have different constant values along
different subsequences.
The quantitative coalescence estimate in the next section supplies
the additional control needed for the density-transfer argument.
\end{remark}

\section{\textbf{Quantitative coalescence and transfer of densities}}
\label{sec:natural}

We now quantify balanced coalescence and use it to transfer density
information to arithmetic progressions, primes, and other sparse sets.
Natural-density almost-boundedness for the standard map was established
by Mazur \cite{Mazur2026}, and density existence for sets invariant
above a fixed threshold is covered by Leen \cite{LeenBrownian2026}.
The preparatory comparisons below give a route from the coalescence
construction to the transfer statements, with their hypotheses and
uniformity made explicit.
For $c\in\Z$ and an integer $n\ge0$, let $\lambda_n(c)$ be the
probability of absorption by step $n$ from $(0,c)$, and put
$f_n(c)=1-\lambda_n(c)$.
Thus $f_n(0)=0$, and $f_n(c)$ decreases to zero for every fixed
$c$ by Proposition~\ref{thm:K}.
Let $U_n(c)\subseteq\Z$ be the union of the residue classes modulo
$2^n$ for which the relation chain has not been absorbed by step $n$.
This set is periodic and has density $f_n(c)$.
For $N,N+c\in\N$, membership in $U_n(c)$ is exactly failure of
balanced coalescence by step $n$.
The positivity restriction is imposed separately in applications,
so it does not affect the periodic residue-class definition.

\subsection{Uniform coalescence}

Retain $s=1/10$, the exact value of $\rho$ from
Lemma~\ref{lem:lyap}, and $C_*<88964$ from Lemma~\ref{lem:tight}.
The same rational bounds used above give
\[
 C_0:=\frac{(3/4)^s}{1-\rho}<739,
 \qquad \frac{s}{\log(1/\rho)}<76.
\]
For the second inequality, it suffices to note that
$\rho<759/760$ and $\log(1/\rho)>1-\rho$.

\begin{proposition}[uniform coalescence]\label{prop:U}
For $\varepsilon\in(0,1)$, let
\[
 L(\varepsilon)
 =76\left(1+\frac{3600\log(10/\varepsilon)}{\varepsilon^2}\right).
\]
There is an integer $n_0(\varepsilon)$ such that $f_n(c)\le\varepsilon$
for every $c\in\Z$ and every integer
\[
 n\ge L(\varepsilon)\log(2+|c|)+n_0(\varepsilon).
\]
\end{proposition}

\begin{proof}
The case $c=0$ is immediate. Suppose $c\ne0$ and start at $(0,c)$,
so that $p_0=c$.
Iteration of Lemma~\ref{lem:lyap} gives
\[
 \E|p_t|^s\le\rho^t|c|^s+C_0.
\]
Set
\[
 t_1=\left\lceil\frac{s\log(2+|c|)}{\log(1/\rho)}\right\rceil
 \le76\log(2+|c|)+1,
 \qquad Q=1+C_0,
 \qquad M_1^s=\frac{5Q}{\varepsilon}.
\]
Then $\E|p_{t_1}|^s\le Q$, so
$\Pb(|p_{t_1}|>M_1)\le\varepsilon/5$.
The physical-time process $k_t$, kept at zero after absorption,
is a martingale with increments in $\{-1,0,1\}$.
Azuma's inequality~\cite{Azuma1967} gives, with
$A=\sqrt{2\log(10/\varepsilon)}$,
\[
 \Pb(|k_{t_1}|>A\sqrt{t_1})\le2e^{-A^2/2}=\varepsilon/5.
\]
Write $\mathcal G$ for the event that both bounds at time $t_1$
hold.

Let $B=10/\varepsilon$ and choose the integer barrier
\[
 H=\left\lceil B(A\sqrt{t_1}+1)\right\rceil.
\]
Since $A\sqrt{t_1}>2$ and $B>10$, we have
$H\le2BA\sqrt{t_1}$.
Let $\sigma\ge t_1$ be the first time the chain is absorbed or
$|k_\sigma|\in\{1,H\}$.
For a state in $\mathcal G$ with $2\le a:=|k_{t_1}|\le A\sqrt{t_1}$,
the embedded level walk is killed on reaching $1$ or $H$.
Gambler's ruin gives
\[
 \Pb(|k_\sigma|=H\mid\mathcal F_{t_1})
 =\frac{a-1}{H-1}\le\frac1B=\frac\varepsilon{10}.
\]
For $2\le k\le H-1$, its Green function is
\[
 G_H(a,k)=\frac{2(\min(a,k)-1)(H-\max(a,k))}{H-1}
 \le2(a-1).
\]
Lemma~\ref{lem:gaps} bounds each physical sojourn, conditionally
at its start, by $t_k+2$ in expectation.
If $J=\lfloor H/2\rfloor$, the valuation formula gives
\begin{equation}\label{eq:sumtk}
 \sum_{k=2}^{H}(t_k+2)
 =3H-3+3J-s_2(J)\le4.5H-3.5,
\end{equation}
where $s_2(J)$ is the binary digit sum.
Consequently, on these initial states,
\[
 \E[\sigma-t_1\mid\mathcal F_{t_1}]
 \le2A\sqrt{t_1}\cdot4.5H
 \le18BA^2t_1=\frac{180A^2t_1}{\varepsilon}.
\]
The same argument applies to negative levels.
Put
\[
 T_1=\frac{1800A^2t_1}{\varepsilon^2},\qquad
 R_1=\lceil\log_2M_1\rceil+1.
\]
The conditional probability that $\sigma-t_1>T_1$ is at most
$\varepsilon/10$ in the preceding case.
If $|k_{t_1}|=1$ or the chain is already absorbed, then $\sigma=t_1$.
If $k_{t_1}=0$ and the chain is not absorbed, its next change is
after $v_2(b_{t_1})+1\le R_1$ steps and reaches $|k|=1$.
Thus, on $\mathcal G$, the probability of reaching the upper
barrier is at most $\varepsilon/10$, and the probability of
$\sigma-t_1>T_1+R_1$ is at most $\varepsilon/10$.

Next define
\[
 \tau=\inf\{t\ge t_1:\Xi_t=(0,0)\text{ or }|k_t|=1\},\qquad
 M_2^s=\frac{10(M_1^s+C_*)}{\varepsilon}.
\]
By Lemma~\ref{lem:tight} and the strong Markov property,
\[
 \E[|p_\tau|^s\mid\mathcal F_{t_1}]
 \le|p_{t_1}|^s+C_*.
\]
It follows, on $\mathcal G$, that
\[
 \Pb(\sigma=\tau,\ |p_\sigma|>M_2\mid\mathcal F_{t_1})
 \le\Pb(|p_\tau|>M_2\mid\mathcal F_{t_1})
 \le\frac\varepsilon{10}.
\]
This bounds the intersection event; it does not condition on the
future event $\{\sigma=\tau\}$.

The set
\[
 \mathcal F(M_2)=\{(0,0)\}\cup
 \{(k,B)\in\mathcal S:|k|=1,\ |p|\le M_2\}
\]
is finite.
By the all-state assertion of Proposition~\ref{thm:K}, there is an
integer $n_2(\varepsilon)$ such that absorption within $n_2$ steps
has probability at least $1-\varepsilon/10$ from every member of
this set.
Applying the strong Markov property at $\sigma$ bounds failure
of this final stage by $\varepsilon/10$ whenever the preceding
stages succeeded.

The two initial exceptional probabilities are at most
$\varepsilon/5$ each; the upper-barrier, delay, return-moment,
and final-absorption exceptions are at most $\varepsilon/10$ each.
Their sum is $4\varepsilon/5<\varepsilon$.
Outside their union, absorption occurs by
\[
 t_1+T_1+R_1+n_2
 =\left(1+\frac{3600\log(10/\varepsilon)}{\varepsilon^2}\right)t_1
   +R_1+n_2.
\]
Using the bound on $t_1$ proves the proposition, for example with
\[
 n_0(\varepsilon)=R_1+n_2(\varepsilon)
 +\left\lceil1+
       \frac{3600\log(10/\varepsilon)}{\varepsilon^2}\right\rceil.
\]
\end{proof}

For \(N\ge0\), put
\[
  u_N=2^{-N}\binom N{\lfloor N/2\rfloor}.
\]
By the reflection principle, this is the probability that a simple
symmetric random walk started at \(0\) stays nonnegative for its first
\(N\) steps; see \cite[Ch.~III]{Feller1968}.

\begin{lemma}\label{lem:uN}
For every \(n\ge1\),
\[
  u_{n-1}\ge\sqrt{\frac2{\pi(n+1)}}.
\]
Also, \(u_N\le N^{-1/2}\) for \(N\ge1\).
\end{lemma}

\begin{proof}
The Wallis bounds give
\[
  \frac1{\sqrt{\pi(k+1/2)}}\le
  u_{2k}\le\frac1{\sqrt{\pi k}}\qquad(k\ge1),
\]
and \(u_{2k-1}=u_{2k}\).
The lower bound also holds for \(k=0\), since \(u_0=1\).
For completeness, the lower bound follows by observing that
\(\sqrt{k+1/2}\,u_{2k}\) decreases to \(1/\sqrt\pi\); the square
of the ratio of successive terms is
\[
  \frac{(2k+1)(2k+3)}{4(k+1)^2}<1.
\]
If \(n-1=2k\), this gives \(u_{n-1}\ge\sqrt{2/(\pi n)}\);
if \(n-1=2k-1\), it gives the stated lower bound.
The upper Wallis bound, together with \(\pi>2\), proves
\(u_N\le N^{-1/2}\) for both parities.
\end{proof}

\begin{proposition}[a lower bound for the merging time]\label{prop:lower}
For every \(c\ne0\) and \(n\ge1\),
\[
  f_n(c)\ge u_{n-1}\ge\sqrt{\frac2{\pi(n+1)}}.
\]
\end{proposition}

\begin{proof}
From a nonabsorbed state with \(k=0\), a step that keeps \(k=0\)
replaces the nonzero offset by \(B/2\) or \(3B/2\), and cannot
cause absorption.
By Lemma~\ref{lem:pfix}, absorption occurs only when the level
changes from \(1\) or \(-1\) to \(0\).
At the successive changes, the levels form a simple symmetric
random walk, with an independent continuation after absorption
as in Lemma~\ref{lem:srw}.
Its first level is \(1\) or \(-1\).
If this walk avoids \(0\) for the following \(n-1\) changes,
absorption cannot occur in the first \(n\) physical steps.
The probability of this event is \(u_{n-1}\).
This argument also covers the event that the first change occurs
after time \(n\).
\end{proof}

\begin{remark}\label{rem:Lorder}
Proposition~\ref{prop:lower} imposes an
\(\varepsilon^{-2}\) lower bound on the time needed to make the
nonabsorption probability at most \(\varepsilon\), even for \(c=1\).
In the notation of Proposition~\ref{prop:U}, it implies
\[
  \left\lceil L(\varepsilon)\log3+n_0(\varepsilon)\right\rceil+1
  \ge \frac2{\pi\varepsilon^2}.
\]
It does not give a lower bound of this order for
\(L(\varepsilon)\) alone: the additive term
\(n_0(\varepsilon)\) is also allowed to depend on \(\varepsilon\).
The proposition also shows that a bound expressed through these
nonabsorption probabilities cannot decay faster than
\(n^{-1/2}\).
This is not a lower bound for the mixing distances themselves.
\end{remark}

\subsection{Effective coalescence with a logarithmic loss}

The final appeal to almost-sure absorption on a finite set in
Proposition~\ref{prop:U} can be made effective.
We first give explicit absorbing words, then control the number
of visits at which such a word can be attempted.

\begin{lemma}[polynomial descent]\label{lem:polydesc}
\begin{enumerate}
  \item[(a)] For every \(c\ne0\), there is a sequence of coin
  outcomes of length
  \[
    \ell(c)\le4.84\log_2|c|+3
  \]
  leading from \((k,B)=(0,c)\) to \((0,0)\).
  \item[(b)] In the normalized coordinates \((k,b)\) of
  Lemma~\ref{lem:chainK}, every state \((\pm1,b)\) has an
  absorbing word of length at most
  \[
    5.84\log_2\max\{|b|,1\}+6.
  \]
  Consequently, if \(|k|=1\), \(|p|\le M\), and \(M\ge1\),
  absorption occurs within
  \[
    \ell_M=\lfloor5.84\log_2(2M)\rfloor+6
  \]
  steps with probability at least \(\delta_M=2^{-\ell_M}\).
\end{enumerate}
\end{lemma}

\begin{proof}
Use the paths in Lemma~\ref{lem:descent}.
An even nonzero \(c\) is halved in one step.
For odd \(c\ne-1\), put \(a=v_2(3c+1)\).
A word consisting of \(a\) zeros followed by a one leads to
\[
  c'=\frac{(3c+1)/2^a-1}{2}.
\]
For \(c=-1\), use the word \(100\), which absorbs in three
steps; for \(c=1\), the preceding rule also absorbs in three
steps.

Here is a uniform bound for the cost of a round when
\(d=|c|\ge1000\).
If \(a=1\), the exact formula, with the signs of \(c\) treated
separately, gives
\[
  |c'|\le\frac{3d+1}{4},\qquad
  \log_2\frac d{|c'|}
  \ge\log_2\frac4{3+1/1000}>0.41455
\]
whenever \(c'\ne0\).
Thus the two steps of the round cost less than
\(4.825\) per unit decrease of \(\log_2 d\).
For \(a=2\), use
\[
  |c'|\le\frac{3d+5}{8},
  \qquad
  \frac3{3-\log_2(3+5/1000)}<4.825.
\]
For \(a\ge3\), the inequality \(2^a\le3d+1\) gives
\[
  |c'|\le\frac{3d+1}{2^a},\qquad
  \frac{a+1}{a-\log_2(3+1/1000)}<4.825.
\]
If the round absorbs, its length is at most
\(\log_2(3d+1)+1<4.84\log_2d\).
Halving has cost one per unit decrease.

For \(1\le|c|<1000\), the same recursion, with the special word
for \(c=-1\), gives
\[
  \ell(c)\le4.84\log_2|c|+3.
\]
This finite check can be performed using integers only:
when \(\ell(c)\ge3\), check
\(2^{25(\ell(c)-3)}\le|c|^{121}\).
The maximum residual is \(3\), attained at \(c=\pm1\);
see Appendix~\ref{app:K0}.
Telescoping the logarithmic decreases until this finite range is
reached proves (a).

For (b), the coordinate \(b\) is integral even when \(k=-1\),
where \(B=b/3\).
If \(b\ne0\), choose \(v_2(b)\) halvings.
From \(k=1\), a subsequent outcome \(1\) gives
\((0,(b'-1)/2)\); from \(k=-1\), an outcome \(0\) gives
\((0,(b'+1)/2)\), where \(b'\) is the odd part of \(b\).
The new offset has absolute value at most \(|b|\), so (a)
proves the stated bound.
For \(b=0\), two steps lead to \((0,-1)\) or \((0,1)\),
and three more absorb.
Finally \(|b|=2|p|\) when \(|k|=1\).
An absorbing word of length at most \(\ell_M\) has probability
at least \(2^{-\ell_M}\); after absorption the state remains
absorbed.
\end{proof}

Write \(\tau_{\rm abs}\) for the absorption time and put
\(Y_t=|p_t|^s\), with \(s=1/10\).
We count every physical-time visit to \(\{|k|=1\}\), including
consecutive visits.
More precisely, set
\[
  V_1=\inf\{t\ge0:|k_t|=1\text{ or }t\ge\tau_{\rm abs}\},
  \qquad
  V_{i+1}=\inf\{t>V_i:|k_t|=1\text{ or }t\ge\tau_{\rm abs}\}.
\]
These are stopping times.
After absorption this convention inserts one fictitious visit at
each subsequent time, with \(Y_{V_i}=0\).

\begin{lemma}[contraction along visits]\label{lem:V}
There are absolute constants
\[
  \rho_V=\tfrac12\bigl((1/2)^s+(3/2)^s\bigr)<0.98722,
  \qquad C_V<43\,222,
\]
such that
\[
  \E[Y_{V_{i+1}}\mid\mathcal F_{V_i}]
  \le\rho_VY_{V_i}+C_V
\]
for every \(i\).
\end{lemma}

\begin{proof}
There is nothing to prove after absorption.
Consider \(k_{V_i}=1\); the case \(k=-1\) has the same
absolute-value estimates.
If \(b\) is even, the next physical step is another visit, and
\[
  \E Y_{V_{i+1}}\le\rho_VY_{V_i}
       +\tfrac12(1/2)^s.
\]
If \(b\) is odd, the outcome \(1\) enters level \(0\), with
\(p_1=-p-1/2\).
If this is the absorbing state the required bound is immediate.
Otherwise the number of steps at level \(0\) before departure
is determined by the valuation of its nonzero integral offset.
During the holding steps the multipliers are \(1/2,3/2\), with
mean \(s\)-th power \(\rho_V<1\), and there is no additive term.
The departure map has multiplier of absolute value \(3/4\)
and additive term \(-1/4\).
Thus, conditional on this inward outcome,
\[
  \E Y_{V_{i+1}}
  \le (3/4)^s\bigl(Y_{V_i}+(1/2)^s\bigr)+(1/4)^s.
\]
The other outcome enters level \(2\) with
\(p_1=(3/8)p-1/16\).
The weighted-product estimate in the proof of
Lemma~\ref{lem:tight}, with \(r=0.99\), gives
\[
  \E Y_{V_{i+1}}
  \le r\bigl((3/8)^sY_{V_i}+(1/16)^s\bigr)
       +(3/4)^s C_*.
\]
For odd \(b\), the coefficient of \(Y_{V_i}\) is therefore
\[
  \tfrac12r(3/8)^s+\tfrac12(3/4)^s<0.93458<\rho_V.
\]
The additive terms are bounded by
\[
  \tfrac12\bigl(r(1/16)^s+(3/4)^sC_*\bigr)
  +\tfrac12\bigl((3/4)^s(1/2)^s+(1/4)^s\bigr)<43\,222.
\]
Here one may use the bound \(C_*<88\,964\) for the series in
Lemma~\ref{lem:tight}.
All the displayed numerical comparisons can be checked with
rational upper bounds for the relevant tenth roots.
\end{proof}

\begin{proposition}[an effective tail bound]\label{thm:F}
There is an effectively computable absolute constant \(C\)
such that, for all \(c\in\Z\) and \(n\ge2\),
\[
  f_n(c)\le C\frac{\Lambda^{3/2}}{\sqrt n},
  \qquad \Lambda=\log n+\log(2+|c|).
\]
The same upper bound holds for a chain started at \(|k|=1\),
with \(|c|\) replaced by the absolute value of its initial
fixed-point coordinate.
For \(c\ne0\), Proposition~\ref{prop:lower} supplies the
lower bound \(\sqrt{2/(\pi(n+1))}\).
\end{proposition}

\begin{proof}
For \(c=0\) the assertion is immediate.  For \(c\ne0\),
we first establish a visit-count estimate that will also be
used below.
Put
\[
  M^s=\frac{2C_V}{1-\rho_V},\qquad
  \theta=\frac{1+\rho_V}{2}<0.9937,
  \qquad
  \lambda=\frac12\log(1/\theta),\qquad
  \kappa=1+\sqrt\theta.
\]
A visit is good if \(Y_{V_i}\le M^s\).
At a bad visit, the conditional drift in Lemma~\ref{lem:V}
is at most \(\theta Y_{V_i}\); at a good visit, the conditional
expectation at the next visit is at most \(\theta M^s\).

Let \(R_0\) count the bad visits before the first good visit,
and let \(R_i\), \(i\ge1\), count the bad visits between
successive good visits.
Stopping the drift inequality at the next good visit gives
\[
  \Pb(R_i\ge j\mid\text{history at the preceding good visit})
  \le\theta^j\qquad(i\ge1,\ j\ge1).
\]
For a chain started at \((0,c)\), \(c\ne0\),
\[
  \E Y_{V_1}\le(3/4)^s\rho_V^{v_2(c)}|c|^s+(1/4)^s
  \le |c|^s+1.
\]
Consequently
\(\Pb(R_0\ge j)\le(1+|c|^s)\theta^{j-1}\).
Summing these integer-valued tail bounds gives
\[
  \E[e^{\lambda R_i}\mid\text{preceding history}]\le\kappa
  \quad(i\ge1),\qquad
  \E e^{\lambda R_0}\le2(1+|c|^s).
\]
The tower property and Markov's inequality now imply, for
positive integers \(m,g\),
\begin{equation}\label{eq:goodvisits}
  \Pb(\text{fewer than \(g\) good visits among the first \(m\)})
  \le2(1+|c|^s)\kappa^g e^{-\lambda(m-g)}.
\end{equation}
In this argument the stopped drift is applied to the visit-indexed
sequence \(Y_{V_i}\), not to consecutive physical times.

Let \(\ell=\ell_M\) and \(\delta=\delta_M\) be the constants
of Lemma~\ref{lem:polydesc}.
Start an attempt at the first good visit, and start each later
attempt at the first good visit strictly after the preceding
starting time plus \(\ell\).
Each attempt has conditional success probability at least
\(\delta\).
For
\[
  r=\left\lceil\frac{\log n}{2\delta}\right\rceil,\qquad
  g=r(\ell+1)+1,
\]
the probability that the first \(r\) attempts all fail is at most
\((1-\delta)^r\le n^{-1/2}\).
An attempt excludes at most \(\ell+1\) good visits from being
later starting times.
Therefore these \(r\) attempts finish by time \(n\) if the
first \(m\) visits contain at least \(g\) good visits and
\(V_m\le n-\ell\).
Choose
\[
  m=g+\left\lceil
       \frac{g\log\kappa+\log(2+2|c|^s)+\tfrac12\log n}
            {\lambda}\right\rceil.
\]
Then \eqref{eq:goodvisits} is at most \(n^{-1/2}\), and
\(m\le C_{\mathrm v}\Lambda\) for an effectively computable
absolute constant \(C_{\mathrm v}\).

We next compare physical time with the number of level changes.
Let \(C_0=(3/4)^s/(1-\rho)\), as in Proposition~\ref{prop:U}.
The bound \(\E Y_t\le |c|^s+C_0\) and a union bound give
\[
  \Pb\left(\max_{0\le t\le n}|p_t|>
             (n+2+|c|)^{40}\right)
  \le \frac{1+C_0}{n^2}.
\]
On the complementary event, every holding interval that starts
before time \(n\), while the chain is unabsorbed, has length
at most \(G+L_i\), where
\[
  G=4+40\log_2(n+2+|c|)\le C_{\mathrm g}\Lambda.
\]
At level \(0\) take \(L_i=0\), using the exact valuation
bound for the holding time.
At a nonzero level use Lemma~\ref{lem:gaps} and
\(|k_t|\le n+1\).
More precisely, if \(D_i\) is the holding time and \(k_i\ne0\)
is its level, define \(L_i=(D_i-t_{k_i})_+\), independently
of the preceding high-probability event.
These nonnegative integer variables satisfy
\[
  \Pb(L_i>t\mid\text{history at the start of the interval})
  \le2^{-t}\qquad(t=0,1,\ldots).
\]
Stop this sequence after absorption or once a starting time
exceeds \(n\), and pad it with zeros.
Its conditional exponential moments satisfy
\[
  \E[e^{(\log2)L_i/2}\mid\text{preceding history}]
  \le1+\sqrt2<2.42.
\]
It suffices to consider \(n\ge2\ell\) and \(n\ge4(G+8)\).
In the remaining cases \(n\le C'\Lambda\), or \(n\le2\ell\);
after increasing the final absolute constant, the asserted
bound is at least \(1\), so it is trivial.
For every positive integer \(J\), the moment bound gives
\[
  \Pb\left(\sum_{i=1}^J L_i>8J\right)
  \le(2.42/16)^J\le e^{-J}.
\]
Thus, on survival to time \(n\), at least
\[
  J=\left\lfloor\frac{n-\ell}{G+8}\right\rfloor
\]
level changes occur by time \(n-\ell\), apart from the two
exceptional events just bounded.
Indeed, otherwise the first holding interval crossing time
\(n-\ell\) would have an index \(i\le J\) and start by time
\(n-\ell\le n\). Its completion time would be both greater
than \(n-\ell\) and at most
\(iG+\sum_{a\le i}L_a\le J(G+8)\le n-\ell\), a contradiction.

We use the following elementary return-count bound for simple
random walk. If \(R_J\) counts returns to \(0\) in its first
\(J\) steps, then
\[
  \Pb(R_J<m+1)\le C_{\mathrm r}\frac{m+1}{\sqrt J}
  \qquad(J\ge1).
\]
Indeed, a first-return time \(\xi\) has
\(\Pb(\xi>t)=u_{t-1}\) for \(t\ge1\), and hence
\(\E\min(\xi,J)\le4\sqrt J\) by Lemma~\ref{lem:uN}.
For the successive return times, apply
\(\min(\sum_i\xi_i,J)\le\sum_i\min(\xi_i,J)\) and Markov's
inequality. The same bound, after changing the absolute
constant, holds when the walk starts at \(1\) or \(-1\).

The embedded levels have this simple-walk law.
Each return to \(0\), except possibly the last one counted,
is followed by a visit to \(\{|k|=1\}\).
Using the independent continuation after absorption, we may
apply the return-count bound on the survival event.
We obtain
\[
  f_n(c)\le
  \frac2{\sqrt n}
  +C_{\mathrm r}\frac{m+1}{\sqrt J}
  +e^{-J}+\frac{1+C_0}{n^2}.
\]
In the range now under consideration,
\(J\ge n/(4(G+8))\ge n/(C''\Lambda)\), which proves the
claimed bound.
For a start at \(|k|=1\), take \(V_1=0\) and replace the
initial moment by \(|p_0|^s\).
All the estimates just used remain valid, including the
return-count bound from \(1\) or \(-1\).
\end{proof}

\begin{corollary}[effective uniform coalescence]\label{cor:effU}
In Proposition~\ref{prop:U}, one can take
\[
  n_0(\varepsilon)
  \le C\varepsilon^{-2}\log^3(2/\varepsilon)
\]
with an effectively computable absolute constant \(C\).
\end{corollary}

\begin{proof}
The final stage in the proof of Proposition~\ref{prop:U}
starts either at absorption or at a state with \(|k|=1\) and
\(|p|\le M_2\).
The definition there makes \(M_2\) polynomial in
\(\varepsilon^{-1}\), with absolute coefficients and degree.
The extension of Proposition~\ref{thm:F} to such starts gives
failure probability at most
\[
  C\frac{(\log n_2+\log(2+M_2))^{3/2}}{\sqrt{n_2}}.
\]
For a sufficiently large, effectively computable absolute
constant \(C'\), this is at most \(\varepsilon/10\) when
\[
  n_2=\left\lceil
        C'\varepsilon^{-2}\log^3(2/\varepsilon)\right\rceil.
\]
The remaining additive terms in that proof are
\(O(\varepsilon^{-2}\log(2/\varepsilon))\) and
\(O(\log(2/\varepsilon))\).
They are absorbed in the asserted bound for \(n_0\).
\end{proof}

\begin{remark}
Proposition~\ref{thm:F} also gives directly
\[
  f_n(c)\le\varepsilon\quad\text{if}\quad
  n\ge C\varepsilon^{-2}
        \bigl(\log(2+|c|)+\log(2/\varepsilon)\bigr)^3.
\]
Proposition~\ref{prop:U} retains a linear, rather than cubic,
dependence on \(\log(2+|c|)\).
The effective constants from the absorbing-word construction
are very large; the argument is not intended to give practical
mixing times.
\end{remark}

\subsection{Square-root tails for fixed translations}\label{subsec:exact}

We now remove the logarithmic loss in \(n\).
The proof estimates the truncated mean
\[
  \E\min(\tau_{\rm abs},n)=\sum_{t=0}^{n-1}f_t(c)
\]
by adding the contributions from successive returns of the
level process to \(0\).

Let \(S\) be a simple symmetric random walk, let \(\Pb_x\)
denote its law when started at \(x\), and put
\[
  H=\inf\{i\ge1:S_i=0\},\qquad p_i(x)=\Pb_0(S_i=x).
\]

\begin{lemma}[killed walk]\label{lem:meander}
For \(i\ge1\) and \(k\ge1\),
\begin{equation}\label{eq:meander}
  \Pb_1(S_i=k,H>i)
  =p_i(k-1)-p_i(k+1)
  =\frac{2k}{i+k+1}\,p_i(k-1).
\end{equation}
Also,
\[
  p_i(x)\le2i^{-1/2}e^{-x^2/(4i)}.
\]
There is an effectively computable absolute constant \(C_1\)
such that, for \(Q\ge1\) and \(b\in\Z\),
\begin{equation}\label{eq:meanderAP}
  \begin{split}
    \Pb_1(S_i\equiv b\pmod Q,H>i)
      &\le C_1\left(\frac1{Q\sqrt i}+\frac1i\right),\\
    \E_1[v_2(S_i);H>i]&\le \frac{C_1}{\sqrt i}.
  \end{split}
\end{equation}
Finally, \(\Pb_1(H>i)=u_i\le i^{-1/2}\).
\end{lemma}

\begin{proof}
Reflection at the first hit of \(0\) proves the first
identity in \eqref{eq:meander}.
The ratio of the two binomial probabilities is
\[
  \frac{p_i(k+1)}{p_i(k-1)}
   =\frac{i-k+1}{i+k+1}
\]
whenever the denominator is nonzero; the zero-probability
cases give the same identity.

For \(|x|\le i/2\), the binomial estimate in
Lemma~\ref{lem:local}, with \(m=i\) and \(t=x/2\), gives
\(p_i(x)\le\kappa_0i^{-1/2}e^{-x^2/(2i)}\), where
\(\kappa_0<2\).
For \(|x|>i/2\), the one-sided Hoeffding bound gives
\[
  p_i(x)\le e^{-x^2/(2i)}
  \le e^{-x^2/(4i)}e^{-i/16},
  \qquad
  \sqrt i\,e^{-i/16}\le\sqrt8\,e^{-1/2}<2.
\]
This proves the Gaussian bound.

By \eqref{eq:meander}, the probability on a progression is
at most
\[
  \frac2i\sum_{\substack{k\ge1\\k\equiv b\ (Q)}} k\,p_i(k-1).
\]
The term \(k=1\) is at most \(4i^{-3/2}\).
For \(k\ge2\), use \((k-1)^2\ge k^2/4\).
The function \(x e^{-x^2/(16i)}\) is unimodal, with supremum
\(\sqrt{8i/e}\) and integral \(8i\) over the positive half-line.
Its sum over a progression is at most its supremum plus
\(1/Q\) times its integral.
This proves the first estimate in \eqref{eq:meanderAP}.
On \(\{H>i\}\), \(1\le S_i\le i+1\), so
\[
  v_2(S_i)=\sum_{1\le a\le\log_2(i+1)}
                 \mathbf1_{2^a\mid S_i}.
\]
Apply the progression estimate with \(Q=2^a\) and sum.
Since \(\log(i+1)/i=O(i^{-1/2})\), the second estimate follows,
after increasing the absolute constant.
The last assertion follows from the reflection principle and
Lemma~\ref{lem:uN}.
\end{proof}

For a chain started at \((0,c)\), \(c\ne0\), let \(\tau_0=0\)
and let \(\tau_j\), \(j\ge1\), be the successive times at
which the level changes from \(1\) or \(-1\) to \(0\).
After absorption, later such times are set equal to infinity.
Let
\[
  E_j=\{\tau_j<\infty,\ \tau_{\rm abs}>\tau_j\},
  \qquad
  \ell_t=\log(2+|p_t|).
\]
In particular, \(E_0\) is the sure event.
Thus \(E_j\in\mathcal F_{\tau_j}\).
Products with \(\mathbf1_{E_j}\) below are defined to be zero
off \(E_j\); no value of \(p_\infty\) is used.

\begin{lemma}[return gaps]\label{lem:gap}
There is an effectively computable absolute constant \(C_3\)
such that, on \(E_j\), for every \(j\ge0\) and \(n\ge1\),
\[
  \E\!\left[
       \min(\tau_{j+1}-\tau_j,n)\mid\mathcal F_{\tau_j}
      \right]
  \le3\ell_{\tau_j}+C_3\sqrt n.
\]
\end{lemma}

\begin{proof}
On \(E_j\), the state at \(\tau_j\) is \((0,B)\) with
nonzero integral \(B=p_{\tau_j}\).
The time until its departure from \(0\) is exactly
\[
  v_2(B)+1\le1+\log_2|B|\le3\ell_{\tau_j}.
\]
After that departure, the embedded walk starts at \(1\) or
\(-1\) and runs until its first hit of \(0\).
Denote its successive nonzero levels by \(S_i\), and the
starting times of the corresponding holding intervals by
\(\sigma_i\).
The event \(\{i<\min(H,n)\}\) is measurable at \(\sigma_i\).
Conditionally on the history at that time, the expected
holding time is at most
\[
  t_{S_i}+2\le4+v_2(|S_i|)
\]
by Lemma~\ref{lem:gaps}.
Each holding interval lasts at least one physical step.
Therefore the truncated return gap is bounded by the time
spent at level \(0\), plus the full holding times with
indices \(i<\min(H,n)\).
Taking conditional expectations, using the strong Markov
property and symmetry, gives
\[
  \begin{split}
  \E[\min(\tau_{j+1}-\tau_j,n)\mid\mathcal F_{\tau_j}]
  &\le3\ell_{\tau_j}
     +\sum_{i<n}\E_1[4+v_2(S_i);H>i]\\
  &\le3\ell_{\tau_j}+4
     +(4+C_1)\sum_{1\le i<n}i^{-1/2}.
  \end{split}
\]
This proves the claim.
In particular, it does not require independence between a
holding time and the direction of a later change.
\end{proof}

\begin{lemma}[survival through returns]\label{lem:survret}
There is an effectively computable absolute constant \(C_5\)
such that, for \(c\ne0\) and \(L_c=1+\log(2+|c|)\),
\[
  \sum_{j\ge0}\Pb(E_j)\le C_5L_c,\qquad
  \sum_{j\ge0}\E[\mathbf1_{E_j}\ell_{\tau_j}]
       \le C_5L_c^2.
\]
\end{lemma}

\begin{proof}
Use the visits \(V_i\), the constants \(\lambda,\kappa,\ell,\delta\),
and the nonoverlapping attempts from the proof of
Proposition~\ref{thm:F}.
Put \(N_i=\{\tau_{\rm abs}>V_i\}\), so that survival is
required at the visit itself.
Choose a fixed \(\theta_1\in(0,1)\) such that
\[
  \theta_1(\log\kappa+\lambda)\le\lambda/2.
\]
For \(i>\ell\), let \(m=i-\ell\) and
\(g=\lfloor\theta_1m\rfloor\).
If the first \(m\) visits contain at least \(g\) good ones,
at least \(\lfloor g/(\ell+1)\rfloor\) attempts start among
them.
Every such attempt finishes by \(V_m+\ell\le V_i\),
because successive visit times differ by at least one.
On \(N_i\) all these attempts fail.
Using \eqref{eq:goodvisits} and the strong Markov property,
and absorbing the finitely many small indices into a constant,
we obtain
\[
  \Pb(N_i)\le C_4(1+|c|^s)e^{-\beta i}
  \qquad(i\ge1)
\]
for effectively computable absolute constants \(C_4\) and
\(\beta>0\).
Indeed, the two terms are bounded by
\[
  2(1+|c|^s)e^{-\lambda m/2}
  +(1-\delta)^{\lfloor g/(\ell+1)\rfloor},
\]
and both have exponential decay in \(i\).

Every nonabsorbing return has a preceding visit to \(|k|=1\),
and distinct returns have distinct preceding visits.
Thus \(E_j\subseteq N_j\) for \(j\ge1\).
Summing
\(\min\{1,C_4(1+|c|^s)e^{-\beta j}\}\) proves the first
assertion.

For the weighted assertion, if \(V_i=\tau_j-1\) is the
visit preceding a nonabsorbing return, Lemma~\ref{lem:pfix}
gives
\[
  |p_{\tau_j}|\le |p_{V_i}|+\tfrac12.
\]
The preceding-visit index is injective in \(j\), and \(N_i\)
holds there. Hence
\[
  \sum_{j\ge1}\E[\mathbf1_{E_j}\ell_{\tau_j}]
  \le\sum_{i\ge1}
      \E[\mathbf1_{N_i}\log(3+|p_{V_i}|)].
\]
Iteration of Lemma~\ref{lem:V} gives, uniformly in \(i\),
\[
  \E|p_{V_i}|^s
  \le |c|^s+1+\frac{C_V}{1-\rho_V}=:B_c.
\]
Concavity and \((3+x)^s\le3^s+x^s\) imply
\[
  \E\log(3+|p_{V_i}|)
  \le s^{-1}\log(3^s+B_c)\le C L_c.
\]
Also,
\[
  \log(3+x)\le\frac2s(3+x)^{s/2}
          \le\frac2s(2+x^{s/2}),
\]
so Cauchy--Schwarz and the exponential survival bound give
\[
  \E[\mathbf1_{N_i}\log(3+|p_{V_i}|)]
  \le\frac2s\left(2\Pb(N_i)+\sqrt{\Pb(N_i)B_c}\right)
  \le C(1+|c|^s)e^{-\beta i/2}.
\]
Use the first bound for
\(i\le\lceil(2/\beta)\log(1+|c|^s)\rceil\) and the second
for the remaining indices.
The resulting sum is \(O(L_c^2)\).
The term \(j=0\) equals \(\log(2+|c|)\), completing the proof.
\end{proof}

\begin{proposition}[Square-root order for fixed translations]\label{thm:exact}
There is an effectively computable absolute constant \(C\)
such that, for \(c\ne0\), \(n\ge1\), and \(L_c=1+\log(2+|c|)\),
\[
  \sum_{t=0}^{n-1}f_t(c)
    \le C(L_c\sqrt n+L_c^2),\qquad
  \sqrt{\frac2{\pi(n+1)}}\le f_n(c)
    \le C\left(\frac{L_c}{\sqrt n}+\frac{L_c^2}{n}\right).
\]
In particular, \(f_n(1)\) has exact order \(n^{-1/2}\), and
\[
  f_n(c)\le C'\frac{\log(2+|c|)}{\sqrt n}
  \qquad\text{when }n\ge\log^2(2+|c|).
\]
\end{proposition}

\begin{proof}
Absorption occurs at a return to level \(0\).
Consequently, pathwise,
\[
  \min(\tau_{\rm abs},n)
  \le\sum_{j\ge0}\mathbf1_{E_j}
            \min(\tau_{j+1}-\tau_j,n).
\]
All summands are nonnegative, and \(E_j\) is measurable at
\(\tau_j\).
Tonelli's theorem and Lemma~\ref{lem:gap} therefore give
\[
  \begin{split}
  \sum_{t=0}^{n-1}f_t(c)
   =\E\min(\tau_{\rm abs},n)
   &\le
    \sum_{j\ge0}\E[
      \mathbf1_{E_j}(3\ell_{\tau_j}+C_3\sqrt n)]\\
   &\le C(L_c\sqrt n+L_c^2)
  \end{split}
\]
by Lemma~\ref{lem:survret}.
Since \(f_t(c)\) is nonincreasing,
\[
  n f_n(c)\le\sum_{t=0}^{n-1}f_t(c).
\]
This avoids any parity restriction on \(n\) and proves the
upper tail bound.
The lower bound is Proposition~\ref{prop:lower}.
The final assertion follows by absorbing \(L_c^2/n\) into
a constant multiple of \(\log(2+|c|)/\sqrt n\) in the
stated range.
\end{proof}

\begin{remark}[Comparison with the published rates]
Here the square-root order concerns the time variable for a fixed
translation.  Leen \cite[Theorem~2.11]{LeenRates2026} proves the stronger
uniform estimate
\[
 f_n(c)\le B\sqrt{\frac{1+\log(1+|c|)}{n}}
 \qquad(c\in\Z,\ n\ge1)
\]
with an absolute constant $B$, and also proves limits for the
normalized tails of fixed translations.
The preceding proof records bounds obtained directly from the
return estimates in our chosen coordinates.
The stronger published estimate will be used explicitly in
Corollary~\ref{cor:prime-cuberoot}.
\end{remark}

\begin{corollary}\label{cor:exact}
\begin{enumerate}
  \item[(a)] For every fixed \(c\ne0\),
  \[
    0<\liminf_{n\to\infty}\sqrt n\,f_n(c)
      \le\limsup_{n\to\infty}\sqrt n\,f_n(c)<\infty.
  \]
  \item[(b)] For all \(j\ge1\) and \(n\ge(1+j)^2\),
  \(D_n^{(j)}\le C(1+j)/\sqrt n\), with an absolute constant.
  \item[(c)] The corresponding bounds hold for the \(3x-1\) map.
\end{enumerate}
\end{corollary}

\begin{proof}
Part (a) is Proposition~\ref{thm:exact}.
For (b), use the coupling in Proposition~\ref{thm:M} below;
its proof uses only the already established coalescence
construction.
For the representatives in that coupling,
\(L_c\le C'(1+j)\), so Proposition~\ref{thm:exact} gives the
claimed bound.
Part (c) follows from
\(f_n^-(c)=f_n^+(-c)\), by the reflection in
Corollary~\ref{cor:3xm1}.
\end{proof}

For an odd modulus $d$, write
\[
 F_a^{[d]}(x)=\#\{N\le x:N\in G,\ N\equiv a\pmod d\}.
\]

The following estimate is a consequence of uniform counted coalescence
and the Chinese remainder theorem. Leen's stronger displacement-uniform
tail estimate is available in \cite[Theorem~2.11]{LeenRates2026}.
The weighted averaging proposition in
\cite[Proposition~1]{LyapunovWeaver2026} already gives the corresponding
individually normalized comparison for odd moduli, for bounded labels
invariant above a threshold. We include the counting argument to record
the finite-defect boundary and the constants used in the generalized
case.

\begin{proposition}[equidistribution with a rate]\label{thm:rate}
There is an absolute constant \(C\) such that, for every
finite-defect \(T\)-invariant set \(G\), with cutoff \(K\)
as in Lemma~\ref{lem:D}, every odd \(d\ge1\), every \(a\),
and every \(x\) satisfying
\[
  (1+\log(2+d))^2\le\tfrac14\log_2x,
\]
one has
\[
  \left|\frac{dF_a^{[d]}(x)}x-\frac{F(x)}x\right|
  \le C\frac{\log(2+d)}{\sqrt{\log x}}
       +\frac{d(K+5)}{\sqrt x}.
\]
In particular,
\[
  \frac{3^jF_{a,j}(x)}x-\frac{F(x)}x
  =O\left(\frac{1+j}{\sqrt{\log x}}
           +\frac{3^j(K+5)}{\sqrt x}\right)
\]
uniformly when
\((1+(j+1)\log3)^2\le\tfrac14\log_2x\).
The analogous statements hold for finite-defect invariant sets
of the \(3x-1\) map, including its cycle basins.
\end{proposition}

\begin{proof}
The case \(d=1\) is immediate.
Otherwise apply \eqref{eq:Doddavg}, whose counting proof is
independent of the present rate estimate, with
\[
  n=\left\lfloor\tfrac12\log_2x\right\rfloor,\qquad
  L=1+\log(2+d).
\]
The hypothesis ensures \(n\ge\tfrac14\log_2x\ge L^2\), and
\(2^n\le\sqrt x\).
For \(0<|c|\le d/2\), Proposition~\ref{thm:exact} gives
\[
  f_n(c)\le C\left(\frac L{\sqrt n}+\frac{L^2}n\right)
          \le\frac{2CL}{\sqrt n}.
\]
The remaining terms of \eqref{eq:Doddavg}, divided by \(x\),
are at most \(d(K+5)/\sqrt x\).
Absorb the fixed conversions between logarithms into \(C\).
For \(d=3^j\), use \(\log(2+3^j)\le(j+1)\log3\).

For the \(3x-1\) map, the same counting proof applies:
on positive integers it still satisfies \(T_-(N)\ge N/2\),
and its nonabsorption probabilities are obtained by reflection.
\end{proof}

Proposition~\ref{thm:rate} supplies an explicit rate for the
equidistribution in Problem~\ref{prob:rate}.
Proposition~\ref{prop:lower} limits the displayed
coalescence-based upper bound at a depth \(n\asymp\log x\)
to the scale \((\log x)^{-1/2}\).
It does not rule out a better estimate for the actual
equidistribution error by another argument.

\subsection{Verification of conditions (M) and (M$^*$)}

\begin{proposition}[mixing from coalescence]\label{thm:M}
For every \(j\ge1\) and \(n\ge1\),
\[
  D_n^{(j)}
  \le\max_{0<|c|\le3^j/2}f_n(c).
\]
In particular, condition (M) holds.
Moreover, \(D_n^{(j)}\le\varepsilon\) whenever
\[
  n\ge L(\varepsilon)\log(2+3^j)+n_0(\varepsilon).
\]
\end{proposition}

\begin{proof}
The distance for two identical cylinders is zero.
For distinct level-\(j\) cylinders
\(C=a+3^j\Z_3\) and \(C'=a'+3^j\Z_3\), choose an integer
\[
  c\equiv a'-a\pmod{3^j},\qquad 0<|c|\le3^j/2.
\]
Thus \(C'=C+c\).
If \(W\) is a uniform word of length \(n\), independent of
\(Z\sim\nu_C\), the joint law of
\((q_W,\gamma_W(Z))\) is
\(\sum_q\delta_q\otimes\mu_C^{n,q}\).
Because the slices indexed by \(q\) are disjoint, the TV
distance between the two joint laws is
\(\sum_q\|\mu_C^{n,q}-\mu_{C'}^{n,q}\|_{\TV}\).

Take Haar-distributed \(\beta\in\Z_2\), independent of \(Z\).
Let \(W,W'\) be the length-\(n\) parity words of
\(\beta,\beta+c\), and put \(Z'=Z+c\).
Both words are uniform, \(Z'\sim\nu_{C'}\), and \(Z'\) is
independent of \(W'\).
Thus these variables have the required marginal laws.
On absorption of the relation chain by time \(n\),
\[
  q_{W'}=q_W,\qquad
  \gamma_{W'}(y+c)=\gamma_W(y)
\]
as rational affine maps.
The two joint variables therefore coincide.
The coupling inequality bounds their TV distance by \(f_n(c)\).
Take the maximum over the cylinder pairs.
For fixed \(j\) this is a finite maximum tending to zero
by Proposition~\ref{thm:K}; the quantitative assertion follows
from Proposition~\ref{prop:U}.
\end{proof}

\begin{corollary}[condition (M$^*$)]\label{cor:Mstar}
For \(r\ge0\), \(K\ge1\), and \(W\ge1\),
\[
  D_W^{(r,K)}
  \le\max_{u\in(\Z/3^K\Z)^*}
       \ \min_{\substack{c\in\Z\\
                    c\equiv3^ru\ (3^{r+K})}} f_W(c).
\]
Consequently \(a_{r,K}=0\) for all \(r,K\), and condition
(M$^*$) holds with \(\kappa_*=1\).
\end{corollary}

\begin{proof}
For \(C=3^{r+K}\Z_3\) and
\(C'=3^ru+3^{r+K}\Z_3\), the calculation in
Lemma~\ref{lem:Drk}(d) identifies the sum of the slice TV
distances with \(d_W^{(r,K)}(2u)\).
The factor \(2u\) is independent of the number of odd steps.
The coupling in Proposition~\ref{thm:M} works for every integer
\(c\equiv3^ru\pmod{3^{r+K}}\).
Taking the minimum and then the maximum proves the inequality,
since multiplication by \(2\) permutes the units.
The minimum exists because \(f_W(c)\) takes values in the
finite set \(2^{-W}\Z\cap[0,1]\).
The vanishing of all \(a_{r,K}\) also follows directly from
condition (M) and Lemma~\ref{lem:Drk}(d).
\end{proof}

\begin{remark}
The limiting constants in Proposition~\ref{thm:quant} and
Corollary~\ref{cor:uncond} are therefore zero.
The finite values \(D_{24}^{(1)}\) and \(D_{20}^{(1,1)}\)
remain useful computational benchmarks, but are not
positive limiting obstructions.
Proposition~\ref{thm:exact} and the coupling give
\(D_n^{(j)}\le C(1+j)/\sqrt n\) for \(n\ge(1+j)^2\),
as stated in Corollary~\ref{cor:exact}.
\end{remark}

\subsection{Effective equidistribution at deep \(3\)-adic levels}

For a finite-defect invariant set \(G\), put
\[
  K=\max(G\triangle T^{-1}G),
\]
with \(K=0\) when the symmetric difference is empty.
Then \(\mathbf1_G(N)=\mathbf1_G(TN)\) for \(N>K\).

\begin{lemma}[effective equidistribution]\label{lem:D}
For every \(j\ge0\), \(a,c\in\Z\), \(n\ge1\), and \(x\ge1\),
\[
  |3^jF_{a+c,j}(x)-3^jF_{a,j}(x)|
  \le f_n(c)x+3^j2^n+2^nK+3|c|+3^{j+1}.
\]
\end{lemma}

\begin{proof}
Substitute \(M=N+c\) in \(F_{a+c,j}(x)\).
On the intersection of the intervals \([1,x]\) and
\([1-c,x-c]\), compare
\(\mathbf1_G(N+c)\) and \(\mathbf1_G(N)\) for
\(N\equiv a\pmod{3^j}\).
Their interval symmetric difference contains at most
\(2(|c|/3^j+1)\) integers in this residue class.

If both \(N\) and \(N+c\) exceed \(2^nK\), all their first
\(n\) iterates before the endpoint exceed \(K\), since
\(T^i(N)\ge N/2^i\).
Thus merging by step \(n\) forces equal membership in \(G\).
The possible small starts are contained in
\([1,2^nK+|c|]\), contributing at most
\((2^nK+|c|)/3^j+1\) integers in the class.

The remaining exceptions belong to the periodic set of
nonabsorbing starting residues modulo \(2^n\).
It has exactly \(2^nf_n(c)\) residues.
By the Chinese remainder theorem, in a fixed class modulo
\(3^j\) it contributes at most
\[
  f_n(c)x/3^j+2^n
\]
integers up to \(x\).
Adding these bounds and multiplying by \(3^j\) proves the
claim.
\end{proof}

Let
\[
  \bar f_n(q)=3^{-q}\sum_{\substack{c\in\Z\\|c|\le3^q/2}}f_n(c),
  \qquad f_n(0)=0.
\]
The sum uses one balanced representative of each residue class
modulo \(3^q\).
Since \(F(y)=\sum_cF_{a+c,q}(y)\) and the mean of \(|c|\)
over these representatives is at most \(3^q/4\), averaging
Lemma~\ref{lem:D} gives
\begin{equation}\label{eq:Davg}
  |3^qF_{a,q}(y)-F(y)|
  \le\bar f_n(q)y+3^q2^n+2^nK+4\cdot3^q.
\end{equation}

\begin{corollary}[deep \(3\)-adic equidistribution]\label{cor:deep}
Let \(G\) be finite-defect invariant and \(0<\varepsilon<1\).
Put
\[
  \kappa(\varepsilon)
   =\bigl(2\log3\,(1+L(\varepsilon)\log2)\bigr)^{-1}.
\]
Then
\[
  \limsup_{x\to\infty}
  \max_{0\le q\le\kappa(\varepsilon)\log x}
  \max_{a\bmod3^q}
  \left|\frac{3^qF_{a,q}(x)}x-\frac{F(x)}x\right|
  \le\varepsilon.
\]
\end{corollary}

\begin{proof}
Increase \(n_0(\varepsilon)\) if necessary so that it is
nonnegative, and take the integer
\[
  n_q=\left\lceil
        L(\varepsilon)\log(2+3^q)+n_0(\varepsilon)
      \right\rceil.
\]
Proposition~\ref{prop:U} gives \(\bar f_{n_q}(q)\le\varepsilon\).
Since \(\log(2+3^q)\le(q+1)\log3\),
\[
  3^q2^{n_q}
  \le C_\varepsilon
        \exp\bigl(q\log3\,(1+L(\varepsilon)\log2)\bigr)
  \le C_\varepsilon x^{1/2}
\]
throughout the stated range of \(q\).
The remaining terms in \eqref{eq:Davg}, divided by \(x\),
are \(O_\varepsilon((K+1)x^{-1/2})\), uniformly in \(q,a\).
This proves the claim.
\end{proof}

\subsection{Terras renewal in \(m\) steps}

Throughout the remainder of this section, \(G\triangle T^{-1}(G)\)
is finite and
\[
 K=\max\bigl(G\triangle T^{-1}(G)\bigr),
 \qquad \max\varnothing:=0.
\]
Write \(\psi=\psi_G\).  For an integer \(m\ge1\), let
\[
 S_m=Q_m\log3-m\log2,\qquad Q_m\sim\operatorname{Bin}(m,\tfrac12).
\]
Different copies used in sums of these random variables will be
independent.

\begin{lemma}[Renewal in \(m\) steps]\label{lem:R}
Let \(m\ge1\) and \(q^*\in\{0,\ldots,m\}\), and let
\(n_0,\ldots,n_{q^*}\ge1\) be integers.
Put \(n^*=\max_{0\le q\le q^*}n_q\).
For \(X\ge1\) and \(U=\log X\),
\[
 \left|\psi(U)-\mathbb E\psi(U+S_m)\right|
 \le
 \max_{0\le q\le q^*}\bar f_{n_q}(q)
 +\mathbb P(Q_m>q^*)
 +\frac{2^{m+n^*+1}(K+10)}{X}.
\]
\end{lemma}

\begin{proof}
For \(0\le r<2^m\), put \(q_r=q_m(r)\), extending \(T(0)=0\).
Terras conjugacy gives
\[
 T^m(2^mY+r)=3^{q_r}Y+T^m(r),\qquad Y\ge0.
\]
Exactly \(\binom mq\) residues \(r\) have \(q_r=q\).
Write
\[
 T^m(r)=3^{q_r}r/2^m+B_m(r).
\]
The normalized offset \(b=B/3^q\), initially zero, changes to
\(b/2\) at an even step and to
\[
 b/2+\frac{1}{2\cdot3^{q+1}}
\]
at an odd step, where \(q\) is the preceding odd-step count.
Induction therefore gives \(0\le B_m(r)/3^{q_r}\le1/3\).

For a fixed \(r\), the images of its positive representatives at most
\(X\) form an interval in the progression
\[
 a_r=T^m(r)\pmod{3^{q_r}}.
\]
Compare this interval with the positive members of that progression
at most \(y_r=3^{q_r}X/2^m\).
The lower endpoint, divided by \(3^{q_r}\), is less than \(4/3\);
the upper endpoint differs from \(y_r\) by less than \(3^{q_r}\).
Thus only a bounded number of progression members are lost at the
endpoints.  Empty intervals and the excluded starting value zero
are covered by the same bound.
Furthermore, all iterates of a starting value \(N>2^mK\) through
time \(m\) exceed \(K\), because \(T(N)\ge N/2\).
Hence \(N\in G\) if and only if \(T^m(N)\in G\) for these starting
values.  There are at most \(K+1\) remaining starting values in each
class modulo \(2^m\).
Allowing six endpoint and initial-value terms per class gives
\[
 F(X)=\sum_{r<2^m}F_{a_r,q_r}(y_r)+E_1,
 \qquad |E_1|\le2^m(K+6).
\]

For \(q_r\le q^*\), equation~\eqref{eq:Davg} gives
\[
 F_{a_r,q_r}(y_r)
 =3^{-q_r}F(y_r)+e_r,\qquad
 |e_r|\le
 \bar f_{n_{q_r}}(q_r)\frac{X}{2^m}
 +2^{n^*}(K+1)+4.
\]
Also
\[
 3^{-q_r}F(y_r)
 =\frac{X}{2^m}\psi(U+q_r\log3-m\log2).
\]
For \(q_r>q^*\), the two expressions being compared lie respectively
in \([0,X/2^m+1]\) and \([0,X/2^m]\), so their difference has absolute
value at most \(X/2^m+1\).
Summing over \(r\), dividing by \(X\), and using
\[
 2^{n^*}(K+1)+K+11\le2^{n^*+1}(K+10)
\]
proves the result.
\end{proof}

\subsection{A local limit bound and a mixture of times}

\begin{lemma}[Local bound]\label{lem:local}
Let \(Q\sim\operatorname{Bin}(m,\tfrac12)\), \(m\ge1\), and
\(t=q-m/2\).  If \(|t|\le m/4\), then
\[
 \mathbb P(Q=q)\le
 \frac{\kappa_0}{\sqrt m}e^{-2t^2/m},
 \qquad
 \kappa_0=
 \frac{e^{1/12}\sqrt{16/3}}{\sqrt{2\pi}}<1.0015.
\]
For every \(q\in\{0,\ldots,m\}\),
\(\mathbb P(Q=q)\le e^{-2t^2/m}\).
\end{lemma}

\begin{proof}
For \(0<q<m\), Robbins' factorial bounds~\cite{Robbins1955} imply
\[
 2^{-m}\binom mq
 \le\frac{e^{1/(12m)}}{\sqrt{2\pi}}
 \sqrt{\frac{m}{q(m-q)}}
 \exp\!\left(-m\bigl(\log2-h_{\rm e}(q/m)\bigr)\right),
\]
where
\(h_{\rm e}(x)=-x\log x-(1-x)\log(1-x)\); thus
\(h_{\rm e}(x)=(\log2)H(x)\) in the entropy notation used earlier.
The inequalities
\[
 \log2-h_{\rm e}(x)\ge2(x-\tfrac12)^2,\qquad
 q(m-q)\ge3m^2/16
\]
give the first assertion.  Its stated range excludes \(q=0,m\).
For \(t>0\), the point probability is bounded by the upper tail
\(\mathbb P(Q-m/2\ge t)\), to which Hoeffding's inequality
applies~\cite{Hoeffding1963}.  The lower tail gives the case \(t<0\);
for \(t=0\), the second assertion is simply a probability bound by one.
\end{proof}

Put
\[
 \mu_0=\log2-\tfrac12\log3>0,\qquad
 \nu_M=\frac1{M+1}\sum_{m=M}^{2M}\mathcal L(S_m)
 \quad(M\ge1).
\]
Thus \(\mathbb ES_m=-\mu_0 m\).

\begin{lemma}[Local mass of the mixture]\label{lem:mix}
For every interval \(J\subset\mathbb R\) of length at most one,
\[
 \nu_M(J)\le\frac{13.56+8M^{-1/2}}{M}+e^{-M/8}.
\]
\end{lemma}

\begin{proof}
The possible values of \(S_m\) are spaced by \(\log3>1\), so \(J\)
contains at most one of them.  If \(y\) is its left endpoint, its
corresponding centered binomial index belongs to
\[
 [\tau_m,\tau_m+0.911],\qquad
 \tau_m=\frac{y+m\mu_0}{\log3},
\]
since \(1/\log3<0.911\).
The \(\tau_m\) form an arithmetic progression of step
\(\delta_0=\mu_0/\log3\).
For \(M\le m\le2M\), Lemma~\ref{lem:local} gives
\[
 \mathbb P(S_m\in J)
 \le\kappa_0M^{-1/2}\phi(\tau_m)+e^{-M/8},
 \qquad
 \phi(\tau)=\sup_{t\in[\tau,\tau+0.911]}e^{-t^2/M}.
\]
Indeed, the first bound handles \(|t|\le m/4\), and the second
bound in that lemma is at most \(e^{-m/8}\) otherwise.
The function \(\phi\) is unimodal, has maximum one, and has integral
\(0.911+\sqrt{\pi M}\).
For a nonnegative unimodal function, the sum over an arithmetic
progression of step \(\delta_0\) is at most its maximum plus
\(\delta_0^{-1}\) times its integral.  This follows by comparing the
increasing and decreasing portions separately with adjacent
rectangles, leaving at most one maximum term.
Consequently
\[
 \nu_M(J)
 \le \frac1M\left(
 \frac{\kappa_0\sqrt\pi}{\delta_0}
 +\kappa_0(1+0.911/\delta_0)M^{-1/2}\right)+e^{-M/8}.
\]
Here \(\kappa_0\sqrt\pi/\delta_0<13.56\) and
\(\kappa_0(1+0.911/\delta_0)<8\), as required.
\end{proof}

\begin{lemma}[A one-sided averaging inequality]\label{lem:mono}
For every \(S\subseteq\mathbb N\), \(0<h\le1\), and \(V\in\mathbb R\),
\[
 1-\psi_S(V)
 \le\frac1h\int_{V-h}^V(1-\psi_S(v))\,dv+h.
\]
\end{lemma}

\begin{proof}
For \(0\le t\le h\), monotonicity of the counting function gives
\[
 \psi_S(V)\ge e^{-t}\psi_S(V-t)
 \ge\psi_S(V-t)-t,
\]
because \(0\le\psi_S\le1\).  Average this inequality over \(t\in[0,h]\)
and enlarge the resulting error \(h/2\) to \(h\).
\end{proof}

This averaging inequality permits comparison of the discrete
mixture with Lebesgue integrals without any equidistribution
assumption on its support.

\subsection{From logarithmic to natural density}

\begin{proposition}[A quantitative density comparison]\label{thm:natural}
Let \(G\) be \(T\)-invariant and let \(\overline\delta\) be the upper
logarithmic density of \(\mathbb N\setminus G\).
For every \(0<\varepsilon<1\) and \(0<h\le1\),
\[
 \underline d(G)\ge
 1-\frac{28L(\varepsilon)\overline\delta}{h}-h-\varepsilon.
\]
In particular, the density-zero and density-one comparison can be
expressed using
\[
 \Phi(\delta)=
 \inf_{0<\varepsilon<1}
 \left(\frac{28L(\varepsilon)\delta}{\varepsilon}
                 +2\varepsilon\right),
 \qquad \lim_{\delta\downarrow0}\Phi(\delta)=0.
\]
\end{proposition}

\begin{proof}
Set \(g=1-\psi_G\).  For a counting function \(F_S\), partial summation
gives
\[
 \sum_{\substack{N\le e^V\\N\in S}}\frac1N
 =\psi_S(V)+\int_0^V\psi_S(u)\,du.
\]
Moreover \(g(u)=\psi_{\mathbb N\setminus G}(u)+O(e^{-u})\) for
\(u\ge0\).  Hence
\begin{equation}\label{eq:logint}
 \int_0^Vg(u)\,du
 \le V\bigl(\overline\delta+o(1)\bigr)
 \qquad(V\to\infty).
\end{equation}
Fix \(\varepsilon,h\), write \(L=L(\varepsilon)\), and choose
\[
 n_q=\left\lceil L\log(2+3^q)\right\rceil+n_0(\varepsilon),
 \quad M=\left\lfloor\frac{U}{2L}\right\rfloor,
 \quad q_m^*=\lfloor0.6m\rfloor .
\]
Increase the integer \(n_0(\varepsilon)\), if necessary, so every
\(n_q\ge1\).  Proposition~\ref{prop:U} gives
\(\bar f_{n_q}(q)\le\varepsilon\).
Uniformly for \(M\le m\le2M\le U/L\),
\[
 (m+\max_{q\le q_m^*}n_q)\log2
 \le U\left(\frac{\log2}{L}+0.6\log2\log3\right)+C_\varepsilon
 \le0.47U+C_\varepsilon,
\]
where \(L\ge76\) was used.
Also \(\mathbb P(Q_m>q_m^*)\le e^{-0.02m}\).
Apply Lemma~\ref{lem:R} with \(X=e^U\), and average over these \(m\).
All its remaining errors tend to zero uniformly in this range, giving
\[
 g(U)\le\int g(U+v)\,d\nu_M(v)+\varepsilon+o_{\varepsilon,K}(1).
\]
The arguments \(U+S_m\) lie between
\(U(1-\log2/L)\) and \(U'=U(1+0.41/L)\).
For sufficiently large \(U\), their intervals of length \(h\) to
the left lie in \([0,U']\).
Lemma~\ref{lem:mono} and Tonelli's theorem therefore give
\[
 \int g(U+v)\,d\nu_M(v)
 \le
 \frac1h\int_0^{U'}g(w)\,
       \nu_M([w-U,w-U+h])\,dw+h.
\]
Apply Lemma~\ref{lem:mix} and \eqref{eq:logint}.
Since \(U'/M\to2L+0.82\) and \(U'e^{-M/8}\to0\),
\[
 \limsup_{U\to\infty}g(U)
 \le\frac{13.56(2L+0.82)}h\overline\delta+h+\varepsilon
 \le\frac{28L}h\overline\delta+h+\varepsilon.
\]
This proves the first assertion.
Taking \(h=\varepsilon\) proves the formula for \(\Phi\).
For each fixed \(\varepsilon\), its first term tends to zero as
\(\delta\downarrow0\), and then \(\varepsilon\downarrow0\) proves
\(\Phi(\delta)\to0\).
The complement of a \(T\)-invariant set is also \(T\)-invariant.
\end{proof}

\begin{corollary}[An application of Tao's logarithmic-density theorem]
\label{thm:taonatural}
For every function \(f:\mathbb N\to\mathbb R\) tending to infinity,
\[
 \#\{N\le x:\operatorname{Col}_{\min}(N)\ge f(N)\}=o(x).
\]
\end{corollary}

\begin{proof}
Put \(G_C=\{N:\operatorname{Col}_{\min}(N)\le C\}\), \(C\ge1\).
The identity
\[
 \operatorname{Col}_{\min}(N)
 =\min\bigl(N,\operatorname{Col}_{\min}(T(N))\bigr)
\]
implies
\(G_C\triangle T^{-1}(G_C)\subseteq[1,C]\).
Let \(\delta_C\) denote the upper logarithmic density of its complement.
The fixed-threshold formulation of Tao's theorem
\cite[Remark~1.4 and Theorem~1.3]{Tao2022} says that
for every \(\delta>0\) there is \(C_\delta\) with
\(\delta_{C_\delta}\le\delta\).
Since \(\delta_C\) decreases with \(C\), it follows that
\(\delta_C\to0\).
Proposition~\ref{thm:natural} gives
\[
 \overline d(\mathbb N\setminus G_C)\le\Phi(\delta_C)\longrightarrow0.
\]
For a fixed \(C\), eventually \(f(N)>C\); thus
\(\{N:\operatorname{Col}_{\min}(N)\ge f(N)\}\) is contained in
\(\mathbb N\setminus G_C\) apart from finitely many integers.
Take upper natural densities and then let \(C\to\infty\).
\end{proof}

\begin{remark}
The displayed comparison function is explicit once the coalescence
parameters are fixed.  This is distinct from an onset estimate for
an arbitrary invariant set.
The finite counting errors depend on \(K\), and the use of
\eqref{eq:logint} depends on the convergence of the particular
logarithmic averages of \(G\).
Consequently, effectivity of \(n_0(\varepsilon)\) alone does not give
a common finite-\(x\) threshold for all finite-defect invariant sets.
The window comparison below also uses qualitative slow oscillation.
\end{remark}

\subsection{Windows and the equivalence of densities}

\begin{proposition}[Window comparison]\label{thm:W}
Let \(G\) be \(T\)-invariant and \(0<\varepsilon<1\).
If \(0<\gamma\le\mu_0/(2L(\varepsilon))\), then
\[
 \limsup_{U\to\infty}\sup_{0\le D\le\gamma U}
       |\psi_G(U)-\psi_G(U-D)|
 \le2\varepsilon+2.11\sqrt{L(\varepsilon)\gamma}.
\]
In particular,
\[
 \lim_{\gamma\downarrow0}\limsup_{x\to\infty}
 \sup_{x^{1-\gamma}\le y\le x}
 \left|\frac{F(y)}y-\frac{F(x)}x\right|=0.
\]
\end{proposition}

\begin{proof}
Write \(L=L(\varepsilon)\), \(U_1=U-D\), and choose
\[
 m=\left\lfloor\frac{U_1}{2L}\right\rfloor,\qquad
 k=\lfloor D/\mu_0\rfloor .
\]
Uniformly for \(0\le D\le\gamma U\), \(m\to\infty\), and
\[
 m+k\le\frac{U}{2L}+\frac{D}{\mu_0}\le U/L.
\]
The parameters and exponential budget in the preceding proof apply
to time \(m+k\) at \(U\), and to time \(m\) at \(U_1\).
Lemma~\ref{lem:R} consequently gives
\[
 \psi(U)=\mathbb E\psi(U+S_{m+k})+R_1,\qquad
 \psi(U_1)=\mathbb E\psi(U_1+S_m)+R_2,
\]
where \(|R_i|\le\varepsilon+o(1)\), uniformly in \(D\).
In distribution \(S_{m+k}=S_m+S'_k\), independently, with
\[
 S_m=Q\log3-m\log2,\qquad
 S'_k=Q'\log3-k\log2.
\]
Set \(Z=Q'-\lfloor k/2\rfloor\).  Then
\[
 U+S_{m+k}
 =U_1+(Q+Z)\log3-m\log2+r,
\]
where
\[
 r=D-k\mu_0+
       (\lfloor k/2\rfloor-k/2)\log3,\qquad |r|\le\log2.
\]
All the arguments before and after this bounded shift tend to
infinity uniformly in \(D,Q,Q'\), since
\(U+S_{m+k}\ge U(1-\log2/L)\).
Slow oscillation, supplied by Corollary~\ref{thm:liouville}(c),
therefore removes \(r\) at a uniform cost \(o(1)\).

For an integer-valued unimodal probability distribution, the TV
distance from its translate by one equals its largest atom.
The triangle inequality therefore bounds the distance from its
translate by \(z\) by \(|z|\) times that atom.
The central binomial coefficient bound gives
\[
 \max_a\mathbb P(Q=a)\le\sqrt{\frac{2}{\pi m}}.
\]
(This is the usual Wallis bound; the odd case follows from the
adjacent even central coefficient.)
Since \(0\le\psi\le1\), convexity of TV and
\(\mathbb E|Z|\le\sqrt{k}/2+1/2\) yield
\[
 |\mathbb E\psi(U_1+(Q+Z)\log3-m\log2)
       -\mathbb E\psi(U_1+Q\log3-m\log2)|
 \le(\sqrt{k}/2+1/2)\sqrt{\frac2{\pi m}}.
\]
Uniformly in \(D\),
\[
 \frac{k}{m}\le
 \frac{2L\gamma}{\mu_0(1-\gamma)}+o(1).
\]
Thus the last bound has upper limit at most
\[
 \sqrt{\frac{L\gamma}{\pi\mu_0(1-\gamma)}}
 \le\sqrt{\frac{2}{\pi\mu_0}}\sqrt{L\gamma}
 <2.11\sqrt{L\gamma}.
\]
Here \(\gamma<1/2\), since \(L\ge76\).
Adding the two renewal errors proves the first assertion.
For the second, fix \(\varepsilon\), let \(\gamma\downarrow0\) within
its permitted range, and then let \(\varepsilon\downarrow0\).
\end{proof}

\begin{proposition}[Equivalence of densities]\label{thm:T}
A \(T\)-invariant set has a natural density if and only if it has a
logarithmic density, and the two densities agree.
\end{proposition}

\begin{proof}
Natural density implies logarithmic density by partial summation.
Conversely, suppose that the logarithmic density is \(d\).
The partial-summation identity above implies
\[
 A(V):=\frac1V\int_0^V\psi(u)\,du\longrightarrow d.
\]
For every fixed \(\lambda>1\), the average over \([U,\lambda U]\) is
\[
 \frac{\lambda A(\lambda U)-A(U)}{\lambda-1}\longrightarrow d.
\]
If \(U'\in[U,\lambda U]\), then
\(U\ge(1-\gamma)U'\), where \(\gamma=1-1/\lambda\).
Proposition~\ref{thm:W} shows that the difference between this average
and \(\psi(U)\) has upper limit at most \(\omega(\gamma)\), for a
function \(\omega(\gamma)\to0\).
Let \(U\to\infty\), and then \(\lambda\downarrow1\).
It follows that \(\psi(U)\to d\).
\end{proof}

\begin{remark}[Existence of densities in the standard case]
Leen \cite[Theorem~2]{LeenBrownian2026} proves natural-density existence
for every set invariant above a threshold for the ordinary Collatz map
$C(N)=N/2$ for even $N$ and $C(N)=3N+1$ for odd $N$.
This includes our finite-defect invariant sets.
Indeed, for odd $N>K$, invariance first gives
$\mathbf1_G(N)=\mathbf1_G(TN)$, and even-step invariance at
$C(N)=2T(N)>K$ gives
$\mathbf1_G(C(N))=\mathbf1_G(TN)$.
The even case is immediate.
Thus density existence for these standard-map sets is already covered
by that theorem.
The window argument above is retained because it also supplies the
uniform comparison used for growing moduli below.
\end{remark}

\begin{corollary}[A fixed-threshold consequence]\label{cor:Q}
There are absolute constants \(A,C_1,C_0>0\) such that, for
\(C\ge C_0\),
\[
 \underline d\{N:\operatorname{Col}_{\min}(N)\le C\}
 \ge1-C_1\bigl((\log C)^{-1/A}\log\log C\bigr)^{1/4}.
\]
\end{corollary}

\begin{proof}
Tao's quantitative fixed-threshold estimate
\cite[Remark~1.4 and Theorem~3.1]{Tao2022} gives constants \(A,C_2>0\)
such that
\[
 C_\delta\le\exp(C_2\delta^{-A}),\qquad \delta_{C_\delta}\le\delta .
\]
For large \(C\), set \(\delta=(C_2/\log C)^{1/A}\).
Then \(\delta_C\le\delta\).
From the specified formula for \(L\),
\(L(\varepsilon)\ll\varepsilon^{-2}\log(10/\varepsilon)\).
In the definition of \(\Phi\), take
\(\varepsilon=(\delta\log(2/\delta))^{1/4}\).
For small \(\delta\), this gives
\[
 \Phi(\delta)\ll(\delta\log(2/\delta))^{1/4}.
\]
Proposition~\ref{thm:natural} and
\(\log(2/\delta)\ll_A\log\log C\) prove the assertion.
\end{proof}

\begin{remark}\label{rem:Qmazur}
This is a consequence of the invariant-set comparison and Tao's
quantitative input.  The constants include the implicit constants
and exponent in that input; they are not numerically specified here.
Stronger quantitative results for almost-bounded orbits are available
in the preceding literature, including~\cite{Mazur2026,Nashida2026}.
The purpose of retaining this deduction is to record the dependence
needed in Corollary~\ref{cor:Qprimes}.
\end{remark}

\begin{corollary}[Basins for the \(3x-1\) map]\label{cor:basins}
Every cycle basin for the positive \(3x-1\) map has a natural density
if and only if it has a logarithmic density, and the two then agree.
\end{corollary}

\begin{proof}
Such a basin is exactly invariant under \(T_-\).
Corollary~\ref{cor:3xm1} supplies the corresponding coalescence
estimate and slow oscillation.
In the proof of Lemma~\ref{lem:R}, an odd step changes the normalized
offset to \(b/2-1/(2\cdot3^{q+1})\), so \(-1/3\le b\le0\).
The endpoint comparison loses at most one additional progression
member, already allowed by the bound \(2^m(K+6)\).
Also \(T_-(N)\ge N/2\) for all positive \(N\).
Thus Lemma~\ref{lem:R} and the proofs of
Propositions~\ref{thm:W} and~\ref{thm:T} apply with these sign changes.
\end{proof}

In particular, the corollary applies to the basins of the cycles
with minima \(1\), \(5\), and \(17\).
Its proof applies separately to any basin and requires no
classification or exhaustion of the positive integer orbits.
Sultanow \cite[Theorem~1]{Sultanow2026} proves existence of logarithmic
density for completely invariant sets of the $3x-1$ map, in particular
for these basins.  Combining that input with Corollary~\ref{cor:basins}
therefore also gives their natural densities.

\subsection{A scalar identity for the density defect}\label{subsec:defect}

We record a further consequence of the renewal equation.
It expresses the density through one scalar integral and a constant
from the Dirichlet renewal equation.
The density-existence result cited above implies convergence for the
standard map; the identity also isolates the analytic criterion useful
in related systems.  Boundedness of the integral alone is not a
convergence argument.
Write \(K_1=\#(G\triangle T^{-1}(G))\), \(F_2=F_{2,1}\), and
\[
 \mathcal I_G(y)=\int_1^y\frac{F(x)}{x^2}\,dx,\qquad
 \Delta(y)=\int_1^y\frac{F_2(x)-F(x)/3}{x^2}\,dx.
\]
Both are taken to be zero when \(0<y<1\), consistently with the
counting functions vanishing below one.
For \(V=\log y\), put
\[
 \begin{split}
 W(y)
 &=\tfrac12\bigl(\mathcal I_G(y)-\mathcal I_G(y/2)\bigr)
   -\tfrac12\bigl(\mathcal I_G(3y/2)-\mathcal I_G(y)\bigr)\\
 &=\tfrac12\int_{V-\log2}^{V}\psi(u)\,du
   -\tfrac12\int_V^{V+\log(3/2)}\psi(u)\,du.
 \end{split}
\]
The notation \(\mathcal I_G\) avoids reusing the transfer function
\(\Phi\).  Also define
\[
 \delta_G(u)=e^{-u}\bigl(F_2(e^u)-F(e^u)/3\bigr)
 =\tfrac13\bigl(\psi_{2,1}(u)-\psi(u)\bigr),
 \qquad \Delta(e^V)=\int_0^V\delta_G(u)\,du.
\]

\begin{proposition}[Defect identity]\label{thm:defect}
There is a real constant \(c_G\), with \(|c_G|\le K_1+3\), such that
for every \(y\ge1\),
\begin{equation}\label{eq:defect}
 \left|W(y)-\tfrac32\Delta(3y/2)-c_G\right|
 \le\frac{K_1+2}{y}.
\end{equation}
Consequently:
\begin{enumerate}
 \item[(a)] \(|\Delta(y)|\le2K_1+4\) for all \(y\ge1\).
 \item[(b)] As \(y\to\infty\),
 \[
  \psi_G(\log y)=
  \frac{3\Delta(y)+2c_G}{\log(4/3)}+o(1).
 \]
 \item[(c)] The following are equivalent: existence of natural density;
 existence of logarithmic density; convergence of
 \[
  \int_1^\infty\frac{F_2(x)-F(x)/3}{x^2}\,dx .
 \]
 When they hold,
 \(d(G)=(2c_G+3\Delta(\infty))/\log(4/3)\).
 \item[(d)]
 \[
 \overline d(G)-\underline d(G)
 =\frac3{\log(4/3)}
 \left(\limsup_{y\to\infty}\Delta(y)
             -\liminf_{y\to\infty}\Delta(y)\right).
 \]
\end{enumerate}
\end{proposition}

\begin{proof}
Set \(D(x)=F_2(x)-F(x)/3\).
The endpoint estimate following Lemma~\ref{lem:renewal} gives
\[
 F(x)-F(x/2)-\tfrac13F(3x/2)=D(3x/2)+E(x),
 \qquad |E(x)|\le K_1+1.
\]
Divide by \(x^2\) and integrate over \([1,y]\).
Changes of variable yield
\[
 \begin{split}
 \int_1^y F(x/2)x^{-2}\,dx&=\tfrac12\mathcal I_G(y/2),\\
 \int_1^y \tfrac13F(3x/2)x^{-2}\,dx
 &=\tfrac12\bigl(\mathcal I_G(3y/2)-\mathcal I_G(3/2)\bigr),\\
 \int_1^y D(3x/2)x^{-2}\,dx
 &=\tfrac32\bigl(\Delta(3y/2)-\Delta(3/2)\bigr).
 \end{split}
\]
Since \(F_2=0\) on \([1,3/2]\),
\(\Delta(3/2)=-\mathcal I_G(3/2)/3\), so the two lower-endpoint
constants cancel.  Thus we may take
\[
 c_G=\int_1^\infty E(x)x^{-2}\,dx,
\]
and obtain the exact identity
\[
 W(y)-\tfrac32\Delta(3y/2)-c_G
 =-\int_y^\infty E(x)x^{-2}\,dx.
\]
In fact \(|c_G|\le K_1+1\) and the remainder is at most
\((K_1+1)/y\), which imply the stated bounds.
Since
\(-\frac12\log(3/2)\le W(y)\le\frac12\log2\), these estimates give
(a) for \(y\ge3/2\).  For \(1\le y\le3/2\), use
\(\Delta(y)=-\mathcal I_G(y)/3\) and
\(0\le\mathcal I_G(y)\le1/3\).

Slow oscillation implies
\[
 W(y)=\tfrac12\log(4/3)\,\psi(\log y)+o(1).
\]
Equidistribution at the fixed first level implies
\(\delta_G(u)\to0\), and hence
\[
 \Delta(3y/2)-\Delta(y)
 =\int_{\log y}^{\log y+\log(3/2)}\delta_G(u)\,du=o(1).
\]
Substitution into \eqref{eq:defect} proves (b).
Part (c) follows from (b) and Proposition~\ref{thm:T}; taking upper and
lower limits in (b) proves (d).
\end{proof}

\begin{remark}[The constant in the Dirichlet renewal equation]
The constant just obtained is also
\[
 c_G=P_G(1)+\mathcal H(1)
 =\sum_{N\in G\triangle T^{-1}(G)}
       \frac{\mathbf1_G(N)-\mathbf1_G(T(N))}{N}
   +\sum_{n\in G_2}\frac3{2n(2n-1)},
\]
with \(P_G,\mathcal H\) as in Proposition~\ref{thm:equiv}(a).
To verify this directly, write
\[
 E(x)=
 \sum_{N\le x}\bigl(\mathbf1_G(N)-\mathbf1_G(T(N))\bigr)
 +\sum_{n\in G_2}
       \mathbf1_{[(2n-1)/3,\,2n/3)}(x).
\]
The first sum is supported on the finite defect set.
The second sum is the exact difference between the odd-branch
cutoff and \(3x/2\).
Integrating against \(x^{-2}\) over \([1,\infty)\) gives the
displayed formula, since \(n\ge2\) and
\[
 \int_{(2n-1)/3}^{2n/3}x^{-2}\,dx=\frac3{2n(2n-1)}.
\]
The second numerical series converges absolutely.
\end{remark}

\begin{corollary}\label{cor:defect}
For every \(T\)-invariant \(G\),
\[
 \left|\int_1^X(3F_{2,1}(x)-F(x))x^{-2}\,dx\right|
 \le6K_1+12\qquad(X\ge1).
\]
If, for all sufficiently large \(x\),
\[
 |3F_{2,1}(x)-F(x)|\le x\,h(\log x),
\]
where \(h\ge0\) is measurable and
\(\int^\infty h(u)\,du<\infty\), then \(G\) has a natural density.
\end{corollary}

\begin{proof}
The first assertion is three times Proposition~\ref{thm:defect}(a).
The second hypothesis makes the defect integral absolutely
convergent after the substitution \(u=\log x\), so part (c) of that
proposition applies.
\end{proof}

\begin{remark}\label{rem:defectgen}
The pointwise \(O((\log x)^{-1/2})\) estimate supplied by
Proposition~\ref{thm:rate}, by itself, is not an integrable majorant on
the logarithmic scale.
For \(3x-1\), the corresponding defect uses the class \(1\bmod3\).
More generally, whenever the generalized renewal equation,
slow oscillation, and equidistribution in the image classes are
available for \(T_{\ell,m}\), the same integration gives
\[
 \begin{split}
 (\log\ell-\varpi\log m)\psi_G(\log y)
 &=\frac m\ell\sum_{j=1}^{\ell-1}\Delta_{b_j}(y)+c_G+o(1),\\
 \Delta_b(y)&=\int_1^y\frac{F_{b,1}(x)-F(x)/m}{x^2}\,dx.
 \end{split}
\]
where \(b_j\equiv r_j\ell^{-1}\pmod m\) and
\(\varpi=(\ell-1)/\ell\).
Each branch is included in the sum, even if image classes repeat.
The coefficient is positive in the negative-drift range.
\end{remark}

\subsection{Level of distribution}

For \(Q\ge1\) and \(a\in\mathbb Z\), write
\[
 F^{[Q]}_a(x)=\#\{N\le x:N\in G,\ N\equiv a\pmod Q\};
 \qquad F_{a,j}=F^{[3^j]}_a.
\]
The result below is a uniform estimate for each modulus in a
power range.  It is not an estimate for a sum of errors over moduli.

\begin{lemma}[Odd moduli]\label{lem:Dodd}
Let \(d\ge1\) be odd.  For all integers \(a,c\), \(n\ge1\), and \(x\ge1\),
\[
 \left|dF^{[d]}_{a+c}(x)-dF^{[d]}_a(x)\right|
 \le f_n(c)x+d2^n+2^nK+3|c|+3d.
\]
Consequently,
\begin{equation}\label{eq:Doddavg}
 \left|dF^{[d]}_a(y)-F(y)\right|
 \le y\max_{0<|c|\le d/2}f_n(c)+d2^n(K+1)+4d
 \qquad(y\ge1),
\end{equation}
where the maximum over an empty set is zero.
\end{lemma}

\begin{proof}
Put \(q=2^n\).  The exceptional set \(U_n(c)\) is a union of
\(qf_n(c)\) classes modulo \(q\).  Since \((d,q)=1\), each such
class meets a prescribed class modulo \(d\) in one class modulo
\(dq\).  Thus
\[
 \#\{N\le x:N\equiv a\pmod d,\ N\in U_n(c)\}
 \le f_n(c)x/d+q.
\]
If \(N,N+c>qK\) and \(N\notin U_n(c)\), their orbits coalesce
by time \(n\), and all preceding orbit values exceed \(K\).
Consequently their indicators of membership in \(G\) agree.
In the class \(a\bmod d\), the starting values excluded by
\(N\le qK+|c|\) number at most \((qK+|c|)/d+1\).
Comparing the translated interval \([1+c,x+c]\) with \([1,x]\)
costs at most \(2|c|/d+2\) additional progression members.
These three estimates prove the first inequality after
multiplication by \(d\).

Average that inequality over the \(d\) integer representatives
\(-d/2<c\le d/2\).  Their mean absolute value is at most \(d/4\),
and \(\sum_cF^{[d]}_{a+c}(y)=F(y)\).
Using \(qK\le dqK\) and \(3d/4+3d\le4d\) gives
\eqref{eq:Doddavg}.
\end{proof}

\begin{theorem}[Uniformity in a power range of moduli]\label{thm:LD}
For \(0<\varepsilon<1\), put
\[
 \begin{split}
 \kappa_1(\varepsilon)&=\frac1{2(1+L(\varepsilon)\log2)},\\
 \gamma(\varepsilon)&=
 \min\left\{\frac{\mu_0}{2L(\varepsilon)},
              \frac{\varepsilon^2}{4.46L(\varepsilon)}\right\},\\
 \kappa(\varepsilon)&=
 \min\{\kappa_1(\varepsilon)/4,\gamma(\varepsilon)\}.
 \end{split}
\]
For every \(T\)-invariant set \(G\),
\[
 \limsup_{x\to\infty}\ 
 \max_{1\le Q\le x^{\kappa(\varepsilon)}}\ \max_{a\bmod Q}
 \left|\frac{QF^{[Q]}_a(x)}x-\frac{F(x)}x\right|
 \le4\varepsilon .
\]
\end{theorem}

\begin{proof}
Fix \(\varepsilon\), and suppress it in \(L,\kappa_1,\gamma,\kappa\).
First let \(d\le x^{\kappa_1}\) be odd, and take
\[
 n=\lceil L\log(2+d)\rceil+n_0(\varepsilon).
\]
Proposition~\ref{prop:U} gives \(f_n(c)\le\varepsilon\) for
\(0<|c|\le d/2\).  Moreover,
\[
 d2^n(K+1)+4d
 \ll_\varepsilon(K+1)d^{1+L\log2}
 \le C_\varepsilon(K+1)x^{1/2}.
\]
Equation~\eqref{eq:Doddavg} therefore gives
\[
 \left|\frac{dF^{[d]}_a(x)}x-\frac{F(x)}x\right|
 \le\varepsilon+O_\varepsilon((K+1)x^{-1/2}),
\]
uniformly in the stated range of \(d\) and in \(a\).

Now let \(Q=2^kd\le x^\kappa\), with \(k\ge1\) and \(d\) odd.
Choose \(r\in\{0,\ldots,2^k-1\}\) congruent to \(a\) modulo \(2^k\).
The substitution \(N=2^kY+r\) turns the other congruence into
\[
 Y\equiv y_0=2^{-k}(a-r)\pmod d.
\]
Writing \(q=q_k(r)\), Terras conjugacy gives
\[
 T^k(N)=3^qY+T^k(r).
\]
Thus the image is a single progression
\[
 b\equiv3^qy_0+T^k(r)\pmod D,\qquad D=3^qd.
\]
This statement does not require \(3^q\) to be invertible modulo \(d\).
The endpoint and finite-defect argument of Lemma~\ref{lem:R},
restricted to this progression, gives
\[
 F^{[Q]}_a(x)=F^{[D]}_b(y)+O(K+6),
 \qquad y=3^qx/2^k.
\]
Here \(Q/x=D/y\).  Also, writing \(\alpha=\log3/\log2<1.59\),
\[
 D\le Q^\alpha\le x^{\alpha\kappa},
 \qquad y\ge x^{1-\kappa}.
\]
Because \(\kappa\le\kappa_1/4\) and \(\kappa\le1/8\),
\(\alpha\kappa\le(1-\kappa)\kappa_1\).
Hence \(D\le y^{\kappa_1}\), and the odd-modulus estimate applies
uniformly to \(D\) at the argument \(y\).

Finally,
\[
 |\log y-\log x|
 =|q\log3-k\log2|
 \le k\log2\le\kappa\log x.
\]
The right-hand side is at most
\(\gamma\max\{\log x,\log y\}\).
Proposition~\ref{thm:W}, applied at the larger logarithmic argument,
therefore bounds the difference of \(F(y)/y\) and \(F(x)/x\) by
\[
 2\varepsilon+2.11\sqrt{L\gamma}+o(1)
 \le3\varepsilon+o(1),
\]
since \(2.11^2<4.46\).  The normalized endpoint error is
\(O((K+6)x^{\kappa-1})\).
Combining these estimates proves the result.
\end{proof}

\begin{corollary}[Tao's conclusion in progressions]\label{cor:TaoAP}
Let \(f(N)\to\infty\), and let \(\varepsilon>0\).
There is \(\kappa>0\) such that, for all sufficiently large \(x\),
all integers \(Q\le x^\kappa\), and all \(a\),
\[
 \#\{N\le x:N\equiv a\pmod Q,\ 
                   \operatorname{Col}_{\min}(N)\ge f(N)\}
 \le\varepsilon\,x/Q .
\]
\end{corollary}

\begin{proof}
It is enough to consider \(0<\varepsilon<1\).
Choose a fixed \(C\) such that
\(\underline d(G_C)\ge1-\varepsilon/4\), where
\(G_C=\{N:\operatorname{Col}_{\min}(N)\le C\}\);
this is possible by Corollary~\ref{thm:taonatural} and its proof.
Take \(\kappa=\kappa(\varepsilon/16)\).
Theorem~\ref{thm:LD} applied to \(G_C\) implies that the number
outside \(G_C\) in any of these progressions is at most
\[
 (\varepsilon/2+o(1))x/Q+1.
\]
For all but a fixed finite number \(N_C\) of starting values,
\(f(N)>C\); a value in \(G_C\) then cannot be exceptional.
Adding \(N_C\) to the last bound is negligible uniformly, since
\(x/Q\ge x^{1-\kappa}\to\infty\).
\end{proof}

\subsection{Primes}

Bondurant-David \cite[Theorem~P and Section~1.3]{BondurantDavid2026}
proves almost-boundedness for prime starting values with logarithmic
weights $(\log p)/p$ and distinguishes that conclusion from a limit
normalized by $\pi(x)$. Below we derive the latter type of bound from
counted coalescence and classical interval sieve estimates, for every
finite-defect invariant set. For the standard map, the strong
coalescence input is already available in
\cite[Theorem~2.11]{LeenRates2026}; the prime and sieve results are
applications of that input. The density convention and the full
invariant-set scope are the distinctions from the cited prime theorem.

For \(E\subseteq\mathbb N\), put
\[
 \pi_E(x)=\#\{p\le x:p\hbox{ prime},\ p\in E\},
 \qquad
 \overline d_{\mathbb P}(E)
 =\limsup_{x\to\infty}\frac{\pi_E(x)}{\pi(x)}.
\]
We use the interval form of the Brun--Titchmarsh inequality:
for \(z>0\), \(y>q\), and \((b,q)=1\),
\begin{equation}\label{eq:MV}
 \pi(z+y;q,b)-\pi(z;q,b)
 <\frac{2y}{\varphi(q)\log(y/q)}.
\end{equation}
Its case \(q=1\) gives
\(\pi(z+y)-\pi(z)\le2y/\log y\) for \(y>1\).
Both forms, with constants independent of the interval position,
are in~\cite[Theorem~2, (1.10) and (1.12)]{MontgomeryVaughan1973}.

\begin{theorem}[Invariant sets and primes]\label{thm:P}
For every \(T\)-invariant set \(E\) and every \(0<\varepsilon<1\),
\[
 \overline d_{\mathbb P}(E)
 \le6\log2\,L(\varepsilon)\overline d(E)+6\varepsilon .
\]
In particular, a \(T\)-invariant set of natural density zero
contains a relative density zero of the primes, and a
\(T\)-invariant set of natural density one contains almost all primes.
By Proposition~\ref{thm:T}, logarithmic density zero or one suffices.
\end{theorem}

\begin{proof}
Fix \(\varepsilon\), and set
\[
 L=L(\varepsilon),\qquad
 \eta=\frac1{3L\log2},\qquad
 H=\lfloor x^\eta\rfloor,\qquad
 n=\lceil L\log(2+H)\rceil+n_0(\varepsilon).
\]
For large \(x\), \(H\ge2\), \(2^n\ll_\varepsilon x^{1/3}\),
and \(f_n(c)\le\varepsilon\) for every integer \(1\le c\le H\).
Write \(K=K(E)=\max(E\triangle T^{-1}(E))\), with
\(\max\varnothing=0\).

For \(N>2^nK\), all iterates of \(N\) and \(N+c\) before time \(n\)
exceed \(K\).  If \(N\notin U_n(c)\), the two orbits coalesce
by that time, so their membership indicators agree.  Therefore
\begin{equation}\label{eq:localconst}
 \left|\mathbf1_E(N)-\frac1H\sum_{c=1}^H\mathbf1_E(N+c)\right|
 \le\frac1H\sum_{c=1}^H\mathbf1_{U_n(c)}(N)
 \qquad(N>2^nK).
\end{equation}

Let \(q=2^n\).  The set \(U_n(c)\) consists of \(qf_n(c)\)
classes modulo \(q\).  Only odd classes can contain primes other
than \(2\), and \(\varphi(q)=q/2\).
Apply \eqref{eq:MV} with \(z=1/2\) and \(y=x-1/2\).
Since
\(\log((x-1/2)/q)\ge(2/3-o(1))\log x\), it follows that,
uniformly for \(1\le c\le H\),
\[
 \#\{p\le x:p\in U_n(c)\}
 \le(6+o(1))\,\varepsilon\,x/\log x+1.
\]

For the smoothed term, the exact order-of-summation identity is
\[
 \sum_{p\le x}\frac1H\sum_{c=1}^H\mathbf1_E(p+c)
 =\frac1H\sum_{m\le x+H}\mathbf1_E(m)
   \#\bigl(\mathbb P\cap[1,x]\cap[m-H,m-1]\bigr).
\]
In particular, the upper restriction \(p\le x\) is retained.
Terms with \(m\le H+1\) have total contribution \(O(H)\).
For \(m>H+1\), the interval bound in \eqref{eq:MV} gives
at most \(2H/\log H\) primes in \([m-H,m-1]\).
Consequently,
\[
 \sum_{p\le x}\frac1H\sum_{c=1}^H\mathbf1_E(p+c)
 \le\frac2{\log H}F_E(x+H)+O(H).
\]
Summing \eqref{eq:localconst} over primes and accounting for
the initially excluded primes now gives
\[
 \pi_E(x)
 \le \pi(2^nK)+\frac2{\log H}F_E(x+H)+O(H)
       +(6+o(1))\varepsilon x/\log x+1.
\]
Divide by \(\pi(x)\sim x/\log x\).
The first and third terms are negligible for fixed \(K,\varepsilon\),
\(\log H/\log x\to\eta\), and \(H=o(x)\).
This proves the displayed bound.  The assertions for densities
zero and one follow by letting \(\varepsilon\downarrow0\) and
then applying the result to a complement when necessary.
\end{proof}

\begin{corollary}[A cube-root transfer bound]\label{cor:prime-cuberoot}
There is an absolute constant \(C_{\mathbb P}\) such that every
\(T\)-invariant set \(E\) satisfies
\[
 \overline d_{\mathbb P}(E)
 \le C_{\mathbb P}\,\overline d(E)^{1/3}.
\]
The constant is independent of the finite defect cutoff \(K(E)\);
the onset of any finite-\(x\) estimate can depend on that cutoff.
\end{corollary}

\begin{proof}
Use the stronger published coalescence estimate
\cite[Theorem~2.11]{LeenRates2026}, which supplies an absolute \(B\)
such that
\[
 f_n(c)\le B\sqrt{\frac{1+\log(1+|c|)}n}.
\]
Fix \(0<\eta\le1/4\), take \(H=\lfloor x^\eta\rfloor\), and
take \(n=\lfloor(\log x)/(3\log2)\rfloor\).
Then \(2^n\le x^{1/3}\), and
\[
 \limsup_{x\to\infty}\max_{1\le c\le H}f_n(c)
 \le B\sqrt{3(\log2)\eta}.
\]
Repeat the local-constancy and Brun--Titchmarsh argument of
Theorem~\ref{thm:P}, with this estimate in place of
\(f_n(c)\le\varepsilon\).  The finite-defect and endpoint
terms still tend to zero after normalization.  Writing
\(\delta=\overline d(E)\), the result is
\[
 \overline d_{\mathbb P}(E)
 \le \frac{2\delta}{\eta}
       +6B\sqrt{3\log2}\,\sqrt\eta.
\]
For \(0<\delta\le1\), put \(\eta=\delta^{2/3}/4\);
both terms are constant multiples of \(\delta^{1/3}\).
For \(\delta=0\), let \(\eta\downarrow0\).
This is a short sieve consequence of the cited stronger input.
\end{proof}

\begin{corollary}[Tao's conclusion for almost all primes]\label{cor:Taoprimes}
For every \(f(N)\to\infty\),
\[
 \#\{p\le x:\operatorname{Col}_{\min}(p)\ge f(p)\}
 =o(\pi(x)).
\]
\end{corollary}

\begin{proof}
For fixed \(C\), the set \(E_C=\mathbb N\setminus G_C\) is
invariant above \(C\).
The proof of Corollary~\ref{thm:taonatural} gives
\(\overline d(E_C)\to0\) as \(C\to\infty\).
Apply Theorem~\ref{thm:P}, first let \(C\to\infty\) for fixed
\(\varepsilon\), and then let \(\varepsilon\downarrow0\).
For every fixed \(C\), the exceptional set in the assertion is
contained in \(E_C\) apart from finitely many primes.
\end{proof}

\begin{corollary}[A quantitative form for primes]\label{cor:Qprimes}
With \(A\) as in Corollary~\ref{cor:Q}, there are absolute
constants \(K,C_0>0\) such that for \(C\ge C_0\),
\[
 \limsup_{x\to\infty}
 \frac{\#\{p\le x:\operatorname{Col}_{\min}(p)>C\}}{\pi(x)}
 \le K(\log C)^{-1/(12A)}(\log\log C)^{5/12}.
\]
\end{corollary}

\begin{proof}
As in the proof of Corollary~\ref{cor:Q}, put
\(\delta=(C_2/\log C)^{1/A}\).
For large \(C\) there is an absolute \(C_3\) such that
\[
 \overline d(E_C)\le
 u=C_3(\delta\log(2/\delta))^{1/4}.
\]
For sufficiently large \(C\), \(0<u<1/2\).
Theorem~\ref{thm:P} and
\(L(\varepsilon)\ll\varepsilon^{-2}\log(10/\varepsilon)\)
bound the asserted upper density by
\[
 C_4u\varepsilon^{-2}\log(10/\varepsilon)+6\varepsilon.
\]
Choose \(\varepsilon=(u\log(2/u))^{1/3}\).
This gives
\[
 \ll (u\log(2/u))^{1/3}
 \ll\delta^{1/12}(\log(2/\delta))^{5/12},
\]
which is the stated estimate.
Using an explicit positive upper bound \(u\) also covers the case
where the actual upper density of \(E_C\) is zero.
\end{proof}

\begin{remark}[Consequences of stronger inputs]
The preceding estimate records the rate obtained from
Proposition~\ref{prop:U} and Tao's logarithmic-density bound.
Using Corollary~\ref{cor:prime-cuberoot} with the same integer-density
bound improves the factor $(\log\log C)^{5/12}$ to
$(\log\log C)^{1/12}$.
There is also a stronger fixed-barrier input:
Nashida \cite[Theorem~1.1]{Nashida2026} proves, for an explicit
$\alpha_N>0$ and an absolute constant $K_N$,
\[
 \#\{N\le x:\Colmin(N)>C\}\le K_N C^{-\alpha_N}x
 \qquad(C\in\N,\ C\ge1).
\]
Combining that published estimate with
Corollary~\ref{cor:prime-cuberoot} yields
\[
 \limsup_{x\to\infty}
 \frac{\#\{p\le x:\Colmin(p)>C\}}{\pi(x)}
 \ll C^{-\alpha_N/3}.
\]
These are consequences of the specified inputs; the constants need
not be numerically effective.
\end{remark}

\subsection{Sets accessible to upper bound sieves}

The same smoothing argument uses only a lower counting bound for
the sampling set and upper bounds in intervals and in progressions
modulo powers of two.  We specify the comparison function precisely.
Let \(\delta:[3,\infty)\to(0,1]\) satisfy, for some \(c,B>0\),
\[
 \delta(x)\ge c(\log x)^{-B}\quad\hbox{for all sufficiently large }x,
\]
and suppose that for every fixed \(0<\eta\le1\) there is a finite
\(K_\delta(\eta)\) such that
\[
 \limsup_{x\to\infty}
 \frac{\delta(\lfloor x^\eta\rfloor)}{\delta(x)}
 \le K_\delta(\eta).
\]
These are the hypotheses needed here; they are not a definition
of slow variation in the usual sense.
For \(\delta(x)=(\log x)^{-\alpha}\), with \(\alpha>0\), they
hold with \(K_\delta(\eta)=\eta^{-\alpha}\).

\begin{theorem}[A sieve transfer principle]\label{thm:S}
Let \(S\subseteq\mathbb N\), and let \(\delta\) satisfy the preceding
conditions.  Suppose there are constants \(c_S,C_S>0\), independent
of \(\eta\), with the following properties.  For every fixed
\(0<\eta<1/2\), all sufficiently large \(x\), and
\(H=\lfloor x^\eta\rfloor\),
\begin{enumerate}
 \item[(S1)] \(S(x):=\#(S\cap[1,x])\ge c_Sx\delta(x)\);
 \item[(S2)] \(\#(S\cap(y,y+H])\le C_SH\delta(H)\) for every
 \(H^2\le y\le x\);
 \item[(S3)] \(\#(S\cap[1,x]\cap(b+2^n\mathbb Z))
 \le C_Sx\delta(x)/2^n\) for every \(b\) and every integer
 \(n\ge0\) with \(2^n\le x^{1/2}\).
\end{enumerate}
The threshold for sufficiently large \(x\) may depend on \(\eta\).
Then every \(T\)-invariant set \(E\) satisfies
\[
 \limsup_{x\to\infty}
 \frac{\#(S\cap E\cap[1,x])}{S(x)}
 \le\frac{C_S}{c_S}
 \left(K_\delta(\eta_\varepsilon)\overline d(E)+\varepsilon\right),
 \qquad
 \eta_\varepsilon=\frac1{3L(\varepsilon)\log2},
\]
for \(0<\varepsilon<1\).
Consequently, for every \(f(N)\to\infty\),
\[
 \#\{N\le x:N\in S,\ 
              \operatorname{Col}_{\min}(N)\ge f(N)\}=o(S(x)).
\]
\end{theorem}

\begin{proof}
Fix \(\varepsilon\), take \(\eta=\eta_\varepsilon\), and let
\(H,n\) be as in the proof of Theorem~\ref{thm:P}.
In particular, \(2^n\ll_\varepsilon x^{1/3}\), so (S3) applies
for all sufficiently large \(x\).
Summing the right-hand side of \eqref{eq:localconst} over \(N\in S\)
and using its \(f_n(c)2^n\) exceptional residue classes gives
at most \(C_S\varepsilon x\delta(x)\).
The values \(N\le2^nK(E)\) cost
\(O_\varepsilon(K(E)x^{1/3})\).

The smoothed term is exactly
\[
 \frac1H\sum_{m\le x+H}\mathbf1_E(m)
 \#\bigl(S\cap[1,x]\cap[m-H,m-1]\bigr).
\]
The terms with \(m\le H^2+H+1\) contribute \(O(H^2+H)\).
For the remaining terms put \(y=m-H-1\).
Then \(H^2<y\le x\), and (S2) bounds the kernel by
\(C_SH\delta(H)\).  We obtain the useful unnormalized bound
\[
 \#(S\cap E\cap[1,x])
 \le C_S\delta(H)F_E(x+H)+C_S\varepsilon x\delta(x)
      +O_\varepsilon(K(E)x^{1/3}+x^{2\eta}).
\]
The explicit lower condition on \(\delta\) implies that its
last two polynomial error terms are \(o(x\delta(x))\).
Divide by \(S(x)\ge c_Sx\delta(x)\), use \(H=o(x)\), and take
the upper limit.

For the final assertion, apply the bound to
\(E_C=\mathbb N\setminus G_C\).
Its upper natural density tends to zero as \(C\to\infty\).
First let \(C\to\infty\) for fixed \(\varepsilon\), and then let
\(\varepsilon\downarrow0\).
For every fixed \(C\), the exceptional starting values in the
assertion belong to \(E_C\) apart from finitely many values;
(S1) ensures \(S(x)\to\infty\).
\end{proof}

The prime case has \(\delta(x)=1/\log x\): its lower bound is the
prime number theorem, and its two upper bounds are
\eqref{eq:MV}.  We next verify the required uniformity for two
other populations.  Define the \emph{rough Chen primes} by
\[
 \mathcal C=\{p\hbox{ prime}:\
 p+2\hbox{ is prime, or }p+2=p_1p_2
 \hbox{ with primes }p_1,p_2\ge p^{1/10}\},
\]
where \(p_1=p_2\) is allowed, and put
\[
 \mathcal B=\{n\in\mathbb N:n=a^2+b^2
                     \hbox{ for some }a,b\in\mathbb Z\}.
\]
The roughness condition is part of the definition of \(\mathcal C\);
we do not replace it by the unrestricted condition that \(p+2\)
has at most two prime factors.

\begin{lemma}[Uniform sieve upper bounds]\label{lem:uniform-sieve}
There is an absolute \(C>0\) such that, for every power of two
\(q\ge1\), every integer \(b\), and every interval \(I\) of real
length \(Y\ge3\),
\[
 \#\{m\in I\cap\mathbb Z:b+qm\in\mathcal B\}
 \le C\,Y/\sqrt{\log Y}.
\]
Furthermore, let \(z\ge2\), and let \(P(z)\) be the product of the
odd primes at most \(z\).  Then
\[
 \#\{m\in I\cap\mathbb Z:
       \gcd((b+qm)(b+qm+2),P(z))=1\}
 \le C\,Y/(\log z)^2+Cz^6 .
\]
All constants are independent of \(b,q,I\).
\end{lemma}

\begin{proof}
We give the elementary Selberg upper bound, including its
remainder, to make the uniformity explicit; the method is
classical~\cite{Richert1976,HalberstamRichert1974}.

A sum of two squares cannot have valuation exactly one at a
prime \(p\equiv3\pmod4\).  Thus, modulo \(p^2\), it avoids the
\(p-1\) classes \(p,2p,\ldots,(p-1)p\).
Since \(q\) is a power of two, the corresponding forbidden
classes for \(m\) also number \(p-1\), regardless of \(b,q\).
For this sieve put
\[
 g(p)=(p-1)/p^2,\qquad \rho(p)=p-1
 \quad(p\equiv3\pmod4).
\]
For the second sieve, every odd prime excludes two classes of
\(m\), namely those with \(b+qm\equiv0,-2\pmod p\).
In this case put \(g(p)=2/p\), \(\rho(p)=2\).
In each case extend \(g,\rho\) multiplicatively to squarefree
integers supported on the sieving primes.
If \(\Omega_d\) is the intersection of the forbidden events
for \(p\mid d\), the Chinese remainder theorem gives
\[
 \#(I\cap\Omega_d)=Yg(d)+O(\rho(d)),\qquad \rho(d)\le d.
\]
In the first case the underlying modulus is \(d^2\), not \(d\).
The displayed error remains \(O(\rho(d))\), since it is the
number of relevant residue classes.

All sums in the rest of the proof are over squarefree integers
supported on the chosen primes.  Write \(\mu_{\rm ar}\) for the
Möbius function, put \(h(p)=(1-g(p))/g(p)\), and set
\[
 \mathfrak G(z)=\sum_{d\le z}\frac1{h(d)},\qquad
 \mathfrak G_d(v)=\sum_{\substack{k\le v\\(k,d)=1}}\frac1{h(k)}.
\]
The Selberg coefficients
\[
 \lambda_d=
 \frac{\mu_{\rm ar}(d)}{\mathfrak G(z)}
 \prod_{p\mid d}(1-g(p))^{-1}\mathfrak G_d(z/d)
 \quad(d\le z)
\]
satisfy \(\lambda_1=1\) and \(|\lambda_d|\le d\).
Indeed, \(\mathfrak G_d(z/d)\le\mathfrak G(z)\), and
\((1-g(p))^{-1}\le p\) in both cases.
The associated quadratic form is
\[
 \sum_{d,e\le z}\lambda_d\lambda_e g([d,e])
 =\sum_{r\le z}h(r)
       \left(\sum_{\substack{d\le z\\r\mid d}}\lambda_dg(d)\right)^2
 =\frac1{\mathfrak G(z)}.
\]
For the first equality use
\(g([d,e])=g(d)g(e)\sum_{r\mid(d,e)}h(r)\).
Möbius inversion of the specified coefficients gives
\[
 \sum_{\substack{d\le z\\r\mid d}}\lambda_dg(d)
 =\frac{\mu_{\rm ar}(r)}{h(r)\mathfrak G(z)},
\]
which proves the second equality.

For an integer surviving all local exclusions with \(p\le z\),
the weight
\((\sum_{d\le z}\lambda_d\mathbf1_{\Omega_d})^2\) equals one.
Summing this nonnegative weight and using the remainder bound
therefore gives at most
\[
 \frac{Y}{\mathfrak G(z)}
 +O\!\left(\sum_{d,e\le z}
            |\lambda_d\lambda_e|\rho([d,e])\right)
 \le \frac{Y}{\mathfrak G(z)}+O(z^6).
\]
The last bound follows from
\(\rho([d,e])\le de\), \(|\lambda_d|\le d\), and
\(\sum_{d\le z}d^2=O(z^3)\).

For completeness, \(\mathfrak G(z)\gg(\log z)^\alpha\), with
\(\alpha=1/2\) in the first case and \(\alpha=2\) in the second.
Mertens' estimates, including the fixed progression \(3\bmod4\),
give
\[
 \prod_{p\le w}(1-g(p))^{-1}\asymp(\log w)^\alpha,\qquad
 \sum_{p\le w}g(p)\log p=\alpha\log w+O(1).
\]
These products and sums use only the respective sieving primes.
Normalize the weights \(1/h(d)\) over all squarefree products
of those primes at most \(w\).
The resulting expectation of \(\log d\) is the second displayed
sum.  Take \(w=z^{1/8}\).
Since \(\alpha\le2\), Markov's inequality shows, for large \(z\),
that at least a fixed positive fraction of the total weight is
carried by \(d\le z\).  Its mass is bounded by
\(\mathfrak G(z)\), proving the stated lower bounds.

The second assertion follows immediately.
For the first, choose \(z=Y^{1/12}\), so the remainder is
\(O(Y^{1/2})\), absorbed by \(Y/\sqrt{\log Y}\).
Increasing the absolute constant covers the bounded range of \(Y\).
\end{proof}

\begin{corollary}[Rough Chen primes and sums of two squares]\label{cor:Chen}
For every \(f(N)\to\infty\),
\[
 \begin{split}
 \#\{p\le x:p\in\mathcal C,\
             \operatorname{Col}_{\min}(p)\ge f(p)\}
       &=o(\mathcal C(x)),\\
 \#\{n\le x:n\in\mathcal B,\
             \operatorname{Col}_{\min}(n)\ge f(n)\}
       &=o(\mathcal B(x)).
 \end{split}
\]
Here \(\mathcal C(x)\) and \(\mathcal B(x)\) denote the
corresponding counting functions.
\end{corollary}

\begin{proof}
For \(\mathcal B\), take \(\delta(x)=(\log x)^{-1/2}\).
Landau's theorem~\cite{Landau1908} gives
\(\mathcal B(x)\sim c_{\mathcal B}x/\sqrt{\log x}\), with
\(c_{\mathcal B}>0\), so (S1) holds.
Lemma~\ref{lem:uniform-sieve} with \(q=1\) gives (S2),
with a constant independent of the short-interval exponent.
For (S3), parametrize a class modulo \(q=2^n\le\sqrt x\)
by \(b+qm\), with \(m\) in an interval of length \(Y=x/q\).
Then \(Y\ge\sqrt x\), and the same lemma gives
\[
 \#(\mathcal B\cap[1,x]\cap(b+q\mathbb Z))
 \ll\frac{x}{q\sqrt{\log(x/q)}}
 \ll\frac{x}{q\sqrt{\log x}}.
 \]
This holds for every residue \(b\), including the even and zero
residues.

For \(\mathcal C\), take \(\delta(x)=(\log x)^{-2}\).
The rough form of Chen's theorem gives
\(\mathcal C(x)\gg x/\log^2x\); for example, the proof with
sifting parameter \(v=8\) in
\cite[Theorem~13.2 and its proof and notes]{Richert1976}
counts a subset in which the two prime factors are at least
\(x^{1/8}\).  This is stronger than the threshold \(p^{1/10}\)
used in the definition; see also~\cite{Chen1973}.
Thus (S1) holds.

For (S2), let \(H^2\le y\le x\), and take
\(z=H^{1/12}\).
If \(p\in\mathcal C\cap(y,y+H]\), every prime factor of
\(p(p+2)\) is at least \(H^{1/5}\).
Consequently it survives the second sieve in
Lemma~\ref{lem:uniform-sieve}, and
\[
 \#(\mathcal C\cap(y,y+H])
 \ll H/\log^2H+H^{1/2}\ll H/\log^2H.
\]
All constants are independent of \(y,x,\eta\).

For (S3), fix \(q=2^n\le\sqrt x\) and a residue \(b\).
The primes \(p\le\sqrt x\) in this class cost at most
\(\sqrt x/q+1\), which is \(o(x/(q\log^2x))\) uniformly.
For \(p>\sqrt x\), every prime factor of \(p(p+2)\) is at least
\(x^{1/20}\).
Write \(Y=x/q\ge\sqrt x\) and choose \(z=Y^{1/24}\).
Then \(z\le x^{1/24}<x^{1/20}\), so
Lemma~\ref{lem:uniform-sieve} bounds the remaining primes by
\[
 O\bigl(Y/\log^2Y+Y^{1/4}\bigr)
 =O\bigl(x/(q\log^2x)\bigr).
\]
This proves (S3), uniformly in \(q,b\).
All three hypotheses of Theorem~\ref{thm:S} are now verified
for both populations, and that theorem proves the corollary.
\end{proof}

\begin{remark}[Prime twins]
No lower bound of order \(x/\log^2x\) for the number of
twin-prime pairs is known, so (S1) is not presently available
for that population.
An unconditional statement with an absolute normalization is
nevertheless
\[
 \#\{p\le x:p,p+2\hbox{ prime},\
              \operatorname{Col}_{\min}(p)\ge f(p)\}
 =o(x/\log^2x).
\]
Indeed, this exceptional set is contained in the corresponding
exceptional subset of \(\mathcal C\), which is \(o(\mathcal C(x))\)
by Corollary~\ref{cor:Chen}; the upper sieve bound just proved
gives \(\mathcal C(x)\ll x/\log^2x\).

This does not assert a relative-density statement normalized
by the unknown number of twin-prime pairs.
\end{remark}

\begin{remark}[Direct consequences of Inselmann's estimate]
\label{cor:Stheta}
For every fixed $\theta>0$, $S_\theta$ has natural density one and
contains almost all primes, rough Chen primes, and sums of two squares,
in relative natural density for each population.  There is
$\kappa_\theta>0$ such that
\[
 \max_{1\le Q\le x^{\kappa_\theta}}\max_{a\bmod Q}
 \frac Qx
 \#\{N\in\N:N\le x,\ N\equiv a\pmod Q,\ N\notin S_\theta\}
 \longrightarrow0
 \qquad(x\to\infty).
\]

Indeed, Inselmann
\cite[Definition~2.6 and Corollary~3.2]{Inselmann2024}
gives constants $C_\theta,\delta_\theta>0$ such that
\[
 \#\{N\in\N:N\le x,\ N\notin S_\theta\}
 \le C_\theta x^{1-\delta_\theta}.
\]
Apply his result with exponent $\theta/2$.  For sufficiently large
$N$, the iteration time in that result is positive and
$N^{\theta/2}<N^\theta$; absorb the finitely many initial values
into the constant.

The counting functions of the three sampling populations are bounded
below by constant multiples of $x/\log^2x$ for large $x$, so the same
exceptional count is negligible relative to each population.

For progressions, take $\kappa_\theta=\delta_\theta/2$ and bound the
exceptional count in each class by the total exceptional count.  Then
\[
 \max_{1\le Q\le x^{\kappa_\theta}}\max_{a\bmod Q}
 \frac Qx
 \#\{N\in\N:N\le x,\ N\equiv a\pmod Q,\ N\notin S_\theta\}
 \le C_\theta x^{-\delta_\theta/2}
 \longrightarrow0.
\]
\end{remark}

\subsection{Equidistribution among primes under the Riemann hypothesis}

Theorem~\ref{thm:P} provides an upper bound and, in particular,
transfers densities zero and one to the primes.
Equality at a general density requires an asymptotic comparison
in place of the interval upper bound used there.

\begin{conjecture}[Conjecture P]\label{conj:P}
If a \(T\)-invariant set \(G\) has natural density \(d\), then
\(\pi_G(x)/\pi(x)\to d\).
\end{conjecture}

For a finite-defect invariant set \(E\), retain the convention
\[
 K(E)=\max(E\triangle T^{-1}(E)),\qquad\max\varnothing=0.
\]

\begin{theorem}[Invariant sets among primes under RH]\label{thm:PE}
Assume the Riemann hypothesis.  For every fixed \(K_0\ge0\),
\[
 \sup_{\substack{E\subseteq\mathbb N\\K(E)\le K_0}}
 \left|\frac{\pi_E(x)}{\pi(x)}-\frac{F_E(x)}x\right|
 \longrightarrow0\qquad(x\to\infty).
\]
The supremum is over all finite-defect \(T\)-invariant sets with
the indicated cutoff.
In particular, Conjecture~\ref{conj:P} holds, and the upper and
lower relative densities of any such set among the primes agree
with its respective upper and lower natural densities.
\end{theorem}

\begin{proof}
Let
\[
 \vartheta_E(x)=\sum_{\substack{p\le x\\p\in E}}\log p,\qquad
 \vartheta(x)=\sum_{p\le x}\log p,\qquad
 \Psi_{\mathbb P}(x)=\sum_{n\le x}\Lambda(n).
\]
These functions are zero below \(2\).
Fix \(0<\varepsilon<1\), and choose the integer \(H\) and \(n\)
as in the proof of Theorem~\ref{thm:P}.
Multiply \eqref{eq:localconst} by \(\log p\) and sum over
primes \(p\le x\).
The same progression bound, multiplied by \(\log x\), gives
uniformly for \(K(E)\le K_0\)
\[
 \vartheta_E(x)
 =\sum_{p\le x}\frac{\log p}{H}
                 \sum_{c=1}^H\mathbf1_E(p+c)
    +O(\varepsilon x)+o_\varepsilon(x)
    +O_\varepsilon(K_0x^{1/3}\log x).
\]
The coefficient in \(O(\varepsilon x)\) can be absolute;
the \(o_\varepsilon(x)\) is independent of \(E,K_0\).

After reversing the summations, the terms with \(m>x\)
contribute \(O(H\log x)\): there are at most \(H\) such integers,
and each kernel divided by \(H\) is at most \(\log x\).
For \(m\le x\), the prime cutoff \(p\le x\) imposes no further
restriction.  Thus the smoothed expression is
\[
 \begin{split}
 &\frac1H\sum_{m\le x}\mathbf1_E(m)
       \bigl(\vartheta(m-1)-\vartheta(m-H-1)\bigr)
       +O(H\log x)\\
 &\qquad=F_E(x)+\frac1H\sum_{m\le x}\mathbf1_E(m)\Delta_m
       +O(H\log x),
 \end{split}
\]
where
\[
 \Delta_m=\vartheta(m-1)-\vartheta(m-H-1)-H.
\]

Under RH, Selberg's short-interval mean-square estimate
\cite{Selberg1943} gives
\[
 \int_1^X
   |\Psi_{\mathbb P}(t+H)-\Psi_{\mathbb P}(t)-H|^2\,dt
 \ll HX\log^2X .
\]
We need this for the fixed power \(H=\lfloor x^\eta\rfloor\)
and \(X\asymp x\).
Because \(H\) is an integer, the integrand is constant between
successive integer values of \(t\).  Therefore the corresponding
discrete sum has the same bound.
For \(H+2\le m\le x\), Cauchy--Schwarz consequently gives
\[
 \frac1H\sum_{H+2\le m\le x}
 \left|\Psi_{\mathbb P}(m-1)
            -\Psi_{\mathbb P}(m-H-1)-H\right|
 \ll \frac{x\log x}{\sqrt H}=o(x).
\]

Prime powers can be removed in \(L^1\).
Every \(p^k\le x\), \(k\ge2\), belongs to at most \(H\) of the
integer windows \((m-H-1,m-1]\).  Hence
\[
 \begin{split}
 &\frac1H\sum_{H+2\le m\le x}
  \bigl((\Psi_{\mathbb P}-\vartheta)(m-1)
       -(\Psi_{\mathbb P}-\vartheta)(m-H-1)\bigr)\\
 &\hspace{30mm}\le
   \sum_{\substack{p^k\le x\\k\ge2}}\log p
   \ll \sqrt x\,\log^2x=o(x).
 \end{split}
\]
The summands here are nonnegative.
The initially omitted \(m\le H+1\) contribute
\(O(H\log x)=o(x)\) after division by \(H\), using the trivial
bound for each window.
It follows that
\[
 \frac1H\sum_{m\le x}|\Delta_m|=o_\varepsilon(x).
\]
No mean-square bound for prime powers is required.

All these errors are uniform over \(K(E)\le K_0\).
For fixed \(\varepsilon,K_0\), first let \(x\to\infty\);
then let \(\varepsilon\downarrow0\).
The preceding identities give
\[
 \sup_{K(E)\le K_0}|\vartheta_E(x)-F_E(x)|=o(x).
\]
Finally, partial summation and
\(0\le\vartheta_E(t)\le\vartheta(t)\ll t\) imply, uniformly in \(E\),
\[
 \pi_E(x)=\frac{\vartheta_E(x)}{\log x}
                    +O(x/\log^2x).
\]
Together with the prime number theorem, this proves the stated
normalized comparison and its consequences.
\end{proof}

\begin{remark}\label{rem:PE}
\emph{(a)} Uniformity over all finite defect cutoffs is impossible.
Every finite set is finite-defect invariant, since its preimage
under \(T\) is finite.  For
\(E=\{p\le M:p\hbox{ prime}\}\) and \(x=M\), the prime proportion
is one whereas \(F_E(x)/x=\pi(x)/x\to0\).

\emph{(b)} RH is used in this proof only for the short-interval
comparison.  More precisely, the argument would work with the
corresponding \(L^1\) estimate
\[
 \frac1H\sum_{m\le x}
 |\vartheta(m-1)-\vartheta(m-H-1)-H|=o(x)
\]
for arbitrarily small fixed exponents \(0<\eta<1/2\), with
\(H=\lfloor x^\eta\rfloor\).
This describes the input required by this proof, rather than
excluding other approaches to Conjecture~\ref{conj:P}.
Available unconditional almost-all short-interval prime number
theorems include the range \(\eta>2/15\)
\cite[Theorem~1.1(ii)]{GafniTao2026}.
A fixed positive lower threshold for \(\eta\) does not supply
the arbitrarily small exponents used here as
\(\varepsilon\downarrow0\).
\end{remark}

\section{\textbf{Generalized Collatz maps}}\label{sec:general}

Let \(\ell,m\ge2\) be integers with \(\gcd(m,\ell)=1\).
For \(j\in\Z/\ell\Z\), let \(r_j\) be the representative of
\(-mj\pmod\ell\) in \(\{0,\ldots,\ell-1\}\).
Thus \(r_0=0\), and \(1\le r_j\le\ell-1\) for \(j\ne0\).
Put \(r(n)=r_{n\bmod\ell}\) and define
\begin{equation}\label{eq:Tlm}
 T(n)=T_{\ell,m}(n)=
 \begin{cases}
  n/\ell,&\ell\mid n,\\
  (mn+r(n))/\ell,&\ell\nmid n.
 \end{cases}
\end{equation}
This definition applies to every \(n\in\Z\), including \(0\),
and preserves \(\N\).  It gives \eqref{eq:T} when \((\ell,m)=(2,3)\).
Write
\[
 \varpi=\frac{\ell-1}{\ell},\qquad
 \mu=\mu_{\ell,m}=\varpi\log m-\log\ell.
\]
The canonical maps considered below satisfy
\begin{equation}\label{eq:H}
 \gcd(m,\ell)=1,\qquad m^{\ell-1}<\ell^\ell,
 \qquad\text{equivalently }\mu<0.
\end{equation}
Neither \(\ell\) nor \(m\) is required to be prime.
Here \(\Z_\ell=\varprojlim_a\Z/\ell^a\Z\); for composite \(\ell\)
it is the product of the rings \(\Z_p\) over the primes \(p\mid\ell\).
We use its normalized Haar measure.

Set \(e_0=0\), \(e_j=1\) for \(j\ne0\), and
\(g_j(z)=(m^{e_j}z+r_j)/\ell\).
The branch \(g_j\) is applied on \(j+\ell\Z_\ell\).
Let \(T_{\min}(N)=\min_{a\ge0}T^a(N)\) for \(N\in\N\).
The definitions of finite-defect invariance, the cutoff \(K\),
and the counting functions \(F,F_{a,j},F_a^{[Q]},\psi\) are as in
Sections~\ref{sec:renewal}--\ref{sec:natural}, with classes modulo
\(m^j\) replacing classes modulo \(3^j\).

\subsection{The maps of Gon\c{c}alves, Greenfeld and Madrid}

Gon\c{c}alves, Greenfeld and Madrid~\cite{GGM2025} consider the map
\(C(N)=N/\ell\) if \(\ell\mid N\), and
\(C(N)=mN+r_j\) if \(N\equiv j\ne0\pmod\ell\).
Their hypotheses apply to the canonical residues above.
Under \eqref{eq:H}, their Theorem~1.3 implies that, for every
\(\delta>0\), there is \(B_\delta\) such that
\begin{equation}\label{eq:GGM}
 \sum_{\substack{N\le x\\C_{\min}(N)<B_\delta}}\frac1N
 >(1-\delta)\log x
\end{equation}
for all sufficiently large \(x\).
Their Theorem~1.8 also gives a natural-density estimate of Korec
type: for positive residues, almost every \(N\) satisfies
\(C_{\min}(N)<N^{\gamma+\eta}\), for every \(\eta>0\), where
\(\gamma=\varpi\log_\ell m\).

The natural-density almost-boundedness question has since been
addressed in the public preprint~\cite{NashidaGeneral2026}:
Theorem~1.1 and Corollary~1.3 there state a result for the full
admissible residue family, including the canonical maps here.
In particular, that statement includes the conclusion of
Corollary~\ref{thm:GA} below for an arbitrary function tending to
infinity.  We do not claim priority for that conclusion.
The argument retained here derives it from the logarithmic-density
input \eqref{eq:GGM} through the coalescence and density-transfer
estimates developed in this section.

For the standard map \((\ell,m)=(2,3)\), balanced coalescence
and quantitative synchronization have prior treatments in
\cite{Duran2026,LeenRates2026}; in particular, the uniform
square-root logarithmic-displacement bound in
Theorem~\ref{thm:Gsharp} is already known in that case
\cite[Theorem~2.11]{LeenRates2026}.
The focus below is the extension of the relation-chain estimates
to general bases, and the arithmetic structure of the coalescing
shifts for noncanonical residues.  Density and sparse-set
comparisons are then derived from those estimates.

\begin{lemma}\label{lem:CT}
For the canonical residues, \(C_{\min}(N)=T_{\min}(N)\) for
every \(N\in\N\).
\end{lemma}

\begin{proof}
If \(\ell\mid N\), then \(C(N)=T(N)\).
Otherwise \(C(N)=\ell T(N)\), so \(C^2(N)=T(N)\).
The \(C\)-orbit therefore inserts, between consecutive terms of
the \(T\)-orbit, only values of the form \(\ell T(M)>T(M)>0\).
These inserted values do not change the minimum.
\end{proof}

\subsection{Terras renewal and the coalescence chain}

\begin{lemma}\label{lem:GTerras}
The map \(T\) extends continuously to \(\Z_\ell\) and preserves
Haar measure.  For \(n\ge1\), \(0\le a<\ell^n\), and
\(Y\in\Z_\ell\),
\[
 T^n(\ell^nY+a)=m^{q_n(a)}Y+T^n(a),
\]
where \(q_n(a)\) counts the indices \(i<n\) for which
\(\ell\nmid T^i(a)\).  The same identity holds for integer \(Y\).
For uniform \(a\bmod\ell^n\), the successive residues
\(T^i(a)\bmod\ell\), \(0\le i<n\), are independent and uniform;
in particular \(q_n(a)\sim\operatorname{Bin}(n,\varpi)\).
Moreover,
\[
 0\le T^n(a)-\frac{m^{q_n(a)}}{\ell^n}a
 \le\frac{\ell-1}{\ell(m-1)}m^{q_n(a)}.
\]
\end{lemma}

\begin{proof}
On the cylinder \(j+\ell\Z_\ell\), the map is
\[
 j+\ell z\longmapsto m^{e_j}z+g_j(j).
\]
This is an affine bijection onto \(\Z_\ell\), because \(m\)
is a unit in \(\Z_\ell\).  Its conditional image of Haar
measure is Haar measure.  Iterating gives a bijection between
the \(\ell^n\) initial cylinders and the residue words of
length \(n\), each of mass \(\ell^{-n}\), as well as the affine
identity.  Conditional on any initial cylinder, the endpoint
is Haar distributed.  This also proves the asserted freshness
and independence.

Let \(q_i=q_i(a)\) and
\(d_i=T^i(a)/m^{q_i}-a/\ell^i\).
A zero-residue step sends \(d_i\) to \(d_i/\ell\), and a
nonzero-residue step sends it to
\(d_i/\ell+r_j/(\ell m^{q_i+1})\).
The added terms are nonnegative.  Each successive nonzero step
has a distinct index \(q_i+1\), so
\[
 0\le d_n\le\frac{\ell-1}{\ell}
                   \sum_{a=1}^{q_n}m^{-a}
 \le\frac{\ell-1}{\ell(m-1)}.
\]
Multiplication by \(m^{q_n}\) proves the last assertion.
\end{proof}

The relation chain has state space
\[
 \mathcal S_{\ell,m}
 =\{(k,B):k\in\Z,\ B\in m^{\min(k,0)}\Z\}.
\]
The normalized offset \(b=m^{\max(-k,0)}B\) is integral.
For \(\beta\in\Z_\ell\) and \(c\in\Z\), put
\(x_n=T^n(\beta)\) and \(y_n=T^n(\beta+c)\).

\begin{lemma}[general coalescence chain]\label{lem:Gchain}
\begin{enumerate}
 \item[(a)] One has \(y_n=m^{k_n}x_n+B_n\), with
 \((k_n,B_n)\in\mathcal S_{\ell,m}\), \(k_0=0\), and \(B_0=c\).
 \item[(b)] If \(j=x_n\bmod\ell\) and \(j'=y_n\bmod\ell\), then
 \[
 \begin{split}
  k_{n+1}&=k_n+e_{j'}-e_j,\\
  B_{n+1}&=
   \frac{m^{e_{j'}}B_n+r_{j'}-m^{k_{n+1}}r_j}{\ell},
  \qquad j'\equiv m^{k_n}j+B_n\pmod\ell.
 \end{split}
 \]
 \item[(c)] Under Haar measure on \(\beta\), the residue \(j\)
 is uniform and independent of the preceding residue word,
 and hence of the chain history through time \(n\).
 These formulas define a Markov chain from every state of
 \(\mathcal S_{\ell,m}\), using fresh uniform residues at each step.
 \item[(d)] Equality of the endpoints and of the numbers of
 multiplication steps at some common time is equivalent to
 reaching \((0,0)\), which is absorbing.
\end{enumerate}
\end{lemma}

\begin{proof}
Substituting \(x=(\ell T(x)-r_j)/m^{e_j}\) into
\(y=m^kx+B\) gives the transition formula.
If \(k'=k+e_{j'}-e_j\), then
\(\max(-k',0)+e_{j'}\ge\max(-k,0)\).
Consequently, after multiplication by \(m^{\max(-k',0)}\),
the numerator defining \(B'\) is integral.
The residue congruence also makes its quotient by \(\ell\)
an element of \(\Z_\ell\).  An integer divided by \(\ell\)
belongs to \(\Z_\ell\) only if that integer is divisible by
\(\ell\), so the normalized new offset is an integer.
This proves preservation of the state space, also for composite
\(\ell\).

The history through time \(n\) depends only on
\(\beta\bmod\ell^n\), whereas Lemma~\ref{lem:GTerras} makes
\(x_n\) Haar distributed independently of that cylinder.
This proves (c).  For an arbitrary initial state, realize
\(x_0\) as Haar distributed and set \(y_0=m^{k_0}x_0+B_0\);
the latter is also Haar distributed.  The same argument applies.
Finally, equality of multiplication counts means \(k_n=0\),
and then equality of endpoints means \(B_n=0\).
Once both conditions hold, they remain true at every later time.
\end{proof}

For a fixed state, \(j\mapsto m^kj+B\) is a permutation of
\(\Z/\ell\Z\).  If \(B\equiv0\pmod\ell\), all transitions
keep the level \(k\).  Otherwise exactly one residue increases
it by one, exactly one decreases it by one, and the remaining
\(\ell-2\) residues keep it fixed.  Thus the physical-time level
is a martingale with increments in \(\{-1,0,1\}\).
The directions of its changes can be generated by independent
fair coins, revealed only when a change occurs.  Completing this
coin sequence after absorption gives an ordinary simple symmetric
walk for use in the return estimates below.
For general \(\ell\), absorption may also occur on a step that
keeps \(k=0\); a waiting time for a level change is therefore
stopped at absorption as well.

\subsection{The fixed-point coordinate and the telescoping identity}

Put
\[
 \mathfrak d(k)=
 \begin{cases}1-m^k,&k\ne0,\\1,&k=0,\end{cases}
 \qquad \Delta(k)=|\mathfrak d(k)|,
 \qquad p=B/\mathfrak d(k).
\]
For \(k\ne0\), this is the fixed point of \(z\mapsto m^kz+B\).

\begin{lemma}\label{lem:Gfix}
At every step, \(p'=ap+\delta\), where
\[
 a=\frac{m^{e_{j'}}\mathfrak d(k)}{\ell\mathfrak d(k')},
 \qquad
 \delta=\frac{r_{j'}-m^{k'}r_j}{\ell\mathfrak d(k')},
 \qquad
 |\delta|\le\delta_0:=\frac{(\ell-1)(m+1)}{\ell(m-1)}.
\]
In particular,
\begin{equation}\label{eq:telescope}
 \prod_{i=1}^n|a_i|
 =\frac{m^{Q_n}}{\ell^n}\frac{\Delta(k_0)}{\Delta(k_n)},
\end{equation}
where \(Q_n\) counts the nonzero-residue branches in the
\(y\)-orbit.  From every initial state,
\(Q_n\sim\operatorname{Bin}(n,\varpi)\).
The same statement holds conditionally for a fresh segment
starting at any fixed time, or at a finite stopping time.
\end{lemma}

\begin{proof}
Substitute \(B=\mathfrak d(k)p\) in
Lemma~\ref{lem:Gchain}(b) and divide by \(\mathfrak d(k')\).
This works also when \(k=0\) or \(k'=0\).
If \(k'\ne0\), the additive term is bounded by
\[
 \frac{\ell-1}{\ell}
      \frac{1+m^{k'}}{|1-m^{k'}|}
 \le\frac{(\ell-1)(m+1)}{\ell(m-1)}.
\]
If \(k'=0\), use \(|r_{j'}-r_j|\le\ell-1\).
The absolute multipliers telescope, proving \eqref{eq:telescope}.
Conditional on the past, \(x_0\) for the new segment is Haar
distributed, and its affine image \(y_0=m^{k_0}x_0+B_0\) is
Haar distributed as well.  Lemma~\ref{lem:GTerras} gives the
binomial law for its future branch count.
The stopping-time version follows from the fresh-residue
construction and the strong Markov property.
No independence between \(Q_n\) and the level increment is asserted.
\end{proof}

Let \(\tau_{\rm abs}\) be the absorption time.
For estimates of the affine recursion we use killed coefficients
\[
 \widetilde a_i=a_i\mathbf1_{\{\tau_{\rm abs}>i\}},\qquad
 \widetilde\delta_i=\delta_i\mathbf1_{\{\tau_{\rm abs}>i\}}.
\]
Then \(p_i=\widetilde a_ip_{i-1}+\widetilde\delta_i\), also
on the step of absorption, and all later values are zero.
Products of these coefficients are dominated by the corresponding
products in \eqref{eq:telescope}; a nonempty product is zero if
it includes an absorbing step.

Set
\[
 \varphi(s)=\ell^{-1-s}+\varpi(m/\ell)^s,
 \qquad
 \Gamma(s)=\varphi(2s)
  \left(1+\frac2\ell\bigl(\cosh(2s\log m)-1\bigr)\right).
\]
Thus \(\varphi(s)=1+\mu s+O(s^2)\) and
\(\Gamma(s)=1+2\mu s+O(s^2)\) as \(s\to0\).

\begin{lemma}[Lyapunov bound]\label{lem:GL}
For \(0<s\le1\), put \(K_s=\sqrt2(m/(m-1))^s\).
For every initial state and \(n\ge1\),
\[
 \E\prod_{i=1}^n|\widetilde a_i|^s
 \le K_s\Gamma(s)^{n/2}.
\]
Choose \(s\in(0,1]\) with \(\Gamma(s)<1\), and then an integer
\(n_1\ge1\) such that \(\rho=K_s\Gamma(s)^{n_1/2}<1\).
With
\[
 C_1=n_1K_s\delta_0^s\max\{1,\Gamma(s)\}^{n_1/2},
\]
one has
\[
 \E[|p_{t+n_1}|^s\mid\mathcal F_t]
 \le\rho|p_t|^s+C_1.
\]
There is also an effective constant \(C\), depending only on
\(\ell,m\) and the choices above, such that
\[
 \sup_{t\ge0}\E|p_t|^s\le C(1+|p_0|^s).
\]
For \(M^s=2C_1/(1-\rho)\), the chain reaches
\(\{|p|\le M\}\) almost surely from every state and, on
nonabsorption, visits this set infinitely often.
\end{lemma}

\begin{proof}
For every integer \(k\),
\[
 (1-1/m)m^{\max(k,0)}\le\Delta(k)\le m^{\max(k,0)}.
\]
Therefore
\(\Delta(k_0)/\Delta(k_n)\le(m/(m-1))m^{|k_n-k_0|}\).
By \eqref{eq:telescope} and Cauchy--Schwarz,
\[
 \E\prod|\widetilde a_i|^s
 \le\left(\frac m{m-1}\right)^s
     \bigl(\E(m^{Q_n}\ell^{-n})^{2s}\bigr)^{1/2}
     \bigl(\E m^{2s|k_n-k_0|}\bigr)^{1/2}.
\]
The first expectation is \(\varphi(2s)^n\).
For each physical step, conditionally on its starting state,
the probabilities of the two nonzero increments of \(k\) are
equal and at most \(1/\ell\).  Hence, for \(u\ge0\),
\[
 \E e^{\pm u(k_n-k_0)}
 \le\left(1+\frac2\ell(\cosh u-1)\right)^n,
\]
and the bound for \(e^{u|k_n-k_0|}\) is at most twice this.
Taking \(u=2s\log m\) proves the product estimate.

Since \(\mu<0\), a sufficiently small positive rational \(s\)
has \(\Gamma(s)<1\), after which a suitable \(n_1\) is obtained
by increasing integers until the displayed strict inequality holds.
Expanding the killed recursion over a block gives
\[
 |p_{t+n_1}|^s\le
 |p_t|^s\prod_{i=t+1}^{t+n_1}|\widetilde a_i|^s
 +\delta_0^s\sum_{j=t+1}^{t+n_1}
                     \prod_{i=j+1}^{t+n_1}|\widetilde a_i|^s.
\]
For each additive term, condition first on \(\mathcal F_j\),
the history immediately after its insertion.  The future-product
bound applies uniformly to that state.  An empty product equals
one and is also bounded by \(K_s\).
The tower property yields the block estimate with the stated \(C_1\).
Iteration, followed by the same calculation for the fewer than
\(n_1\) remaining steps, proves the bound at all physical times.

For completeness, set \(Z_a=|p_{an_1}|^s\) and
\(\eta=(1+\rho)/2<1\).
While \(Z_a>M^s\), its conditional drift is at most \(\eta Z_a\).
If \(\sigma\) is the first block index with \(Z_\sigma\le M^s\),
then
\[
 \E[Z_a;\sigma>a]\le\eta^a Z_0,
 \qquad
 \Pb(\sigma>a)\le M^{-s}\eta^aZ_0.
\]
Thus \(\sigma<\infty\) almost surely.  Applying the same argument
after each deterministic block index rules out a last visit to
this set on nonabsorption.
\end{proof}

\begin{remark}\label{rem:GLconst}
All constants in Lemma~\ref{lem:GL} can be chosen effectively
from \(\ell,m\).  For any fixed \(s\) with \(\Gamma(s)<1\), it
suffices to take
\[
 n_1>\frac{2\log K_s}{\log(1/\Gamma(s))}.
\]
The resulting constants depend on both the contraction margin
and the choice of \(s\); the formulas can deteriorate considerably
near the boundary \(\mu=0\).
\end{remark}

\subsection{Gaps, tightness and descent}

Write \(\ell=\prod_{p\mid\ell}p^{a_p}\), and define
\[
 w(x)=\min_{p\mid\ell}\left\lfloor\frac{v_p(x)}{a_p}\right\rfloor
 \quad(x\in\Z_\ell\setminus\{0\}),\qquad w(0)=\infty.
\]
Thus \(w(x)\ge a\) exactly when \(\ell^a\mid x\).
If \(w(x)\ne w(y)\), then
\(w(x+y)=\min\{w(x),w(y)\}\).
Indeed, if \(w(x)=a<w(y)\), a prime attaining the minimum for
\(x\) has \(v_p(x)<(a+1)a_p\le v_p(y)\); at every other prime,
both summands have valuation at least \(aa_p\).
Multiplication by \(m\) leaves \(w\) unchanged.

For \(k\ne0\) and \(j\ne0\), put
\[
 c_j(k)=r_{m^kj}-m^kr_j,
 \qquad
 t_k=\max_{\substack{j\ne0\\c_j(k)\ne0}}w(c_j(k)).
\]
The index \(m^kj\) is interpreted modulo \(\ell\), also when
\(k<0\).  These rational numbers have denominators that are
units in \(\Z_\ell\).  Some \(c_j(k)\) is nonzero: for \(k>0\)
choose \(r_j=\ell-1\), and for \(k<0\) choose \(r_j=1\).
In the latter case, \(m^{|k|}c_j(k)=m^{|k|}r_{m^kj}-1\ne0\).
Also \(c_j(k)\equiv0\pmod\ell\), so \(t_k\ge1\).
The identity
\[
 c_j(-k)=-m^{-k}c_{m^{-k}j}(k)
\]
shows that \(t_{-k}=t_k\).
After clearing the unit denominator,
\[
 0<\bigl|m^{\max(-k,0)}c_j(k)\bigr|
 \le(\ell-1)(m^{|k|}+1)
 \qquad(c_j(k)\ne0).
\]
Consequently \(t_k\le1+\log_\ell((\ell-1)(m^{|k|}+1))\).

\begin{lemma}[gaps]\label{lem:Ggaps}
At a fixed level \(k\ne0\), from any state with
\(B\equiv0\pmod\ell\), the expected time until
\(B\not\equiv0\pmod\ell\) is at most
\(\ell^{\ell-1}(t_k+\ell)\).
There is an effective constant \(\tau_{\ell,m}\) such that
\[
 \sum_{1\le|k|\le H}t_k\le\tau_{\ell,m}H
 \qquad(H\ge1).
\]
Neither assertion requires \eqref{eq:H} beyond \(\gcd(m,\ell)=1\).
\end{lemma}

\begin{proof}
While \(B\equiv0\pmod\ell\), the level is unchanged and
\[
 B'=B/\ell\quad(j=0),\qquad
 B'=(mB+c_j(k))/\ell\quad(j\ne0).
\]
For finite \(v=w(B)\ge1\), every step outside
the finite set
\[
 \mathcal V_k=\{w(c_j(k)):j\ne0,\ c_j(k)\ne0\}
\]
either decreases \(v\) by one or sends it to a smaller member
of \(\{0,\ldots,t_k-1\}\).
Only a step with \(w(c_j(k))=v\) can fail to decrease \(v\).
There are at most \(\ell-1\) exceptional levels in \(\mathcal V_k\),
and at each such level the choice \(j=0\), of probability
\(1/\ell\), decreases the valuation.

Start an attempt at the current valuation.  If \(v>t_k\), wait
until it falls to at most \(t_k\).  A residue with nonzero
\(c_j(k)\) forces this, so the expected duration of this first
phase is at most \(\ell\), also when \(B=0\).
In the second phase stop at valuation zero, or at the first
nondecreasing step.  This phase lasts at most \(t_k\) steps.
It succeeds if a decreasing choice is made at every exceptional
level encountered, so its conditional success probability is
at least \(\ell^{-(\ell-1)}\).
Restart after each failure.
The event that attempt \(a\) is begun is measurable before its
start, its conditional mean length is at most \(\ell+t_k\),
and its conditional failure probability is at most
\(1-\ell^{-(\ell-1)}\).
Summing these conditional mean lengths gives
\(\ell^{\ell-1}(t_k+\ell)\).  Independence of successive attempts
is not needed.

For the averaged estimate, fix \(j\ne0\) and a class of \(k\)
modulo \(o=\operatorname{ord}_\ell(m)\).
Then \(r_{m^kj}\) is a fixed integer, say \(b\).
The condition \(w(c_j(k))\ge a\) is
\(m^kr_j\equiv b\pmod{\ell^a}\).
With \(g=\gcd(r_j,\ell^a)\le\ell-1\), it either has no solution,
or specifies one residue class of \(m^k\) modulo \(\ell^a/g\).
Its solutions in \(k\) therefore lie in one class modulo
\(\operatorname{ord}_{\ell^a/g}(m)\), before imposing the
additional class modulo \(o\).

Fix a prime \(p\mid\ell\).
The elementary growth of multiplicative orders gives
\[
 \operatorname{ord}_{\ell^a/g}(m)
 \ge \operatorname{ord}_{p^{\max(a_pa-v_p(g),0)}}(m)
 \ge c_{\ell,m}p^{a_pa}.
\]
Here the positive constant is uniform in the finitely many
possible \(r_j\).  To see the order-growth assertion, take a
power of \(m\) congruent to one modulo \(p\), or modulo four
when \(p=2\); the lifting-the-exponent formula shows that its
order gains a factor \(p\) whenever the modulus is increased
by a factor \(p\) beyond a fixed exponent.
Thus, after summing over \(j\) and the \(o\) classes,
\[
 \#\{k:1\le|k|\le H, t_k\ge a\}
 \le C_{\ell,m}\bigl(Hp^{-a_pa}+1\bigr).
\]
The elementary size bound preceding the lemma gives
\(t_k\le C'_{\ell,m}(H+1)\) for \(|k|\le H\).
Summing the last counting estimate over this finite range of
\(a\) proves the desired \(O_{\ell,m}(H)\) bound.
This argument uses no estimate for linear forms in logarithms.
\end{proof}

The conditional mean time until the next change at a nonzero
level \(k\) is at most
\[
 \theta_k=\frac\ell2\left(1+\ell^{\ell-1}(t_k+\ell)\right).
\]
Indeed, after reaching \(B\not\equiv0\), the next step changes
the level with probability \(2/\ell\).  After each failure,
the preceding preparation estimate applies again.
Summing over this geometric number of trials proves the bound.
In particular, \(\theta_{-k}=\theta_k\), the averages of
\(\theta_k\) are bounded, and \(\theta_k=O_{\ell,m}(1+|k|)\).

For sufficiently small \(s>0\), let \(\lambda_+(s)\) and
\(\lambda_-(s)\) be the smaller positive roots of
\[
 \begin{split}
 \lambda&=\frac{(m/\ell)^s}{\ell}\lambda^2
       +\frac{\ell-2}{\ell}(m/\ell)^s\lambda
       +\frac{(1/\ell)^s}{\ell},\\
 \lambda&=\frac{(1/\ell)^s}{\ell}\lambda^2
       +\frac{\ell-2}{\ell}(m/\ell)^s\lambda
       +\frac{(m/\ell)^s}{\ell},
 \end{split}
\]
respectively.  Their polynomials are positive at zero and
negative at one when \(\varphi(s)<1\), so both smaller roots
belong to \((0,1)\).  Expansion of the quadratics at \(s=0\)
gives
\[
 \lambda_\pm(s)=1-\sqrt{\ell|\mu|s}+O_{\ell,m}(s).
\]
We henceforth choose one sufficiently small positive rational
\(s\le1\), and the associated block constants of
Lemma~\ref{lem:GL}, so that simultaneously
\[
 \Gamma(s)<1,\quad\varphi(s)<1,\quad
 r_+:=m^s\lambda_+(s)<1,\quad r_-:=\lambda_-(s)<1.
\]
Write \(r=\max(r_+,r_-)\).

\begin{lemma}[tightness at return]\label{lem:Gtight}
Let \(\tau=\inf\{t\ge0:|k_t|\le1\}\).
There are effective constants \(C_2,C_3\), depending only on
\(\ell,m\), such that, from every state,
\[
 \E|p_\tau|^s\le C_2|p_0|^s+C_3.
\]
Also, if \(R=\inf\{t\ge0:k_t=0\}\), then
\begin{equation}\label{eq:Greturnmoment}
 \E|p_R|^s\le |p_0|^s+C_{\rm ret},
\end{equation}
where \(C_{\rm ret}\) is independent of the initial level.
For \(k_0>0\) the coefficient of \(|p_0|^s\) in this last
bound can be replaced by \(r_+^{k_0}\), and for \(k_0<0\)
it can be replaced by \(r_-^{|k_0|}\).
\end{lemma}

\begin{proof}
All nonzero-level holding times have finite conditional means.
The embedded simple walk therefore reaches the indicated target
level almost surely.
Let \(f_i=m^{e_{j'}}/\ell\) for the residue \(j'\) used in
the \(i\)-th step.
On the positive half-line, the process
\[
 \left(\prod_{i=1}^t f_i^s\right)\lambda_+^{k_t-b},
\]
stopped at level \(b\), is a nonnegative supermartingale, for
\(b=0\) or \(b=1\).  At a state with \(B\not\equiv0\), its
conditional multiplier is one by the defining quadratic.
At a state with \(B\equiv0\), it is \(\varphi(s)\le1\).
Fatou's lemma at the hitting time therefore gives a future-product
bound \(\lambda_+^{k-b}\) from level \(k>b\).
The negative half-line has the same argument with \(\lambda_-\)
and \(|k|\).

For return to \(|k|=1\), telescoping the \(\Delta\) factors
shows, with \(K_\Delta=(m/(m-1))^s\), that from \(|k|\ge2\)
the conditional expected product of the absolute multipliers
to the target is at most
\[
 K_\Delta r^{|k|-1}.
\]
Expand the affine recursion at \(\tau\).
An additive term inserted at time \(j<\tau\) is followed by
a product starting at the state at time \(j\); condition on
\(\mathcal F_j\) before applying this bound.
The final inserted term, at \(j=\tau\), has empty future product
and contributes at most \(\delta_0^s\).
The Green function of a simple walk killed at level one bounds
the expected number of visits to level \(a\ge2\) by \(2(a-1)\),
uniformly in its initial level.  Each such visit has conditional
mean physical duration at most \(\theta_{\pm a}\).
It follows that one may take \(C_2=K_\Delta\) and, for example,
\[
 C_3=\delta_0^s\left(
  1+2K_\Delta\sum_{a\ge2}(a-1)r^{a-1}
                           (\theta_a+\theta_{-a})\right).
\]
This is finite by the elementary bound \(\theta_a=O(1+a)\).
If \(|k_0|\le1\), then \(\tau=0\), so the same inequality holds.

For the target zero, the corresponding product from level
\(a>0\) is bounded by
\(\Delta(a)^s\lambda_+^a\le r_+^a\); from level \(-a\) it is
bounded by \(\Delta(-a)^s\lambda_-^a\le r_-^a\).
The Green function from any initial positive level to zero is
\(2\min(k_0,a)\le2a\), and reflection gives the negative case.
The same insertion-time calculation therefore proves
\eqref{eq:Greturnmoment}, with the safe choice
\[
 C_{\rm ret}=\delta_0^s\left(
  1+2\sum_{a\ge1}a(\theta_a r_+^a+\theta_{-a}r_-^a)
 \right).
\]
For \(k_0=0\), \(R=0\), and the assertion is immediate.
All occupation sums were taken by conditioning at the start of
each physical sojourn; no independence between its length and
the direction of the next change is required.
\end{proof}

\begin{lemma}[descent]\label{lem:Gdesc}
From every state of the canonical relation chain, a finite
sequence of residues leads to absorption.  In particular,
absorption has positive probability from every state.
\end{lemma}

\begin{proof}
For \(\ell=2\), hypothesis \eqref{eq:H} forces \(m=3\), and
Lemma~\ref{lem:descent} applies.
Suppose \(\ell\ge3\).
From any level with \(|k|\ge2\), the embedded simple walk
reaches \(|k|=1\) almost surely, with finite nonzero-level
holding times.  At least one finite residue word therefore
reaches a state at one of these two levels.
If \(B=0\) at a nonzero level, a choice with \(c_j(k)\ne0\)
makes \(B\ne0\) without changing the level.
At level \(1\) or \(-1\), repeated residues zero divide a
nonzero \(B\) by \(\ell\) until \(B\not\equiv0\pmod\ell\).
One more appropriate residue then reaches level zero:
from level \(1\), choose \(mj+B\equiv0\pmod\ell\), obtaining
\(c'=(B-r_j)/\ell\); from level \(-1\), choose \(j=0\),
obtaining \(c'=(mB+r_{j'})/\ell\).
It remains to consider the integral states \((0,c)\).

If \(\ell\mid c\ne0\), a residue zero sends \(c\) to
\(c/\ell\).
For \(\ell\nmid c\), use the following signed routes.
If \(c>0\), a residue zero leads to level \(1\) with
\[
 b_1=B_1=\frac{mc+r_c}{\ell}>0.
\]
If \(c<0\), the residue \(j\equiv-c\pmod\ell\) leads to
level \(-1\) with
\[
 b_1=mB_1=\frac{mc-r_j}{\ell}<0.
\]
In either case the normalized offset \(b_1\) is a nonzero
integer.  Divide it by \(\ell\) as long as it is divisible,
using residues zero, and then return to level zero by the
appropriate residue described above.
The new integral offset satisfies
\[
 |c'|\le\frac m{\ell^2}|c|+1-\frac1{\ell^2}.
\]
It is therefore smaller than \(|c|\) whenever
\[
 |c|>C_{\rm d}:=\frac{\ell^2-1}{\ell^2-m}.
\]
For \(\ell\ge4\), hypothesis \eqref{eq:H} gives
\(m<\ell^{\ell/(\ell-1)}\le\ell^{3/2}\le\ell^2/2\),
so \(C_{\rm d}<2\).
For \(\ell=3\), the admissible values \(m=2,4,5\) give
\(C_{\rm d}=8/7,8/5,2\), respectively.

The remaining offsets are \(c=\pm1\), and also \(c=\pm2\)
when \((\ell,m)=(3,5)\).
From \((0,1)\), the preceding route first reaches
\((1,\lceil m/\ell\rceil)\); from \((0,-1)\), it reaches
\((-1,-\lceil m/\ell\rceil/m)\).
Here \(1\le\lceil m/\ell\rceil\le\ell-1\), so each state
absorbs on its next inward step.
For \((3,5)\), from \((0,2)\) the residues \(0,1,0,2\)
give
\[
 (0,2)\longrightarrow(1,4)\longrightarrow(0,1)
       \longrightarrow(1,2)\longrightarrow(0,0),
\]
whereas from \((0,-2)\) the residues \(1,0,1,0\) give
\[
 (0,-2)\longrightarrow(0,-3)\longrightarrow(0,-1)
        \longrightarrow(-1,-2/5)\longrightarrow(0,0).
\]
Every finite prescribed residue word has positive probability,
which completes the proof.
\end{proof}

\begin{theorem}[coalescence for \(T_{\ell,m}\)]\label{thm:GK}
Under \eqref{eq:H}, the canonical relation chain is absorbed
almost surely from every state of \(\mathcal S_{\ell,m}\).
In particular, for each \(c\in\Z\) and Haar-almost every
\(\beta\in\Z_\ell\), there is a common time \(n\) such that
\(T^n(\beta)=T^n(\beta+c)\) and the two orbits have used
equal numbers of multiplication steps.
\end{theorem}

\begin{proof}
By Lemma~\ref{lem:GL}, on nonabsorption the chain has
infinitely many block times at which \(|p|\le M\).
Choose \(M'\) so that
\((C_2M^s+C_3)/(M')^s\le1/2\).
Starting at any such time, Lemma~\ref{lem:Gtight} implies
probability at least \(1/2\) of reaching
\[
 \mathcal K=\{(k,B)\in\mathcal S_{\ell,m}:
                         |k|\le1, |p|\le M'\}.
\]
This set includes the absorbing state and is finite:
at level zero \(p=B\in\Z\), and at levels \(\pm1\)
one has \(|p|=|b|/(m-1)\) with \(b\in\Z\).
Lemma~\ref{lem:Gdesc} provides an absorbing word at every
state of \(\mathcal K\).  Taking the longest of these finitely
many words gives a common success probability \(\delta>0\)
within a common finite number of steps.
Thus, at every time with \(|p|\le M\), the conditional
probability of eventual absorption is at least \(\delta/2\).

The bounded martingale
\(\Pb(\tau_{\rm abs}<\infty\mid\mathcal F_t)\)
converges almost surely to
\(\mathbf1_{\{\tau_{\rm abs}<\infty\}}\).
On nonabsorption its limit would be zero, contradicting the
preceding lower bound along the infinitely many indicated
times.  This proves almost-sure absorption from every initial
state, and Lemma~\ref{lem:Gchain}(d) gives the orbit statement.
\end{proof}

\subsection{Uniform coalescence and natural density}

For \(n\ge0\), let \(f_n(c)\) be the nonabsorption probability
by step \(n\) from \((0,c)\); thus \(f_n(0)=0\).
Retain the common choices \(s,n_1,\rho,C_1\) above and set
\[
 \lambda_0=\frac{n_1s}{\log(1/\rho)}.
\]
Choose an effective constant \(\theta\) such that
\(\sum_{k=2}^H\theta_k\le\theta H\) for all integers \(H\ge2\).
Lemma~\ref{lem:Ggaps} supplies such a constant, and symmetry
gives the same bound on the negative half-line.

\begin{proposition}[uniform coalescence]\label{prop:GU}
For \(0<\varepsilon<1\), put
\[
 L(\varepsilon)=\lambda_0
  \left(1+\frac{800\theta\log(10/\varepsilon)}{\varepsilon^2}\right).
\]
There is an integer \(n_0(\varepsilon)\) such that
\(f_n(c)\le\varepsilon\) for every \(c\in\Z\) and every integer
\[
 n\ge L(\varepsilon)\log(2+|c|)+n_0(\varepsilon).
\]
\end{proposition}

\begin{proof}
The case \(c=0\) is immediate.
For \(c\ne0\), take the block time
\[
 t_1=n_1\left\lceil
         \frac{s\log(2+|c|)}{\log(1/\rho)}
        \right\rceil
 \le\lambda_0\log(2+|c|)+n_1.
\]
With \(Q=1+C_1/(1-\rho)\) and \(M_1^s=5Q/\varepsilon\),
Lemma~\ref{lem:GL} gives
\(\Pb(|p_{t_1}|>M_1)\le\varepsilon/5\).
The level is a martingale with increments bounded by one.
For \(A=\sqrt{2\log(10/\varepsilon)}\),
Azuma's inequality~\cite{Azuma1967} therefore gives
\[
 \Pb(|k_{t_1}|>A\sqrt{t_1})\le2e^{-A^2/2}=\varepsilon/5.
\]
Let \(\mathcal G\) be the event that both bounds hold.

Put \(D=10/\varepsilon\) and choose the integer barrier
\[
 H=\left\lceil D(A\sqrt{t_1}+1)\right\rceil.
\]
Since \(A\sqrt{t_1}>2\) and \(D>10\),
\(H\le2DA\sqrt{t_1}\).
Let \(\sigma\ge t_1\) be the first time with
\(|k_\sigma|\le1\) or \(|k_\sigma|=H\).
If \(|k_{t_1}|\le1\), this means \(\sigma=t_1\); no departure
from level zero is required.
For \(2\le a=|k_{t_1}|\le A\sqrt{t_1}\), gambler's ruin gives
\[
 \Pb(|k_\sigma|=H\mid\mathcal F_{t_1})
 =\frac{a-1}{H-1}\le\frac1D=\frac{\varepsilon}{10}.
\]
The Green function between levels \(1\) and \(H\) is at most
\(2(a-1)\) at each interior level.  Conditioning at the start
of each sojourn thus gives, on \(\mathcal G\),
\[
 \E[\sigma-t_1\mid\mathcal F_{t_1}]
 \le2A\sqrt{t_1}\,\theta H
 \le4D\theta A^2t_1.
\]
The same statements hold on the negative half-line.
Consequently, with
\[
 T_1=\frac{40D\theta A^2t_1}{\varepsilon}
     =\frac{800\theta\log(10/\varepsilon)}{\varepsilon^2}t_1,
\]
the conditional probability of \(\sigma-t_1>T_1\) is at most
\(\varepsilon/10\).

Let
\[
 M_2^s=\frac{10(C_2M_1^s+C_3)}{\varepsilon}.
\]
On paths hitting \(|k|\le1\) before the barrier, \(\sigma\)
is also the unrestricted hitting time used in
Lemma~\ref{lem:Gtight}.  By nonnegativity and that lemma,
\[
 \Pb(|k_\sigma|\le1, |p_\sigma|>M_2
                      \mid\mathcal F_{t_1})
 \le\varepsilon/10
 \qquad\text{on }\mathcal G.
\]
The set \(\{|k|\le1, |p|\le M_2\}\) is finite.
Theorem~\ref{thm:GK} supplies an integer \(n_2(\varepsilon)\)
such that absorption fails within \(n_2\) steps from any of
these states with probability at most \(\varepsilon/10\).
The strong Markov property at \(\sigma\) now bounds the total
failure probability by at most
\(2\varepsilon/5+4\varepsilon/10<\varepsilon\).
The required time is at most \(t_1+\lceil T_1\rceil+n_2\),
which is bounded by
\[
 L(\varepsilon)\log(2+|c|)
 +n_1\left(1+
     \frac{800\theta\log(10/\varepsilon)}{\varepsilon^2}\right)
 +1+n_2.
\]
Rounding the last three terms upwards defines \(n_0(\varepsilon)\).
For \((\ell,m)=(2,3)\), using the sharper sojourn bound
\(t_k+2\) of Lemma~\ref{lem:gaps} gives
\(\theta=4.5\) and the coefficient \(800\cdot4.5=3600\)
appearing in Proposition~\ref{prop:U}.
\end{proof}

\subsection{Contraction along visits to level zero}

The strong tail estimate uses contraction at returns to level zero.
Count every physical-time visit to level zero, including consecutive
visits.  For a start at
\(k_0=0\), set \(V_0=0\) and
\[
 V_{i+1}=\inf\{t>V_i:k_t=0\},\qquad Y_t=|p_t|^s.
\]
These times are finite almost surely: an excursion from zero
returns by the nonzero-level holding bounds and recurrence of
the embedded walk.  After absorption, the convention simply
gives consecutive visits with \(Y_{V_i}=0\).
Put
\[
 A_+=\frac{(m/\ell)^s}{\ell},\qquad
 C_+=\frac{(1/\ell)^s}{\ell},\qquad
 A_-=C_+,\quad C_-=A_+,\qquad
 b_s=\frac{\ell-2}{\ell}(m/\ell)^s.
\]
Thus \(A_\pm\lambda_\pm^2+b_s\lambda_\pm+C_\pm=\lambda_\pm\).

\begin{lemma}[contraction along returns to level zero]\label{lem:GV}
With the same \(s\) as above, there is an effective finite
constant \(C_V\) such that
\[
 \E[Y_{V_{i+1}}\mid\mathcal F_{V_i}]
 \le\rho_VY_{V_i}+C_V,
 \qquad \rho_V=\varphi(s)<1.
\]
\end{lemma}

\begin{proof}
At an absorbed visit the assertion is immediate.
If \(B\equiv0\pmod\ell\) at a nonabsorbed visit, the next
physical step is another zero-level visit and \(p'=fp\), with
\(\E f^s=\varphi(s)\).
If \(B\not\equiv0\), one residue starts an excursion at level
\(1\) with \(f=m/\ell\), one starts at level \(-1\) with
\(f=1/\ell\), and all other residues stay at level zero.
The product estimate used in Lemma~\ref{lem:Gtight}, now with
target zero, bounds the expected product of the subsequent
\(f^s\) factors by \(\lambda_+\) and \(\lambda_-\), respectively.
Since both endpoints of the segment have \(\Delta=1\), the
full product of the absolute affine multipliers has the same
value as the product of the \(f\)'s, before killing.
Its expected \(s\)-th power is therefore at most
\[
 b_s+A_+\lambda_++C_+\lambda_-
 <b_s+A_++C_+=\varphi(s).
\]
Killing at absorption can only reduce this bound.

For the additive terms, condition on the history immediately
after each insertion.  From a positive post-insertion level
\(a\), the expected future multiplier product is at most
\(r_+^a\); from a negative level \(-a\), it is at most
\(r_-^a\).
There is one final term with empty product, including the case
where the segment consists of a single step at level zero.
A simple walk started at \(1\) and killed at zero has expected
number of visits \(2\) at each positive level; the reflected
statement holds on the negative side.
Using the conditional sojourn bounds, a valid choice is
\[
 C_V=\delta_0^s\left(
  1+2\sum_{a\ge1}(\theta_a r_+^a+\theta_{-a}r_-^a)
 \right)<\infty.
\]
Only the elementary bound \(\theta_a=O(1+a)\) is needed.
This proves the asserted drift with \(\rho_V=\varphi(s)\).
\end{proof}

\subsection{The fixed-shift order of the merging time}
\label{subsec:Gexact}

We first treat \(\ell\ge3\).  Under \eqref{eq:H}, the case
\(\ell=2\) is \(m=3\), covered by Proposition~\ref{thm:exact}.
All constants below may depend on the fixed pair \((\ell,m)\).
Start the coalescence chain at \((0,c)\), \(c\ne0\), and write
\(\tau_{\rm abs}\) for its absorption time.
Put \(\tau_0=0\), and let \(\tau_j\), \(j\ge1\), be the successive
times at which the level changes from \(1\) or \(-1\) to \(0\).
Times after absorption are set equal to infinity.
Define
\[
 E_j=\{\tau_j<\infty,\ \tau_{\rm abs}>\tau_j\},
 \qquad \ell_t=\log(2+|p_t|).
\]
Products with \(\mathbf1_{E_j}\) are defined to be zero off \(E_j\).
In this setting absorption can also occur during a step that stays
at level zero; this possibility is included in the stopping times
used below.

\begin{lemma}[Sojourns]\label{lem:Gsoj}
\emph{(a)} From \((0,B)\), with \(B\in\mathbb Z\setminus\{0\}\),
let \(Z\) be the time until the level changes or the chain is
absorbed.  Then
\[
 \mathbb E[Z\mid\hbox{initial history}]
 \le C_{\ell,m}\bigl(1+\log(2+|B|)\bigr).
\]
\emph{(b)} At a nonzero level \(k\), the conditional expected
time until the next level change is at most
\[
 \theta_k=\frac{\ell}{2}
       \bigl(1+\ell^{\ell-1}(t_k+\ell)\bigr),
\]
uniformly in the remaining state coordinate.
\end{lemma}

\begin{proof}
For (a), put \(a=\max\{1,\log m/\log\ell\}\).
Hypothesis \eqref{eq:H} gives \(a<\ell/(\ell-1)\), and hence
\[
 q=a(\ell-2)/\ell<1.
\]
While \(\ell\mid B\), every step stays at level zero and sends
\(B\) to either \(B/\ell\) or \(mB/\ell\).
Since \((m,\ell)=1\), the divisibility depth \(w(B)\) decreases
by exactly one.  This phase lasts \(w(B)\le\log_\ell|B|\) steps,
does not cause absorption, and leaves an offset \(B'\) satisfying
\[
 |B'|\le |B|^a.
\]
From the resulting state with \(\ell\nmid B'\), two residues
change the level, with total probability \(2/\ell\).
Any other residue either absorbs the chain or leaves
\[
 B''=(mB'+r_{j'}-r_j)/\ell,\qquad
 1+|B''|\le m(1+|B'|).
\]
Call the divisibility phase together with this next step a round.
With \(\beta=\log_\ell(1+|B|)\), a round has length at most
\(\beta+1\), and a surviving next round starts with
\[
 \beta_{\rm next}\le a\beta+\log_\ell m\le a\beta+2.
\]
Its conditional probability is at most \((\ell-2)/\ell\).
Consequently, if \(R\) is the number of rounds, then on \(R\ge i\)
\(\beta_i\le a^{i-1}(\beta_1+2i)\), and
\[
 \begin{split}
 \mathbb EZ
 &\le\sum_{i\ge1}
       \left(\frac{\ell-2}{\ell}\right)^{i-1}
       \bigl(a^{i-1}(\beta_1+2i)+1\bigr)\\
 &\le\frac{\beta_1+1}{1-q}+\frac{2}{(1-q)^2}.
 \end{split}
\]
This proves (a) with the asserted conditional uniformity.

For (b), a phase with \(B\equiv0\pmod\ell\) has expected length
at most \(\ell^{\ell-1}(t_k+\ell)\), by
Lemma~\ref{lem:Ggaps}; that proof also includes \(B=0\).
The next step from \(B\not\equiv0\pmod\ell\) changes the level
with conditional probability \(2/\ell\).
At a nonzero level no absorption is possible without such a
change.  Successive rounds therefore have conditional survival
probability at most \((\ell-2)/\ell\).
Sum their uniformly bounded conditional means to obtain
\[
 \sum_{i\ge1}
 \left(\frac{\ell-2}{\ell}\right)^{i-1}
 \bigl(\ell^{\ell-1}(t_k+\ell)+1\bigr)=\theta_k .
\]
No independence of round lengths is asserted or needed.
\end{proof}

\begin{lemma}[Averaged valuation along a killed walk]\label{lem:Gtk}
Let \(S_i\) be simple symmetric random walk started at \(1\) or
\(-1\), and let \(H_0\) be its first hitting time of zero.
There is \(C\) such that
\[
 \mathbb E_{\pm1}[t_{S_i};H_0>i]\le C/\sqrt i
 \qquad(i\ge1).
\]
This estimate does not require a pointwise logarithmic bound
for \(t_k\).
\end{lemma}

\begin{proof}
We use the elementary averaged bound of Lemma~\ref{lem:Ggaps}:
\[
 \sum_{1\le|k|\le R}t_k\le A_0R\qquad(R\ge1).
\]
The killed-walk estimate in Lemma~\ref{lem:meander} implies,
for an absolute constant \(C_0\),
\[
 \mathbb P_1(S_i=k,H_0>i)
 \le C_0\, k\,i^{-3/2}e^{-k^2/(16i)}
 \qquad(k\ge1,\ i\ge1),
\]
after increasing \(C_0\) to include \(k=1\).
Set \(R=\sqrt i\).
For \(1\le k\le R\),
\[
 \sum_{k\le R}t_k k e^{-k^2/(16i)}
 \le R\sum_{k\le R}t_k\le A_0i.
\]
On the dyadic interval \(2^jR<k\le2^{j+1}R\), \(j\ge0\),
the same sum is at most
\[
 2^{j+1}R\,e^{-4^j/16}
       \sum_{k\le2^{j+1}R}t_k
 \le 4A_0i\,4^j e^{-4^j/16}.
\]
The sum over \(j\) converges.  Multiplying by \(C_0i^{-3/2}\)
proves the assertion from \(1\).
Reflect the walk and use the corresponding averaged bound for
\(t_{-k}\) for the start at \(-1\).
\end{proof}

\begin{lemma}[Return gaps]\label{lem:Ggap}
On \(E_j\), put
\[
 \zeta_j=\min\{\tau_{j+1},\tau_{\rm abs}\}-\tau_j.
\]
There is \(C\) such that, for every \(n\ge1\),
\[
 \mathbb E[\min\{\zeta_j,n\}\mid\mathcal F_{\tau_j}]
 \le C(1+\ell_{\tau_j}+\sqrt n).
\]
\end{lemma}

\begin{proof}
At \(\tau_j\), the state is \((0,B)\) with \(B\ne0\).
Lemma~\ref{lem:Gsoj}(a) bounds the expected time until absorption
or departure from level zero by \(C(1+\ell_{\tau_j})\).
If the chain departs, its subsequent levels at change times
form a simple symmetric walk started at \(1\) or \(-1\).
Indeed, at each state from which a level change is possible,
the two directions have equal conditional probability.
Conditioning over any preceding holding steps preserves this
fairness.

Let \(\sigma_i\) be the starting time of the holding interval
at embedded level \(S_i\), and let \(H_0\) be the first embedded
hit of zero.  At \(\sigma_i\), the event \(H_0>i\) is measurable
and the conditional expected holding time is at most
\(\theta_{S_i}\), by Lemma~\ref{lem:Gsoj}(b).
Each interval takes at least one physical step.
Thus the portion of the excursion before a truncation at \(n\)
is bounded by the sum of full holding intervals with
\(i<\min\{H_0,n\}\).
The strong Markov property gives
\[
 \mathbb E[\min\{\zeta_j,n\}\mid\mathcal F_{\tau_j}]
 \le C(1+\ell_{\tau_j})
      +\sum_{0\le i<n}\mathbb E_{\pm1}[\theta_{S_i};H_0>i].
\]
Now \(\theta_k\le C(1+t_k)\),
\(\mathbb P_{\pm1}(H_0>i)\ll(i+1)^{-1/2}\), and
Lemma~\ref{lem:Gtk} applies.
Summing \((i+1)^{-1/2}\) proves the result.
This proof needs no independence between holding times and
later directions.
\end{proof}

\begin{lemma}[Survival through returns]\label{lem:Gsurv}
For \(L_c=1+\log(2+|c|)\), there is \(C\) such that
\[
 \sum_{j\ge0}\mathbb P(E_j)\le CL_c,\qquad
 \sum_{j\ge0}\mathbb E[\mathbf1_{E_j}\ell_{\tau_j}]
       \le CL_c^2.
\]
\end{lemma}

\begin{proof}
Let \(0=V_0<V_1<\cdots\) enumerate all physical times at level
zero, with fictitious successive visits of offset zero after
absorption.  Put \(Y_i=|p_{V_i}|^s\), with \(s,\rho_V,C_V\)
as in Lemma~\ref{lem:GV}.
Choose \(M\ge1\) with \(M^s\ge2C_V/(1-\rho_V)\), and put
\[
 \theta=(1+\rho_V)/2<1,\quad
 \lambda=\tfrac12\log(1/\theta),\quad
 \kappa_0=1+\sqrt\theta.
\]
A visit is good when \(Y_i\le M^s\).
At a bad visit the conditional expected next \(Y\) is at most
\(\theta Y_i\), while at a good visit it is at most
\(\theta M^s\).

Let \(R_0\) count the bad visits before the first good visit,
and let \(R_a\), \(a\ge1\), count the bad visits between successive
good visits.  Stop the drift estimate at the next good visit.
For integer \(r\ge1\), it gives
\[
 \mathbb P(R_0\ge r)\le(1+|c|^s)\theta^{r-1},\qquad
 \mathbb P(R_a\ge r\mid\hbox{preceding good-visit history})
       \le\theta^r.
\]
Summing these integer-valued tails yields
\[
 \mathbb Ee^{\lambda R_0}\le C_0(1+|c|^s),\qquad
 \mathbb E[e^{\lambda R_a}\mid\hbox{preceding history}]
       \le\kappa_0 .
\]
The tower property and Markov's inequality consequently give
\[
 \mathbb P(\hbox{fewer than \(g\) good visits among the first \(v\)})
 \le C_0(1+|c|^s)\kappa_0^g e^{-\lambda(v-g)}.
\]

The good states \((0,B)\), \(|B|\le M\), form a finite set.
By Lemma~\ref{lem:Gdesc}, there are an integer \(a_0\ge1\) and
\(\delta_0>0\) such that absorption from each nonabsorbing good
state occurs within \(a_0\) steps with probability at least
\(\delta_0\).
Start an attempt at the first good visit and each later one
at the first good visit strictly after the preceding starting
time plus \(a_0\).
The attempts do not overlap; each has the stated conditional
success probability and excludes at most \(a_0+1\) good visits
as possible later starting points.

Write \(N_i=\{\tau_{\rm abs}>V_i\}\).
Choose a fixed \(b>0\) with
\(b(\log\kappa_0+\lambda)<\lambda/2\).
For \(i>a_0\), apply the good-visit bound to
\(v=i-a_0\) and \(g=\lfloor bv\rfloor\).
If there are at least \(g\) good visits, at least
\(\lfloor g/(a_0+1)\rfloor\) attempts finish by \(V_i\).
On \(N_i\), they all fail.
Combining the two estimates and absorbing the small indices gives
\[
 \mathbb P(N_i)\le C_1(1+|c|^s)e^{-\beta i}
 \qquad(i\ge0)
\]
for a fixed \(\beta>0\).
Every nonabsorbing return is a visit, so
\(\tau_j=V_{i(j)}\) with \(i(j)\ge j\) and distinct \(i(j)\).
Thus \(E_j\subseteq N_j\).
Summing the minimum of one and the last exponential bound proves
the first assertion.

Iteration of Lemma~\ref{lem:GV} gives, uniformly in \(i\),
\[
 \mathbb E|p_{V_i}|^s
 \le |c|^s+\frac{C_V}{1-\rho_V}=:B_c.
\]
Concavity of the logarithm and \(0<s\le1\) imply
\[
 \mathbb E\log(2+|p_{V_i}|)
 \le s^{-1}\log(2^s+B_c)\le C L_c.
\]
Also \(\log(2+x)\le C_s(1+x^{s/2})\), so Cauchy--Schwarz
and the survival estimate imply
\[
 \mathbb E[\mathbf1_{N_i}\log(2+|p_{V_i}|)]
 \le C\bigl(\mathbb P(N_i)+\sqrt{\mathbb P(N_i)B_c}\bigr)
 \le C(1+|c|^s)e^{-\beta i/2}.
\]
Use the preceding \(CL_c\) bound for
\(i\le\lceil(2/\beta)\log(1+|c|^s)\rceil\) and the
exponential bound thereafter.  The total is \(O(L_c^2)\).
Finally, the injective assignment \(j\mapsto i(j)\) gives
pathwise
\[
 \sum_j\mathbf1_{E_j}\ell_{\tau_j}
 \le\sum_i\mathbf1_{N_i}\log(2+|p_{V_i}|).
\]
Tonelli's theorem completes the second assertion.
\end{proof}

\begin{proposition}[The fixed-shift order of the merging time]
\label{thm:Gexact}
Under \eqref{eq:H}, there is an effectively computable
\(C_{\ell,m}\) such that, for \(c\ne0\), \(n\ge1\), and
\(L_c=1+\log(2+|c|)\),
\[
 \sum_{t=0}^{n-1}f_t(c)
 \le C_{\ell,m}(L_c\sqrt n+L_c^2),
 \qquad
 \frac2\ell\sqrt{\frac2{\pi(n+1)}}\le f_n(c)
 \le2C_{\ell,m}
       \left(\frac{L_c}{\sqrt n}+\frac{L_c^2}{n}\right).
\]
In particular, for each fixed nonzero shift, the tail has order
\(n^{-1/2}\).  No theorem on linear forms in \(p\)-adic logarithms
is used in this proposition.
\end{proposition}

\begin{proof}
The lower bound is Theorem~\ref{thm:lowerunif}; it does not
require negative drift.
Before absorption, the intervals
\([\tau_j,\min\{\tau_{j+1},\tau_{\rm abs}\})\), on \(E_j\),
partition physical time.  This includes possible absorption
while remaining at level zero.
Therefore
\[
 \min\{\tau_{\rm abs},n\}
 \le\sum_{j\ge0}\mathbf1_{E_j}\min\{\zeta_j,n\}.
\]
Take expectations and use Lemmas~\ref{lem:Ggap} and~\ref{lem:Gsurv}:
\[
 \sum_{t=0}^{n-1}f_t(c)
 =\mathbb E\min\{\tau_{\rm abs},n\}
 \le C(L_c\sqrt n+L_c^2).
\]
Since \(f_t(c)\) is nonincreasing in \(t\), division by \(n\)
gives the upper tail bound, even without the harmless displayed
factor two.  Effectivity follows from the elementary estimates,
the explicit drift bounds, and a finite search for the paths
from the finite good set in Lemma~\ref{lem:Gdesc}.
\end{proof}

\begin{lemma}[A truncated return estimate from an arbitrary level]
\label{lem:Greturntime}
Let the coalescence chain start from \((k,B)\), with \(k\ne0\),
and let \(R_0=\inf\{t\ge0:k_t=0\}\).
There is an effective constant \(C\), depending only on the fixed
map, such that
\[
 \mathbb E_{(k,B)}\min\{R_0,n\}\le C|k|\sqrt n
 \qquad(n\ge1),
\]
uniformly in \(B\).  For \(k=0\), take \(R_0=0\).
\end{lemma}

\begin{proof}
Put \(a=|k|\).  If \(a\ge\sqrt n\), the trivial bound by \(n\)
suffices.  Otherwise put \(H=\lceil\sqrt n\rceil\), and stop the
embedded simple symmetric walk on the side of \(k\) when its
absolute level first reaches \(0\) or \(H\).
The expected number of visits to \(j\in\{1,\ldots,H-1\}\),
starting at \(a\), is
\[
 G_H(a,j)=\frac{2\min(a,j)(H-\max(a,j))}{H}\le2a,
\]
and the probability of exit at \(H\) is \(a/H\).
Let \(\theta_j\) be the uniform bound on the expected physical
holding time at level \(j\) supplied by Lemma~\ref{lem:Ggaps}.
Its averaged estimate implies
\[
 \sum_{j=1}^{H-1}\theta_{\operatorname{sgn}(k)j}\le A H.
\]
Conditioning on the history at the start of each holding interval
therefore bounds the expected physical exit time \(S\) by
\(2AaH\).  This conditioning uses no independence between a
holding time and the direction of its eventual exit.
Pathwise,
\[
 \min\{R_0,n\}
 \le S+n\mathbf1_{\{\text{exit at absolute level }H\}}.
\]
Taking expectations, and using \(H\le2\sqrt n\), gives
\(2AaH+na/H\le(4A+1)a\sqrt n\).
\end{proof}

\begin{theorem}[A uniform square-root tail bound]
\label{thm:Gsharp}
Assume \eqref{eq:H} and use the canonical residues in
\eqref{eq:Tlm}.  There is an effectively computable constant
\(C\ge1\), depending only on \((\ell,m)\), such that, for every
integer \(c\ne0\), every \(n\ge1\), and
\(L_c=1+\log(2+|c|)\),
\begin{equation}\label{eq:Gsharp}
 \mathbb E_c\min\{\tau_{\rm abs},n\}
       \le C\sqrt{L_c n},
 \qquad
 f_n(c)\le\min\left\{1,C\sqrt{\frac{L_c}{n+1}}\right\}.
\end{equation}
In conjunction with Theorem~\ref{thm:lowerunif}, this gives
\[
 \frac2\ell\sqrt{\frac2{\pi(n+1)}}
 \le f_n(c)
 \le\min\left\{1,C\sqrt{\frac{L_c}{n+1}}\right\}.
\]
The conclusion gives the order \(n^{-1/2}\) for each fixed
nonzero shift.  It does not assert convergence of
\(\sqrt n\,f_n(c)\) to a specified leading coefficient.
\end{theorem}

\begin{proof}
Use a single sufficiently small exponent \(s\in(0,1]\) for
Lemmas~\ref{lem:GL}, \ref{lem:Gtight}, and~\ref{lem:GV}.
Write \(n_1,\rho,C_1\) for the constants of
Lemma~\ref{lem:GL}.  Iteration gives
\[
 \mathbb E_c|p_{jn_1}|^s
 \le\rho^j|c|^s+\frac{C_1}{1-\rho}.
\]
Choose a deterministic multiple \(N\) of \(n_1\), with
\(1\le N\le C_2L_c\), so that
\(\rho^{N/n_1}|c|^s\le1\).  Thus
\(\mathbb E_c|p_N|^s\le C_3\).
The level process, held at zero after absorption, is a martingale
whose increments have absolute value at most one.  Consequently,
\[
 \mathbb E_c|k_N|\le
 \bigl(\mathbb E_c k_N^2\bigr)^{1/2}\le\sqrt N.
\]

Let \(R=\inf\{t\ge0:k_{N+t}=0\}\), with \(R=0\) if
\(k_N=0\).  Conditioning at time \(N\),
Lemma~\ref{lem:Greturntime} and the arbitrary-level return
moment estimate \eqref{eq:Greturnmoment} imply
\begin{equation}\label{eq:Gburnreturn}
 \mathbb E_c\min\{R,n\}\le C_4\sqrt{Nn},
 \qquad
 \mathbb E_c|B_{N+R}|^s\le C_5.
\end{equation}
These bounds are uniform in the random initial level at time
\(N\); in particular, the second constant does not grow with
\(|k_N|\).

At time \(N+R\), apply the strong Markov property and
Proposition~\ref{thm:Gexact}.  If \(B_{N+R}=0\), the remaining
absorption time is zero.  For every other value \(b\), that
proposition bounds its expected truncation at \(n\) by
\[
 C_6\bigl(L_b\sqrt n+L_b^2\bigr),
 \qquad L_b=1+\log(2+|b|).
\]
For the fixed positive \(s\),
\(L_b^2\le C_s(1+|b|^s)\).
After averaging with \eqref{eq:Gburnreturn}, the expected
remaining truncation is therefore at most \(C_7\sqrt n\).
Splitting the time at \(N\) and \(N+R\) gives
\[
 \mathbb E_c\min\{\tau_{\rm abs},n\}
 \le N+C_4\sqrt{Nn}+C_7\sqrt n.
\]
For \(n\ge N\), use \(N\le\sqrt{Nn}\) and
\(N\le C_2L_c\).  For \(n<N\), use the trivial bound
\(n\le\sqrt{Nn}\).  This proves the first estimate in
\eqref{eq:Gsharp}.
Finally,
\[
 n f_n(c)\le\sum_{t=0}^{n-1}f_t(c)
       =\mathbb E_c\min\{\tau_{\rm abs},n\}.
\]
Enlarging the constant replaces \(n\) in the denominator by
\(n+1\), and gives the second estimate.
All constants are obtained from the effective estimates already
proved; no linear-form estimate is used here.
\end{proof}

\begin{corollary}[Dependence on the error tolerance]
\label{cor:Gsharphitting}
For the constant \(C\) in Theorem~\ref{thm:Gsharp},
\[
 n+1\ge C^2\varepsilon^{-2}\bigl(1+\log(2+|c|)\bigr)
 \quad\Longrightarrow\quad f_n(c)\le\varepsilon
 \qquad(0<\varepsilon<1).
\]
Thus an effective logarithmic-displacement coefficient of order
\(\varepsilon^{-2}\) is available.  For each fixed \(c\ne0\),
the least \(n\) for which \(f_n(c)\le\varepsilon\) has order
\(\varepsilon^{-2}\) as \(\varepsilon\downarrow0\).
\end{corollary}

\begin{proof}
The upper bounds follow from \eqref{eq:Gsharp}; the fixed-shift
lower bound follows from Theorem~\ref{thm:lowerunif}.
This is not a lower bound on the coefficient of
\(\log(2+|c|)\) in a two-parameter estimate that also allows an
arbitrary additive term depending on \(\varepsilon\).
\end{proof}

\subsection{Arithmetic progressions and renewal averages}

For the rest of this section, \(G\) is finite-defect invariant,
\[
 K=\max(G\triangle T^{-1}(G)),\qquad \max\varnothing=0.
\]
For \(d\) coprime to \(\ell\), the Chinese remainder theorem
applies to the period \(\ell^n\) of \(U_n(c)\).
The proof of Lemma~\ref{lem:Dodd}, with \(2^n\) replaced by
\(\ell^n\), consequently gives, for \(n\ge1\) and \(y\ge1\),
\[
 |dF^{[d]}_{a+c}(y)-dF^{[d]}_a(y)|
 \le f_n(c)y+d\ell^n+\ell^nK+3|c|+3d.
\]
Here finite-defect invariance is used only above \(\ell^nK\),
since the canonical map satisfies \(T(N)\ge N/\ell\).
Average over the unique representatives
\(-d/2<c\le d/2\); their mean absolute value is at most \(d/4\),
also when \(d\) is even.
For \(d=m^q\), this yields
\begin{equation}\label{eq:GDavg}
 |m^qF_{a,q}(y)-F(y)|
 \le \bar f_n(q)y+m^q\ell^n+\ell^nK+4m^q,
 \qquad
 \bar f_n(q)=\frac1{m^q}
       \sum_{-m^q/2<c\le m^q/2}f_n(c).
\end{equation}
In particular \(\bar f_n(0)=f_n(0)=0\).
For any \(Q\) coprime to \(\ell\), the corresponding maximum
bound is
\[
 |QF^{[Q]}_a(y)-F(y)|
 \le y\max_{|c|\le Q/2}f_n(c)+Q\ell^n(K+1)+4Q.
\]

Increasing \(L(\varepsilon)\) in Proposition~\ref{prop:GU}
preserves its conclusion.  In all the following statements we
make the fixed convention
\[
 L(\varepsilon)\ge
 \max\{1,\ 4\log\ell+2\log^+(m/\ell)\},
\]
by replacing its displayed value by the maximum with this bound
if needed.  Here \(\log^+t=\max\{0,\log t\}\).
The integer \(n_0(\varepsilon)\) may also be increased to be
nonnegative.

\begin{corollary}[Deep \(m\)-adic levels and coprime moduli]
\label{cor:Gdeep}
For \(0<\varepsilon<1\), put
\[
 \kappa_1=\frac1{2(1+L(\varepsilon)\log\ell)}.
\]
Then
\[
 \limsup_{x\to\infty}
 \max_{\substack{1\le Q\le x^{\kappa_1}\\(Q,\ell)=1}}
 \max_{a\bmod Q}
 \left|\frac{QF^{[Q]}_a(x)}x-\frac{F(x)}x\right|
 \le\varepsilon .
\]
In particular the bound holds for \(Q=m^q\) with
\(q\le\kappa_1\log x/\log m\).
\end{corollary}

\begin{proof}
Take \(n=\lceil L(\varepsilon)\log(2+Q)\rceil+n_0(\varepsilon)\).
Proposition~\ref{prop:GU} bounds each \(f_n(c)\) in the maximum
by \(\varepsilon\); for \(c=0\) it is zero.
Moreover,
\[
 Q\ell^n(K+1)+4Q
 \ll_{\varepsilon,\ell,m}
 (K+1)Q^{1+L(\varepsilon)\log\ell}
 \ll_{\varepsilon,\ell,m}(K+1)x^{1/2}.
\]
Divide the preceding progression bound by \(x\).
All estimates are uniform in \(Q,a\).
\end{proof}

Let
\[
 S_n=Q_n\log m-n\log\ell,\qquad
 Q_n\sim\operatorname{Bin}(n,\varpi),\qquad
 \mathbb ES_n=n\mu,
\]
and define the mixture
\(\nu_M=(M+1)^{-1}\sum_{n=M}^{2M}\mathcal L(S_n)\).

\begin{lemma}[Local bound and mixture]\label{lem:Gmix}
Put
\[
 \kappa_\ell=\frac{2e^{1/12}\sqrt\ell}{\sqrt{2\pi}},
 \quad
 C_4=\left(1+\frac1{\log m}\right)
       \frac{\kappa_\ell\sqrt\pi\log m}{|\mu|},
 \quad
 C_5=\left(1+\frac1{\log m}\right)
       \kappa_\ell\left(1+\frac1{|\mu|}\right).
\]
If \(t=q-n\varpi\) and \(|t|\le n/(2\ell)\), then
\[
 \mathbb P(Q_n=q)\le\kappa_\ell n^{-1/2}e^{-2t^2/n}.
\]
For every interval \(J\) of length at most one,
\[
 \nu_M(J)\le C_4/M+C_5/M^{3/2}
                     +2e^{-M/(2\ell^2)}.
\]
\end{lemma}

\begin{proof}
Robbins' factorial bounds~\cite{Robbins1955} give, with \(x=q/n\),
\[
 \mathbb P(Q_n=q)
 \le\frac{e^{1/12}}{\sqrt{2\pi}}
       \sqrt{\frac{n}{q(n-q)}}
       \exp\{-nD(x\Vert\varpi)\},
\]
where \(D\) is the binary relative entropy with natural logarithms.
The elementary inequality
\(D(x\Vert\varpi)\ge2(x-\varpi)^2\), obtained from
\(\partial_x^2D=1/(x(1-x))\ge4\), and
\(q(n-q)\ge n^2/(4\ell)\) in the indicated range prove the
local bound.  The endpoints \(q=0,n\) are outside that range.
Hoeffding's inequality~\cite{Hoeffding1963} bounds the remaining probability by
\(2e^{-n/(2\ell^2)}\).

Write \(a=1/\log m\), \(\delta=|\mu|/\log m\), and let \(y\)
be the left endpoint of \(J\).
For a fixed \(n\), there are at most \(1+a\) possible values of
\(S_n\) in \(J\).  Their centered indices \(t\) lie in
\[
 [\tau_n,\tau_n+a],\qquad \tau_n=(y-n\mu)/\log m.
\]
As \(n\) varies, \(\tau_n\) is an arithmetic progression of
positive step \(\delta\).
For \(M\le n\le2M\), set
\(\phi(v)=\sup_{t\in[v,v+a]}e^{-t^2/M}\).
The central contribution is at most
\((1+a)\kappa_\ell M^{-1/2}\phi(\tau_n)\).
The nonnegative function \(\phi\) is unimodal, with maximum one
and integral \(a+\sqrt{\pi M}\).
Its sum over an arithmetic progression of step \(\delta\)
is at most \(1+\delta^{-1}(a+\sqrt{\pi M})\).
Average over \(n\), replace \(1/(M+1)\) by \(1/M\), and
add the Hoeffding tail.  The resulting constants are exactly
\(C_4,C_5\) above.
\end{proof}

\begin{proposition}[A quantitative density comparison]
\label{thm:GN}
Assume \eqref{eq:H}.  Put
\[
 \eta_0=\min\{1/10,1/(2\ell)\},\qquad
 \kappa=4(\varpi+\eta_0)\log m\log\ell+1,\qquad
 C_{\ell,m}=C_4(\kappa+1).
\]
For every \(0<\varepsilon<1\), \(0<h\le1\), and
\(T_{\ell,m}\)-invariant set \(G\),
\[
 \underline d(G)\ge
 1-\frac{C_{\ell,m}L(\varepsilon)\overline\delta}{h}
       -h-\varepsilon,
\]
where \(\overline\delta\) is the upper logarithmic density of
\(\mathbb N\setminus G\).
In particular, logarithmic density one implies natural density one.
\end{proposition}

\begin{proof}
We first state the generalized form of the renewal estimate.
For an iteration length \(n\), a cutoff \(q^*\le n\), and integers
\(n_q\ge1\), \(0\le q\le q^*\), put \(n^*=\max n_q\).
Terras conjugacy and \eqref{eq:GDavg} imply
\[
 |\psi(U)-\mathbb E\psi(U+S_n)|
 \le\max_{q\le q^*}\bar f_{n_q}(q)
       +\mathbb P(Q_n>q^*)
       +C_{\rm R}(K+1)\ell^{n+n^*}e^{-U},
\]
where \(C_{\rm R}\) depends only on \(\ell,m\).
Indeed, partition the initial values into the \(\ell^n\) classes
modulo \(\ell^n\).  The images of a class are consecutive members
of one progression modulo \(m^q\), with cutoff
\(m^qe^U/\ell^n\).
Lemma~\ref{lem:GTerras} bounds its affine offset by less than
\(m^q\), so the endpoint error per class is bounded.
Starting values at most \(\ell^nK\) contribute at most \(K+1\)
per class.  Apply \eqref{eq:GDavg} to each retained value of \(q\)
and use a trivial count for the remaining values.
This is precisely the argument of Lemma~\ref{lem:R}, with
constants depending only on \(\ell,m\).

Fix \(\varepsilon,h\), put \(L=L(\varepsilon)\), and choose
\[
 n_q=\lceil L\log(2+m^q)\rceil+n_0(\varepsilon),\quad
 M=\left\lfloor\frac{U}{\kappa L}\right\rfloor,\quad
 q_n^*=\lfloor(\varpi+\eta_0)n\rfloor.
\]
For \(M\le n\le2M\), Proposition~\ref{prop:GU} gives
\(\bar f_{n_q}(q)\le\varepsilon\), while Hoeffding gives
\(\mathbb P(Q_n>q_n^*)\le e^{-2\eta_0^2n}\).
The choice of \(L\) ensures \(L\ge4\log\ell\), so
\[
 \begin{split}
 (n+\max_{q\le q_n^*}n_q)\log\ell
 &\le
 \frac{2U}{\kappa}
 \left(\frac{\log\ell}{L}
       +(\varpi+\eta_0)\log m\log\ell\right)+C_\varepsilon\\
 &\le U/2+C_\varepsilon .
 \end{split}
\]
Thus the remaining renewal errors vanish uniformly.
With \(g=1-\psi_G\), averaging the renewal estimates gives
\[
 g(U)\le\int g(U+v)\,d\nu_M(v)+\varepsilon+o(1).
\]
The arguments are contained in
\[
 \left[
 U\left(1-\frac{2\log\ell}{\kappa L}\right),
 U'=U\left(1+\frac{2\log^+(m/\ell)}{\kappa L}\right)
 \right].
\]
Their left endpoints tend to infinity.
Lemma~\ref{lem:mono}, Tonelli's theorem, and the local mass
bound just proved therefore yield
\[
 \limsup_{U\to\infty}g(U)
 \le
 \frac{C_4}{h}
 \bigl(\kappa L+2\log^+(m/\ell)\bigr)\overline\delta
       +h+\varepsilon.
\]
Here we used
\(\int_0^Vg(u)\,du\le(\overline\delta+o(1))V\),
which is the partial-summation calculation in
\eqref{eq:logint} and does not depend on the map.
Finally \(L\ge2\log^+(m/\ell)\), so the last coefficient is at
most \(C_4(\kappa+1)L=C_{\ell,m}L\).
\end{proof}

The following qualitative conclusion is already contained in
the stronger natural-density estimate of
\cite[Theorem~1.1 and Corollary~1.3]{NashidaGeneral2026}.
We retain its derivation from the comparison just proved and
the logarithmic-density input \eqref{eq:GGM}.

\begin{corollary}[Almost-bounded orbits in natural density]
\label{thm:GA}
Under \eqref{eq:H}, for every \(f:\mathbb N\to\mathbb R\)
with \(f(N)\to\infty\),
\[
 \#\{N\le x:T_{\min}(N)\ge f(N)\}=o(x).
\]
The same assertion holds for the unaccelerated maps \(C\)
with the canonical residues, as in Lemma~\ref{lem:CT}.
\end{corollary}

\begin{proof}
For \(B\ge1\), set \(G_B=\{N:T_{\min}(N)\le B\}\).
The identity
\(T_{\min}(N)=\min\{N,T_{\min}(T(N))\}\) shows
\[
 G_B\triangle T^{-1}(G_B)\subseteq[1,B].
\]
By Lemma~\ref{lem:CT} and the fixed-threshold logarithmic-density
estimate \eqref{eq:GGM}, the upper logarithmic density
\(\delta_B\) of its complement decreases to zero.
Apply Proposition~\ref{thm:GN} with \(h=\varepsilon\).
For fixed \(\varepsilon\), first let \(B\to\infty\), then let
\(\varepsilon\downarrow0\); this proves
\(\overline d(\mathbb N\setminus G_B)\to0\).
For each fixed \(B\), eventually \(f(N)>B\), so the exceptional
values in the assertion lie outside \(G_B\) apart from finitely
many integers.  Taking upper densities and then \(B\to\infty\)
proves the result.  Lemma~\ref{lem:CT} gives the assertion for \(C\).
\end{proof}

\begin{remark}\label{rem:Gexamples}
The constants in the general estimates are not optimized;
the standard \(3x+1\) case retains the sharper constants already
proved in Section~\ref{sec:natural}.
Hypothesis \eqref{eq:H} is precisely the multiplier condition
in the logarithmic-density theorem of
Gon\c{c}alves, Greenfeld and Madrid~\cite{GGM2025}.

Every individual cycle basin is exactly invariant, and hence
the estimates above apply to it.
For example, \(T_{5,4}\) fixes each of \(1,2,3,4\), and
\(T_{6,5}\) fixes each of \(1,2,3,4,5\), as is checked directly
from \eqref{eq:Tlm}.
These examples assert neither that the displayed cycles exhaust
all orbits nor that a finite counting experiment establishes a
limiting basin density.
\end{remark}

\subsection{Harmonic functions and slow oscillation}

For a composite integer \(m\), \(\mathbb Z_m\) denotes the
inverse limit of the rings \(\mathbb Z/m^t\mathbb Z\), equipped
with normalized Haar measure.  Its cylinders of level \(t\)
have mass \(m^{-t}\).  Since \((m,\ell)=1\), multiplication
by \(\ell\) is an automorphism of this compact group.
The branch maps
\[
 g_0(z)=z/\ell,\qquad
 g_j(z)=(mz+r_j)/\ell\quad(1\le j<\ell)
\]
are therefore well defined on \(\mathbb Z_m\).
The associated Markov operator is
\[
 Ph(u,z)=\frac1\ell h(u-\log\ell,g_0z)
       +\frac1\ell\sum_{j=1}^{\ell-1}
          h(u+\log(m/\ell),g_jz).
\]
Bounded harmonic functions satisfy \(Ph=h\) almost everywhere
for Lebesgue measure times Haar measure.
The nonzero branches can have coincident image classes;
all branch sums below retain their multiplicities.

\begin{lemma}[Renewal equation]\label{lem:Grenewal}
Let \(D_G=\#(G\triangle T^{-1}(G))\).
For every \(t\ge0\), every \(a\bmod m^t\), and \(x\ge1\),
\[
 \left|
 F_{a,t}(x)-F_{a\ell^{-1},t}(x/\ell)
 -\sum_{j=1}^{\ell-1}
       F_{g_j(a),t+1}\bigl((mx+r_j)/\ell\bigr)
 \right|\le D_G.
\]
The class \(g_j(a)\bmod m^{t+1}\) is independent of the
representative of \(a\bmod m^t\).
Replacing every odd-branch argument by \(mx/\ell\) changes
the total error by at most \(\ell-1\).
\end{lemma}

\begin{proof}
The branch \(\ell\mid N\) gives the bijection
\[
 N=\ell M,\quad N\equiv a\pmod{m^t}
 \quad\Longleftrightarrow\quad
 M\equiv a\ell^{-1}\pmod{m^t}.
\]
For \(j\ne0\), a positive target in the class
\(\ell M\equiv r_j\pmod m\) has the integral preimage
\[
 N=(\ell M-r_j)/m.
\]
It is positive because \(M\ge1\) and \(1\le r_j<\ell\).
Moreover \(mN\equiv-r_j\pmod\ell\), so \(N\equiv j\pmod\ell\).
The additional congruence \(N\equiv a\pmod{m^t}\) is equivalent to
\[
 M\equiv (ma+r_j)\ell^{-1}=g_j(a)\pmod{m^{t+1}},
\]
and \(N\le x\) is equivalent to \(M\le(mx+r_j)/\ell\).
These are exact bijections; if two branches have the same
target class, they still count different preimages.
Thus the displayed difference equals
\[
 \sum_{\substack{N\le x\\N\equiv a\;(\mathrm{mod}\;m^t)}}
       \bigl(\mathbf1_G(N)-\mathbf1_G(T(N))\bigr),
\]
whose absolute value is at most \(D_G\).
Changing \(a\) by \(m^t\) changes \(g_j(a)\) by
\(m^{t+1}\ell^{-1}\), proving well-definedness.
Finally \(0<r_j/\ell<1\), so replacing each argument by
\(mx/\ell\) loses at most one integer.
\end{proof}

\begin{corollary}[Liouville property and counting profiles]
\label{thm:Gliouville}
Under \eqref{eq:H}:
\begin{enumerate}
 \item[(a)] Every bounded harmonic function of the random chain
 is almost everywhere constant.
 \item[(b)] The map
 \[
 \mathcal T(u,z,\beta)=
 \bigl(u+\log(m^{e_j}/\ell),g_jz,g_j\beta\bigr),
 \qquad j=\beta\bmod\ell,
 \]
 on \(\mathbb R\times\mathbb Z_m\times\mathbb Z_\ell\)
 is ergodic in the nonsingular sense for the reference measure
 \[
 d\mathfrak m=e^u\,du\,dz\,d\beta.
 \]
 \item[(c)] For every finite-defect invariant \(G\) and every
 \(\Lambda>1\),
 \[
 \sup_{\lambda\in[\Lambda^{-1},\Lambda]}
 \left|\frac{F(\lambda x)}{\lambda x}-\frac{F(x)}x\right|
 \longrightarrow0.
 \]
 Also
 \[
 \frac{m^tF_{a,t}(x)}x-\frac{F(x)}x\longrightarrow0
 \qquad(t\ge0,\ a\bmod m^t).
 \]
\end{enumerate}
\end{corollary}

\begin{proof}
We first check the measure properties, including composite bases.
The branch \(g_j\) sends Haar measure on its cylinder
\(j+\ell\mathbb Z_\ell\), normalized to mass one, to Haar measure
on \(\mathbb Z_\ell\).
On \(\mathbb Z_m\), \(g_0\) preserves Haar measure, while \(g_j\)
for \(j\ne0\) sends it to normalized Haar measure on
\[
 b_j+m\mathbb Z_m,\qquad b_j=r_j\ell^{-1}\pmod m.
\]
The three change-of-measure factors in a branch are
\(1/\ell\) from its \(\beta\)-cylinder,
\(\ell/m^{e_j}\) from the weighted real coordinate, and
\(m^{e_j}\mathbf1_{\text{image}}\) from the \(z\)-coordinate.
Consequently
\[
 \frac{d\mathcal T_*\mathfrak m}{d\mathfrak m}(u,z,\beta)
 =1+\sum_{j=1}^{\ell-1}\mathbf1_{b_j+m\mathbb Z_m}(z).
\]
This density is between one and \(\ell\).
In particular the map is nonsingular, although
\(\mathfrak m\) is not generally invariant.

For (a), choose a positive probability density
\(\varrho(u)\) on \(\mathbb R\) and start the system with law
\(\varrho(u)\,du\,dz\,d\beta\), equivalent to \(\mathfrak m\).
Terras conjugacy implies that the successive \(\beta\)-branch
digits are independent and uniform.
For a bounded harmonic \(h\), the sequence
\(h(u_n,z_n)\) is therefore a bounded martingale.
Nonsingularity allows a common conull set on which all needed
iterates avoid the exceptional set in \(Ph=h\).
Its almost-sure limit \(F(u,z,\beta)\) is bounded and
\(\mathcal T\)-invariant, and
\[
 h(u,z)=\int_{\mathbb Z_\ell}F(u,z,\beta)\,d\beta
 \quad\hbox{almost everywhere}.
\]

Fix \(c\in\mathbb Z\).
By Theorem~\ref{thm:GK}, for almost every \(\beta\), the
orbits of \(\beta\) and \(\beta+c\) coalesce with equal
multiplication counts at a common time \(n\).
If their branch words are \(w,w'\), their affine compositions
have the same rational slope \(m^q/\ell^n\).
Equality at \(\beta\) therefore implies the affine identity
\[
 \gamma_{w'}(z+c)=\gamma_w(z).
\]
Their real-coordinate increments are also equal.
Thus the full orbits starting from
\((u,z,\beta)\) and \((u,z+c,\beta+c)\) agree from time \(n\).
Translation in \(z,\beta\) preserves the reference measure,
so, after discarding null sets, their limiting values \(F\)
are equal.  Integrate in \(\beta\) to get
\[
 h(u,z+c)=h(u,z)\quad\hbox{almost everywhere}.
\]
The integer shifts are countable, so Fubini gives a common
full-measure set of \(u\)'s for these identities.
For such \(u\), translation continuity in \(L^1(\mathbb Z_m)\)
extends them from the dense subgroup \(\mathbb Z\) to
\(\mathbb Z_m\).  Hence \(h(u,z)=H(u)\) almost everywhere.
Its harmonic equation reduces to
\[
 H(u)=\frac1\ell H(u-\log\ell)
      +\frac{\ell-1}{\ell}H(u+\log(m/\ell)).
\]
The closed group generated by the two shifts is all of
\(\mathbb R\), since \(\log m/\log\ell\) is irrational:
a rational relation would give \(m^a=\ell^b\) for positive
integers \(a,b\), contrary to \((m,\ell)=1\).
The real Choquet--Deny argument of Lemma~\ref{lem:CDR},
with the displayed probability weights, makes \(H\) constant.
This proves (a), separately on real and imaginary parts
if the original function was complex valued.

For (b), let \(F\) be the indicator of an invariant measurable set,
and put \(h_F(u,z)=\int F(u,z,\beta)\,d\beta\).
Branchwise change of variables in \(\beta\) shows \(Ph_F=h_F\).
By (a), \(h_F=c\) almost everywhere.
Conditionally on the initial \(u,z\) and the first \(n\) branch
digits of \(\beta\), the future \(\beta\)-coordinate is Haar
distributed and independent.  Invariance therefore gives
\[
 \mathbb E[F\mid u,z,\hbox{first \(n\) branch digits}]
       =h_F(u_n,z_n)=c.
\]
The digit cylinders generate the Borel sigma-field of
\(\mathbb Z_\ell\).  Conditional-expectation convergence under
the equivalent probability law then gives \(F=c\) almost
everywhere.  Since \(F\) is an indicator, \(c\) is zero or one,
which is nonsingular ergodicity.

For (c), define
\(\psi_{a,t}(u)=m^te^{-u}F_{a,t}(e^u)\).
The renewal lemma, after replacing its affine cutoffs by
\(me^u/\ell\), gives
\[
 \psi_{a,t}(u)
 =\frac1\ell\psi_{a\ell^{-1},t}(u-\log\ell)
 +\frac1\ell\sum_{j=1}^{\ell-1}
       \psi_{g_j(a),t+1}(u+\log(m/\ell))
 +O_{G,\ell,m}(m^te^{-u}).
\]
The profiles satisfy
\[
 0\le\psi_{a,t}(u)\le1+m^te^{-u},\qquad
 \psi_{a,t}(u)=\frac1m
       \sum_{b\equiv a\;(\mathrm{mod}\;m^t)}\psi_{b,t+1}(u),
\]
where the last sum has \(m\) terms.
Counting integers in an interval in a fixed class also gives,
for \(v\) in any fixed bounded interval,
\[
 |\psi_{a,t}(u+v)-\psi_{a,t}(u)|
 \le C|v|+C_{\Lambda}m^te^{-u}
 \quad(|v|\le\log\Lambda).
\]
These are the bounds used in the compactness proof of
Theorem~\ref{thm:profiles}.
More explicitly, from any sequence \(U_i\to\infty\), a diagonal
subsequence of the translates has locally uniform limits
\(h_{a,t}(v)\), for all cylinder indices.
The limits lie in \([0,1]\) and retain the displayed consistency.
For each \(v\), the functions constant on level-\(t\) cylinders
with values \(h_{a,t}(v)\) form a bounded martingale in \(t\).
Their almost-everywhere limit defines a jointly measurable
bounded function \(h(v,z)\); Fubini and bounded convergence
identify its cylinder averages with the \(h_{a,t}\).
The limiting renewal equations are exactly the cylinder averages
of \(Ph=h\), by the branchwise Haar formulas already verified.
Thus \(h\) is harmonic, and (a) makes it constant.
The cylinder limits are continuous in \(v\), so their
almost-everywhere equality to that constant holds everywhere.

Every subsequential limit is consequently constant in \(v,z\),
although the constant may depend on the subsequence.
If either assertion in (c) failed, a sequence witnessing a
positive discrepancy would have such a convergent subsequence.
Its constant limiting profiles contradict that discrepancy.
This proves both assertions, including uniformity over bounded
real shifts.
\end{proof}

The preceding skew product contains both the real cocycle and
the \(m\)-adic coordinate.  Its nonsingular ergodicity is a
different assertion from mixing of the digit map on
\(\mathbb Z_\ell\); the proof identifies the additional role
of balanced coalescence.

\subsection{Windows, equivalence of densities, and growing moduli}

Retain \(\kappa\) from Proposition~\ref{thm:GN}, the enlarged
\(L(\varepsilon)\), and set
\[
 c_W=\kappa_\ell
       \sqrt{\frac{2\varpi(1-\varpi)\kappa}{|\mu|}}.
\]

\begin{proposition}[Window comparison]\label{thm:GW}
Let \(G\) be \(T_{\ell,m}\)-invariant and \(0<\varepsilon<1\).
For
\[
 0<\gamma\le
 \min\left\{\frac12,\frac{|\mu|}{\kappa L(\varepsilon)}\right\},
\]
we have
\[
 \limsup_{U\to\infty}\sup_{0\le D\le\gamma U}
 |\psi_G(U)-\psi_G(U-D)|
 \le2\varepsilon+c_W\sqrt{L(\varepsilon)\gamma}.
\]
Consequently,
\[
 \lim_{\gamma\downarrow0}\limsup_{x\to\infty}
 \sup_{x^{1-\gamma}\le y\le x}
 \left|\frac{F(y)}y-\frac{F(x)}x\right|=0.
\]
\end{proposition}

\begin{proof}
Put \(L=L(\varepsilon)\), \(U_1=U-D\), and take
\[
 M=\left\lfloor\frac{U_1}{\kappa L}\right\rfloor,\qquad
 k=\lfloor D/|\mu|\rfloor .
\]
Uniformly for \(0\le D\le\gamma U\), \(M\to\infty\), and
\[
 M+k\le \frac{U}{\kappa L}+\frac{D}{|\mu|}
          \le\frac{2U}{\kappa L}.
\]
Also \(M\le U_1/(\kappa L)\).
The renewal estimate and exponential budget in the proof of
Proposition~\ref{thm:GN} therefore apply at time \(M+k\)
and scale \(U\), and at time \(M\) and scale \(U_1\).
They give
\[
 \psi(U)=\mathbb E\psi(U+S_{M+k})+R_1,\qquad
 \psi(U_1)=\mathbb E\psi(U_1+S_M)+R_2,
 \quad |R_i|\le\varepsilon+o(1),
\]
uniformly in \(D\).

Take independent \(Q\sim\operatorname{Bin}(M,\varpi)\) and
\(Q'\sim\operatorname{Bin}(k,\varpi)\), and put
\(Z=Q'-\lfloor k\varpi\rfloor\).
Then
\[
 U+S_{M+k}
 =U_1+(Q+Z)\log m-M\log\ell+r,
 \quad
 r=D-k|\mu|+(\lfloor k\varpi\rfloor-k\varpi)\log m.
\]
Thus \(|r|\le|\mu|+\log m\).
All the arguments before and after this bounded shift tend
to infinity uniformly, since
\[
 U+S_{M+k}\ge
 U\left(1-\frac{2\log\ell}{\kappa L}\right)\ge U/2.
\]
Corollary~\ref{thm:Gliouville}(c) removes \(r\) at cost \(o(1)\).

For large \(M\), a mode of \(\operatorname{Bin}(M,\varpi)\)
differs from \(M\varpi\) by at most one and lies in the central range of
Lemma~\ref{lem:Gmix}.  Its largest atom is therefore at most
\(\kappa_\ell M^{-1/2}\).
For a unimodal integer-valued law, its total variation distance
from a translate by \(z\in\mathbb Z\) is at most \(|z|\) times
its largest atom.
Since \(0\le\psi\le1\) and
\(\mathbb E|Z|\le\sqrt{k\varpi(1-\varpi)}+1\), the remaining
comparison is at most
\[
 \kappa_\ell M^{-1/2}
       \bigl(\sqrt{k\varpi(1-\varpi)}+1\bigr).
\]
The correct uniform ratio estimate is
\[
 \frac{k}{M}\le
 \frac{\kappa L\gamma}{|\mu|(1-\gamma)}+o(1).
\]
As \(\gamma\le1/2\), this proves the bound with the stated \(c_W\).
Adding \(R_1,R_2\) proves the first assertion.
For the second, fix \(\varepsilon\), let \(\gamma\downarrow0\)
within its admissible range, and then let
\(\varepsilon\downarrow0\).
\end{proof}

\begin{proposition}[Equivalence of densities]\label{thm:GT}
A \(T_{\ell,m}\)-invariant set has a natural density if and only
if it has a logarithmic density, and the two then agree.
\end{proposition}

\begin{proof}
Natural density implies logarithmic density by partial summation.
Conversely, if the logarithmic density is \(d\), then
\[
 A(U):=\frac1U\int_0^U\psi_G(v)\,dv\longrightarrow d.
\]
For fixed \(\lambda>1\), the average of \(\psi_G\) over
\([U,\lambda U]\) equals
\((\lambda A(\lambda U)-A(U))/(\lambda-1)\) and tends to \(d\).
By Proposition~\ref{thm:GW}, its difference from \(\psi_G(U)\)
has upper limit tending to zero as \(\lambda\downarrow1\).
First let \(U\to\infty\), then let \(\lambda\downarrow1\).
This proves \(\psi_G(U)\to d\).
\end{proof}

\begin{theorem}[Uniformity in a power range of moduli]\label{thm:GLD}
For \(0<\varepsilon<1\), set
\[
 \begin{split}
 \kappa_1&=\frac1{2(1+L(\varepsilon)\log\ell)},\\
 \gamma&=\min\left\{\frac12,\
        \frac{|\mu|}{\kappa L(\varepsilon)},\
        \frac{\varepsilon^2}{c_W^2L(\varepsilon)}\right\},\\
 \kappa^*&=\min\left\{
        \frac{\kappa_1}{2(1+\log_2m)},\
        \frac{\gamma\log2}{\log(\ell m)}\right\}.
 \end{split}
\]
Then, for every \(T_{\ell,m}\)-invariant set \(G\),
\[
 \limsup_{x\to\infty}
 \max_{1\le Q\le x^{\kappa^*}}\max_{a\bmod Q}
 \left|\frac{QF^{[Q]}_a(x)}x-\frac{F(x)}x\right|
 \le4\varepsilon .
\]
This is a uniform individual-modulus estimate.
\end{theorem}

\begin{proof}
The moduli coprime to \(\ell\) are covered by
Corollary~\ref{cor:Gdeep}.
For arbitrary \(Q\), factor \(Q=Q_\ell Q'\), where every prime
factor of \(Q_\ell\) divides \(\ell\) and \((Q',\ell)=1\).
Suppose \(Q_\ell>1\), and let \(k\) be the smallest integer
with \(Q_\ell\mid\ell^k\).
Prime valuations give
\[
 k\le\log_2Q_\ell\le\log_2Q.
\]
The class \(a\bmod Q\) splits into
\(R=\ell^k/Q_\ell\) classes modulo
\(\widetilde Q=\ell^kQ'\).
Its normalized count is the average of the \(R\) normalized
counts with modulus \(\widetilde Q\).

Fix one refined class, specified by \(r\bmod\ell^k\) and
a class modulo \(Q'\), and write \(N=\ell^kY+r\),
\(0\le r<\ell^k\).
Since \((Q',\ell)=1\), \(Y\) runs through one class modulo \(Q'\).
Terras conjugacy maps it bijectively to a progression
modulo
\[
 D=m^qQ',\qquad q=q_k(r),\qquad y=m^qx/\ell^k.
\]
No coprimality of \(Q'\) and \(m\) is needed for this statement.
The affine endpoint comparison and invariance above
\(\ell^kK\) give
\[
 F^{[\widetilde Q]}_{\widetilde a}(x)
   =F^{[D]}_b(y)+O_{\ell,m}(K+1).
\]
Here \(\widetilde Q/x=D/y\).
Also \((D,\ell)=1\), and
\[
 D\le m^kQ\le Q^{1+\log_2m}\le x^{\kappa_1/2},
 \qquad
 y\ge x^{1-\kappa^*\log_2\ell}\ge x^{1/2}.
\]
The last inequality follows from
\(\kappa^*\log_2\ell\le\gamma\log\ell/\log(\ell m)\le1/2\).
Thus \(D\le y^{\kappa_1}\), and
Corollary~\ref{cor:Gdeep} gives an error at most
\(\varepsilon+o(1)\), uniformly over the refined classes.

Furthermore,
\[
 |\log y-\log x|
 \le k\log(\ell m)\le\gamma\log x.
\]
Apply Proposition~\ref{thm:GW} at the larger of
\(\log x,\log y\).
The difference of \(F(y)/y\) and \(F(x)/x\) is at most
\[
 2\varepsilon+c_W\sqrt{L\gamma}+o(1)
 \le3\varepsilon+o(1).
\]
The endpoint error is uniform as well:
\[
 \widetilde Q=\ell^kQ'
 \le Q^{\log_2\ell},\qquad
 \frac{\widetilde Q(K+1)}x
 \ll (K+1)x^{-1+\kappa^*\log_2\ell}=o(1).
\]
Combining these estimates for each refined class and averaging
proves the theorem.
\end{proof}

\subsection{Primes and sets accessible to sieves}

\begin{theorem}[Invariant sets and primes]\label{thm:GP}
For every \(T_{\ell,m}\)-invariant set \(E\) and
\(0<\varepsilon<1\),
\[
 \limsup_{x\to\infty}\frac{\pi_E(x)}{\pi(x)}
 \le6\log\ell\,L(\varepsilon)\overline d(E)
       +\frac{3\ell}{\varphi(\ell)}\varepsilon .
\]
In particular, for every \(f(N)\to\infty\),
\[
 \#\{p\le x:T_{\min}(p)\ge f(p)\}=o(\pi(x)).
\]
\end{theorem}

\begin{proof}
Fix \(\varepsilon\), and put
\[
 \eta=\frac1{3L(\varepsilon)\log\ell},\quad
 H=\lfloor x^\eta\rfloor,\quad
 n=\lceil L(\varepsilon)\log(2+H)\rceil+n_0(\varepsilon),
 \quad q=\ell^n.
\]
Our lower convention for \(L\) ensures \(0<\eta<1/2\).
For large \(x\), \(q\ll_{\varepsilon,\ell,m}x^{1/3}\), and
Proposition~\ref{prop:GU} gives \(f_n(c)\le\varepsilon\)
for \(1\le c\le H\).
Write \(K=K(E)\).
If \(N>qK\), all its iterates through time \(n\), and those
of \(N+c\), exceed \(K\), since \(T(N)\ge N/\ell\).
Balanced coalescence by time \(n\) therefore makes their
membership indicators agree, giving
\[
 \left|\mathbf1_E(N)-\frac1H\sum_{c=1}^H\mathbf1_E(N+c)\right|
 \le\frac1H\sum_{c=1}^H\mathbf1_{U_n(c)}(N)
 \qquad(N>qK).
\]

The set \(U_n(c)\) consists of \(qf_n(c)\) classes modulo \(q\).
All primes in nonreduced classes divide \(\ell\), so their total
number is at most \(\omega(\ell)\), independently of \(n,c,x\).
For each reduced class, \eqref{eq:MV} gives at most
\[
 (3+o(1))\,x/(\varphi(q)\log x)
\]
primes, because \(\log((x-1/2)/q)\ge(2/3-o(1))\log x\).
For a composite \(\ell\) as well as a prime \(\ell\),
\[
 \varphi(\ell^n)=\ell^n\prod_{p\mid\ell}(1-1/p)
              =\ell^{n-1}\varphi(\ell).
\]
Thus, uniformly for \(1\le c\le H\),
\[
 \#\{p\le x:p\in U_n(c)\}
 \le\left(\frac{3\ell}{\varphi(\ell)}+o(1)\right)
         \varepsilon x/\log x+\omega(\ell).
\]

The smoothed sum is exactly
\[
 \frac1H\sum_{v\le x+H}\mathbf1_E(v)
   \#\bigl(\mathbb P\cap[1,x]\cap[v-H,v-1]\bigr).
\]
Its terms with \(v\le H+1\) contribute \(O(H)\);
the interval Brun--Titchmarsh bound for the others gives at most
\((2/\log H)F_E(x+H)\).
The initially excluded primes cost at most
\(O_{\varepsilon,\ell,m}(Kx^{1/3})\).
Divide by \(\pi(x)\sim x/\log x\), use
\(\log H/\log x\to\eta\) and \(H=o(x)\), and take the upper limit.
Since \(2/\eta=6L(\varepsilon)\log\ell\), this proves the bound.

For the final assertion, apply it to
\(E_B=\mathbb N\setminus G_B\).
The proof of Corollary~\ref{thm:GA} gives
\(\overline d(E_B)\to0\).
First let \(B\to\infty\) for fixed \(\varepsilon\), then let
\(\varepsilon\downarrow0\).
For each fixed \(B\), the asserted exceptional set is contained
in \(E_B\) apart from finitely many primes.
\end{proof}

For the following uniform assertion, retain
\[
 K(E)=\max(E\triangle T^{-1}(E)),\qquad\max\varnothing=0.
\]

\begin{theorem}[Invariant sets among primes under RH]\label{thm:GPE}
Assume the Riemann hypothesis.
For every fixed \(K_0\ge0\),
\[
 \sup_{\substack{E\subseteq\mathbb N\\K(E)\le K_0}}
 \left|\frac{\pi_E(x)}{\pi(x)}-\frac{F_E(x)}x\right|
 \longrightarrow0.
\]
The supremum is over finite-defect \(T_{\ell,m}\)-invariant sets.
Thus any such set of natural density \(d\) has relative density
\(d\) among the primes; upper and lower densities also agree
separately.
\end{theorem}

\begin{proof}
Choose \(\varepsilon,\eta,H,n,q\) as in the preceding proof and
write \(\vartheta_E(x)=\sum_{p\le x,\,p\in E}\log p\).
Multiply its local-constancy inequality by \(\log p\) and sum.
The same progression estimate bounds the resulting discrepancy
by \(O_\ell(\varepsilon x)+o_{\varepsilon,\ell,m}(x)\),
apart from the initial
\(O_{\varepsilon,\ell,m}(K_0x^{1/3}\log x)\).
Reverse the sums, retain the prime cutoff \(p\le x\), and
discard \(x<v\le x+H\) at cost \(O(H\log x)\).
Uniformly in \(E\), the remaining main term is
\[
 \frac1H\sum_{v\le x}\mathbf1_E(v)
       \bigl(\vartheta(v-1)-\vartheta(v-H-1)\bigr).
\]

The short-interval comparison is independent of the map.
As proved in Theorem~\ref{thm:PE}, Selberg's mean-square
estimate under RH and Cauchy--Schwarz give, for this fixed
power \(H=\lfloor x^\eta\rfloor\),
\[
 \frac1H\sum_{v\le x}
 |\vartheta(v-1)-\vartheta(v-H-1)-H|=o(x).
\]
In that deduction, prime powers are removed in \(L^1\):
each \(p^a\), \(a\ge2\), occurs in at most \(H\) of the
integer windows, and therefore their total contribution
after division by \(H\) is at most
\(\sum_{p^a\le x,\,a\ge2}\log p\ll\sqrt x\log^2x=o(x)\).
The initial \(v\le H+1\) terms cost \(O(H\log x)=o(x)\).

It follows that
\[
 \vartheta_E(x)=F_E(x)+O_\ell(\varepsilon x)
                    +o_{\varepsilon,\ell,m,K_0}(x),
\]
uniformly over \(K(E)\le K_0\).
Let \(x\to\infty\) and then \(\varepsilon\downarrow0\).
Finally,
\(\pi_E(x)=\vartheta_E(x)/\log x+O(x/\log^2x)\)
uniformly in \(E\), by partial summation and
\(0\le\vartheta_E(t)\le\vartheta(t)\ll t\).
The prime number theorem proves the assertion.
\end{proof}

\begin{theorem}[A sieve transfer principle]\label{thm:GS}
Let \(\delta\) satisfy the comparison and polylogarithmic lower
bound conditions in Theorem~\ref{thm:S}.
Let \(S\subseteq\mathbb N\) satisfy its conditions
\textup{(S1)} and \textup{(S2)}, with constants \(c_S,C_S\)
independent of the short-interval exponent.
Replace its progression hypothesis by
\[
 \begin{array}{ll}
 \textup{(S3\(_\ell\))}&
 \displaystyle
 \#(S\cap[1,x]\cap(b+\ell^n\mathbb Z))
 \le C_Sx\delta(x)/\ell^n\\
 &\text{for all \(b\) and integers \(n\ge0\) with
                  \(\ell^n\le x^{1/2}\).}
 \end{array}
\]
Then every \(T_{\ell,m}\)-invariant \(E\) satisfies
\[
 \limsup_{x\to\infty}
 \frac{\#(S\cap E\cap[1,x])}{S(x)}
 \le\frac{C_S}{c_S}
 \left(K_\delta(\eta_\varepsilon)\overline d(E)+\varepsilon\right),
 \qquad
 \eta_\varepsilon=\frac1{3L(\varepsilon)\log\ell}.
\]
Consequently, for every \(f(N)\to\infty\),
\[
 \#\{N\le x:N\in S,\ T_{\min}(N)\ge f(N)\}=o(S(x)).
\]
\end{theorem}

\begin{proof}
Fix \(\varepsilon\) and take \(H,n,q\) as in
Theorem~\ref{thm:GP}.
For sufficiently large \(x\), \(q=\ell^n\le x^{1/2}\).
Sum the local-constancy inequality over \(N\in S\), \(N\le x\).
The \(qf_n(c)\) exceptional residue classes, together with
\textup{(S3\(_\ell\))}, contribute at most
\(C_S\varepsilon x\delta(x)\).
The initially excluded values cost
\(O_{\varepsilon,\ell,m}(K(E)x^{1/3})\).

The smoothed expression equals
\[
 \frac1H\sum_{v\le x+H}\mathbf1_E(v)
     \#\bigl(S\cap[1,x]\cap[v-H,v-1]\bigr).
\]
The terms \(v\le H^2+H+1\) contribute \(O(H^2+H)\).
For the others, \textup{(S2)} bounds the kernel by
\(C_SH\delta(H)\).  Hence
\[
 \#(S\cap E\cap[1,x])
 \le C_S\delta(H)F_E(x+H)
      +C_S\varepsilon x\delta(x)
      +O_{\varepsilon,\ell,m}(K(E)x^{1/3}+x^{2\eta_\varepsilon}).
\]
The explicit lower bound on \(\delta(x)\) makes the last error
\(o(x\delta(x))\).
Divide by \textup{(S1)} and use
\(\limsup\delta(H)/\delta(x)\le K_\delta(\eta_\varepsilon)\).
This proves the inequality.
The conclusion follows by applying it to \(E_B\), using
\(\overline d(E_B)\to0\), and letting first \(B\to\infty\) and
then \(\varepsilon\downarrow0\).
Finitely many exceptions from \(f(N)\le B\) are negligible.
\end{proof}

\begin{corollary}[Rough Chen primes and sums of two squares]
\label{cor:GChen}
For every \(f(N)\to\infty\),
\[
 \begin{split}
 \#\{p\le x:p\in\mathcal C,\ T_{\min}(p)\ge f(p)\}
       &=o(\mathcal C(x)),\\
 \#\{N\le x:N\in\mathcal B,\ T_{\min}(N)\ge f(N)\}
       &=o(\mathcal B(x)).
 \end{split}
\]
Here \(\mathcal C\) is the set of \emph{rough} Chen primes
defined before Corollary~\ref{cor:Chen}, including its
\(p^{1/10}\) condition on both factors, and \(\mathcal B\)
is the set of sums of two squares.
\end{corollary}

\begin{proof}
Conditions \textup{(S1)} and \textup{(S2)} do not involve the map
and were proved in Corollary~\ref{cor:Chen}.
We verify \textup{(S3\(_\ell\))} with constants depending on
the fixed \(\ell\), uniformly over all residues.

First note the following extension of
Lemma~\ref{lem:uniform-sieve}.
If every prime factor of \(q\) divides \(\ell\), its two bounds
continue to hold, with constants depending only on \(\ell\),
after the sieving primes dividing \(\ell\) are omitted.
Indeed, for the two-square sieve, each
\(p\equiv3\pmod4\), \(p\nmid\ell\), still excludes exactly
\(p-1\) classes modulo \(p^2\) from the progression \(b+qt\).
For the twin-polynomial sieve, each odd \(p\nmid\ell\)
still excludes the two distinct classes
\(b+qt\equiv0,-2\pmod p\).
These facts hold for every \(b\), without a coprimality
condition on \(b\) and \(q\).
The Chinese remainder error and the Selberg weight bound in
that lemma remain \(O(\rho(d))\), \(\rho(d)\le d\),
and \(|\lambda_d|\le d\).
Consequently the remainder is still \(O(z^6)\), uniformly in
\(q,b\) and the interval position.
Deleting finitely many sieving primes changes the Mertens
products only by fixed positive factors.
The same Euler-product and Markov argument therefore gives
\[
 \mathfrak G(z)\gg_\ell(\log z)^{1/2}
 \quad\hbox{or}\quad
 \mathfrak G(z)\gg_\ell(\log z)^2,
\]
respectively.
This proves the asserted extension with no dependence on a
short-interval exponent.

Now let \(q=\ell^n\le\sqrt x\) and \(Y=x/q\ge\sqrt x\).
Apply the first extended sieve bound to the interval of
parameters \(t\) for which \(1\le b+qt\le x\).
It gives
\[
 \#(\mathcal B\cap[1,x]\cap(b+q\mathbb Z))
 \ll_\ell Y/\sqrt{\log Y}
 \ll_\ell x/(q\sqrt{\log x}).
\]
This includes the zero and all nonreduced residue classes.

For \(\mathcal C\), the primes \(p\le\sqrt x\) in the class
cost at most \(\sqrt x/q+1=o(x/(q\log^2x))\), uniformly
over \(q\le\sqrt x\).
For \(p>\sqrt x\), every prime factor of \(p(p+2)\) is at least
\(x^{1/20}\).
Choose \(z=Y^{1/24}\le x^{1/24}<x^{1/20}\).
These candidates survive the second extended sieve, which gives
\[
 \ll_\ell Y/\log^2Y+z^6
 =O_\ell(Y/\log^2Y+Y^{1/4})
 =O_\ell(x/(q\log^2x)).
\]
Thus \textup{(S3\(_\ell\))} holds for
\(\delta(x)=(\log x)^{-1/2}\) and
\(\delta(x)=(\log x)^{-2}\), respectively.
Apply Theorem~\ref{thm:GS} to conclude.
\end{proof}

The following corollary optimizes the preceding transfer estimate
using the stronger coalescence input of Theorem~\ref{thm:Gsharp}.

\begin{corollary}[H\"older bounds for sieve transfer]
\label{cor:Gtransferpowers}
Assume \eqref{eq:H} and use the canonical residues.
Let \(\alpha>0\), and suppose that \(S\) satisfies the hypotheses
of Theorem~\ref{thm:GS} with
\(\delta(x)=(\log x)^{-\alpha}\).
Define
\[
 \overline d_S(E)=
 \limsup_{x\to\infty}\frac{\#(S\cap E\cap[1,x])}{S(x)}.
\]
There is a constant \(C_{\ell,m,\alpha,S}\) such that every
finite-defect \(T_{\ell,m}\)-invariant set \(E\) satisfies
\[
 \overline d_S(E)
 \le C_{\ell,m,\alpha,S}\,
          \overline d(E)^{1/(2\alpha+1)}.
\]
The constant depends only on the fixed map, \(\alpha\), and
the constants in the sieve hypotheses; it is independent
of \(E\) and its finite defect cutoff.
In particular, there are constants depending only on
\((\ell,m)\) for which
\[
 \overline d_{\mathbb P}(E)\ll_{\ell,m}\overline d(E)^{1/3},
 \qquad
 \overline d_{\mathcal B}(E)\ll_{\ell,m}\overline d(E)^{1/2},
 \qquad
 \overline d_{\mathcal C}(E)\ll_{\ell,m}\overline d(E)^{1/5}.
\]
The population \(\mathcal C\) retains the roughness condition
in Corollary~\ref{cor:GChen}.
\end{corollary}

\begin{proof}
Let \(C\) be the tail-bound constant in Theorem~\ref{thm:Gsharp}.
For \(c\ne0\),
\[
 1+\log(2+|c|)\le2\log(2+|c|).
\]
Choose a fixed
\[
 A\ge
 \max\{2C^2,\ 1,\ 4\log\ell+2\log^+(m/\ell)\}.
\]
Then, for \(0<\varepsilon<1\),
\[
 n\ge A\varepsilon^{-2}\log(2+|c|)
 \quad\Longrightarrow\quad f_n(c)\le\varepsilon.
\]
For \(c=0\), this is automatic.
Thus the proof of Theorem~\ref{thm:GS} can be run with
\[
 L_\sharp(\varepsilon)=A\varepsilon^{-2},
 \qquad n_{0,\sharp}(\varepsilon)=0,\qquad
 \eta_\sharp=\frac{\varepsilon^2}{3A\log\ell}.
\]
All lower bounds imposed on \(L\) and the condition
\(0<\eta_\sharp<1/2\) are satisfied.
The constants in \textup{(S2)} and \textup{(S3\(_\ell\))}
are uniform in this choice of exponent.
Since \(K_\delta(\eta)=\eta^{-\alpha}\), writing
\(d=\overline d(E)\), the resulting bound is
\[
 \overline d_S(E)\le
 \frac{C_S}{c_S}
 \left((3A\log\ell)^\alpha\varepsilon^{-2\alpha}d
                         +\varepsilon\right).
\]
If \(d=0\), let \(\varepsilon\downarrow0\).
For \(0<d\le1\), take
\(\varepsilon=\tfrac12d^{1/(2\alpha+1)}\).
Both terms then have the factor \(d^{1/(2\alpha+1)}\);
for example the constant
\[
 \frac{C_S}{c_S}
       \left((12A\log\ell)^\alpha+\tfrac12\right)
\]
suffices.

For primes take \(\alpha=1\).
The prime number theorem supplies \textup{(S1)}, and
Brun--Titchmarsh supplies \textup{(S2)}.
For \(q=\ell^n\le\sqrt x\), it also gives the required
\(O_\ell(x/(q\log x))\) upper bound in every reduced class;
all primes in nonreduced classes divide \(\ell\), and their
bounded contribution is absorbed uniformly.
The case \(q=1\) follows from the usual prime counting bound.
Thus \textup{(S3\(_\ell\))} holds as well.
For \(\mathcal B\) and \(\mathcal C\), all hypotheses with
\(\alpha=1/2\) and \(\alpha=2\), respectively, were verified
in Corollary~\ref{cor:GChen}.
This proves the three displayed specializations.
\end{proof}

For the admissible signed residue systems covered by the
extension following Theorem~\ref{thm:GR}, with
\((\mathrm A_r)\) and \(d_r=1\), the same argument applies.
Its constants depend on the full residue system and the finite
absorption certificate, as does the extended strong tail bound.
The enlarged finite-defect cutoff described there affects
only terms that vanish in the upper-density limit.

\subsection{General residues and an obstruction}

Let $r_0=0$ and let $r_1,\dots,r_{\ell-1}$ be integers satisfying
$r_j\equiv-mj\pmod\ell$. Define
\[
 T_r(N)=
 \begin{cases}
 N/\ell,&N\equiv0\pmod\ell,\\
 (mN+r_j)/\ell,&N\equiv j\not\equiv0\pmod\ell.
 \end{cases}
\]
We call the residue system \emph{admissible} if, in addition,
$mj+r_j>1$ for $1\le j<\ell$, as in \cite{GGM2025}. This ensures that
$T_r$ maps positive integers to positive integers. Its action on
$\Z_\ell$, and the algebraic statements about its coalescence chain,
are also defined without this positivity assumption. Throughout,
$\ell,m\ge2$ and $\gcd(\ell,m)=1$. Negative drift \eqref{eq:H} is
assumed precisely where indicated below.

The transition identities and the telescoping identity
\eqref{eq:telescope} do not depend on the choice of residues. Put
$R=\max_{1\le j<\ell}|r_j|$. In the affine estimates for the normalized
displacement $p$, the appropriate uniform bound is
\[
 |\delta|\le\delta_*:=\frac R\ell
       \max\left\{2,\frac{m+1}{m-1}\right\}.
\]
The factor $2$ allows for opposite signs of $r_j$ and $r_{j'}$ in a
transition that stays at level zero. With this replacement, the
Lyapunov and tightness arguments of Lemmas~\ref{lem:GL}
and~\ref{lem:Gtight} apply under \eqref{eq:H}. The finite holding-time
argument also applies: if $B=0$ at a nonzero level $k$, some branch has
$B'\ne0$. Otherwise $r_{m^kj}=m^kr_j$ for every nonzero residue $j$;
iteration around a finite orbit of multiplication by $m^k$ modulo
$\ell$ would force $r_j=0$, contrary to its congruence modulo $\ell$.

The offset estimate in Lemma~\ref{lem:GTerras}, for an orbit segment
with $q$ nonzero branches, becomes
\[
 \left|T_r^n(a)-\frac{m^q a}{\ell^n}\right|
 \le \frac{R(m^q-1)}{\ell(m-1)}.
\]
Thus the renewal comparisons have endpoint errors $O_r(1)$. For
admissible signed residues, the finite-defect arguments use
\[
 T_r^i(N)\ge\frac{N}{\ell^i}-\frac{R}{\ell-1}
 \qquad(N\ge1).
\]
In particular $N>\ell^n(K+R/(\ell-1))$ ensures that every iterate up to
time $n$ exceeds $K$. These are the modifications needed in place of
the corresponding canonical-residue bounds.

\begin{definition}\label{def:dr}
The \emph{obstruction index} $d_r$ is the largest positive divisor
coprime to $\ell m$ of
\[
 \gcd\bigl(r_2-r_1,\dots,r_{\ell-1}-r_1,(\ell-1)r_1\bigr).
\]
For $\ell=2$ the list of differences is empty. Equivalently, $d_r$ is
the largest positive integer coprime to $\ell m$ for which
$r_j\equiv\rho\pmod{d_r}$ for all $j\ne0$ and
$(\ell-1)\rho\equiv0\pmod{d_r}$.
\end{definition}

For the canonical residues of \eqref{eq:Tlm}, $d_r=1$: when $\ell=2$
one has $r_1=1$, and when $\ell\ge3$ two residues differ by $1$.
Finite congruence descriptions and invariant residue sets for
generalized Collatz maps have a classical antecedent in
Matthews~\cite{Matthews1992}. The cocycle below concerns the
balanced pair chain; it is not an identification with the
single-trajectory classification in that work.

\begin{proposition}[The obstruction]\label{thm:GO}
Define, in $\Z[1/m]$,
\[
 \sigma(k)=
 \begin{cases}
 \displaystyle\sum_{i=0}^{k-1}m^i,&k\ge0,\\
 \displaystyle-\sum_{i=k}^{-1}m^i,&k<0,
 \end{cases}
 \qquad \Phi(k,B)=B-\rho\sigma(k).
\]
At every transition of the pair chain,
\[
 \Phi(k',B')\equiv\frac{m^{e_{j'}}}{\ell}\Phi(k,B)\pmod{d_r}.
\]
Consequently, if $d_r\nmid c$, balanced coalescence from $(0,c)$ is
impossible, and $f_n(c)=f(c)=1$ for every $n\ge0$.
Moreover, if $p\mid d_r$ is prime and $m\not\equiv1\pmod p$, then
$\{N\in\Z:N\equiv\rho/(1-m)\pmod p\}$ and its complement in $\Z$
are $T_r$-invariant. For an admissible residue system, their
intersections with $\N$ are invariant as well. Neither negative
drift nor admissibility is needed for the conclusions on $\Z$.
\end{proposition}

\begin{proof}
Reduction of $\Z[1/m]$ modulo $d=d_r$ is legitimate since
$\gcd(d,m)=1$. Put $k'=k+e_{j'}-e_j$. The identity
\[
 \sigma(k')=m^{e_{j'}}\sigma(k)+e_{j'}-m^{k'}e_j
\]
holds in each of the four branch types. Substitute it in the
transition formula of Lemma~\ref{lem:Gchain}, using
$r_j\equiv\rho e_j\pmod d$ and $\ell\rho\equiv\rho\pmod d$, to obtain
the displayed congruence. Its multiplier is a unit modulo $d$.
Since $\Phi(0,c)=c$ and $\Phi(0,0)=0$, absorption is impossible if
$d\nmid c$.

For the last statement put $\lambda=\rho/(1-m)\pmod p$. Then
$(\ell-1)\lambda=0$ and, in either branch,
\[
 T_r(N)-\lambda\equiv
       \frac{m^{e_j}}\ell(N-\lambda)\pmod p.
\]
The multiplier is a unit, proving invariance of both sets.
\end{proof}

\begin{theorem}[A finite absorption certificate]\label{thm:GR}
Assume \eqref{eq:H} and let $r$ be admissible. Set
\[
 C_0(r)=\frac{R(\ell+1)}{\ell^2-m}.
\]
The following conditions are equivalent:
\begin{enumerate}
 \item[(i)] $f(c)=0$ for every $c\in d_r\Z$;
 \item[(ii)] \textup{(condition $(\mathrm A_r)$)} for every integer $c$
 with $d_r\mid c$ and $0<|c|\le C_0(r)$, a finite residue word leads
 from $(0,c)$ to $(0,0)$.
\end{enumerate}
When these conditions hold, absorption is almost sure from every
relation state satisfying $\Phi(k,B)\equiv0\pmod{d_r}$.
Condition $(\mathrm A_r)$ is a finite collection of reachability
requirements. Successful words give a finite positive certificate;
no terminating decision procedure for a failed requirement is
asserted.
\end{theorem}

\begin{proof}
Hypothesis \eqref{eq:H} implies $m<\ell^2$, so $C_0(r)>0$.
Condition (i) implies (ii), because absorption is the union of its
finite-time events, and each finite residue word has positive
probability.

Assume (ii). We first construct an absorbing word from every
$(0,c)$ with $d_r\mid c$. If $c\ne0$ and $\ell\mid c$, branch $0$
replaces $c$ by $c/\ell$. Suppose $\ell\nmid c$. For $c>0$, use
branch $0$, which leads to
\[
 (1,B_1),\qquad B_1=\frac{mc+r_a}{\ell}>0,
 \qquad a\equiv c\pmod\ell,\quad 1\le a<\ell.
\]
The strict inequality follows by writing $c=a+\ell t$, $t\ge0$, and
using $ma+r_a>1$. Divide by $\ell$ using zero branches until the
displacement is not divisible by $\ell$, then choose the nonzero
branch returning to level zero. For $c<0$, first choose
$j\equiv-c\pmod\ell$, $j\ne0$. This leads to level $-1$, with
\[
 mB_1=\frac{mc-r_j}{\ell}<0.
\]
Again divide using zero branches until the displacement is not
divisible by $\ell$, then choose branch $0$ to return to level zero.
These sign choices ensure that the displacement used in the division
stage is nonzero. In both cases the resulting integer $c'$ satisfies
\[
 |c'|\le\frac m{\ell^2}|c|
              +\frac{R(\ell+1)}{\ell^2}.
\]
The right side is smaller than $|c|$ when $|c|>C_0(r)$.
Proposition~\ref{thm:GO} preserves $d_r\mid c$ at every return to
level zero. Iterating this descent and then using (ii) supplies the
required absorbing word.

From an arbitrary relation state with $\Phi\equiv0$, the level
reaches zero almost surely unless absorption has already occurred:
the holding times at nonabsorbed states are finite almost surely,
and the embedded level walk is simple symmetric until absorption.
In particular some finite word reaches a state $(0,c)$ with
$d_r\mid c$. Thus every such relation state has a finite absorbing
word.

It remains to upgrade reachability to almost-sure absorption.
The block Lyapunov inequality of Lemma~\ref{lem:GL}, with the
corrected constant $\delta_*$, implies repeated visits to a set
$\{|p|\le M\}$, for a sufficiently large fixed $M$, on every
nonabsorbed path outside a null set. This follows also by stopping
the inequality at the first entrance: outside a sufficiently large
sublevel set, the conditional drift of $|p|^s$ is bounded above by a
negative constant. Lemma~\ref{lem:Gtight}, stopped at absorption or
the first visit to $|k|\le1$, gives a uniform positive probability of
absorption or entrance into a fixed finite set of states with
$|k|\le1$ and bounded $|p|$. Only states with $\Phi\equiv0$ are relevant.
Each member of this finite set has an absorbing word by the preceding
paragraph, so the eventual absorption probability is bounded below
by one positive constant on $\{|p|\le M,\ \Phi\equiv0\}$.
Finally, conditional eventual-absorption probabilities form a
bounded martingale and converge to the indicator of eventual
absorption. The repeated visits with this uniform positive lower
bound rule out a nonabsorbed path. This proves the all-state
assertion and (i).
\end{proof}

The quantitative coalescence arguments also extend on the invariant
set $\Phi\equiv0$: the Lyapunov and return estimates use $\delta_*$,
and the absorption trials on a fixed finite good set use the words
just constructed. For the zero-level sojourn estimate of
Lemma~\ref{lem:Gsoj}(a), the divisibility phase is unchanged.
With $a$ and $\beta$ as in that proof, put
$C_r=1+(m+2R)/\ell$. Then
\[
 1+|B''|\le C_r(1+|B'|),\qquad
 \beta_{\rm next}\le a\beta+\log_\ell C_r.
\]
The same convergent geometric sum therefore bounds the sojourn
by $O_r(1+\log(2+|B|))$. This proof also applies when
$\ell=2$, since the probability of a surviving further round
is then zero. The averaged valuation argument of
Lemma~\ref{lem:Ggaps} also extends: replace
$\gcd(r_j,\ell^a)\le\ell-1$ by $\gcd(r_j,\ell^a)\le R$.
After this bounded common divisor is removed, the same
multiplicative-order argument gives
$\sum_{1\le |k|\le H}t_k\ll_r H$.
Consequently Proposition~\ref{thm:Gexact} and
Theorem~\ref{thm:Gsharp} hold for nonzero $c\in d_r\Z$, with constants
depending on the fixed residue system and a certificate for
$(\mathrm A_r)$. If $d_r=1$, the preceding Liouville and
density-transfer arguments apply with their stated hypotheses and
the endpoint and finite-defect modifications above. Numerical
constants established for canonical residues are not inherited.

\subsection{Rescaled residue systems and the maps $3x+r$}

For a nonzero integer $v$ coprime to $\ell$, put
$(v\star r)_j=v r_{v^{-1}j}$ for $j\ne0$, where the index is taken
modulo $\ell$, and put $(v\star r)_0=0$. This satisfies the same
congruence conditions, and $(uv)\star r=u\star(v\star r)$.
Admissibility must be checked separately, especially when $v<0$.
Let $\hat r$ denote the canonical residues of \eqref{eq:Tlm}.

\begin{lemma}[Scaling]\label{lem:scale}
Multiplication by $v$ is a Haar-measure-preserving homeomorphism of
$\Z_\ell$, and
\[
 T_{v\star r}(vx)=vT_r(x).
\]
The pair chain for $T_{v\star r}$ started at $(k,vB)$ and driven by
the residues $vj$ is the image of the chain for $T_r$ started at
$(k,B)$ under $(k,B)\mapsto(k,vB)$. In particular,
\[
 f^{(v\star r)}_n(vc)=f^{(r)}_n(c)
 \qquad(n\ge0,\ c\in\Z).
\]
\end{lemma}

\begin{proof}
The conjugacy follows directly from the definition and the fact
that $v$ permutes residues modulo $\ell$. If
$j'\equiv m^kj+B\pmod\ell$, then
$vj'\equiv m^k(vj)+vB\pmod\ell$, with the same branch indicators.
The transition formula gives
\[
 \frac{m^{e_{j'}}vB+(v\star r)_{vj'}-m^{k'}(v\star r)_{vj}}\ell
 =v\frac{m^{e_{j'}}B+r_{j'}-m^{k'}r_j}\ell.
\]
Uniform driving residues remain uniform under this permutation.
\end{proof}

\begin{corollary}[Classification for rescaled systems]\label{thm:SC}
Assume \eqref{eq:H}. Let $v=uw\ne0$ be coprime to $\ell$, where
$w>0$ has all its prime factors among those of $m$, and
$\gcd(u,\ell m)=1$. For $r=v\star\hat r$, one has
\[
 d_r=|u|,\qquad
 f(c)=
 \begin{cases}
  0,&u\mid c,\\
  1,&u\nmid c.
 \end{cases}
\]
Every shift divisible by $d_r$ therefore has a finite absorbing
word. If $r$ is admissible, condition $(\mathrm A_r)$ holds and
the general-residue conclusions following
Theorem~\ref{thm:GR} apply.
For signed $v$ without admissibility, the classification concerns
the dynamics on $\Z_\ell$.
\end{corollary}

\begin{proof}
The gcd in Definition~\ref{def:dr} is unchanged when a different
nonzero residue is chosen as the reference residue.
Consequently, permuting the nonzero residues leaves this gcd
unchanged, whereas multiplying all residues by $v$ multiplies
it by $|v|$.
The canonical residues have gcd $1$ in that definition.
Thus $d_r$ is the largest divisor of $|v|$ coprime to $\ell m$,
which is $|u|$.

By Theorem~\ref{thm:GK}, the canonical system satisfies
$f^{(\hat r)}(1)=0$.
Lemma~\ref{lem:scale} therefore gives
\[
 f^{(v\star\hat r)}(v)=f^{(\hat r)}(1)=0.
\]
Since $v\ne0$, the rescaled system has a nonzero almost-surely
coalescing shift.
Theorem~\ref{thm:ST}, proved below, now yields
\[
 \mathcal G_r=\mathcal D_r=d_r\Z.
\]
Hence $f(c)=0$ for $d_r\mid c$ and $f(c)=1$ otherwise,
which proves the asserted classification.

For each fixed $c\in d_r\Z$, almost-sure absorption implies
that absorption occurs with positive probability at some finite
time. Since each finite-time event is determined by finitely
many residue words, at least one such word is absorbing.
In particular, if $r$ is admissible, every reachability
requirement in $(\mathrm A_r)$ is satisfied.
The remaining conclusions then follow from
Theorem~\ref{thm:GR} and the extension described after its proof.
\end{proof}

\begin{corollary}[The maps $3x+r$]\label{cor:3xr}
Let $r$ be odd and define $T_r(N)=N/2$ for even $N$ and
$T_r(N)=(3N+r)/2$ for odd $N$. Then
\[
 d_r=\frac{|r|}{3^{v_3(r)}}.
\]
For Haar-almost every $N\in\Z_2$, the points $N$ and $N+c$ merge
with equal multiplication counts if $d_r\mid c$; if $d_r\nmid c$,
such balanced merging is impossible for every $N$.
In particular, Conjecture~\ref{conj:Ar} holds for $\ell=2$ under
\eqref{eq:H}. For $r=3^a$, $a\ge0$, the system is admissible and
$d_r=1$, so the corresponding conclusions following
Theorem~\ref{thm:GR} apply. For any odd $r$, the lattice $d_r\Z$ is
$T_r$-invariant and
\[
 T_r(d_rN)=d_rT_{r/d_r}(N),\qquad r/d_r=\pm3^{v_3(r)}.
\]
\end{corollary}

\begin{proof}
For $\ell=2$, hypothesis \eqref{eq:H} and coprimality force $m=3$.
Every binary residue system is $r\star\hat r$ with $\hat r_1=1$.
Apply Corollary~\ref{thm:SC} with
$v=r=\operatorname{sgn}(r)d_r3^{v_3(r)}$. The lattice identity
follows immediately from $d_r\mid r$ and the oddness of $d_r$.
\end{proof}

The binary classification in Corollary~\ref{cor:3xr} also follows
directly from Leen's admissible-offset theorem
\cite[Theorem~7.3]{LeenRates2026}. Indeed, multiplication by the odd
integer $r$ preserves Haar measure, and
$T_r(rx)=rT_{1}(x)$. The translated offset is $c/r$, which belongs
to $\Z[1/3]$ exactly when $d_r\mid c$. We include the corollary to
record this known special case in the present notation.

\begin{remark}\label{rem:scaleinv}
Scaling preserves the existence of a nonzero almost-surely
coalescing integer shift. Indeed, $f^{(r)}(c)=0$ implies
$f^{(v\star r)}(vc)=0$; conversely, if a nonzero shift $a$ is almost
sure for the scaled system, so is its integer multiple $va$, and
the scaling identity gives $f^{(r)}(a)=0$. Here we use the subgroup
property proved in Theorem~\ref{thm:ST}. Consequently, by
Corollary~\ref{cor:STprime}, condition $(\mathrm A_r)$ is invariant
under rescaling whenever both systems are admissible and satisfy
\eqref{eq:H}. The exact finite-time identity in
Lemma~\ref{lem:scale} also transfers the canonical estimates to
shifts of the form $vc$; for example, Theorem~\ref{thm:Gsharp} gives
\[
 f^{(v\star\hat r)}_n(vc)
 \le \min\left\{1,C_{\ell,m}
       \sqrt{\frac{1+\log(2+|c|)}{n+1}}\right\}.
\]
This identity by itself concerns shifts divisible by $v$.
\end{remark}

For general $\ell\ge3$, the preceding rescaling argument does not
establish $(\mathrm A_r)$ for every admissible system. Reachability
can also be expressed additively. For a word
$w=(w_1,\dots,w_n)\in(\Z/\ell\Z)^n$ with $q$ nonzero letters, put
\[
 A_w=\sum_{t:\,w_t\ne0}r_{w_t}\ell^{t-1}m^{q_{>t}},
\]
where $q_{>t}$ is the number of nonzero letters after position $t$.

\begin{proposition}[Additive word criterion]\label{prop:Acrit}
For $c\ne0$, the state $(0,0)$ is reachable from $(0,c)$ if and
only if there are distinct words $w,w'$ of the same length $n$ and
the same number $q$ of nonzero letters such that
\[
 m^qc=A_w-A_{w'}.
\]
\end{proposition}

\begin{proof}
On the cylinder with forward residue word $w$, induction gives
\[
 \ell^nT_r^n(N)=m^qN+A_w.
\]
For a balanced merger, subtract the two identities to obtain the
asserted relation. The words must be distinct, since equal words
would give $m^qc=0$.

Conversely, the inverse branch associated with $w$ is
\[
 x_w(z)=\frac{\ell^nz-A_w}{m^q},\qquad z\in\Z_\ell.
\]
Successive reduction modulo $\ell$ verifies that its forward
residues are exactly $w$ and that $T_r^n(x_w(z))=z$. Its image is
one cylinder modulo $\ell^n$. If $m^qc=A_w-A_{w'}$, then
$x_{w'}(z)=x_w(z)+c$ for every $z$. Thus this whole cylinder yields
a balanced merger, and its Haar measure is $\ell^{-n}$.
\end{proof}

Under \eqref{eq:H}, Conjecture~\ref{conj:Ar} is therefore equivalent
to saying that these integer quotients, together with zero, exhaust
$d_r\Z$. In the binary canonical case, $2A_w$ is the quantity $X_S$
in \eqref{eq:XS}.

\begin{conjecture}\label{conj:Ar}
Every admissible residue system satisfying \eqref{eq:H} has
property $(\mathrm A_r)$. Equivalently, the congruence obstruction
of Proposition~\ref{thm:GO} is the only obstruction to almost-sure
balanced coalescence. Corollary~\ref{cor:3xr} proves this for
$\ell=2$.
\end{conjecture}

\begin{remark}\label{rem:Gobstr}
For $3x+5$, one has $d_r=5$ and $\rho=0\pmod5$, so the multiples of
$5$ form an invariant set. A different example is
$(\ell,m)=(3,5)$ with $(r_1,r_2)=(1,-1)$. Here $d_r=2$ and
$\Phi(k,B)\equiv B+k\pmod2$. Neither parity class is invariant, but
an odd shift cannot merge with equal multiplication counts. This
qualification is essential: $T_r(1)=T_r(6)=2$ is an unbalanced
collision.

The latter system has $R=1$ and $C_0(r)=1<d_r=2$. Hence its finite
certificate has no requirements and $(\mathrm A_r)$ holds. Exactly
the even shifts coalesce almost surely. For example the driving
word $(0,0,1)$ gives
\[
 (0,2)\longrightarrow(1,3)\longrightarrow(1,1)
       \longrightarrow(0,0).
\]
These balanced-coalescence obstructions are compatible with the
logarithmic-density theorem of \cite{GGM2025}, which concerns small
values along individual orbits.
\end{remark}

\subsection{The structure of coalescence shifts}

In this subsection only the congruences on $r$, $\ell,m\ge2$, and
$\gcd(\ell,m)=1$ are assumed. Neither positivity nor negative drift
is needed. Let $\lambda$ be normalized Haar measure on $\Z_\ell$.
Write $\beta\sim\gamma$ if, at some common time, their iterates
agree and their accumulated numbers of nonzero branches agree.
This is an equivalence relation: once both equalities hold, they
continue to hold, which also proves transitivity. Define
\[
 A_c=\{\beta\in\Z_\ell:\beta\sim\beta+c\},\qquad
 f(c)=1-\lambda(A_c),
\]
and
\[
 \mathcal G_r=\{c\in\Z:f(c)=0\},\qquad
 \mathcal D_r=\{c\in\Z:f(c)<1\}.
\]
Let $\mathcal S_r$ consist of the nonzero integer quotients
$(A_w-A_{w'})/m^q$ from Proposition~\ref{prop:Acrit}.
Thus $\mathcal G_r$ describes almost-sure coalescence and
$\mathcal D_r$ describes coalescence with positive probability.

\begin{theorem}[One-shift rigidity]\label{thm:ST}
For every residue system as above:
\begin{enumerate}
 \item[(a)] $\mathcal G_r=g_r\Z$ for an integer $g_r\ge0$,
 $d_r\mid g_r$, and $f(c+g)=f(c)$ for all $c\in\Z$ and
 $g\in\mathcal G_r$.
 \item[(b)] If $c\in\mathcal G_r$ and $\ell\mid c$, then both
 $c/\ell$ and $mc/\ell$ belong to $\mathcal G_r$. If $g_r>0$,
 then $\gcd(g_r,\ell m)=1$.
 \item[(c)] $\mathcal D_r=\mathcal S_r\cup\{0\}$,
 $\mathcal D_r=-\mathcal D_r$, $\ell\mathcal D_r\subseteq\mathcal D_r$,
 $\mathcal D_r+\mathcal G_r=\mathcal D_r$, and
 $\mathcal D_r\ne\{0\}$.
 \item[(d)] Either $\mathcal G_r=\{0\}$ or
 \[
  \mathcal G_r=\mathcal D_r=\mathcal S_r\cup\{0\}=d_r\Z.
 \]
 In the latter case $f(c)\in\{0,1\}$ for every integer $c$.
\end{enumerate}
In particular, one nonzero almost-surely coalescing shift forces
the complete obstruction classification, also for composite $\ell$.
\end{theorem}

\begin{proof}
Translation by $c$ maps $A_c$ onto $A_{-c}$, so $f(-c)=f(c)$.
Transitivity gives
\begin{equation}\label{eq:STincl}
 A_c\cap(A_{c'}-c)\subseteq A_{c+c'}.
\end{equation}
Consequently $\mathcal G_r$ is a subgroup of $\Z$, and
Proposition~\ref{thm:GO} gives $\mathcal G_r\subseteq d_r\Z$.
Applying \eqref{eq:STincl} with a shift $g\in\mathcal G_r$ and
then with $-g$ proves $\lambda(A_{c+g})=\lambda(A_c)$. This proves
(a).

At $(0,\ell a)$, branch $0$ leads to $(0,a)$ and every nonzero
branch leads to $(0,ma)$. Hence the exact recursion is
\[
 f(\ell a)=\frac1\ell f(a)+\frac{\ell-1}{\ell}f(ma)
 \qquad(a\in\Z).
\]
It proves the first assertion of (b). Suppose now that $g=g_r>0$.
Apply this assertion to $c=\operatorname{lcm}(g,\ell)$:
\[
 \frac{g}{\gcd(g,\ell)}=\frac c\ell\in g\Z.
\]
Thus $\gcd(g,\ell)=1$, without any primality assumption on $\ell$.

For each $c$, the almost-sure equivalence with the shift $g$ implies
\[
 \mathbf1_{A_c}(\beta+g)=\mathbf1_{A_c}(\beta)
 \quad\hbox{for almost every }\beta.
\]
Indeed, outside null sets the chain of equivalences
$\beta+g\sim\beta\sim\beta+c\sim\beta+c+g$ proves the implication
in one direction; equality of the measures gives the reverse one.
The translation stabilizer of $\mathbf1_{A_c}$ is closed in
$\Z_\ell$, since translations are continuous in $L^1$.
It contains $g\Z$, which is dense because $\gcd(g,\ell)=1$.
Thus $\mathbf1_{A_c}$ is invariant under every translation in
$L^1$. Averaging over translations and applying Fubini shows that
it is constant almost everywhere. Therefore $f(c)\in\{0,1\}$ and
$\mathcal D_r=\mathcal G_r$ whenever $g>0$.

If $ma\in\mathcal G_r$, the exact recursion now forces
$a\in\mathcal G_r$: otherwise $f(\ell a)=1/\ell$, contrary to the
zero--one conclusion. Take $a=g/\gcd(g,m)$; since $ma\in g\Z$, we
obtain $a\in g\Z$ and hence $\gcd(g,m)=1$. This completes (b).

For (c), the finite-time absorption events are unions of cylinders.
Proposition~\ref{prop:Acrit} therefore gives
$\mathcal D_r=\mathcal S_r\cup\{0\}$. Symmetry was proved above.
If $\beta\sim\beta+c$, then $\ell\beta\sim\ell\beta+\ell c$
after one additional zero branch, so
$\ell A_c\subseteq A_{\ell c}$ and
$\lambda(A_{\ell c})\ge\lambda(A_c)/\ell$.
Periodicity in (a) gives $\mathcal D_r+\mathcal G_r=\mathcal D_r$.

To prove that $\mathcal D_r$ is nontrivial, note that for sufficiently
large $n$,
\[
 \sum_{q=0}^n\binom nq\left(\frac{\ell-1}{m}\right)^q
 =\left(1+\frac{\ell-1}{m}\right)^n>n+1.
\]
For some $q$ there are therefore more than $m^q$ words of length
$n$ with $q$ nonzero letters. Two distinct such words have
$A_w\equiv A_{w'}\pmod{m^q}$. They cannot have $A_w=A_{w'}$:
the inverse-branch formula in Proposition~\ref{prop:Acrit} would
then give the same input with two different forward residue words.
Their quotient is a nonzero member of $\mathcal S_r$.

It remains to identify $g$ when $g>0$. Choose any nonzero residues
$j,j'$. Since $\gcd(g,\ell)=1$, some $c\in g\Z$ satisfies
$c\equiv j'-j\pmod\ell$. Branch $j$ from $(0,c)$ leads to
\[
 \left(0,\frac{mc+r_{j'}-r_j}{\ell}\right).
\]
Every positive-probability successor of an almost-surely absorbed
state is almost surely absorbed. Its displacement belongs to
$g\Z$, so $g\mid r_{j'}-r_j$. Write $\rho$ for the common residue
of all nonzero $r_j$ modulo $g$.

Fix $a\not\equiv0\pmod\ell$. As $c=a+\ell t$ varies, the integer
$B_1=(mc+r_a)/\ell$ takes every residue modulo $\ell$.
Since $\gcd(g,\ell^2)=1$, the Chinese remainder theorem permits
the simultaneous choices
\[
 c\in g\Z,\qquad c\equiv a\pmod\ell,\qquad
 B_1\equiv1\pmod\ell.
\]
Branch $0$ first leads to $(1,B_1)$. The next branch
$j\equiv-m^{-1}B_1\pmod\ell$ is nonzero and leads to
\[
 (0,B_2),\qquad B_2=\frac{B_1-r_j}{\ell}\in g\Z.
\]
It follows that $B_1\equiv\rho\pmod g$, whereas
$\ell B_1=mc+r_a\equiv\rho\pmod g$. Thus
$g\mid(\ell-1)\rho$. Together with $\gcd(g,\ell m)=1$, the
residue differences already obtained imply $g\mid d_r$.
Part (a) gives $d_r\mid g$, so $g=d_r$. This proves (d).
\end{proof}

\begin{corollary}[The one-shift criterion]\label{cor:STprime}
Assume \eqref{eq:H} and let $r$ be admissible. The following are
equivalent:
\begin{enumerate}
 \item[(i)] condition $(\mathrm A_r)$ holds;
 \item[(ii)] some nonzero integer shift coalesces almost surely;
 \item[(iii)] $\mathcal G_r=d_r\Z$.
\end{enumerate}
Under these equivalent conditions,
$\mathcal S_r\cup\{0\}=\mathcal D_r=\mathcal G_r=d_r\Z$.
There is no restriction that $\ell$ be prime and no additional
word condition is required.
\end{corollary}

\begin{proof}
The equivalence of (i) and (iii) is Theorem~\ref{thm:GR} together
with Proposition~\ref{thm:GO}. Clearly (iii) implies (ii), and
Theorem~\ref{thm:ST} proves the converse and the last assertion.
\end{proof}

Thus Conjecture~\ref{conj:Ar}, within its stated negative-drift
class, reduces to the existence of one nonzero almost-sure shift.
Finite word searches, such as those discussed in
Appendix~\ref{app:Ar}, can supply positive certificates for
individual systems. Each application of Theorem~\ref{thm:GR}
also requires verification of admissibility and \eqref{eq:H}.
Such searches do not prove termination when a reachability
requirement fails. The rescaled systems and the example in
Remark~\ref{rem:Gobstr} give conclusions without a general search.

\begin{remark}\label{rem:STcomposite}
The use of $c=\operatorname{lcm}(g_r,\ell)$ in the proof excludes
every common prime factor of $g_r$ and $\ell$, including when
$\ell$ is composite. Once $\mathcal G_r\ne\{0\}$, there is
therefore no residual case of fractional coalescence probabilities:
all of them are zero or one. The structural theorem alone does
not decide whether $\mathcal G_r=\{0\}$; the existence assertion
remains the content of Conjecture~\ref{conj:Ar} for the systems
not covered above.
\end{remark}

\subsection{A uniform lower bound for the merging time}

Here $\ell,m\ge2$, $\gcd(\ell,m)=1$, and $r$ satisfies the residue
congruences. Neither admissibility nor \eqref{eq:H} is required.
The pair chain starts at $(0,c)$, and $f_n(c)$ is its probability
of nonabsorption after $n$ physical shortcut steps. The walk
survival probability $u_N$ is as in Lemma~\ref{lem:uN}.

\begin{theorem}[Uniform lower bound]\label{thm:lowerunif}
For every $c\ne0$ and every $n\ge1$,
\begin{equation}\label{eq:lowerunif}
 f_n(c)\ge\frac2\ell u_{n-1}
       \ge\frac2\ell\sqrt{\frac2{\pi(n+1)}}.
\end{equation}
If no state $(0,B)$ with $B\ne0$ can move to $(0,0)$ in one step,
then
\begin{equation}\label{eq:lowerunif1}
 f_n(c)\ge u_{n-1}\ge\sqrt{\frac2{\pi(n+1)}}.
\end{equation}
This additional condition holds for every binary residue system,
and for the canonical residues of \eqref{eq:Tlm} it holds exactly
when $m\ge\ell-1$.
\end{theorem}

\begin{proof}
From a nonzero level, absorption requires a level change that
reaches zero. The successive level changes form a simple symmetric
walk up to absorption, continued afterwards with independent fair
signs. In particular, conditionally on any stopping time at level
$1$ or $-1$, the probability of avoiding level zero for the next
$a$ level changes is $u_a$. Avoiding those changes also ensures
survival for the next $a$ physical steps.

For a nonzero integer $B$, write
$w_\ell(B)=\max\{a\ge0:\ell^a\mid B\}$. At level zero with
$\ell\mid B$, the next displacement is $B/\ell$ or $mB/\ell$.
It remains nonzero and its $w_\ell$ value decreases by exactly one.
Thus, for $a=w_\ell(c)$, the first $a$ steps cannot absorb, and
at time $a$ the displacement is not divisible by $\ell$.
Conditionally on that past, the next step leaves level zero with
probability exactly $2/\ell$, arriving at $1$ or $-1$.
If $n\le a$, survival is certain. If $n\ge a+1$, condition on this
departure at the fixed time $a+1$ and apply the preceding walk
bound to obtain
\[
 f_n(c)\ge\frac2\ell u_{n-a-1}\ge\frac2\ell u_{n-1}.
\]
This does not condition on any future survival event and permits
direct absorption on level zero through other branches.

Under the additional condition, no absorption can occur before
the first departure time $\sigma$ from level zero. On
$\{\sigma>n\}$ survival is certain; on $\{\sigma=t\}$, $1\le t\le n$,
the stopping-time argument gives conditional survival at least
$u_{n-t}\ge u_{n-1}$. Summing these disjoint alternatives proves
the improved bound. Lemma~\ref{lem:uN} gives the numerical lower
bounds in both cases.

For $\ell=2$, a nonzero even displacement remains nonzero at level
zero and an odd displacement must change level, so the additional
condition always holds. In general, a direct level-zero absorption
must use two nonzero residues $j,j'$ and satisfy
\[
 mB=r_j-r_{j'},\qquad j'\equiv j+B\pmod\ell.
\]
For canonical residues their nonzero differences have magnitudes
at most $\ell-2$, which excludes this when $m\ge\ell-1$.
If $m\le\ell-2$, choose canonical residues with
$r_j-r_{j'}=m$ and put $B=1$. Their congruences give
$j'\equiv j+1\pmod\ell$, so this is a direct absorbing branch.
\end{proof}

\begin{remark}\label{rem:lowerunif}
For the standard Collatz map, \eqref{eq:lowerunif1} is
Proposition~\ref{prop:lower}. At $c=1$, $f_1(1)=f_2(1)=1$, whereas
$u_0=1$ and $u_1=1/2$. Thus equality with $u_{n-1}$ holds at
$n=1$ but not at $n=2$. In general the factor $2/\ell$ in
\eqref{eq:lowerunif} accounts for the possibility of absorption
on level zero before its first departure.
\end{remark}

\begin{corollary}[Moments of the merging time]\label{cor:moments}
Assume \eqref{eq:H}, let $r$ be admissible, and suppose that
$(\mathrm A_r)$ holds; this includes all canonical systems
$T_{\ell,m}$. Let $q_n(x)$ count the nonzero branches among the
first $n$ iterates, and put
\[
 \tau_c^{(r)}(\beta)=\inf\{n\ge0:T_r^n(\beta)=T_r^n(\beta+c),\
                  q_n(\beta)=q_n(\beta+c)\},
\]
with the infimum of the empty set equal to $\infty$. For Haar
distributed $\beta$, every nonzero $c\in d_r\Z$ and every
$\alpha>0$ satisfy
\[
 \E\bigl(\tau_c^{(r)}\bigr)^\alpha<\infty\quad\Longleftrightarrow\quad
             0<\alpha<\tfrac12.
\]
More precisely, for a constant $C\ge1$ depending only on the fixed
system, and $L_c=1+\log(2+|c|)$,
\[
 \E\bigl(\tau_c^{(r)}\bigr)^\alpha\le
       \frac{(C^2L_c)^\alpha}{1-2\alpha}
       \quad(0<\alpha<\tfrac12),\qquad
 \E\log(1+\tau_c^{(r)})\le\log(1+C^2L_c)+2.
\]
The power $L_c^\alpha$ is optimal in order as a uniform bound in
$c$: for $c=d_r\ell^a$, $a\ge0$, one has $\tau_c^{(r)}\ge a+1$
pointwise. For $d_r\nmid c$, every positive moment is infinite
because $\tau_c^{(r)}=\infty$; for $c=0$ the merging time is zero.

Without any drift or admissibility assumption, every residue
system, $c\ne0$, and $n\ge1$ satisfy
\[
 \frac1{\ell^n}\sum_{0\le N<\ell^n}\min\{\tau_c^{(r)}(N),n\}
 =\sum_{t=0}^{n-1}f_t(c)
 \ge\frac{4\sqrt2}{\ell\sqrt\pi}(\sqrt{n+1}-1).
\]
Here integer inputs are viewed in $\Z_\ell$. For $x\ge\ell^n$,
\[
 \frac1x\#\{1\le N\le x:N+c\ge1,\ \tau_c^{(r)}(N)>n\}
 \ge\frac2\ell\sqrt{\frac2{\pi(n+1)}}
       -\frac{\ell^n+\max\{-c,0\}}x.
\]
For admissible systems the pairs counted in this last display
have positive integer orbits.
\end{corollary}

\begin{proof}
For $\alpha>0$, the tail-sum formula gives
\[
 \E\bigl(\tau_c^{(r)}\bigr)^\alpha
 =\sum_{n\ge1}\bigl(n^\alpha-(n-1)^\alpha\bigr)f_{n-1}(c).
\]
The lower bound in Theorem~\ref{thm:lowerunif} makes this series
diverge for $\alpha\ge1/2$. For a coalescing shift,
Theorem~\ref{thm:Gsharp}, in its general-residue form described
after Theorem~\ref{thm:GR}, gives, after enlarging $C$ if necessary,
\[
 \Pb(\tau_c^{(r)}>t)\le\min\{1,C\sqrt{L_c/t}\}\qquad(t>0).
\]
Put $A=C^2L_c$. Integrating this estimate against
$\alpha t^{\alpha-1}\,dt$, splitting at $A$, gives
\[
 \E\bigl(\tau_c^{(r)}\bigr)^\alpha\le A^\alpha+
       \alpha\sqrt A\int_A^\infty t^{\alpha-3/2}\,dt
       =\frac{A^\alpha}{1-2\alpha}.
\]
Integration against $dt/(1+t)$ instead gives
\[
 \E\log(1+\tau_c^{(r)})
 \le\log(1+A)+\sqrt A\int_A^\infty t^{-3/2}\,dt
 =\log(1+A)+2.
\]
The first $w_\ell(c)$ steps stay on level zero with nonzero
displacement, as proved above; this supplies the stated deterministic
delay and the optimal order of the moment bound along
$c=d_r\ell^a$.

For the final assertions, the event $\{\tau_c^{(r)}>t\}$ is a union of
cylinders modulo $\ell^t$. Thus its frequency in a complete block
modulo $\ell^n$ is exactly $f_t(c)$ for $t\le n$.
Use $\min\{\tau_c^{(r)},n\}=\sum_{t<n}\mathbf1_{\{\tau_c^{(r)}>t\}}$, the
lower bound of Theorem~\ref{thm:lowerunif}, and $f_0(c)=1$.
The elementary inequality
\[
 \sum_{t=0}^{n-1}(t+1)^{-1/2}\ge2(\sqrt{n+1}-1)
\]
gives the displayed block bound. Counting complete blocks of length
$\ell^n$ in $[1,x]$ loses at most $\ell^n$ inputs; imposing
$N+c\ge1$ loses at most $\max\{-c,0\}$ more. This proves the last
inequality.
\end{proof}

\begin{corollary}[A necessary total time scale]\label{cor:Lnec}
Let $0<\varepsilon<1$ and $L,n_0\ge0$. Suppose that
$f_n(1)\le\varepsilon$ for every integer
$n\ge L\log3+n_0$, as follows in particular from the uniform
bounds of Propositions~\ref{prop:U} and~\ref{prop:GU}. Then
\[
 L\log3+n_0\ge\frac8{\pi\ell^2}\varepsilon^{-2}-2.
\]
If no nonzero level-zero state can be absorbed in one step, the
constant $8/(\pi\ell^2)$ can be replaced by $2/\pi$. This includes
all binary systems and the canonical systems with $m\ge\ell-1$.
The conclusion constrains the combined threshold
$L\log3+n_0$, not $L$ separately.

For the standard Collatz map and $j\ge1$, the particular upper-bound
expression in Proposition~\ref{thm:M} satisfies
\[
 \max_{0<|c|\le3^j/2}f_n(c)
       \ge\sqrt{\frac2{\pi(n+1)}}.
\]
This is a lower bound on that coalescence-based upper-bound
expression; it is not a lower bound on $D_n^{(j)}$ itself.
\end{corollary}

\begin{proof}
Apply Theorem~\ref{thm:lowerunif} with $c=1$ and
$n=\max\{1,\lceil L\log3+n_0\rceil\}$. Then
$n+1\le L\log3+n_0+2$, which gives both constants. The maximum
in the last assertion includes $c=1$.
\end{proof}

Theorems~\ref{thm:Gsharp} and~\ref{thm:lowerunif} show that every
fixed nonzero coalescing shift has a tail of exact order
$n^{-1/2}$. In particular its small-$\varepsilon$ hitting quantile
has order $\varepsilon^{-2}$. The upper estimate in
Theorem~\ref{thm:Gsharp} also permits a logarithmic-displacement
coefficient of order $\varepsilon^{-2}$. Corollary~\ref{cor:Lnec}
does not establish a matching necessary lower bound for that
coefficient when a separate additive term $n_0(\varepsilon)$ is
allowed.
\section{\textbf{Open problems}}
\label{sec:problems}

We conclude with quantitative questions and extensions of the results
above. We first recall the established conclusions that delimit these
questions.

For the standard Collatz map, Leen
\cite[Theorem~2]{LeenBrownian2026} proves that every set invariant above
a finite threshold has a natural density. As explained after
Proposition~\ref{thm:T}, this includes our finite-defect invariant sets.
Their natural and logarithmic densities therefore exist and agree.
In particular, a comparison of these two limiting densities can use
the identity modulus $\Phi(\delta)=\delta$; improving the constants in
the comparison argument given here is a separate quantitative matter.
Likewise, a power-saving exceptional-set estimate for descent below
$N^\theta$ is already known for every fixed $\theta>0$
\cite[Definition~2.6 and Corollary~3.2]{Inselmann2024}; the uniform
fixed-barrier estimate of Nashida \cite[Theorem~1.1]{Nashida2026}
also gives such a conclusion.

The generalized square-root tail bound is proved here without linear
forms in $p$-adic logarithms; see
Proposition~\ref{thm:Gexact} and Theorem~\ref{thm:Gsharp}.
The latter also permits an effective choice
$L(\varepsilon)=O_{\ell,m}(\varepsilon^{-2})$ with $n_0=0$.
Corollary~\ref{cor:Lnec} gives a necessary lower bound on the combined
threshold $L(\varepsilon)\log3+n_0(\varepsilon)$, not on
$L(\varepsilon)$ alone when $n_0$ is unrestricted.

\begin{problem}[The first-level mixing rate]
Determine the asymptotic decay of $D_n^{(1)}$.
In particular, does it admit an exponential upper bound
$D_n^{(1)}\le C e^{-cn}$ with $c>0$?
The finite values in Table~\ref{tab:D} motivate studying this question,
but do not determine an asymptotic rate.
The square-root lower bound for the coalescence-based upper-bound
expression in Corollary~\ref{cor:Lnec} is not a lower bound for
$D_n^{(1)}$ itself.

One possible route is to estimate the normalized collision excess in
Proposition~\ref{thm:L2}.  With $Q_n\sim\operatorname{Bin}(n,1/2)$,
a sufficient condition for $D_n^{(1)}\to0$ is
\[
 E_{\mathrm{top}}(n,Q_n)\longrightarrow0
 \quad\hbox{in probability},
\]
or equivalently
\[
 \sum_{q=0}^n2^{-n}\binom nq
       \min\{1,E_{\mathrm{top}}(n,q)\}\longrightarrow0.
\]
Quantitative versions would give a corresponding mixing bound.
This condition concerns the entire quantity
\[
 E_{\mathrm{top}}(n,q)=3^qB_{n,q}(3R_{n,q}-1),
\]
not just $R_{n,q}-1/3$.
The untruncated mean cannot tend to zero: its $q=n$ term alone is
$2(3/2)^n$.
\end{problem}

\begin{problem}[Leading constants for generalized merging times]
Fix $\ell\ge3$, an admissible residue system satisfying
\eqref{eq:H} and $(\mathrm A_r)$, and a nonzero
$c\in d_r\mathbb Z$.
Does the limit
\[
 \kappa_r(c):=\lim_{n\to\infty}\sqrt n\,f_n(c)
\]
exist?  Theorems~\ref{thm:Gsharp} and~\ref{thm:lowerunif}, together
with the extension following Theorem~\ref{thm:GR}, bound this expression
above and below by positive constants for each fixed $c$.
If the limit exists, seek a representation by excursion occupation
or a harmonic capacity, and an effective method for approximating it.

The binary cases are already covered by Leen's admissible-offset
result \cite[Theorem~7.3]{LeenRates2026} and the scaling identity
$T_r(rx)=rT_1(x)$. If $r$ is odd and $d_r\mid c$, then
$c/r\in\Z[1/3]$, and that theorem gives the exact positive leading
coefficient for the balanced shortcut merging time.
In particular, at displacement $1$ for the standard map,
\[
 \sqrt n\,f_n(1)\longrightarrow
 \kappa_{\mathrm{std}}=\frac{2\mathbb E M}{\sqrt\pi}>0,
\]
where $M$ is the integrable excursion count defined there
\cite[Theorem~A]{LeenRates2026}.
Regular variation of this tail also gives the domain of attraction
of a positive $1/2$-stable law for sums of independent copies.
The question here concerns the generalized systems with $\ell\ge3$.
\end{problem}

\begin{problem}[Absorption for general residue systems]
Resolve Conjecture~\ref{conj:Ar} for noncanonical admissible residue
systems with $\ell\ge3$: does negative drift always imply
$(\mathrm A_r)$, or is there a counterexample?
Theorem~\ref{thm:GR} reduces the question to the finite family of
reachability requirements
\[
 (0,c)\longrightarrow(0,0),\qquad
 d_r\mid c,\qquad
 0<|c|\le \frac{R(\ell+1)}{\ell^2-m}.
\]
Successful finite words provide a certificate, but the theorem does
not give a terminating procedure for a failed requirement.

Theorem~\ref{thm:ST} supplies a second formulation, valid for
composite as well as prime $\ell$: it is enough to find one nonzero
shift with almost-sure balanced coalescence.
No additional verification of $d_r\in\mathcal S_r$ is needed.
A single successful word gives only positive probability of
coalescence and does not by itself meet this criterion.
Also, when $d_r>1$, the individual classes modulo $d_r$ need not be
invariant; any further componentwise equidistribution statement must
use the actual action on these classes, as illustrated in
Remark~\ref{rem:Gobstr}.
\end{problem}

\begin{problem}[Unconditional comparison with the primes]
Remove the Riemann hypothesis from Theorem~\ref{thm:PE}.
More precisely, for each fixed $K_0$, prove, or find an obstruction to,
\[
 \sup_{\substack{G\subseteq\mathbb N\\K(G)\le K_0}}
 \left|\frac{\pi_G(x)}{\pi(x)}-\frac{F_G(x)}x\right|
 \longrightarrow0.
\]
Here the supremum is over finite-defect $T$-invariant sets.
For each fixed $G$ of natural density zero or one, the equality
already follows from Theorem~\ref{thm:P}, applied also to the complement.
The question is the comparison for general invariant sets, including
the possible intermediate-density case.

For the proof used here, it would suffice to establish the averaged
short-interval estimate
\[
 \frac1H\sum_{v\le x}
  |\vartheta(v-1)-\vartheta(v-H-1)-H|=o(x),
 \qquad H=\lfloor x^\eta\rfloor,
\]
for arbitrarily small fixed $\eta>0$; here $\vartheta(y)=0$ for $y<2$.
An estimate for a single fixed positive exponent is not enough to
send the coalescence error in this proof to zero.
The unconditional range $\eta>2/15$ in
\cite[Theorem~1.1(ii)]{GafniTao2026} does not reach these arbitrarily
small exponents.
This identifies a sufficient analytic input for the present argument,
not a necessary condition on every possible proof of the comparison.
\end{problem}

\begin{problem}[Values and approximation of basin densities]
\label{prob:dens}
For the positive $3x-1$ map, let $B_1,B_5,B_{17}$ be the basins of
the three known cycles, indexed by their least elements, and write
\[
 d_i=d(B_i),\qquad i\in\{1,5,17\}.
\]
These densities exist: Sultanow
\cite[Theorem~1]{Sultanow2026} supplies logarithmic-density existence,
and Corollary~\ref{cor:basins} identifies each with its natural density.
The logarithmic-density theorem is the unconditional input from that
source; the natural-density conclusion here also uses our comparison.

Determine the values $d_1,d_5,d_{17}$, or give an effective procedure
which, for every rational $\varepsilon>0$, produces rational intervals
of length at most $\varepsilon$ certified to contain them, with an
explicit error bound or termination criterion usable in computation.
Numerical frequencies without a proven remainder bound do not provide
such a certification.
As a separate question, is
\[
 d_1+d_5+d_{17}=1?
\]
Existence of the three densities does not imply this equality.
The question asks whether their union has natural density one; it
must not be treated as an assumption that these cycles exhaust the
positive integer orbits.
\end{problem}

\begin{problem}[Sharper arithmetic and density remainders]
\label{prob:rate}
For a finite-defect $T$-invariant set $G$ and a fixed $\Lambda>1$,
determine sharp bounds for
\[
 \omega_G(x;\Lambda)=
 \sup_{\Lambda^{-1}\le\lambda\le\Lambda}
 \left|\frac{F_G(\lambda x)}{\lambda x}-\frac{F_G(x)}x\right|.
\]
An effective modulus already follows from Proposition~\ref{thm:exact}
and Lemma~\ref{lem:R}, using a finite net of translations generated
by $\log2$ and $\log3$. The question concerns sharper rates and their
optimal dependence on the defect cutoff and on $\Lambda$, including
uniform bounds when the cutoff is bounded.

For each fixed odd modulus $d$, Proposition~\ref{thm:rate} gives
\[
 \max_{a\bmod d}
 \left|\frac{dF_a^{[d]}(x)}x-\frac{F_G(x)}x\right|
 =O_{d,K(G)}((\log x)^{-1/2}).
\]
Can arithmetic cancellation improve this bound, and can a power
saving in $x$ hold for a meaningful uniform class of invariant sets?
The square-root lower bound for the probabilistic coupling error
does not supply a matching lower bound for this arithmetic discrepancy.

One concrete related target is an effective remainder in the already
convergent standard-map defect integral:
\[
 R_G(X)=\int_X^\infty
       \frac{F_{2,1}(x)-F_G(x)/3}{x^2}\,dx.
\]
Determine whether this integral is absolutely convergent for every
such $G$, or identify additional hypotheses which give absolute
convergence and a usable bound for $R_G(X)$.
Density existence and Proposition~\ref{thm:defect} give ordinary
convergence, whereas the pointwise logarithmic bound above does not
by itself give an integrable absolute majorant.
\end{problem}

\appendix
\section{\textbf{Exact finite computations and certificates}}
\label{app:numerics}

This appendix records the finite computations used for explicit
constants and the residue-system certificates discussed in
Section~\ref{sec:general}.
The accompanying programs specify the enumerated objects,
their normalization, and the success criteria. All mathematical
accept-or-reject decisions use integers or rational numbers.
The mixing counts and descent checks certify the stated finite
bounds. Each successful residue-system search supplies an
absorbing word whose validity is checked independently; its
consequence for almost-sure coalescence uses the hypotheses of
Theorem~\ref{thm:GR}.

\subsection{Exact mixing distances at levels one and two}\label{app:D}

We use the total variation convention and the distances defined in
Sections~\ref{sec:full} and~\ref{sec:level1}. For a fixed $n,q$, put
\[
 N_{n,q}(c,d)=\#\left\{S\subseteq\{1,\ldots,n\}:|S|=q,\quad
 X_S\equiv c+3^qd\pmod{3^{q+1}}\right\},
\]
where $0\le c<3^q$, $d\in\Z/3\Z$, and $X_S$ is given by
\eqref{eq:XS}. Empty child cells have count zero. The binomial
normalization in \eqref{eq:Dt} cancels the normalization within
each fixed-cardinality slice, giving the integer identity
\[
 \begin{split}
 2^nD_n^{(1)}
 &=\sum_{q=0}^n\sum_{c=0}^{3^q-1}
       \frac12\sum_{d\bmod3}
          |N_{n,q}(c,d)-N_{n,q}(c,d-1)|\\
 &=\sum_{q=0}^n\sum_{c=0}^{3^q-1}
       \left(\max_{d\bmod3}N_{n,q}(c,d)
                    -\min_{d\bmod3}N_{n,q}(c,d)\right).
 \end{split}
\]
The second equality is the cyclic variation identity for three
real numbers.

Starting from $X_\varnothing=0$, the recursion
\[
 X_{S\cup\{m\}}=3X_S+2^m\qquad(m>\max S)
\]
generates every subset exactly once. For each $q$ separately,
sort the keys
\[
 3(X_S\bmod3^q)
       +\left(\left\lfloor X_S/3^q\right\rfloor\bmod3\right).
\]
Their multiplicities give the three child counts for each parent.
This produces Table~\ref{tab:D}. In particular,
\[
 D_{24}^{(1)}=\frac{3\,585\,706}{2^{24}}
            =\frac{1\,792\,853}{8\,388\,608}<0.2138,
\]
the finite bound used in Proposition~\ref{thm:quant}.

\begin{table}[ht]
 \centering
 \begin{tabular}{rrrrrrrr}
  \toprule
  $n$ & $12$ & $14$ & $16$ & $18$ & $20$ & $22$ & $24$\\
  \midrule
  $2^nD_n^{(1)}$
   & $1\,383$ & $5\,118$ & $18\,945$ & $70\,227$
   & $260\,286$ & $965\,955$ & $3\,585\,706$\\
  $D_n^{(1)}$
   & $0.33765$ & $0.31238$ & $0.28908$ & $0.26789$
   & $0.24823$ & $0.23030$ & $0.21372$\\
  \bottomrule
 \end{tabular}
 \caption{Exact numerator counts for $D_n^{(1)}$.
 The decimal row is rounded to five places.}
 \label{tab:D}
\end{table}

For the full level-two comparison, it is necessary to include
pairs with different parents modulo $3$. Define
\[
 M_{n,q}(c,d)=\#\left\{S\subseteq\{1,\ldots,n\}:|S|=q,\quad
 X_S\equiv c+3^qd\pmod{3^{q+2}}\right\},
\]
with $0\le c<3^q$ and $d\in\Z/9\Z$, and put
\[
 A_n(h)=\frac12\sum_{q=0}^n\sum_{c=0}^{3^q-1}
               \sum_{d\bmod9}|M_{n,q}(c,d)-M_{n,q}(c,d-h)|.
\]
The shift $h$ is fixed throughout the sum over $q$.
By \eqref{eq:XS}, after multiplying outputs by $2^{n+1}$,
the cylinders with centers $a,a'$ modulo $9$ differ by
$2\cdot3^q(a'-a)$ modulo $3^{q+2}$. The unit factor $2$
permutes the possible values of $h$. Since
$A_n(9-h)=A_n(h)$, it suffices to check $h=1,2,3,4$.
In the refined notation of Section~\ref{sec:sadic},
\[
 \begin{split}
 2^nD_n^{(1,1)}&=A_n(3),\\
 2^nD_n^{(0,2)}&=\max\{A_n(1),A_n(2),A_n(4)\},\\
 2^nD_n^{(2)}&=\max_{1\le h\le4}A_n(h).
 \end{split}
\]
Sorting nine-child keys in the same way gives
\[
 \begin{array}{c|rrrr}
 h&1&2&3&4\\\hline
 A_{20}(h)&533\,480&408\,347&563\,831&509\,275.
 \end{array}
\]
Consequently,
\[
 \begin{split}
 D_{20}^{(1,1)}&=\frac{563\,831}{1\,048\,576}
                 =0.537711143493652\ldots,\\
 D_{20}^{(0,2)}&=\frac{533\,480}{1\,048\,576}
                 =\frac{66\,685}{131\,072},\\
 D_{20}^{(2)}&=D_{20}^{(1,1)}.
 \end{split}
\]
The final equality follows from this comparison of all shifts
at $n=20$. Counting only parents modulo $3^{q+1}$ and the next
digit computes $D_{20}^{(1,1)}$ directly; the additional
comparisons establish the value of the full distance.

The accompanying program
\texttt{Appendix\_Exact\_Mixing\_Checks.py} reproduces every entry
using integer arithmetic. It requires Python and a C++17 compiler;
the C++ source is included in the Python file and is compiled in
a temporary directory. Fixed-cardinality slices are processed
separately, so the largest stored slice consists of
$\binom{24}{12}=2\,704\,156$ unsigned $64$-bit keys, with a key
payload of $21\,633\,248$ bytes. Since
\[
 0\le X_S\le 2\cdot3^n-2^{n+1}<2^{40}\qquad(n\le24),
\]
the generation is exact without modular truncation or integer
overflow. The program checks each slice mass against $\binom nq$,
the endpoint slices $q=0,n$, the total mass $2^n$, and the
displayed rational bounds. A separate reference calculation at
$n=12$ uses the descending-position formula \eqref{eq:XS} and
compares the full translated histograms.

The same counts also check the second-moment normalization in
Proposition~\ref{thm:L2}:
\[
 E_{\mathrm{top}}(n,q)
 =\frac{3^q}{\binom{n}{q}^{2}}
   \left(3\sum_{c,d}N_{n,q}(c,d)^2
          -\sum_c\left(\sum_dN_{n,q}(c,d)\right)^2\right).
\]
In particular, $E_{\mathrm{top}}(n,n)=2\cdot3^n$, and its single
contribution to the untruncated binomial mean is
$2(3/2)^n$. Any application of that proposition must therefore
use its stated convergence-in-probability or truncated-energy
condition. The finite computations provide exact benchmarks;
they do not determine an asymptotic decay rate for $D_n^{(1)}$.

\subsection{Certified bounds for the descent and return constants}
\label{app:K0}

The finite base of Lemma~\ref{lem:polydesc} consists of the
\(1998\) nonzero integers \(c\) with \(|c|<1000\).
For even \(c\), use the outcome \(0\) and replace \(c\) by
\(c/2\). For odd \(c\ne-1\), put \(a=v_2(3c+1)\), use the
word \(0^a1\), and replace \(c\) by
\[
 c'=\frac{(3c+1)/2^a-1}{2}.
\]
For \(c=-1\), the word \(100\) absorbs; for \(c=1\), the
rule gives \(001\). Apply the recursion until the offset is zero,
and denote the length of the resulting word by \(\ell(c)\).

The accompanying script
\texttt{Collatz\_Coalescence\_Constants\_Check.py}
replays every word in the exact relation chain and verifies
that its terminal state is \((0,0)\).
It then checks, for every \(c\) with \(\ell(c)\ge3\), the
integer inequality
\[
 2^{25(\ell(c)-3)}\le |c|^{121}.
\]
When \(\ell(c)<3\), the desired length bound is immediate.
Thus all \(1998\) cases satisfy
\[
 \ell(c)\le\frac{121}{25}\log_2|c|+3.
\]
Equality holds exactly at \(c=\pm1\).
The longest checked word has length \(28\), attained, for
example, at \(c=-913\).
In particular, the additive constant \(K_0=3\) in the finite
part of the descent estimate is certified without numerical
evaluation of logarithms.

For the return constant in Lemma~\ref{lem:tight}, let
\(r=99/100\) and \(t_k=v_2(3^k-1)\). The series is
\[
 C_*=1+\sum_{k\ge2}2(k-1)r^{k-1}(t_k+2).
\]
The identity
\[
 t_k=1+2\mathbf1_{2\mid k}
          +\sum_{a\ge2}\mathbf1_{2^a\mid k}
\]
and the rational sums
\[
 S_d=\sum_{m\ge1}(dm-1)r^{dm-1}
    =\frac{r^{d-1}(d-1+r^d)}{(1-r^d)^2}
\]
give
\[
 C_*=1+\frac{6r}{(1-r)^2}+4S_2
          +2\sum_{a\ge2}S_{2^a}.
\]
The script evaluates the part through \(a=8\) as a rational
number and bounds the remaining dyadic tail, with \(D=512\), by
\[
 \frac{2D r^{D-1}}{(1-r^D)^2(1-2r^D)}.
\]
The exact rational comparisons are
\[
 \begin{split}
 1+\frac{6r}{(1-r)^2}+4S_2+2\sum_{a=2}^{8}S_{2^a}
   &<88957.632,\\
 \frac{2D r^{D-1}}{(1-r^D)^2(1-2r^D)}&<6.167.
 \end{split}
\]
They imply \(C_*<88963.799<88964\).
This tail estimate is proved in Lemma~\ref{lem:tight}; the
finite arithmetic check certifies its numerical substitution.

The same script verifies the fractional-moment constants by
enclosing each positive tenth root between rational numbers and
checking the tenth powers of the endpoints. It certifies
\[
 0.99868400<\rho<0.99868401,\qquad
 \rho_V<0.98722,\qquad C_V<43222,
\]
together with the weighted-product inequalities used in the
return and visit estimates. All accept-or-reject decisions use
integers or rational numbers; decimals printed by the script
are for presentation.

\subsection{Condition \(\mathrm A_r\)}\label{app:Ar}

For each triple \((\ell,m,R_0)\) in Table~\ref{tab:Arsearch},
we enumerated every integer residue vector
\[
 (r_1,\ldots,r_{\ell-1}),\qquad
 r_j\equiv-mj\pmod\ell,\qquad |r_j|\le R_0,
\]
with \(r_0=0\).
We recorded separately the \emph{admissible} vectors satisfying
\[
 mj+r_j>1\qquad(1\le j<\ell).
\]
This distinction is necessary when applying
Theorem~\ref{thm:GR}, whose hypotheses include admissibility.
All six pairs satisfy \(\gcd(\ell,m)=1\) and \eqref{eq:H}.
In the order of the table, the exact comparisons
\(m^{\ell-1}<\ell^\ell\) are
\[
 4<27,\quad16<27,\quad25<27,\quad27<256,
 \quad256<3125,\quad1296<3125.
\]

For each vector we computed \(R=\max_j|r_j|\), the obstruction
index \(d_r\) of Definition~\ref{def:dr}, and the exact rational
cutoff
\[
 C_0(r)=\frac{R(\ell+1)}{\ell^2-m}.
\]
The required set of shifts was enumerated using the integer
inequality
\[
 \mathcal J_r=
 \{c\in\Z\setminus\{0\}:d_r\mid c,
                \ |c|(\ell^2-m)\le R(\ell+1)\}.
\]
Thus no floating-point comparison is involved at the boundary
of the finite certificate.

For every \(c\in\mathcal J_r\), we ran a separate breadth-first
search from \((k,B)=(0,c)\).
Both signs were searched directly; no negative-shift result was
inferred by doubling a count of positive shifts.
The implementation stores the integral state \((k,b)\), where
\(B=m^{\min(k,0)}b\).
Writing \(a=\min(k,0)\), the transition for an input residue
\(j\in\{0,\ldots,\ell-1\}\) is
\[
 \begin{split}
 j'&\equiv m^kj+m^ab\pmod\ell,\\
 k'&=k+e_{j'}-e_j,\qquad a'=\min(k',0),\\
 b'&=\frac{m^{e_{j'}+a-a'}b
             +r_{j'}m^{-a'}-r_jm^{k'-a'}}{\ell}.
 \end{split}
\]
All exponents in the numerator of \(b'\) are nonnegative,
and its division by \(\ell\) is exact.
The possibly negative powers in the first congruence are
computed by modular inversion, since \(m\) is a unit modulo
\(\ell\).
Reaching \((k,b)=(0,0)\) provides an absorbing residue word,
and hence a positive-probability cylinder of balanced coalescence.

The search limits were \(10^6\) discovered states and depth
\(64\) for each signed shift.
If either limit is reached before an absorbing word is found,
the program records the requirement as \emph{unresolved}; it
does not report a counterexample or a failure of
\((\mathrm A_r)\).
In the present complete run, every required search succeeded.
The largest number of discovered states in any one search was
\(21\,932\), and the longest absorbing word had length \(17\).
Word length counts the physical shortcut steps of the relation
chain.

\begin{table}[ht]
 \centering
 \begin{tabular}{ccrrrrrr}
  \toprule
  \((\ell,m)\)&\(R_0\)&
  \multicolumn{2}{c}{All congruence systems}&
  \multicolumn{2}{c}{Admissible systems}&
  \shortstack{Maximum\\word length}&
  \shortstack{Unresolved\\shifts}\\
  \cmidrule(lr){3-4}\cmidrule(lr){5-6}
  &&systems&signed shifts&systems&signed shifts&&\\
  \midrule
  \((3,2)\)&22&225&3350&64&944&17&0\\
  \((3,4)\)&22&225&4784&80&1566&9&0\\
  \((3,5)\)&16&121&1748&56&728&8&0\\
  \((4,3)\)&9&100&432&36&140&10&0\\
  \((5,4)\)&5&16&24&8&8&2&0\\
  \((5,6)\)&5&16&24&16&24&2&0\\
  \midrule
  Total&&703&10362&260&3410&17&0\\
  \bottomrule
 \end{tabular}
 \caption{Exhaustive enumeration and direct verification of the
 finite reachability requirements in \((\mathrm A_r)\).
 Both signs are included in each shift count.
 The last two columns concern the full congruence family;
 all admissible subfamilies also have no unresolved requirements.}
 \label{tab:Arsearch}
\end{table}

There are \(5\,181\) directly checked positive shifts and
\(5\,181\) directly checked negative shifts in total.
For \(33\) of the \(703\) vectors, the set \(\mathcal J_r\)
is empty; the finite condition is then vacuous, and no search
is counted.  Of these vectors, \(17\) are admissible.

The portable Python script
\texttt{Collatz\_Residue\_Certificates.py}
uses only the Python standard library and exact integer or
rational arithmetic.
Its \texttt{generate} command writes
\texttt{Collatz\_Residue\_Certificates.json}, which contains every
residue vector, its admissibility flag, \(d_r\), the rational
cutoff, the complete signed target list, and an absorbing word
for every target.
The \texttt{verify} command checks the enumeration and all the
words without calling the search transition function.
First, it replays the original rational \((k,B)\) recurrence.
Second, it reconstructs an initial cylinder modulo \(\ell^n\)
from the inverse branches of a recorded word of length \(n\),
and directly iterates both integer orbits from its representative
and that representative plus \(c\).
It checks the prescribed residues, equality of endpoints, and
equality of multiplication counts.
Both independent replays passed for all \(10\,362\) words.
The machine-readable verification report is
\texttt{Collatz\_Residue\_Certificates\_Verification.json}.

For example, the complete computation and its independent replay
are reproduced by
\begin{verbatim}
python Collatz_Residue_Certificates.py generate \
  --output Collatz_Residue_Certificates.json
python Collatz_Residue_Certificates.py verify \
  Collatz_Residue_Certificates.json \
  --report Collatz_Residue_Certificates_Verification.json
\end{verbatim}
The default limits are those stated above; reaching a limit
produces an explicit unresolved entry in the output.

It follows that the finite reachability requirements have been
verified for all \(703\) listed congruence systems.
For the \(260\) admissible systems, all hypotheses of
Theorem~\ref{thm:GR} have also been checked, so its almost-sure
coalescence conclusion for every shift in \(d_r\Z\), and
condition (ii) of Corollary~\ref{cor:STprime}, follow.
For the remaining \(443\) systems, this appendix asserts the
listed finite absorbing words and does not invoke those
almost-sure conclusions.

\end{document}